\documentclass[a4paper,10pt]{amsart}
\usepackage{tikz}
\usetikzlibrary{matrix,arrows,decorations.pathmorphing}
\usepackage{tikz-cd}
\usepackage{amsfonts, amssymb, amsthm, amsmath}  
\usepackage{todonotes}
\usepackage{hyperref}
\usepackage{graphicx}
\usepackage{psfrag}
\newtheorem{thm}{Theorem}  

\newtheorem{cor}[thm]{Corollary}  
\newtheorem{lemma}[thm]{Lemma}  
\newtheorem{remark}[thm]{Remark}  
\newtheorem{defn}[thm]{Definition}  
\newtheorem{prop}[thm]{Proposition}  
\newtheorem{claim}[thm]{Claim}  
\newtheorem{example}[thm]{Example}  
\newtheorem{construction}[thm]{Construction}
\numberwithin{thm}{section}  
\numberwithin{equation}{page}
\def\pf{\noindent\emph{Proof: }}  
\def\stop{\hfill$\square$}  
\def\x{^{\sharp}}
\def\co{\colon\thinspace}
\providecommand{\fun}[1]{\boldsymbol {#1}}

\providecommand{\ov}[1]{\hspace{-.1cm}\downarrow_{#1}}

\providecommand{\Y}{\mathcal Y}

\providecommand{\totl}[1]{\ensuremath{\lceil #1\rceil }}
\providecommand{\totb}[1]{\ensuremath{\underline{ #1}}}

\newcommand{\ex}{\bold}
\providecommand {\e}[1]{\mathfrak t^{#1}}
\providecommand{\bd}[1]{^{{}^{\rtimes#1}}}
\providecommand{\sfp}[3]{\ensuremath {#1}_{#2}^{#3}}
\providecommand{\C}[2]{\ensuremath {C^{#1,\underline{#2}}}}

\newcommand{\Mod}{\mathcal M}
\newcommand{\Ms}{\mathcal M^{\infty,\underline 1}}
\newcommand{\Msw}{\mathcal M^{st}}
\newcommand{\dmsw}{\mathcal M^{st}_{\bullet}}
\newcommand{\dmod}{\mathcal M_{\bullet}}
\newcommand{\M}{\mathcal M(pt)}  
\newcommand{\T}{\mathbf{T}}

\newcommand{\exte}{\subset_{e}}

\providecommand{\fp}[2]{{}_{\hspace{3pt}#1\hspace{-2pt}}\times_{#2}}

\DeclareMathOperator{\dist}{dist}  
\DeclareMathOperator{\id}{id}
\DeclareMathOperator{\expl}{Expl}
\DeclareMathOperator{\coker}{coker}
\newcommand{\dbar}{\bar{\partial}}  

\newcommand{\dvert}{d_{\text vert}}

\providecommand{\et}[2]{\ensuremath{\bold T^{#1}_{#2}}}
\providecommand{\lrb}[1]{\ensuremath{\left(#1\right)}}
\providecommand{\abs}[1]{\left\lvert #1\right\rvert}

\author{Brett Parker   }
  
\title{Holomorphic curves in  exploded manifolds  ---  Kuranishi structure}

\begin{document}

\begin{abstract}This paper constructs a Kuranishi structure for the moduli stack of holomorphic curves in  exploded manifolds. To avoid some technicalities of abstract Kuranishi structures, we embed our  Kuranishi structure inside an ambient moduli stack of not-necessarily-holomorphic curves.  The construction also works for the moduli stack of holomorphic curves in any compact symplectic manifold.
\end{abstract}

\maketitle

\tableofcontents

\thanks{This work was supported by ARC grant DP1093094.}
\newpage

\section{Introduction}

\

In this paper, we construct a species of Kuranishi structure on the moduli stack of holomorphic curves in exploded manifolds.\footnote{For a brief introduction to exploded manifolds, see \cite{scgp}; a more thorough introduction is  \cite{iec}, and \cite{elc} is a dictionary comparing log schemes and exploded manifolds. 
} Using \cite{vfc} we can then construct a virtual fundamental class from such a Kuranishi structure. Because the category of exploded manifolds extends that of smooth manifolds, our construction includes a new construction of Kuranishi structures on the moduli stack of holomorphic curves in any compact symplectic manifold.

\

Fukaya and Ono used Kuranishi structures in their definition of Gromov--Witten invariants in \cite{FO}. A similar approach to defining Gromov--Witten invariants was independently taken by Li and Tian in \cite{Tian-Li}, and much subsequent  work has since refined the construction and use of Kuranishi structures; see \cite{KFOOO,KMW, pardon, MW2, MW3, CLW, joycebook}.  Our Kuranishi structures differ from Fukaya and Ono's: our local models use exploded manifolds in place of manifolds, and they are embedded in an infinite dimensional moduli stack of not-necessarily-holomorphic curves, with a level of regularity, $\C\infty1$, appropriate for exploded manifolds.
 This ambient moduli stack plays the role of an ambient space containing the holomorphic curves; embedding our Kuranishi structures this way avoids some  technicalities involved in using Kuranishi structures to define Gromov--Witten invariants.\footnote{ See the preprint  \cite{KMW} by McDuff and Werheim for some discussion of issues that must be overcome when using abstract Kuranishi structures, and  see the recent preprint \cite{KFOOO} of Fukaya, Oh, Ohta, and Ono for some improvements on Fukaya and Ono's original definitions and a much more detailed version of their construction of Gromov--Witten invariants.} 

In \cite{FO},  Fukaya and Ono use a homotopy of the linearization, $D\dbar$, of the $\dbar$ operator to a complex map for orienting their Kuranishi structures. Moreover, in \cite{FOinteger}, they explain how this homotopy gives a stably-almost-complex structure, at least away points in the moduli space representing nodal curves.  They also sketch how, if their stably-almost-complex structure was globally defined, then they could define invariant integer counts of holomorphic curves. Similarly,  in the preprint \cite{kuranishihomology}, Dominic Joyce sketches a construction of integer invariants using  similar almost-complex information and a modified version of Kuranishi structures.  
In the category of smooth manifolds, the construction of a stably-almost-complex structure must be treated with great care around singular curves. In contrast,   using exploded manifolds, no curves are singular in this way, and we encounter no problems constructing  a globally-defined stably-almost-complex structure using a homotopy of $D\dbar$ to a complex map; Fukaya and Ono's suggested construction of integer invariants is carried out using exploded manifolds in \cite{icc}. A separate approach to realising Fukaya and Ono's construction  using smooth manifolds and global Kuranishi charts is \cite{AMS2021,Bai2025,Bai2022}.

\

Our Kuranishi structures consist of a collection of embedded Kuranishi charts. We identify holomorphic curves using  a section, $\dbar$, of a sheaf $\mathcal Y$ over an ambient moduli stack, $\Msw$, of stable not-necessarily-holomorphic $\C\infty1$ curves. In an open neighborhood $\mathcal O$ around a holomorphic curve $f$, we can choose a nice, finite-dimensional subsheaf $V\subset \mathcal Y$. If $D\dbar$ is transverse to $V$ at $f$, then in an open neighborhood $\mathcal U$ of $f$,  the moduli stack of solutions $ h$ to $\dbar  h\in V$ is  represented by the quotient of some $\C\infty1$ family of curves $\hat f$ by a group $G$ of automorphisms.  Our embedded Kuranishi charts have the information $(\mathcal U,V,\hat f/G)$. By putting further assumptions on $V$, we construct embedded Kuranishi charts with nice properties such as being compatible with chosen evaluation maps, or having an equivalent of Fukaya and Ono's stably-almost-complex structure.

To define an embedded Kuranishi structure, a collection of Kuranishi charts must obey compatibility and extendibility conditions. For two embedded Kuranishi charts $(\mathcal U_{i},V_{i},\hat f_{i}/G_{i})$  to be compatible, we require that, on the intersection of $\mathcal U_{1}$ with $\mathcal U_{2}$, one of $V_{1}$ or $V_{2}$ must be a subsheaf of the other. There are weaker versions of compatibility  sufficient for constructing Gromov--Witten invariants, but subject to compactness assumptions\footnote{Appropriate compactness for the moduli stack of holomorphic curves is proved in \cite{cem}.} on the moduli stack of holomorphic curves, a covering by such compatible embedded Kuranishi charts always exists, and any two such embedded Kuranishi structures are homotopic, so a weaker version of compatibility is unnecessary. 
To construct something on a Kuranishi structure, (such as the weighted branched perturbation of $\dbar$ used by Fukaya and Ono to define Gromov--Witten invariants), we often proceed chart by chart, shrinking the prior domain of definition slightly at each stage; see for example, Proposition 2.3 of \cite{vfc}. To facilitate such a procedure, we require that our Kuranishi charts have compatible extensions. 

\

In Section \ref{structure section}, we describe the ambient moduli stack, $\Ms$, of $\C\infty1$ not-necessarily-holomorphic curves, and an open substack, $\Msw$, of well-behaved stable curves. Reading Section \ref{structure section} is essential for understanding the results of this paper. Section \ref{stack section} contains an introduction to stacks over the category of exploded manifolds.   The tangent space of this ambient moduli stack $\Msw$ is defined and demystified in Section \ref{tangent section}, then, in Section \ref{Ddbar section} we define the linearization, $D\dbar$, of the $\dbar$ operator on $\Msw$. Finally, 
embedded Kuranishi structures are defined in Section \ref{K section}.

Section \ref{tpo section} contains a quick summary of necessary regularity results  from \cite{reg}. Sections  \ref{ev0 section} and \ref{+n section} construct some evaluation maps from the ambient moduli stack of stable, not-necessarily-holomorphic curves $\Msw$. Section \ref{core section} defines  a core family giving a concrete local model for this ambient moduli stack $\Msw$. Such concrete models for the infinite-dimensional ambient moduli stack are essential, and used  throughout the rest of the paper. Core families are constructed in Theorem \ref{core criteria} and Proposition \ref{smooth model family}. Section \ref{bump section} uses core families to study the topology of the ambient moduli stack $\Msw$. For example, in Lemma \ref{hausdorff}, we show that  the topology on $\Msw$ is pulled back from a Hausdorff topological space.  Section \ref{piV section} is dedicated to locally analyzing  the moduli stack of solutions to the weakened $\dbar$ equation, $\dbar f\in V$. The main theorem is Theorem \ref{V moduli stack}, where we show that $\dbar^{-1}V$ is locally represented by $\hat f/G$ for a family of curves $\hat f$ and finite group $G$.

Section \ref{K construction} constructs an embedded Kuranishi structure 
 for the moduli stack of holomorphic curves in any family of exploded manifolds  $\hat{\ex B}\longrightarrow \ex B_{0}$, subject to a compactness condition on the moduli stack of holomorphic curves.  See Theorem \ref{K existence} for a precise statement. In Corollary \ref{K homotopy}, we also prove that any two such Kuranishi structures are homotopic. We end the paper with a construction of a stably-almost-complex structure for our embedded Kuranishi structures.

\section{Structure of the ambient moduli stack of stable curves }\label{structure section}

In this section, we work towards describing  embedded Kuranishi structures  by first describing basic properties of the ambient moduli stack $\Ms$ of $\C\infty1$ not-necessarily-holomorphic curves, concentrating on a well-behaved open substack $\Msw\subset\Ms$ of stable $\C\infty1$ curves. We  use  $\Mod\subset \Msw$ for the moduli stack of stable holomorphic curves.

\subsection{The functors $\ex F$ and $\ex C$}

\

This paper studies  families of  holomorphic curves in a smooth family of targets in the exploded category, \[\pi_{\ex B_0}\co(\hat {\ex B},J)\longrightarrow\ex B_0\] where $\pi_{\ex B_0}\co \hat {\ex B}\longrightarrow \ex B_0$ is a smooth map of exploded manifolds\footnote{See definitions 3.1, 3.2, and 3.13 of \cite{iec} for the definition  of exploded manifolds and maps between them. For a family of exploded manifolds, see Definition 10.1, which relies on the definition of the tangent space of an exploded manifold, discussed in Section 6 of \cite{iec}.}, and each fiber of $\pi_{\ex B_{0}}$ is a complete\footnote{For `complete', see Definition 3.15 of \cite{iec}, relying on Definition 3.14 and Example 3.8.}, basic\footnote{For `basic', see Definition 4.6 of \cite{iec}, relying on definitions 4.4. and 4.3.} exploded manifold with an  almost-complex structure\footnote{Definition 8.1 of \cite{iec}.} $J$. For this paper, $\hat{\ex B}\longrightarrow \ex B_0$ will always refer to  a family of exploded manifolds with such structure.  

We will often talk about $\C\infty1$ families of curves\footnote{We use `curve' in a sense analogous to `prestable map' elsewhere in the literature; see Definition 8.3 of \cite{iec}.} $\hat f$ in $\hat {\ex B}\longrightarrow \ex B_0$. These are commutative diagrams
\[\begin{tikzcd}\ex C(\hat f)\arrow{r}{\hat f}\dar{\pi_{\ex F(\hat f)}}&\hat{\ex B}\dar{\pi_{\ex B_0}}
\\\ex F(\hat f)\arrow{r} &\ex B_0\end{tikzcd}\]
where $\hat f$, and all other maps in the above diagram are $\C\infty1$ maps of exploded manifolds\footnote{The regularity $\C\infty1$ is defined in Section 7 of \cite{iec}, and in particular Definition 7.5. Although $\C\infty1$ regularity is slightly weaker than smooth, is as good as smooth for our purposes. Note that smooth manifolds are a full subcategory of the category of exploded manifolds, and a map of exploded manifolds with domain a smooth manifold is $\C\infty1$ if and only if it is smooth. We will use that $\C\infty1$ exploded manifolds form a nice category with notion of transversality and fiber products analogous the case of smooth manifolds; see Section 9 of \cite{iec}. }, and
 $\pi_{\ex F(\hat f)}\co\ex C(\hat f)\longrightarrow \ex F(\hat f)$ is a family of curves\footnote{Definitions 11.1 and 11.3 of \cite{iec}} with a fiberwise complex structure\footnote{Definition 8.1 of \cite{iec}} $j$. In particular, a curve shall, by default,  refer to a map to an exploded manifold, from a domain $\ex C$ that is a complete 2-dimensional  exploded manifold with an almost complex structure $j$.  We will generally refer to a family of curves by the map $\hat f\co \ex C(\hat f)\longrightarrow \hat{\ex B}$, with the rest of the diagram above understood as part of the data. 
 
 We say that the exploded manifold $\ex F(\hat f)$ parametrises the family of curves, and the exploded manifold $\ex C(\hat f)$ is the domain of the family of curves.  A particular case is when the exploded manifold $\ex F(\hat f)$ is a single point. In this case, we call $\hat f$ a curve, and generally use the notation $f$ instead of $\hat f$ to emphasise that $f$ is a single curve instead of a family containing many curves. If we restrict a family $\hat f$ to the fiber over any point in $\ex F(\hat f)$, we get a curve $f$, and we say that such curves $f$ are the curves in the family $\hat f$, and use the notation $f\in\hat f$. For notational convenience, we will sometimes conflate such a curve $f$ in $\hat f$ with the corresponding point in $\ex F(\hat f)$ and write $f\in \ex F(\hat f)$ instead of $\ex F(f)\in \ex F(\hat f)$.
 
So, for a family of curves $\hat f$, $\ex C(\hat f)$ and $\ex F(\hat f)$ are exploded manifolds. Actually,  $\ex F$ and $\ex C$ define functors from the category of $\C\infty1$ families of curves to the category of $\C\infty1$ exploded manifolds, and $\pi_\ex F\co \ex C\Rightarrow\ex F$ is a natural transformation. Concretely, as defined in  \cite[Section 11]{iec}, a morphism  $\psi\co \hat f\longrightarrow \hat g$ between $\C\infty1$ families of curves is a commutative diagram
\begin{equation}\label{morphism}\begin{tikzcd}\ex C(\hat f)\rar{\ex C(\psi)}\dar{\pi_{\ex F(\hat f)}}\ar[bend left]{rr}{\hat f}&\ex C(\hat g)\rar{\hat g}\dar{\pi_{\ex F(\hat g)}}&\hat{\ex B}
\\\ex F(\hat f)\rar{\ex F(\psi)}&\ex F(\hat g)\end{tikzcd}\end{equation}
such that restricted to each fiber of $\pi_{\ex F(\hat f)}$ and $\pi_{\ex F(\hat g)}$, the map $\ex C(\psi)\co \ex C(\hat f)\longrightarrow\ex C(\hat g)$ is a holomorphic isomorphism.
 So,  $\ex C$ and $\ex F$ are functors from the category, $\Ms$, of $\C\infty1$ families of curves to the category of $\C\infty1$ exploded manifolds. Note that the morphism $\psi\co \hat f\longrightarrow \hat g$ is entirely determined by the map $\ex C(\psi)\co \ex C(\hat f)\longrightarrow \ex C(\hat g)$ of exploded manifolds. We will generally use the terminology of a `map' $\psi\co \hat f\longrightarrow \hat g$ instead of a `morphism'. 

\subsection{The ambient moduli stack of  curves}

\

The category, $\Ms$, of $\C\infty1$ families of curves together with the functor $\ex F$ to the category of exploded manifolds is a stack\footnote{For an approachable introduction to stacks over various categories, see \cite{stacks}.} over the category of exploded manifolds. To distinguish this infinite-dimensional moduli stack from the moduli stack of holomorphic curves, we will refer to $\Ms$ as an ambient moduli stack.  In this section, we verify that $\Ms$ is a stack. The reader unfamiliar with stacks should keep in mind that abstract properties defining stacks simply formalise  geometric intuition about how we can pull back families of curves, and also how we can glue together locally defined families of curves.   First, we discuss pullbacks of families of curves, verifying that $(\Ms, \ex F)$ is a category fibered in groupoids.

 The square in (\ref{morphism}) is a fiber product\footnote{Definition 9.1 of \cite{iec}.} diagram of $\C\infty 1$ exploded manifolds, moreover, Lemma 10.4 of \cite{iec} implies that given any $\C\infty1$ map  $h\co\ex A\longrightarrow \ex F(\hat g)$ of exploded manifolds, there exists a pulled back family of curves, $h^*\hat g$ parametrised by $\ex A$ with a map $\psi_h\co h^*\hat g\longrightarrow \hat g$ such that $\ex F(\psi_h)$ is $h$.
 \[\begin{tikzcd}\ex C(h^*\hat g)\rar{\ex C(\psi_h)}\dar \ar[bend left]{rr}{h^*\hat g}&\ex C(\hat g)\rar{\hat g}\dar &\hat{\ex B}
\\\ex F(h^*\hat g)=\ex A\rar{\ex F(\psi_h)=h}&\ex F(\hat g)\end{tikzcd}\]
Note that we can also pull back morphisms: the pullback of a map $\alpha\co \hat g\longrightarrow \hat f$ is the map $h^*\alpha\co h^*\hat g\longrightarrow \hat f$ defined as
\[h^*\alpha:=\alpha\circ \psi_h \ .\]
This pullback $(h^*\hat g,\psi_h)$ is unique up to a canonical isomorphism arising from the following:

\noindent{\bf Universal  property of the pullback of a family of curves}:  Given any other family of curves $\hat f$ with  maps $\psi'\co \hat f\longrightarrow \hat g$ and  $x\co \ex F(\hat f)\longrightarrow \ex A$ such that $\ex F(\psi')=h\circ x$, there exists a unique map $\phi\co \hat f\longrightarrow h^*\hat g$ such that $F(\phi)=x$ and the following diagram commutes.
\[\begin{tikzcd}\ex C(\hat f)\ar[bend left]{rr}{\ex C(\psi')}\rar[dotted]{\ex C(\phi)}\dar&  \ex C(h^*\hat g)\rar{\ex C(\psi_h)}\dar &\ex C(\hat g)\dar &
\\ \ex F(\hat f) \rar{x} & \ex F(h^*\hat g)=\ex A\rar{h}&\ex F(\hat g)\end{tikzcd}\]

Let $\Ms(\hat {\ex B})$ be the category of $\C\infty1$ families of curves in $\hat{\ex B}$. A category such as $\Ms(\hat{\ex B})$ with such a functor $\ex F$ and pullbacks satisfying the above universal property is called a category fibered in groupoids over the category of exploded manifolds; see  \cite[Definition 3.1]{stacks}. We show below that $(\Ms(\hat{\ex B}),\ex F)$ has further nice properties making it a stack  over the category of exploded manifolds, so we call $\Ms(\hat{\ex B})$ the ambient moduli stack of $\C\infty1$ curves in $\hat{\ex B}$.   When no ambiguity is present, we use $\Ms$ in place of $\Ms(\hat {\ex B})$ to refer to this as the ambient moduli stack of $\C\infty 1$ curves. 

The next lemma encapsulates  how we can glue together locally defined families of curves. 

\begin{lemma} Every descent datum for $\Ms$ is effective. In other words, given
\begin{itemize}
\item an exploded manifold $\ex A$ with an open\footnote{ See \cite[Definition 3.1]{iec} for the the standard topology on an exploded manifold $\ex A$. } cover $\{\iota_i\co U_i\longrightarrow \ex A\}$ with inclusions 
\[\iota_{ij}\co  U_i\cap U_j\longrightarrow U_i\ ,\]
\item families $\hat f_i$ in $\Ms$ parametrised by $U_i$,
\item and transition maps that are isomorphisms 
\[\alpha_{ij}\co \iota_{ji}^*\hat f_j\longrightarrow \iota_{ij}^*\hat f_i\]
satisfying the cocycle condtion that over $U_i\cap U_j\cap U_k$,  
\[\alpha_{ij}\circ\alpha_{jk}=\alpha_{ik}\ ,\]
\end{itemize}
there exists a family of curves $\hat f$ parametrized by $\ex A$ with isomorphisms 
\[\alpha_i\co \iota_i^*\hat f\longrightarrow \hat f_i\] such that 
\[\alpha_{ij}= \iota_{ji}^*\alpha_i \circ\lrb{\iota_{ij}^*\alpha_j}^{-1}  \]
\end{lemma}

\pf
Exploded manifolds are defined as topological spaces with a sheaf of functions satisfying local conditions (see definitions 3.1 and 3.13 of \cite{iec}), and similarly, families of curves can be defined as exploded manifolds, $\ex C(\hat f)$, with further structure ($\hat f$, $\pi_{\ex F(\hat f)}$, $j$) satisfying local conditions using the topology from $\ex F(\hat f)$. Accordingly, we can glue together families of curves over local charts on $\ex F(\hat f)$ as usual. The statement that every descent datum is effective is a formalisation of this familiar gluing procedure. 

\stop 

The next lemma is a formalisation of the obvious geometric fact that maps between families of curves are locally determined. 

\begin{lemma} Isomorphisms are a sheaf for $\Ms$. In other words, given
\begin{itemize}\item families $\hat f$ and $\hat g$ in $\Ms$;
\item an open cover $\{\iota_i\co U_i\hookrightarrow \ex F(\hat f)\}$ with intersections $\iota_{ij}\co U_i\cap U_j\longrightarrow U_i$;
\item and maps $\alpha_i\co \iota_i^*\hat f\longrightarrow \hat g$ such that
\[\iota_{ji}^*\alpha_j=\iota_{ij}^*\alpha_i\]
\end{itemize}
there exists a unique map $\alpha\co \hat f\longrightarrow \hat g$ such that
\[\iota_i^*\alpha=\alpha_i \ .\]
\end{lemma}

In particular, because the other sheaf axioms follow immediately from the existence of pullbacks, we get a sheaf of sets over $\ex F(\hat f)$ by assigning to the open set $U$ the set of maps of $\hat f\rvert_{U}$ to $\hat g$.

\pf

Morphisms $\psi$ of families of curves are maps, $\ex C(\psi)$, of sets satisfying conditions which are local in $\ex F(\hat f)$. It follows that we can uniquely glue together locally defined morphisms, so long as their restrictions are identical.

\stop

\subsection{Stacks over the category of exploded manifolds}

\label{stack section}

\

We have now verified that the ambient moduli stack  $\Ms$ satisfies the following standard definition of a stack over the category of $\C\infty1$ exploded manifolds; compare  \cite[Definition 4.3]{stacks}. This section spells out some standard notions of stacks in the case of stacks over the category of $\C\infty1$ manifolds. The reader comfortable with standard stack terminology can skip to Section \ref{Y section}.

\begin{defn} A {\bf stack} over the category of exploded manifolds is a category $\mathcal X$ with a functor $\ex F$ to the category of exploded manifolds, such that
\begin{itemize}\item $(\mathcal X,\ex F)$ is a category fibered in groupoids;
\item isomorphisms in $(\mathcal X,\ex F)$  form a sheaf;
\item and every descent datum is effective.
\end{itemize}\end{defn}

We shall think of an object $\hat f$ in the category $\mathcal X$ as a $\mathcal X$--family parametrised by the exploded manifold $\ex F(\hat f)$. For this reason, we will avoid using the singular word `object' to describe such an $\hat f$ in $\mathcal X$ unless $\ex F(\hat f)$ is a point. Usually our stacks will be moduli stacks of curves, whose objects are families of curves. 
 
 \begin{defn} An {\bf individual object}, or ${\bf curve}$ in a stack $\mathcal X$ over the category of exploded manifolds is an object $f$ of the category $\mathcal X$ with $\ex F(f)$ a point. A {\bf family} $\hat f$ in $\mathcal X$ is an object in the category  $\mathcal X$, and this family is said to be {\bf parametrized} by $\ex F(\hat f)$. The {\bf dimension} of the family $\hat f$ is $\dim\ex F(\hat f)$.
\end{defn}

\begin{example}\label{vector-bundle stack} The stack $\mathcal V$ of rank $n$ complex vector-bundles\footnote{As for smooth manifolds, a complex vector-bundle over an exploded manifold is a submersion $\ex V\longrightarrow \ex A$ of exploded manifolds with a $\mathbb C$--linear structure on each fiber so that there exist  local charts modeled on $\mathbb C^n\times \ex U\longrightarrow \ex U$ and transition maps that are $\mathbb C$--linear isomorphisms on fibers. } has families $\hat f$ rank $n$ complex  vector-bundles $V(\hat f)\longrightarrow \ex F(\hat f)$, and morphisms $\psi\co \hat f\longrightarrow \hat f'$ maps of vector-bundles
\[\begin{tikzcd} V(\hat f)\dar \rar{V(\psi)} &V(\hat f')\dar
\\ \ex F(\hat f) \rar{\ex F(\psi)} & \ex F(\hat f')\end{tikzcd}\]
  so that resticted to fibers,  $V(\psi)$ is a $\mathbb C$--linear isomorphism. In this case, the individual objects are $n$--dimensional complex vector-spaces. 
  
\end{example}

\begin{example}\label{S(A)} For an exploded manifold $\ex A$, there is a natural stack $\mathcal S(\ex A)$ whose families are $\C\infty 1$ maps of exploded manifolds  $\hat f\co \ex F(\hat f)\longrightarrow \ex A$ and whose morphisms are commutative diagrams of exploded manifolds
\[\begin{tikzcd}\ex F(\hat f)\dar \rar{\hat f} & \ex A
\\ \ex F(\hat f')\ar{ur}{\hat f'}\end{tikzcd}\]
Given any $\C\infty1$ map $x\co \ex A_1\longrightarrow \ex A_2$ of exploded manifolds, there is an induced functor $\mathcal S(x)\co \mathcal S(\ex A_1)\longrightarrow \mathcal S(\ex A_2)$ sending $\hat f$ to $x\circ \hat f$. An individual object in $\mathcal S(\ex A)$ is a map of a point into $\ex A$, so  isomorphism classes of individual objects correspond to points in $\ex A$.
\end{example}

\begin{example}\label{over f} For a family $\hat f$ in the stack $\mathcal X$, we can associate a new stack $\mathcal X/\hat f$, defined as follows:
\begin{itemize}\item families  in $\mathcal X/\hat f$ are morphisms $\hat h\longrightarrow \hat f$ within $\mathcal X$, 
\item   morphisms  in $\mathcal X/\hat f$ are commutative diagrams within $\mathcal X$,
\[\begin{tikzcd}\hat h\rar\dar{\psi} & \hat f
\\ \hat h'\ar{ur}\end{tikzcd}\]
\item the functor $\ex F$ from $\mathcal X/\hat f$ is  induced from the functor $\ex F$ from $\mathcal X$, so that $\ex F(\hat h\longrightarrow \hat f):=\ex F(\hat h)$, and $\ex F$ of the above morphism is $\ex F(\psi)$.
\end{itemize}

\end{example}

We can generalise Example \ref{S(A)} by considering the stack quotient of an exploded manifold $\ex A$ under the action of a finite group $G$. To describe this stack quotient, we need the following standard notions of an etale map and  a $G$--bundle. For this definition, recall that each exploded manifold $\ex E$ has a standard topology; \cite[Definition 3.1]{iec}. 

\begin{defn}An etale map $x\co \ex E'\longrightarrow \ex E$ of exploded manifolds is a map such that every point in $\ex E'$ has some open neighborhood $U'$  so that $x$ restricted to $U'$ is an isomorphism onto an open subset of $\ex E$. \end{defn}

Unlike the case of smooth manifolds, a map of exploded manifolds with bijective derivative is not necessarily etale. We use equidimensional submersion to refer to a map with bijective derivative. The following remark characterises when a equidimensional submersion is etale.

\begin{remark}\label{submersion etale} Recall the tropical structure of an exploded manifold,  \cite[Definition 4.4]{iec}. In particular, to each point in $x$ in an exploded manifold, we have an associated $\mathbb Z$--affine polytope $\mathcal P(x)$ such that $x$ is contained in a coordinate chart $U$ with tropical part $\totb U=\mathcal P(x)$. Maps of exploded manifolds induce functorial maps of tropical structure.  An equidimensional submersion $\pi\co \ex X\longrightarrow \ex Y$ of exploded manifolds is etale if and only if $\mathcal P\pi(x)\co \mathcal P(x)\longrightarrow \mathcal P(\pi(x))$ is a $\mathbb Z$--affine isomorphism for all $x\in \ex X_T$. Equivalently, there are local charts $U$ on $\ex X$ and $V$ on $\ex Y$ with $\pi(U)\subset V$ such that $\totb{\pi}\co \totb{U}\longrightarrow \totb{V}$ is a $\mathbb Z$--affine isomorphism. 
\end{remark}

\begin{defn}\label{G bundle} For $G$ a finite group, a {\bf $G$--bundle over an exploded manifold} $\ex X$ is an etale  map $\ex X\bd G\longrightarrow \ex X$ of exploded manifolds with a free $G$--action on $\ex X\bd G$ whose orbits are the  fibers of $\ex X\bd G\longrightarrow \ex X$.  A morphism of $G$--bundles  is a $G$ equivariant map $\ex X_1\bd G\longrightarrow \ex X_2\bd G$.

A {\bf $G$--bundle internal to $\mathcal X$} is  a morphism $\hat f\bd G\longrightarrow \hat f$ in $\mathcal X$ with a $G$--action on $\hat f\bd G$ so that $\ex F(\hat f\bd G)\longrightarrow \ex F(\hat f)$ is a $G$--bundle over the exploded manifold $\ex F(\hat f)$. The family $\hat f$ is called the {\bf base} of the $G$--bundle $\hat f\bd G$. A morphism of $G$--bundles internal to $\mathcal X$ is a $G$--equivariant morphism $\hat f\bd G\longrightarrow \hat g\bd G$.
\end{defn}

\begin{example} For $G$ a finite group, define the stack $\mathcal B G$ to be the category of $G$--bundles over exploded manifolds. The functor $\ex F$ applied to a $G$--bundle is the base of the bundle. 
\end{example}

\begin{defn}\label{stack quotient} For $\ex A$ an exploded manifold with an action of a finite group $G$, define the {\bf stack quotient} $\ex A/G$ to be a stack whose families are $G$--equivariant maps, $\hat f$, of $G$--bundles into $\ex A$
\[ \begin{tikzcd}\ex F\bd G(\hat f)\dar\rar{\hat f}& \ex A 
\\ \ex F(\hat f)\end{tikzcd}\]
and whose morphisms $\psi\co\hat f_1\longrightarrow \hat f_2$  are $G$--equivariant morphisms $\ex F\bd G(\psi)\co \ex F\bd G(\hat f_1)\longrightarrow \ex F\bd G(\hat f_2)$  of $G$--bundles such that the following diagram commutes
\[\begin{tikzcd} \ex F\bd G(\hat f_1)\dar\rar{\ex F\bd G(\psi)}\ar[bend left]{rr}{\hat f_1} & \dar \ex F\bd G(\hat f_2)\rar{\hat f_2} &\ex A
\\ \ex F(\hat f_1)\rar{\ex F(\psi)} & \ex F(\hat f_2)\end{tikzcd}\]
\end{defn}

\begin{defn}\label{stack map} A {\bf map}, or  {\bf morphism of stacks} $(\mathcal X_1,\ex F_1)\longrightarrow (\mathcal X_2,\ex F_2)$ over the category of exploded manifolds is a covariant functor 
\[\fun\Phi\co \mathcal X_1\longrightarrow \mathcal X_2\]
such that $\ex F_2\circ \fun\Phi=\ex F_1$.
\end{defn}

\begin{example}\label{m stack map} For $\ex A$ an exploded manifold, a  map $\ex A\longrightarrow \mathcal X$ means a morphism of stacks $\fun\Phi\co (\mathcal S(\ex A),\ex F)\longrightarrow (\mathcal X,\ex F)$. Such a morphism is equivalent to the information of a family $\hat f:=\fun\Phi(\text{id}_{\ex A})$ parametrized by $\ex A$ in $\mathcal X$, and a  choice\footnote{ Remember that pullbacks are only determined up to canonical isomorphism. The particular choice of pullback is irrelevant for all practical purposes. We have described the functor $\fun\Phi$ on families and morphisms to $\text{id}_{\ex A}$. On all other morphisms, $\psi$,  $\fun\Phi(\psi)$ is uniquely determined by the universal property of pullbacks.} of pullback $h^{*}\hat f$ for every map $h\co \ex A'\longrightarrow \ex A$. 
\end{example}

\begin{example} \label{stack m map} A  map from a stack $\mathcal X$ to an exploded manifold $\ex A$ is a map of stacks $\fun\Phi\co \mathcal X\longrightarrow \mathcal S(\ex A)$. We also refer to such a map as a $\C\infty1$ map.

In particular, the functor $\fun\Phi$ applied to $\hat f$ is a $\C\infty1$ map
\[\fun\Phi(\hat f)\co \ex F(\hat f)\longrightarrow \ex A\]
and these maps must be compatible: given any map $\psi\co\hat g\longrightarrow \hat f$, the following diagram commutes:
\[\begin{tikzcd}\ex F(\hat g)\rar[swap]{\ex F(\psi)}\ar[bend left]{rr}{\fun\Phi(\hat g)}&\ex F(\hat f)\rar[swap]{\fun\Phi(\hat f)}&\ex A\end{tikzcd}\]
For an individual object $f$ of $\mathcal  X$, $\ex F(f)$ is a point, so $\fun\Phi(f)$ has image a point, $\fun\Phi(f):=\fun\Phi(f)(\ex F(f))\in \ex A$. This point $\fun\Phi(f)$ only depends on the isomorphism class of $f$, Moreover, $\fun\Phi$ is determined by its values $\fun\Phi(f)$ on isomorphism classes of individual objects.  
\end{example}

\begin{example}\label{vector-bundle over stack} A rank $n$ complex  {\bf vector-bundle} over a stack $\mathcal X$ is a map $\mathcal X\longrightarrow \mathcal V$ to the stack of rank $n$ complex  vector-bundles from Example \ref{vector-bundle stack}.  

For example, a vector-bundle over $\mathcal S(\ex A)$ consists of a vector-bundle over $\ex A$ and a choice of pullback of this vector-bundle for each map of an exploded manifold into $\ex A$.
\end{example}

\begin{example}\label{quotient map} Let $\ex A$ be an exploded manifold with a finite group action $G$. There is a {\bf quotient map} $\fun\pi\co \ex A\longrightarrow \ex A/G$ defined as follows: for a family $\hat f\co \ex F(\hat f)\longrightarrow \ex A$ in $\mathcal S(\ex A)$,  $\fun\pi(\hat f)$ is the induced $G$--equivariant map of  the trival $G$--bundle to $\ex A$

\[\begin{tikzcd}\ex F(\hat f)\times G\dar\ar{rr}{(x,g)\to g*\hat f(x) }& & \ex A
\\ \ex F(\hat f)\end{tikzcd}\]
and for a morphism $\psi\co \hat f\longrightarrow \hat f'$ in $S(A)$, the morphism $\fun\pi(\psi)$ is the induced map of trivial $G$--bundles $\fun\pi(\psi)(f,g)=(\ex F(\psi)(f),g)$.

The family $\fun\pi(\text{id}_{\ex A})$ in the stack $\ex A/G$ has automorphism group $G$, which acts on the total space of the $G$--bundle $\ex A\times G$ diagonally as $g*(a,g')=(g*a,g'g^{-1})$. For a family $\hat f$ in $\ex A/G$, a morphism $\psi\co\hat f\longrightarrow \pi(\text{id}_{\ex A})$ is equivalent to a section of the $G$--bundle $\ex F\bd G(\hat f)\longrightarrow \ex F(\hat f)$; this section is the section sent to $\ex A\times \text{id}\subset \ex A\times G$.

\end{example}

\begin{example} There is a  map $\fun\Phi_{\mathcal X}\co \mathcal X/\hat f\longrightarrow \mathcal X$ sending $\hat h\longrightarrow \hat f$ to $\hat h$. Given a morphism $\psi\co \hat f\longrightarrow \hat g$ within a stack $\mathcal X$, there is an obvious induced map $\psi_*\co \mathcal X/\hat f\longrightarrow \mathcal X/\hat g$, given by composition with $\psi$. Moreover, the following diagram commutes strictly.
\[\begin{tikzcd}\mathcal X/\hat f\rar{\psi_*}\dar{\fun\Phi_{\mathcal X}} &\mathcal X/\hat g\ar{dl}{\fun\Phi_{\mathcal X}}
\\ \mathcal X
\end{tikzcd}\]
 \end{example}

\begin{defn} \label{family quotient stack}For $\hat f$ a family of curves in $\mathcal X$ with a finite group $G$ of automorphisms, define the stack $\hat f/G$ as follows.
\begin{itemize}
\item A  family $\hat h$ in $\hat f/G$  is  $G$--bundle $\hat h\bd G$ internal to $\mathcal X$, Definition \ref{G bundle}, together with a $G$--equivariant map $\psi_{\hat h}\co \hat h\bd G\longrightarrow \hat f$ within $\mathcal X$.
\item A morphism $\alpha\co \hat g\longrightarrow \hat h$ between families in $\hat f/G$ is given by a $G$--equivariant diagram within $\mathcal X$
\[\begin{tikzcd}\hat g\bd G\rar{\alpha\bd G} \ar[bend left]{rr}{\psi_{\hat g}} & \hat h\bd G\rar{\psi_{\hat h}} & \hat f
\end{tikzcd}\] 
\end{itemize}
 
There is a canonical map of stacks $\fun\Phi_{\mathcal X}\co \hat f/G\longrightarrow \mathcal X$, defined as follows. Each family $\hat h$ involves  a $G$--bundle $\hat h\bd G$ internal to $\mathcal X$. We define $\fun\Phi_{\mathcal X}(\hat h)$ to be the base of this $G$--bundle, so our $G$--bundle internal to $\mathcal X$ is the following map. 
\[\begin{tikzcd}\hat h\bd G \dar \\ \fun\Phi_{\mathcal X}(\hat h)\end{tikzcd} \]
Moreover, $\fun\Phi_{\mathcal X}(\alpha)$ is determined by the following commutative diagram.
\[\begin{tikzcd}\hat g\bd G\rar{\alpha\bd G} \dar\ar[bend left]{rr}{\psi_{\hat g}} & \hat h\bd G\rar{\psi_{\hat h}}\dar & \hat f
\\ \fun\Phi_{\mathcal X}(\hat g) \rar{\fun\Phi_{\mathcal X}(\alpha)} &\fun\Phi_{\mathcal X}(\hat h)
\end{tikzcd}\] 

\end{defn}

\begin{remark} \label{domain in quotient}
In the case that $\mathcal X$ is a moduli stack of curves, $\hat f/G$ is also a moduli stack of curves. The functors $\ex F$ and $\ex C$ are pulled back from $\mathcal X$ using $\fun\Phi_{\mathcal X}$ so, for $\hat h$ a family in the stack $\hat f/G$, the domain of $\hat h$ is  $\ex C(\hat h)=\ex C(\fun\Phi_{\mathcal X}(\hat h))$, and there is a map 
\[\ex C(\hat h)\longrightarrow \ex C(\hat f)/G\]
given by the $G$--equivariant map
\[\begin{tikzcd}\ex C(\hat h\bd G)\rar{\psi_{\hat h}}\dar &\ex C(\hat f)
\\ \ex C(\hat h)\end{tikzcd}\]
\end{remark}

\begin{example}\label{family map} For families $\hat h$ and $\hat f$ in $\mathcal X$,  a map $\fun\Phi\co \hat h\longrightarrow \hat f/G$ means a map of stacks $\fun\Phi\co \mathcal X/\hat h\longrightarrow \hat f/G$ such that the following diagram commutes strictly.
\[\begin{tikzcd}\mathcal X/\hat h\rar{\fun\Phi}\dar{\fun\Phi_{\mathcal X}} &\hat f/G\ar{dl}{\fun\Phi_{\mathcal X}}
\\ \mathcal X
\end{tikzcd}\]
Such a map is determined up to unique $2$-isomorphism by $\fun\Phi$ applied to $\text{id}_{\hat h}\co\hat h\longrightarrow \hat h$. Unpacking definitions, $\fun\Phi(\text{id}_{\hat h})$ is a $G$--bundle $\hat h\bd G\longrightarrow \hat h$ internal to $\mathcal X$, with a $G$--equivariant map 
\[\begin{tikzcd}\hat h\bd G\dar \ar{rr}{\psi_{\fun\Phi(\text{id}_{\hat h})}} && \hat f
\\ \hat h\end{tikzcd}\]
so, up to unique $2$-isomorphism, a map $\hat h\longrightarrow \hat f/G$ is a $G$--bundle over $\hat h$ with a $G$--equivariant map to $\hat f$.
\end{example}

\begin{example}\label{family quotient} There is a quotient map 
 \[\fun \pi\co \hat f\longrightarrow \hat f/G\]
  which sends the family $\hat h\longrightarrow \hat f$ in $\mathcal X/{\hat f}$ to the family
\[\begin{tikzcd}\hat h\times G\dar\ar{rr}{(h,g)\to g*h}&&\hat f
\\ \hat h\end{tikzcd}\]
in the stack $\hat f/G$. 

\end{example}

\begin{example}\label{continuous map} Let $X$ be a toplogical space. Recall, from \cite[Definition 3.1]{iec},  that each exploded manifold $\ex A$ is a topological space with topology pulled back from a Hausdorff topological space $\totl {\ex A}$. In particular, there is a notion of continuous map from an exploded manifold $\ex A$ to $X$. Define $\mathcal S(X)$ to be the stack over the category of $\C\infty1$ exploded manifolds 
\begin{itemize}
\item whose families are continuous maps \[\hat f\co \ex F(\hat f)\longrightarrow X\]
\item and whose morphisms $\psi\co \hat f\longrightarrow \hat f'$ are commutative diagrams of continuous maps
\[\begin{tikzcd} \ex F(\hat f)\dar[swap]{\ex F(\psi)}\rar{\hat f} & X
\\ \ex F(\hat f')\ar{ur}[swap]{\hat f'}\end{tikzcd}\]
where $\ex F(\psi)$ is a $\C\infty1$ map of exploded manifolds.
\end{itemize}
If $\mathcal X$ is a stack over the category of $\C\infty1$ exploded manifolds, a {\bf continuous map} $\fun\Phi\co \mathcal X\longrightarrow X$ is a morphism of stacks $\fun\Phi\co \mathcal X\longrightarrow S(X)$.

In particular, the morphism $\fun\Phi$ applied to a family $\hat f$ is a continuous map
\[\fun\Phi(\hat f)\co \ex F(\hat f)\longrightarrow X\]
and these maps must be compatible: given any map $\psi\co \hat g\longrightarrow \hat f$, the following diagram commutes:
\[\begin{tikzcd}\ex F(\hat g)\rar[swap]{\ex F(\psi)}\ar[bend left]{rr}{\fun\Phi(\hat g)}&\ex F(\hat f)\rar[swap]{\fun\Phi(\hat f)} &X\end{tikzcd}\]
For an individual object $f$,  $\fun\Phi(f)$ has image a point in $X$. This point $\fun\Phi(f)$ only depends on the isomorphism class of $f$, and $\fun\Phi$ is determined by its values $\fun\Phi(f)$ on isomorphism classes of individual objects. 

A {\bf continuous function } on a stack $\mathcal X$ is a continuous map from  $\mathcal X$ to the topological space $\mathbb R$.
\end{example}

\begin{defn} Given two morphisms $\fun\Phi_1, \fun\Phi_2\co \mathcal X_1\longrightarrow \mathcal X_2$ of stacks over the category of exploded manifolds, a {\bf $2$--morphism} between $\fun\Phi_1$ and $\fun\Phi_2$ is a natural transformation $\eta\co \fun\Phi_1 \Rightarrow \fun\Phi_2$  such that for any $\hat f$ in $\mathcal X_1$, $\ex F(\eta_{\hat f})$ is $\text{id}_{\ex F(\hat f)}$. 
\end{defn}

Every $2$--morphism is actually a $2$--isomorphism, because $\eta_{\hat f}\co \fun\Phi_1(\hat f)\longrightarrow \fun\Phi_2(\hat f)$ lifts the identity on $\ex F(\hat f)$, and is hence an isomorphism because stacks are categories fibered in groupoids.

\begin{defn}An {\bf isomorphism} of stacks is a morphism $\fun\Phi\co \mathcal X_1\longrightarrow \mathcal X_2$ such that there exists a  morphism $\fun\Phi^{-1}\co\mathcal X_2\longrightarrow \mathcal X_1$ and $2$--morphisms  $\fun\Phi\circ\fun\Phi^{-1}\Rightarrow \text{\bf id}_{\mathcal X_2}$ and $\fun\Phi^{-1}\circ\fun\Phi\Rightarrow \text{\bf id}_{\mathcal X_1}$.  Two stacks are called {\bf equivalent} or isomorphic  if there exists an isomorphism between them. 

Two maps of stacks $\fun\Phi\co \mathcal X_1\longrightarrow \mathcal X_2$ and $\fun\Phi'\co \mathcal X'_1\longrightarrow \mathcal X'_2$ are {\bf equivalent} if there is diagram which commutes up to $2$--isomorphism
\[\begin{tikzcd}\mathcal X_1\dar{\fun\Phi} \rar & \mathcal X_1'\dar{\fun\Phi'}
\\ \mathcal X_2 &\lar\mathcal X_2' \end{tikzcd}\]
and whose horizontal arrows are isomorphisms.

\end{defn}

\begin{example} The stack $\hat f/G$ from Definition \ref{family quotient} is equivalent to the quotient stack $\ex F(\hat f)/G$,  and the quotient maps $\fun\pi\co\ex F(\hat f)\longrightarrow \ex F(\hat f)/G$ from Example \ref{quotient map} and $\fun\pi\co \mathcal X/f\longrightarrow \hat f/G$ from Example \ref{family quotient} are equivalent. 
\end{example}

\begin{defn}\label{represented by} We say that a stack is {\bf represented} by an exploded manifold $\ex A$ if it is equivalent to $\mathcal S(\ex A)$. A morphism of stacks $\fun\Phi\co  \mathcal X_1\longrightarrow \mathcal X_2$ is {\bf represented} by a map $x\co \ex A_1\longrightarrow \ex A_2$ of exploded manifolds if there is a commutative diagram
\[\begin{tikzcd}\mathcal X_1\dar{\fun\Phi} &\lar \mathcal S(\ex A_1)\dar{\mathcal S(x)}
\\ \mathcal X_2 \rar& \mathcal S(\ex A_2) \end{tikzcd}\]
whose  horizontal arrows are equivalences. 

\end{defn}

\begin{example} The stack $\mathcal X/\hat f$ from Example \ref{over f} is represented by $\ex F(\hat f)$.
\end{example}

For a family of curves $\hat f$ with a finite group of automorphisms $G$, we have the following notion of a stack $\mathcal X$ being represented by $\hat f/G$. This is slightly more specific  than $\mathcal X$ being equivalent to the quotient stack $\ex F(\hat f)/G$.

\begin{defn}\label{represented by family} Let $\hat f$  be a family in a stack $\mathcal X$, and let $G$ be a finite group of automorphisms of $\hat f$.  We say that $\mathcal X$ is {\bf represented by} $\hat f/G$ if the map $\fun\Phi_{\mathcal X}\co \hat f/G\longrightarrow \mathcal X$ from Definition \ref{family quotient stack} is an isomorphism.  
\end{defn}

\begin{example}The quotient stack $\ex A/G$ is represented by $\pi(\text{id}_{\ex A})/G$.\end{example}

When $G$ is the trivial group, this definition reduces to the following: $\mathcal X$ is represented by $\hat f$ if the map $\mathcal X/\hat f\longrightarrow \mathcal X$ is an isomorphism. This happens if and only if every family in $\mathcal X$ has a unique map to $\hat f$.

We can always take the fiber product of stacks (compare \cite[Definition 6.1]{stacks}).

\begin{defn}\label{stack fiber product} If $\fun\Phi_i\co \mathcal X_i\longrightarrow \mathcal Z$ are morphisms of stacks, the {\bf fiber product} $\mathcal X_1\times_{\mathcal Z}\mathcal X_2$ is a stack whose
\begin{itemize}
\item families are triples  $(\hat f_1,\psi,\hat f_2)$ where $\hat f_i$ are families in $\mathcal X_i$ with $\ex F(\hat f_1)=\ex F(\hat f_2)$ and 
\[\psi\co \fun\Phi_1(\hat f_1)\longrightarrow \fun\Phi_2(\hat f_2)\]
is a morphism in $\mathcal Z$ with $\ex F(\psi)=\text{id}$; and
\item morphisms $(\hat f_1,\psi,\hat f_2)\longrightarrow (\hat f_1',\psi',\hat f_2')$ are comprised of a pair of morphisms $\alpha_i\co \hat f_i\longrightarrow \hat f'_i$ such that $\ex F(\alpha_1)=\ex F(\alpha_2)$ and 
\[\psi'\circ \fun\Phi_1(\alpha_1)=\fun\Phi_2(\alpha_2)\circ \psi \ .\]
\end{itemize}
As a special case, let $\mathcal Z$ be the trivial stack $\mathcal S(\text{pt})$, equivalent as a category to the category of $\C\infty 1$ exploded manifolds. Each stack $\mathcal X_i$ admits a unique morphism $\fun\Phi_i=\ex F$ to $\mathcal S(\text{pt})$. Define the {\bf product} $\mathcal X_1\times \mathcal X_2$ of $\mathcal X_1$ and $\mathcal X_2$ to be $\mathcal X_1\times_{\mathcal S(\text{pt})}\mathcal X_2$. In this case, the morphisms $\psi$ above are always the identity, so a family in $\mathcal X_1\times \mathcal X_2$ consists of a pair $(\hat f_1,\hat f_2)$ of families in $\mathcal X_i$ with $\ex F(\hat f_1)=\ex F(\hat f_2)$, and a morphism in $\mathcal X_1\times \mathcal X_2$ consists of a pair of morphisms in $\mathcal X_i$, send by $\ex F$ to the same morphism.
\end{defn}
The stack fiber product comes with obvious morphisms $\fun\pi_{i}\co \mathcal X_{1}\times_{\mathcal Z}\mathcal X_{2}\longrightarrow \mathcal X_{i}$ and a $2$--isomorphism $\eta\co  \fun\Phi_{1}\circ\fun\pi_1\Rightarrow  \fun\Phi_{2}\circ \fun\pi_2$ defined by $\eta_{(\hat f_{1},\psi,\hat f_{2})}=\psi$. 
\[\begin{tikzcd}\mathcal X_{1}\times_{\mathcal Z}\mathcal X_{2}\rar{\fun\pi_1}\dar[swap]{\fun\pi_2} & \mathcal X_1\dar{\fun\Phi_1}\ar[double]{dl}{\eta}
\\ \mathcal X_2\rar[swap]{\fun\Phi_2} & \mathcal Z\end{tikzcd}\]
This data satisfies the following universal property: Given maps $\fun\Psi_{i}\co \mathcal X'\longrightarrow \mathcal X_{i}$ and a $2$--morphism $\eta'\co \fun\Phi_{1}\circ \fun\Psi_{1}\Rightarrow \fun\Phi_{2}\circ \fun\Psi_{2}$, there exists a unique map $\fun\beta\co \mathcal X'\longrightarrow \mathcal X_{1}\times_{\mathcal Z}\mathcal X_{2}$ such that the following diagram strictly commutes
\[\begin{tikzcd} \mathcal X'\ar{dr}{\fun\beta}\rar{\fun\Psi_{1}}\dar{\fun\Psi_{2}} & \mathcal X_{1}
\\ \mathcal X_{2}& \mathcal X_{1}\times_{\mathcal Z}\mathcal X_{2}\uar{\fun\pi_{1}}\lar{\fun\pi_{2}} \end{tikzcd}\]
and such that 
\[\eta'=\fun\beta^{*}\eta\co \fun\Phi_1\circ \fun\pi_{1}\circ \fun\beta\Rightarrow \fun\Phi_2\circ\fun \pi_{2}\circ \fun\beta  \ .\]
In particular, $\fun\beta(\hat f)=(\fun\Psi_{1}(\hat f), \eta'_{\hat f} ,\fun\Psi_{2}(\hat f))$, and $\fun\beta(\alpha)=(\fun\Psi_{1}(\alpha),\fun\Psi_{2}(\alpha))$.

A special case of the stack fiber product is the fiber product of families $\hat f_{i}$ in a stack $\mathcal X$. Recall that  the stacks $\mathcal X/f_{i}$ from Example \ref{over f} come with natural maps $\fun\Phi_{\mathcal X}\co \mathcal X/\hat f_{i}\longrightarrow \mathcal X$. 

\begin{example}\label{family fp example}
If $\hat f_1$ and $\hat f_2$ are families in the stack $\mathcal X$, then $\hat f_1\times_\mathcal X\hat f_2$ means $\mathcal X/\hat f_{1}\times _\mathcal X \mathcal X/\hat f_{2}$.  This stack has
\begin{itemize}
\item families  $(\psi_1, \alpha ,\psi_2)$, where  $\psi_i\co \hat h_{i}\longrightarrow \hat f_i$ are morphisms in $\mathcal X$; and $\alpha\co\hat h_{1}\longrightarrow \hat h_{2}$ is an isomorphism such that $\ex F(\alpha)$ the identity; and 
\item morphisms commutative diagrams in $\mathcal X$:
\[\begin{tikzcd}\hat f_{1}  &\lar[swap]{\psi_{1}}\dar \hat h_{1}\rar{\alpha}&\hat h_{2}\dar \rar{\psi_{2}}\dar& \hat f_{2}
\\ & \ar{ul}{\psi_{1}'}\hat h_{1}'\rar{\alpha'}& \hat h'_{2}\ar{ur}[swap]{\psi_{2}'}
\end{tikzcd}\]
and the obvious functor with $\ex F(\psi_{1},\alpha,\psi_{2}):=\ex F(\hat h_{1})=\ex F(\hat h_{2})$. 
\end{itemize}

This stack is equivalent to the stack of families $\hat h$ in $\mathcal X$ with a pair  of maps $\psi_i\co \hat h\longrightarrow \hat f_i$. In this equivalent stack, 
\begin{itemize}
\item families are triples $(\psi_1,\hat h,\psi_2)$ comprised of a family $\hat h$ in $\mathcal X$ and a pair of morphisms $\psi_i\co \hat h\longrightarrow \hat f_i$ in $\mathcal X$.
\item morphisms $(\psi_1,\hat h,\psi_2)\longrightarrow (\hat h',\psi_1',\psi_2')$ are morphisms $\hat h\longrightarrow \hat h'$ such that the following diagram  commutes. 
\[\begin{tikzcd}\hat f_1 & \lar[swap]{\psi_1} \hat h\dar\rar{\psi_2} &\hat f_2
\\ & \hat h' \ar{ul}{\psi_1'}\ar{ur}[swap]{\psi_2'}\end{tikzcd}\]
\end{itemize}

\end{example}

For the fiber product of families $\hat f_i$ in $\mathcal X$, we have the following notion of $\hat f_1\times_{\mathcal X}\hat f_2$ being represented by a family $\hat h$ in $\mathcal X$. This notion is slightly stronger than the notion of the fiber product (and associated maps) being represented by exploded manifolds. 

\begin{defn}\label{family fiber product} Given families $\hat f_i$ in a stack $\mathcal X$, and maps $\pi_i\co\hat h\longrightarrow \hat f_i$,  we say that the fiber product $\hat f_1\times_\mathcal X\hat f_1$ is {\bf represented by} $(\pi_1,\hat h,\pi_2)$ if  this family $(\pi_1,\hat h,\pi_2)$ in $\hat f_1\times_\mathcal X\hat f_2$ is a final object: in other words, any other family in $\hat f_1\times_\mathcal X\hat f_2$ has a unique morphism to $(\pi_1,\hat h,\pi_2)$.

 In this case, we use the notation $\hat h=\hat f_1\times_\mathcal X\hat f_2$, and say that $\hat f_1\times_\mathcal X \hat f_2$ is represented by the family $\hat h$.
\end{defn}

Equivalently, $(\pi_1,\hat h,\pi_2)$ represents $\hat f_1\times_\mathcal X\hat f_2$ if it satisfies the following universal property: Given any family $\hat h'$ in $\mathcal X$ and maps $\pi'_i\co \hat h'\longrightarrow \hat f_i$, there exists a unique morphism $\hat h'\longrightarrow \hat h$ such that the following diagram commutes.
\[\begin{tikzcd}\hat h'\rar{\pi'_1}\dar[swap]{\pi'_2}\ar{dr}{\exists !} & \hat f_1
\\ \hat f_2  & \hat h=\hat f_1\times_{\mathcal X}\hat f_2\lar{\pi_2}\uar[swap]{\pi_1}\end{tikzcd}\]

\begin{defn}\label{representable map}

A map $\mathcal X\longrightarrow \mathcal Z$ of stacks is called {\bf etale}, or a {\bf representable submersion} if,  for all maps $\ex A\longrightarrow \mathcal Z$,  the map  $\ex A\times_{\mathcal Z}\mathcal X\longrightarrow \ex A$ is represented by a map of exploded manifolds that is respectively etale, or  a submersion. A representable submersion or etale map $\mathcal X\longrightarrow \mathcal Z$ is {\bf injective, closed, proper, complete, a refinement,  or a degree $n$ cover} if all the maps  $\ex A\times_{\mathcal Z}\mathcal X\longrightarrow \ex A$ are represented by maps of exploded manifolds that are respectively injective, closed, proper, complete\footnote{Definition 3.15 of \cite{iec}.},  refinements\footnote{Definition 10.5 of \cite{iec}} or degree $n$ covers.
\end{defn}

Often, representable submersions are simply called submersions, however we also need a notion of a submersion with infinite dimensional fibers; see Definition \ref{holomorphic submersion}.

\begin{example} Given a finite group $G$ acting on an exploded manifold $\ex A$, the  quotient map $\fun\pi\co \ex A\longrightarrow \ex A/G$ from Example \ref{quotient map} is a etale map which is a degree $\abs{G}$ cover. In particular, for a map $\ex X\longrightarrow \ex A/G$, the  map  $\ex X\times_{\ex A/G}\ex A\longrightarrow \ex X$ is represented by the $G$--bundle $\ex X\bd G\longrightarrow \ex X$ parametrized by $\ex X$.  

For a family $\hat f$ in $\ex A/G$, the fiber product $\hat f\times_{\ex A/G}\fun\pi(\text{id}_{\ex A})$ is represented the family $\hat f\bd G$. In particular  $\hat f\bd G$ is the pullback of $\hat f$ under the $G$--bundle map $\ex F\bd G(\hat f)\longrightarrow \ex F(\hat f)$, so $\hat f\bd G\longrightarrow \hat f$ is a $G$--bundle internal to $\ex A/G$, with $\ex F(\hat f\bd G)\longrightarrow \ex F(\hat f)$  the $G$--bundle $\ex F\bd G(\hat f)\longrightarrow \ex F(\hat f)$.

Let us check that $\hat f\bd G$ represents the fiber product $\hat f\times_{\ex A/G}\fun\pi(\text{id}_{\ex A})$ as in Definition \ref{family fiber product}. We have that $\ex F\bd G(\hat f\bd G)=\ex F\bd G(\hat f)\times_{\ex F(\hat f)}\ex F\bd G(\hat f) $, so the diagonal section gives a canonical section of $\ex F\bd G(\hat f\bd G)\longrightarrow \ex F(\hat f\bd G)$, and hence a map $\hat f\bd G\longrightarrow \fun\pi(\id_{\ex A})$. The maps $\hat f\longleftarrow \hat f\bd G\longrightarrow \fun\pi(\text{id}_{\ex A})$ determine a family in the stack $\hat f\times_{\ex A/G}\fun\pi(\text{id}_{\ex A})$. Given a map $\psi\co \hat h\longrightarrow \hat f$ with a section $s$ of $\ex F(\hat h\bd G)\longrightarrow \ex F(\hat h)$, the universal property of the pullback $\hat f\bd G$ determines a unique map $\psi'\co\hat h\longrightarrow \hat f\bd G$ lifting $\psi$ so that $\ex F(\psi\bd G)=\ex F\bd G(\psi)\circ s$.  This section $s$ is equivalent to a map $\hat h\longrightarrow \fun\pi(\text{id}_{\ex A})$, and the map $\psi'$ is the unique map such that the following diagram commutes:
\[\begin{tikzcd} & \hat h\dar{\psi'}\ar{dl}[swap]{\psi}\ar{dr}
\\ \hat f & \lar \hat f\bd G \rar & \fun\pi(\text{id}_{\ex A})
\end{tikzcd}\]
Accordingly, $\hat f\longleftarrow \hat f\bd G\longrightarrow \fun\pi(\text{id}_{\ex A})$ is a final object, and $\hat f\bd G$ represents the fiber product $\hat f\times_{\ex A/G}\fun\pi(\text{id}_{\ex A})$, as required.

\end{example}


To a stack $\mathcal X$ over the category of exploded manifolds, we can associate the set $\mathcal X_{top}$ of isomorphism classes of individual objects in $\mathcal X$ (assuming the category $\mathcal X$ is equivalent to a small category) --- for example if $\mathcal X$ is a moduli space of curves, $\mathcal  X_{top}$ is the coarse moduli space.  We want our substacks to act like subsets of $\mathcal X_{top}$, (with an obvious notion of union,  intersection and complement) so we make a relatively strong definition of substack, using the following notion.

\begin{defn} Let $\mathcal U\subset \mathcal X$ be a subcategory of a stack $\mathcal X$ over the category of exploded manifolds.  For  each family $\hat f$ in $\mathcal X$, define the subset $\mathcal U(\hat f)\subset \ex F(\hat f)$ to be the set of points $p\in \ex F$ such that $\hat f$ restricted to $p$ is  isomorphic to an object in $\mathcal U$. 
\end{defn}

\begin{defn}\label{substack} A {\bf substack} $\mathcal U\subset \mathcal X$ of a stack over the category of exploded manifolds is a full subcategory $\mathcal U$ such that a family $\hat f$ in $\mathcal X$ is in $\mathcal U$  if and only if $\mathcal U(\hat f)=\ex F(\hat f)$.
\end{defn}

 In particular, a substack  $\mathcal U\subset \Ms(\hat{\ex B})$   means a full subcategory of $\Ms(\hat {\ex B})$ such that a family of curves $\hat f$ in $\Ms(\hat {\ex B})$ is in $\mathcal U$ if and only if each individual curve in $\hat f$ is isomorphic to a curve in $\mathcal U$. Accordingly, we can describe such a substack $\mathcal U$  by describing properties of the individual curves in $\mathcal U$. 

\begin{defn}\label{open substack} An {\bf open substack} $\mathcal U\subset \mathcal X$ of a stack over the category of exploded manifolds is a substack $\mathcal U$ such that $\mathcal U(\hat f)\subset \ex F(\hat f)$ is open (using the standard topology on the exploded manifold $\ex F(\hat f)$ from Definition 3.1 of \cite{iec}) for all families $\hat f$ in $\mathcal X$.

A {\bf closed substack} $\mathcal C\subset \mathcal X$ is a substack such that $\mathcal C(\hat f)\subset \ex F(\hat f)$ is closed for all families $\hat f$. 

A {\bf compact substack} $\mathcal K\subset \mathcal X$ is one with the property that any cover of $\mathcal K$ by open substacks admits a finite sub-cover. 
\end{defn}

As usual for topological subsets, the union of open substacks is an open substack, a finite intersection of open substacks is an open substack, and the complement of any open substack is a closed substack. For example,  the closure of a substack is defined to be the intersection of all closed substacks containing it.

Although our definition of a substack is strong, our definition of an open substack above is equivalent to the standard definition of an injective etale map (or open embedding) $\mathcal U\longrightarrow \mathcal X$; compare   \cite[Definition C.6]{joycebook}. In fact, $\mathcal U(\hat f)\longrightarrow \ex F(\hat f)$ is the injective etale map representing $\mathcal U\times_\mathcal X\hat f\longrightarrow \mathcal X/f$.

\begin{remark}\label{extra open}Given a substack $\mathcal C\subset \mathcal X$, the intersection of $\mathcal C$ with an open or closed substack of $\mathcal X$ will always be an open or closed substack of $\mathcal C$, but $\mathcal C$ might have many more open or closed substacks, because exploded manifolds may not be degenerate enough to detect the subspace topology of $\mathcal C$. An example is given by the stack $\mathcal S(X)\subset\mathcal S(\mathbb R^2)$ corresponding to a subset $X\subset \mathbb R^{2}$ that is connected but not path connected. Relevant to this paper, the moduli stack of stable holomorphic curves $\Mod$ is a closed substack of the ambient moduli stack $\Msw$ of stable curves. In general, $\Mod$ is too degenerate to be analysed by itself  as a stack over the category of exploded manifolds, and we need the extra structure of its inclusion into the ambient stack $\Msw$.  \end{remark}

\begin{example}\label{open inverse image} Given any map $\fun\Phi\co \mathcal X_1\longrightarrow \mathcal X_2$,  the inverse image of any (open) substack  $\mathcal O$ of $\mathcal X_2$ is an (open) substack  of $\fun\Phi^{-1}(\mathcal O)$ of $\mathcal X_1$. In particular, a family $\hat f$ is in $\fun\Phi^{-1}(\mathcal O)$ if and only if $\fun\Phi(\hat f)$ is in $\mathcal O$, so, when $\mathcal O$ is open the subfamily of $\hat f$ contained in $\fun\Phi^{-1}\mathcal O$ is always open. Moreover, for $\ex A$ any exploded manifold and $X$ any topological space, there is an obvious correspondence between (open) subsets of $\ex A$ or $X$ and (open) substacks of $\mathcal S(\ex A)$ or $\mathcal S(X)$ respectively. (A  caveat to this is that, as in Remark \ref{extra open}, when $X$ is not a nice topological space, $\mathcal S(X)$ can have more open substacks than just the inverse image of open subsets of $X$.)
Accordingly, if $x\co \mathcal X\longrightarrow X$ is a continous map from a stack to a toplogical space $X$, as in Example \ref{continuous map}, the inverse image of any open subset of $X$ is an open substack of $\mathcal X$.
\end{example}  

\begin{example}\label{Xtop} Assuming $\mathcal X$ is equivalent to a small category, define $\mathcal X_{top}$ to be the set whose points are isomorphism classes of individual objects in $\mathcal X$. Every substack $\mathcal U\subset \mathcal X$ has a corresponding subset $\mathcal U_{top}\subset \mathcal X_{top} $ such that $f$ is in $\mathcal U$ if and only if the isomorphism class, $f_{top}$  of $f$, is in $\mathcal U_{top}$. Moreover $\mathcal X_{top}$ has a natural topology, with $\mathcal U_{top}\subset \mathcal X_{top}$ an open subset if and only if $\mathcal U$ is an open substack. To any map of stacks  $\fun\Phi\co \mathcal X\longrightarrow \mathcal X'$ over the category of exploded manifolds, there is an induced continuous map $\fun\Phi_{top}\co\mathcal X_{top}\longrightarrow\mathcal X'_{top}$ sending the isomorphism class $f_{top}$ of $f$ to the isomorphism class $\fun\Phi(f)_{top}$ of $\fun\Phi (f)$. Clearly $(\fun\Phi\circ\fun\Phi')_{top}=\fun\Phi_{top}\circ\fun\Phi'_{top}$, and $2$--isomorphic maps of stacks are sent to identical continuous maps,  so the construction defines a functor from  the $2$--category of stacks over exploded manifolds to the more familiar $1$--category of topological spaces. 
\end{example}

\begin{example}\label{totlX} Just as the topology of an exploded manifold $\ex A$ is usually non-Hausdorff, but pulled back from a map to a Hausdorff topological space $\totl{\ex A}$, the topology on a stack $\mathcal X$ and $\mathcal X_{top}$ is pulled back from the topology on a simpler toplogical space $\totl{\mathcal X_{top}}$, defined as follows: Say that individual objects $f$ and $f'$ in $\mathcal X$ are {\bf topologically indistinguishable} if every closed substack of $\mathcal X$ contains $f$ if and only if it contains $f'$. Define $\totl{\mathcal X_{top}}$ to be the set equivalence classes of topologically indistinguishable individual objects. There is a topology on $\totl{\mathcal X_{top}}$ such that  each open or closed substack  $\mathcal U\subset \mathcal X $ corresponds to a (respectively open or closed) subset $\totl{\mathcal U_{top}}$,  where $f$ is in $\mathcal U$ if and only if its equivalence class, $\totl{f_{top}}$ is in $\totl{\mathcal U_{top}}$. As in Example \ref{Xtop}, this construction  defines a functor from the $2$--category of stacks over exploded manifolds to the $1$--category of topological spaces.  

For $\ex A$ an exploded manifold, $\mathcal S(\ex A)_{top}$ is canonically isomorphic to the topological space $\ex A$, and $\totl{\mathcal S(\ex A)_{top}}$ is canonically isomorphic to $\totl{\ex A}$.
\end{example}

In many stacks of interest,  curves $f_1$ and $f_2$ will be topologically indistinguishable if and only if they admit maps $f_i\longrightarrow \hat f$ to some family such that their image in $\totl{\ex F(\hat f)}$ is the same point. This property need not hold in a general stack, but it is inherited by any closed or open substack. In Lemma \ref{hausdorff}, we show that an open substack $\Msw$ of $\Ms$ enjoys this property, and we show that $\totl{\Msw_{top}}$ is Hausdorff. In contrast, $\totl{\Ms_{top}}$ is not Hausdorff.


\begin{remark}\label{orbifold remark} As explained in \cite{diffstacks, orbifoldstack,joycebook},  a good formalization of smooth orbifolds is   as Deligne--Mumford stacks over the category of smooth manifolds. In the category of smooth manifolds, a Deligne--Mumford stack is a stack $\mathcal X$ over the category of smooth manifolds,  with a second-countable and Hausdorff topology, and locally equivalent to $U/G$ for some manifold $U$ with an action of a finite group $G$. Equivalently,  there exists a surjective etale map $M\longrightarrow \mathcal X$ from a manifold $M$ such that the map $M\times_{\mathcal X} M\longrightarrow M\times M$ is proper. 

Another equivalent condition is as follows: $\mathcal X$ is a Deligne--Mumford stack if it is equivalent to the quotient stack of a etale proper groupoid in the category of smooth manifolds. Given a surjective etale map $U\longrightarrow \mathcal X$, the Deligne--Mumford stack $\mathcal X$ is equivalent to the quotient stack of the  etale proper groupoid whose objects are parametrized by $U$, and whose morphisms are parametrized by $U\times_{\mathcal X} U$.  Composition of these morphisms is  encoded in the three maps from $U\times_{\mathcal X}U\times_{\mathcal X}U$ to $U$.  
\end{remark}

\begin{defn}\label{orbifold} An {\bf exploded orbifold} is a stack $\mathcal X$ over the category of exploded manifolds  which is  {\bf Deligne--Mumford} in the following sense: $\mathcal X$ is locally equivalent to $\ex A/G$ for $\ex A$ an exploded manifold with an action of a finite group $G$, and $\totl{\mathcal X_{top}}$ is second countable and Hausdorff. Equivalently,  there exists a surjective etale map \footnote{Definition \ref{representable map}.} $\ex U\longrightarrow \mathcal X$ such that $\ex U\times_\mathcal X \ex U\longrightarrow \ex U\times\ex  U$ is represented by proper map of exploded manifolds. 
\end{defn}

An alternate, (and equivalent!) definition is that an exploded orbifold is a stack equivalent to the quotient stack of a proper\footnote{The reader might recall, from Definition 3.15 of \cite{iec}, that the correct analogue of `proper' in the category of exploded manifolds is usually `complete'. In this case, these maps will be complete if and only if they are proper, so we can continue to use the terminology of a etale proper groupoid.} etale groupoid in the category of exploded manifolds. To obtain an etale proper groupoid from a Deligne--Mumford stack $\mathcal X$, simply choose a surjective etale map $\ex U\longrightarrow \mathcal X$.   Then, there is an etale proper groupoid $\ex U\times_{\mathcal X}\ex U\rightrightarrows \ex U$ with objects parametrized by $\ex U$, morphisms parametrized by $\ex U\times_{\mathcal X}\ex U$, and composition encoded by the three natural maps $\ex U\times_{\mathcal X}\ex U\times_{\mathcal X}\ex U\longrightarrow \ex U\times_{\mathcal X}\ex U$. Moreover, the quotient stack of this groupoid is equivalent to $\mathcal X$. See Section 4 of  \cite{orbifoldstack} or Section 2.4 of \cite{diffstacks} for how to define the quotient stack of  an etale proper groupoid.

The following is a standard definition of a sheaf over a stack; compare Section 3.1 of \cite{diffstacks}.

\begin{defn}\label{sheaf} Let $\mathcal X$ be a stack over the category of exploded manifolds. A {\bf sheaf} $\mathcal Y$ (of sets) over $\mathcal X$  is a contravariant functor $\mathcal Y$ with domain $\mathcal X$ (and codomain the category of sets) such that $\mathcal Y$ is a sheaf when restricted to the subcategory of open subfamilies\footnote{On a family $\hat f$, we use the standard topology of the exploded manifold $\ex F(\hat f)$, so an open subfamily is the restiction of $\hat f$ to an open subset of $\ex F(\hat f)$.} of $\hat f$. 

A map of sheaves is a natural transformation $\eta\co \mathcal Y_1\Rightarrow\mathcal Y_2$, and a global section of a sheaf $\mathcal Y$ is a map from the trivial sheaf into $\mathcal Y$. The pullback of $\mathcal Y$  by a  map $\fun\Phi\co \mathcal X'\longrightarrow \mathcal X$ is $\fun\Phi^*\mathcal Y:= \mathcal Y\circ \fun\Phi$.

\end{defn}

For example:
\begin{itemize}\item For a fixed $\hat g$, there is a sheaf over $\mathcal X$ which  assigns to $\hat f$ the set of maps from $\hat f$ to $\hat g$. Similarly, $\C\infty1$ maps to an exploded manifold, or continuous maps to a topological space form a sheaf.

\item There is a natural sheaf of rings over $\mathcal X$ which assigns to every family $\hat f$ the ring $\C\infty 1(\ex F(\hat f))$ of 
$\mathbb R$--valued $\C\infty1$--functions on the exploded manifold $\ex F(\hat f)$, and which assigns to every morphism $\psi\co\hat f\longrightarrow \hat g$, the homomorphism $\C\infty1(\ex F(\hat g))\longrightarrow \C\infty1(\ex F(\hat f))$ induced from composition with the map $\ex F(\psi)\co \ex F(\hat f)\longrightarrow \ex F(\hat g)$. A global section of this sheaf is a $\C\infty1$ map $\mathcal X\longrightarrow \mathbb R$.
\item There is a sheaf which assigns to $\ex F(\hat f)$ the ring of $\C\infty1$ differential forms on $\ex F(\hat f)$.
\item For $\mathcal V$ a stack of vector-bundles (like in Example \ref{vector-bundle stack}), there is a sheaf, $\Gamma$, which assigns to a family $\hat f$ the $\C\infty1$ sections of the vector-bundle $\mathcal  V(\hat f)\longrightarrow \ex F(\hat f)$. These sections form a locally free, finite rank module over the ring,  $\C\infty 1(\ex F(\hat f))$, of functions on $\ex F(\hat f)$. There is a similar sheaf of sections of a vector-bundle over a stack; see Example \ref{vector-bundle over stack}.

\end{itemize}

\subsection{The sheaf $\Y$ over the ambient moduli stack of curves}
\label{Y section}

\

In this paper, we consider holomorphic curves in the ambient stack  $\Ms(\hat {\ex B})$ as the solution of an equation $\dbar\hat f=0$, where $\dbar$ is a section of a sheaf $\Y$ over the $\Ms(\hat {\ex B})$.

We shall often need  the vertical (co)tangent space of families or submersions. Given a submersion between exploded manifolds,
$x\co\ex A\longrightarrow \ex B$,
use  $T\ex A\ov{\ex B}$ or $T_{vert}\ex A$ to indicate the vertical tangent bundle of $\ex A$ over $\ex B$. In other words,  $T\ex A\ov{\ex B}$ is the sub bundle of $T\ex A$ consisting of the kernel of the derivative of $x$.
Use the notation $T^{*}\ex A\ov{\ex B}$ or $T^{*}_{vert}\ex A$ to indicate the vertical cotangent bundle of $\ex A$ over $\ex B$. Define  $T^{*}\ex A\ov{\ex B}$ to be equal to the dual of $T\ex A\ov{\ex B}$. We shall use the notation $T_{vert}$ when no ambiguity shall arise, and when it is less cumbersome. For example, given a family of curves 
\[\begin{tikzcd}\ex C(\hat f)\dar \rar{\hat f} & \hat{\ex B}\dar
\\ \ex F(\hat f) \rar & \ex B_0\end{tikzcd}\]
 the notation $T_{vert}\ex C(\hat f)$ shall always mean $T\ex C(\hat f)\ov{\ex F(\hat f)}$, and  the notation $T_{vert}\hat{\ex B}$ shall always mean $T\hat{\ex B}\ov{\ex B_0}$.

\begin{defn}
Given a $\C\infty1$ family of curves $\hat f$ in $\hat{\ex B}$,   
 define
\[\dvert \hat f\co T_{vert}\ex C(\hat f)\longrightarrow T_{vert}\hat{\ex B}\]
to be the derivative\footnote{See \cite[Section 6]{iec}} of $\hat f$ restricted to the vertical tangent space, $T_{vert}\ex C(\hat f)\subset T\ex C(\hat f)$.

Define
\[\dbar \hat f\co T_{vert}\ex C(\hat f)\longrightarrow T_{vert}\hat{\ex B}\]
as
\[\dbar \hat f:=\frac 12\lrb{ \dvert f +J\circ \dvert f\circ j}\ .\]

\end{defn}

 This $\dbar f$ is a section of the vector-bundle $Y(\hat f):=\lrb{T^{*}_{vert}\ex C(\hat f) \otimes \hat f^{*}T_{vert}\hat{\ex B}}^{0,1}$ over $\ex C(\hat f)$, and is also as a section of a corresponding sheaf $\Y(\hat f)$ over $\ex F(\hat f)$.

 Recall the notion of an integral vector from  \cite[Definition 6.8]{iec}. The above section $\dbar \hat f$ has regularity $\C\infty1$, and automatically vanishes on integral vectors. We can add such sections and multiply them by $\mathbb R$--valued $\C\infty1$ functions pulled back from  $\ex F(\hat f)$, so  the set of such sections has the structure of a $\C\infty1(\ex F(\hat f))$--module over $\ex F(\hat f)$.
 
 \begin{defn}\label{Ydef} Let  $\Y(\hat f)$  be the sheaf of $\C\infty1(\ex F(\hat f))$--modules over $\ex F(\hat f)$ consisting of $\C\infty1$ sections of $Y(\hat f):=\lrb{T^{*}_{vert}\ex C(\hat f) \otimes \hat f^{*}T_{vert}\hat{\ex B}}^{0,1}$  vanishing on integral vectors within $T_{vert}\ex C(\hat f)$.
 
 \end{defn}

Given any map \[\hat f\longrightarrow \hat g\] there is a corresponding pullback diagram of vector-bundles \[\begin{tikzcd}Y(\hat f)\rar\dar & Y(\hat g)\dar
\\ \ex C(\hat f)\dar\rar&\ex C(\hat g)\dar
\\ \ex F(\hat f)\rar&\ex F(\hat g)\end{tikzcd}\] and a functorial map of sheaves
\[\Y(\hat f)\longleftarrow \Y(\hat g)\ \]
compatible with the respective $\C\infty1(\ex F(\hat f))$ and $\C\infty1(\ex F(\hat g))$ module structures using the pullback of functions $\C\infty1(\ex F(\hat g))\longrightarrow \C\infty1 (\ex F(\hat f))$.
Accordingly,  $\Y$ is a sheaf over $\Ms$; see Definition \ref{sheaf}.


\subsection{Fiberwise holomorphic maps from moduli stacks of exploded curves} \label{fiberwise holomorphi map}

\

\

We have defined some general notions of maps from stacks over the category of exploded manifolds (see Definition \ref{stack map} and examples \ref{m stack map}, \ref{stack m map}, and \ref{quotient map}) but we will also need the special notion of a fiberwise holomorphic map from a moduli stack of exploded curves.

\begin{defn}\label{stack of curves} A stack of exploded curves is a stack $\mathcal U$ over the category of exploded manifolds with a map $\mathcal U\longrightarrow \Ms$. \end{defn}

Any stack of exploded curves has functors $\ex F$ and $\ex C$ and a natural transformation $\pi_{\ex F}\co \ex C\Rightarrow \ex F$ pulled back from $\Ms$ so that $\ex C(\hat f)\longrightarrow \ex F(\hat f)$ is a family of exploded curves.

Given any moduli stack of exploded curves $\mathcal U$, there is a universal curve $\mathcal U^{+1}\longrightarrow \mathcal U$ such that given any family $\hat f$ in $\mathcal U$, the following diagram
\[\begin{tikzcd}\ex C(\hat f)\dar \rar & \mathcal U^{+1}\dar
\\ \ex F(\hat f)\rar & \mathcal U\end{tikzcd}\]
is a fibre product diagram in the sense that $\ex C(\hat f)\longrightarrow \ex F(\hat f)$  represents $\hat f\times_{\mathcal U}\mathcal U^{+1}\longrightarrow \ex F(\hat f)$; see Definition \ref{represented by}. Below we define this universal curve as a stack. In Section \ref{+n section}, we also construct $\mathcal U^{+1}$ as a moduli stack of exploded curves; see Definition \ref{+1 curve stack}.

\begin{defn}\label{+1 stack} If $\mathcal U$ is a moduli stack of exploded curves, define the stack $\mathcal U^{+1}$ as follows:
\begin{itemize} \item A family $(\hat f,s)$ in $\mathcal U^{+1}$ is a family $\hat f$ in $\mathcal U$ and a $\C\infty 1$ section 
\[s\co \ex F(\hat f)\longrightarrow \ex C(\hat f)\]
\item A morphism $(\hat f,s)\longrightarrow (\hat f',s')$ in $\mathcal U^{+1}$ is a morphism $\psi\co \hat f\longrightarrow \hat f'$ in $\mathcal U$ such that the following diagram commutes.
\[\begin{tikzcd} \ex C(\hat f)\rar{\ex C(\psi)} & \ex C(\hat f')
\\ \uar{s}\ex F(\hat f)\rar{\ex F(\psi)} & \ex F(\hat f')\uar{s'}  \end{tikzcd}\]
\end{itemize}

\end{defn}

For any family of curves $\hat f$ in $\mathcal U$, there exists a family  $\hat f^{+1}$ in $\mathcal U^{+1}$ representing $\hat f\times_{\mathcal U}\mathcal U^{+1}$. This family is parametrised by $\ex C(\hat f)$ instead of $\ex F(\hat f)$, so the assignment $\hat f\mapsto \hat f^{+1}$ is not a map of stacks, although it is a functor.  We upgrade $\hat f^{+1}$ to a family of curves in Section \ref{+n section}.

\begin{defn} \label{basic+1def}
For $\hat f$ a family of curves in $\mathcal U$, define the family $\hat f^{+1}$ in $\mathcal U^{+1}$ to be $(\pi_{\ex F(\hat f)}^*\hat f, \Delta)$, where $\Delta$ is the diagonal section in the following diagram.

\[\begin{tikzcd} \ex C(\pi_{\ex F(\hat f)}^*\hat f)=\ex C(\hat f)\times_{\ex F(\hat f)}\ex C(\hat f)\rar\dar & \ex C(\hat f)\dar{\pi_{\ex F(\hat f)}}
\\ \ex F(\pi_{\ex F(\hat f)}^*\hat f)=\ex C(\hat f)\rar{\pi_{\ex F(\hat f)}}\uar[bend right, swap]{\Delta} & \ex F(\hat f)\end{tikzcd}\]

For $\psi\co \hat f\longrightarrow \hat g$ a morphism in $\mathcal U$, define $\psi^{+1}\co \hat f^{+1}\longrightarrow \hat g^{+1}$ to be the unique morphism in $\mathcal U^{+1}$ induced by the  morphism $\pi_{\ex F(\hat f)}^*\hat f\longrightarrow \pi_{\ex F(\hat g)}^*\hat g$ from the following diagram.

\[\begin{tikzcd} \ex C(\hat f)\dar{\pi_{\ex F(\hat f)}}\rar{\ex C(\psi)} & \ex C(\hat g)\dar{\pi_{\ex F(\hat g)}}
\\ \ex F(\hat f)\rar & \ex F(\hat g)\end{tikzcd}\]
\end{defn}

\begin{defn}\label{curve fhm} A fiberwise holomorphic map $(\fun\Phi,\fun\Phi^{+1})$ between stacks of exploded curves 
is a strictly commuting diagram map of stacks 
\[\begin{tikzcd} \mathcal U^{+1}\dar \rar{\fun\Phi^{+1}} & \mathcal X^{+1}\dar
\\ \mathcal U\rar{\fun\Phi} & \mathcal X\end{tikzcd}\]
such that, for all families $\hat f$ in $\mathcal U$, there exists a fiberwise holomorphic map  
 $\fun\Phi^{+1}_{\hat f}\co \ex C(\hat f)\longrightarrow \ex C(\fun\Phi(\hat f))$ 
\[\begin{tikzcd}\ex C(\hat f)\dar \rar{\fun\Phi^{+1}_{\hat f}}&\ex C(\fun\Phi(\hat f))\dar
\\ \ex F(\hat f)\rar{\text{id}} & \ex F(\fun\Phi(\hat f))\end{tikzcd}\]
such that 
\[\fun\Phi^{+1}(\hat f,s)=(\fun\Phi(\hat f), \fun\Phi_{\hat f}^{+1} \circ s)\]


\end{defn}

Note that $\fun\Phi_{\hat f}^{+1}$ is completely determined by $(\fun \Phi,\fun\Phi^{+1})$.
%
 Moreover, the assignment $\hat f\mapsto \fun\Phi_{\hat f}^{+1}$ defines a natural transformation $\ex C\Rightarrow \ex C\circ\fun\Phi$, which determines the map of stacks $\fun\Phi^{+1}$.
 
 In general, given any map $\fun\Phi^{+1}\co \mathcal U^{+1}\longrightarrow \mathcal X$ of stacks, and given any family $\hat f$ in $\mathcal U$, there exists a unique map of stacks
 \begin{equation}\label{phi_f}\fun\Phi^{+1}_{\hat f}\co \ex C(\hat f)\longrightarrow \mathcal X\end{equation}
 such that 
 \[\fun\Phi^{+1}(\hat f,s)=\fun\Phi^{+1}_{\hat f}\circ s\ .\]
 In particular, a family $\hat h$ in $\mathcal S(\ex C(\hat f))$ consists of a map $\hat h:\ex F(\hat h)\longrightarrow \ex C(\hat f)=\ex F(\hat f^{+1})$. Then we can define $\fun\Phi^{+1}_{\hat f}(\hat h):=\fun\Phi^{+1}(\hat h^*\hat f^{+1})$.

As well as fiberwise holomorphic maps between moduli stacks of curves, we also need a notion of a fiberwise holomorphic map to a family of almost complex exploded manifolds, and an orbifold generalisation of this. Instead of a general definition, we only need the following concrete case. 

\begin{defn}\label{fhm}Let $\mathcal U\longrightarrow \Ms(\hat{\ex  B})$ be  a stack of exploded curves, and let $\hat {\ex A}\longrightarrow {\ex X}$ be a family of almost-complex\footnote{See Definition 8.1 of \cite{iec} for almost complex structures.} exploded manifolds with a choice of a finite group $G$ of automorphisms. Define a {\bf fiberwise-holomorphic} map to be a strictly commuting diagram of exploded stacks 
\[\begin{tikzcd}\mathcal U^{+1}\rar{\fun\Phi^{+1}}\dar &\hat {\ex A}/G\dar
\\ \mathcal U\rar{\fun\Phi} &{\ex X}/G\end{tikzcd}\]
such that for all families $\hat h$ in $\mathcal U$, the map of stacks
\[\fun\Phi^{+1}_{\hat h}\co \ex C(\hat h)\longrightarrow \hat {\ex A}/G\]
 from (\ref{phi_f}) is fiberwise holomorphic in the following sense: Such a map to the stack quotient $\hat{\ex A}/G$ is determined by\footnote{In the notation of Definition \ref{stack quotient}, $\fun\Phi^{+1}(\hat h)=\fun\Phi^{+1}_{\hat h}(\text{id}_{\ex C(\hat h)})$.} a $G$--bundle\footnote{Definition \ref{G bundle}.} $\hat h\bd G\longrightarrow \hat h$, and a $G$--equivariant,  map of exploded manifolds 
\[\fun\Phi^{+1}(\hat h)\co \ex C(\hat h\bd G)\longrightarrow \hat{\ex A}\ .\]
This map $\fun\Phi^{+1}(\hat h)$ must be fiberwise holomorphic, so  the following is a family of  of holomorphic curves in $\hat{\ex A}\longrightarrow \ex X$.
 \[\begin{tikzcd}\ex C(\hat h\bd G)\ar{rr}{\fun\Phi^{+1}(\hat h)} \dar &&\hat{\ex A}\dar
\\ \ex F(\hat h\bd G)\ar{rr}{\fun\Phi(\hat h)}&&\ex X \end{tikzcd}\]

Say that this fiberwise holomorphic map has an {\bf effective $G$--action } if the map $\fun\Phi^{+1}(\hat h)\co\ex C(\hat h\bd G)\longrightarrow \hat {\ex A}$ is never preserved by  the action of a nontrivial element of $G$ on $\hat{\ex A}$. 

%
\end{defn}

To define such a fiberwise holomorphic map $(\fun\Phi,\fun\Phi^{+1})$, from $\mathcal U$ it suffices to specify the following:
\begin{itemize}
\item for each family $\hat h$ in $\mathcal U$, a $G$--bundle $\hat h\bd G\longrightarrow \hat h$ internal to $\mathcal U$;
\item for each morphism $\alpha\co \hat h_1\longrightarrow \hat h_2 $ in $\mathcal U$, a $G$--equivariant morphism \[\alpha\bd G\co \hat h_1\bd G\longrightarrow\hat h_2\bd G\] such that $(\alpha\circ \beta)\bd G=\alpha\bd G\circ \beta\bd G$; and
\item for each family $\hat h$ in $\mathcal U$ a fiberwise holomorphic $G$--equivariant map
 \[\begin{tikzcd}\ex C(\hat h\bd G)\ar{rr}{\fun\Phi^{+1}(\hat h)} \dar &&\hat{\ex A}\dar
\\ \ex F(\hat h\bd G)\ar{rr}{\fun\Phi(\hat h)}&&\ex X \end{tikzcd}\]
 such that, given $\alpha\co\hat h_1\longrightarrow \hat h_2$, we have
\[\fun\Phi^{+1}(\hat h_1)=\fun\Phi^{+1}(\hat h_2)\circ \ex C(\alpha\bd G)\ .\] 
\end{itemize}

It is immediate that the above data defines a map of stacks $\fun\Phi$ sending the family $\hat h$ to the family in $\ex X/G$ represented by the above $G$--bundle and  $G$--equivariant map $\fun\Phi(\hat h)$. Moreover, the $G$--equivariant map $\fun\Phi^{+1}(\hat h)$ induces a map of stacks $\fun\Phi^{+1}_{\hat h}\co \ex C(\hat h)\longrightarrow \hat{\ex A}/G$, and we can define the map of stacks $\fun\Phi^{+1}\co \mathcal U\longrightarrow \hat{\ex A}/G$ so that $\fun\Phi^{+1}(\hat h,s)=\fun\Phi^{+1}_{\hat h}\circ s$.
%
%
%



We use fiberwise-holomorphic maps study the structure of the moduli stack of curves and to parametrise obstruction bundles. Section \ref{ev0 section} constructs a fiberwise-holomorphic map: 
\[\begin{tikzcd}\Ms(\hat {\ex B})^{+1}\rar\dar &\M^{+1}\dar
\\ \Ms(\hat{\ex B})\rar&\M\end{tikzcd}\]
This map may be insufficient for our purposes because it collapses bubbles in the domain of curves. In Section \ref{core section}, we  construct `core families' $\hat f$ with a group of automorphisms $G$ (of finite order, $\abs G$) so that there is a neighborhood $\mathcal U$ of $\hat f$ with a fiberwise-holomorphic map
\[\begin{tikzcd}\mathcal U^{+1}\rar{\fun\Phi^{+1}}\dar &\ex C(\hat f)/G\dar
\\\mathcal U\rar{\fun\Phi}&\ex F(\hat f)/G \end{tikzcd}\]
where,  for any curve $f$ in $U$, the corresponding  $\abs G$ connected components of the map $\fun\Phi^{+1}(f)\co \ex C( f)\times G\longrightarrow\ex C(\hat f)$ are degree-1 holomorphic maps onto curves in  $\ex C(\hat f)$. So long as $\mathcal U$ is small enough, we can also ensure that these maps are holomorphic isomorphisms.   Note that  $\hat f/G$ is also a moduli stack of curves, so we can also apply Definition \ref{curve fhm}, to define a  fibrewise holomorphic map of the following type.
\[\begin{tikzcd} \mathcal U^{+1}\dar \rar{\fun\Phi^{+1}} & (\hat f/G)^{+1}\dar
\\ \mathcal U \rar{\fun\Phi} & \hat f/G\end{tikzcd} \]
Noting that the stack $\ex C(\hat f)/G$ is equivalent to $(\hat f/G)^{+1}$, it is not hard to check that these two notions of a fiberwise holomorphic map are equivalent. 

\subsection{Topology of the ambient moduli stack of curves}\label{topology section}

\

A family of curves $\hat f$ comes with a natural  topology on the exploded manifold $\ex F(\hat f)$, so there is a natural notion of an open substack of the ambient moduli stack of curves, given by Definition \ref{open substack}. On the other hand,  Definition 11.4 of \cite{iec}, gives a notion of convergence of a sequence of $\C\infty1$ curves. Accordingly, we could define an open substack of $ \Ms$ as a substack $\mathcal U\subset \Ms$ such that every sequence of curves in $\Ms$ converging to a curve in $\mathcal U$ is eventually contained in $\mathcal U$. It is not immediately clear that these notions of open substack agree. So,  in this section, we show that such open substacks $\mathcal U$ of $\Ms$ are the substacks intersecting any family of curves in an open subset; and hence satisfying Definition \ref{open substack}.

\begin{lemma}\label{sequence family}Given any   sequence of curves in $\Ms$ converging to $f$ in $\C\infty1$,  there exists a family $\hat f$ of curves containing  $f$ and a subsequence $\{f_{i}\}$ of the given sequence such that each $f_{i}\in \hat f$, and within $\ex F(\hat f)$,  $f_{i}$ converges to $f$.
\end{lemma}

\pf

The definition of the $\C\infty1$ topology on $\Ms(\hat{\ex B})$ (from sections  7 and 11 of \cite{iec}) states that there exists 
\begin{itemize}
\item a family of curves $\hat g$ containing $f$,
\item
a sequence of curves $f_{i}'$ in $\hat g$ converging to $f$, 
\item a sequence of fiberwise almost-complex structures $j_{i}$ on $\ex C(\hat g)$ converging in $\C\infty1$ to the given almost-complex structure on $\ex C(\hat g)$,
\item and a sequence of sections $\psi_{i}$ of $\hat g^{*}T_{vert}\hat {\ex B}$ converging in $\C\infty1$ to $0$
\end{itemize}
so that there is an identification of $\ex C(f_{i})$ with  $\ex C(f_{i}')$ with the almost-complex structure $j'$, and so that the map $f_{i}$ is $f'_{i}$ followed by exponentiation of $\psi_{i}$ in some (fixed) metric on $\hat{\ex B}$.\footnote{See remark 6.6 of \cite{iec} for metrics on an exploded manifold.}

Let $\hat f_{0}$ be $\hat g\times \mathbb R$. The section $\psi_{i}$ being small in $\C\infty1$ implies that there exists a section $\psi_{i}'$ of $\hat f_{0}^{*}T_{vert}\hat{\ex B}\subset \hat g\times \mathbb R$  equal to $\psi$ at $\hat g\times\{i^{-1}\}$,  supported within the region $\hat g\times ((i+1)^{-1},(i-1)^{-1})$, and small in $\C\infty1$. As noted before  \cite[Definition 7.6]{iec},  convergence in $\C\infty1$
for fiberwise complex structures on $\ex C(\hat f_{0})$ or sections of $\hat f_{0}^{*}T_{vert}\hat{\ex B}$ is equivalent to convergence in some countable sequence of norms. By passing to a subsequence, we can assume that this $\psi_{i}'$ has size less than $2^{-i}$ in the first $i$ norms. Then $\sum_{i}\psi_{i}'$ is a $\C\infty1$ section of $\hat f_{0}^{*}T_{vert}\hat{\ex B}$, which restricts to be $\psi_{i}$ on $\hat g\times\{i^{-1}\}$ and which is zero on $\hat g\times\{0\}$. Define the map $\hat f$ to be $\hat f_{0}$ followed by exponentiation of $\sum_{i}\psi_{i}'$.

Similarly, by passing to a subsequence we can construct a $\C\infty1$ fiberwise complex structure $j$ on $\ex C(\hat f)$ which restricts to $\ex C(\hat g)\times \{i^{-1}\}$ to be $j_{i}$, and which is the original fiberwise complex structure on $\ex C(\hat g)\times\{0\}$. Now there is a sequence of inclusions $ f_{i}\longrightarrow \hat f$ converging to $f\longrightarrow \hat f$, as required.

\stop

\begin{lemma} \label{open pullback} For a substack $\mathcal U$ of $\Ms$, the following two conditions are equivalent.
\begin{enumerate}
\item Every sequence of curves in $\Ms$ that converges in $\C\infty1$ to a curve  in $\mathcal U$ is eventually contained in $\mathcal U$.
\item $\mathcal U$ is open in the sense of Definition \ref{open substack}: for all families $\hat f$ in $\Ms$, the subset  $\mathcal U(\hat f)\subset \ex F(\hat f)$ consisting of curves in $\mathcal U$ is open.
\end{enumerate}

\end{lemma}

\pf 

 Denote by $\mathcal U(\hat f)$ the subset of $\ex F(\hat f)$ consisting of curves in $\mathcal U$. 

Suppose that  any sequence of curves converging to a curve in $\mathcal U$ is eventually contained in $\mathcal U$. We must show that $\mathcal U(\hat f)\subset \ex F(\hat f)$ is open for all families of curves $\hat f$. Convergence within $\hat f$ is at least as strong as convergence within $\Ms$, so any sequence of points in $\ex F(\hat f)$ converging to a point in $\mathcal U(\hat f)\subset\ex F(\hat f)$ is eventually contained in $\mathcal U(\hat f)$. The topology on $\ex F$ is pulled back from a metrisable topology on $\totl{\ex F}$, so this sequential condition implies that  $\mathcal U(\hat f)\subset \ex F(\hat f)$ is open as required.

Alternately, suppose that $\mathcal U(\hat f)\subset \ex F(\hat f)$ is open for all families $\hat f$ in $\Ms$. Let $f_{i}$ be a sequence of curves converging to a curve  $f$ in $\mathcal U$. It remains to show that $f_i$ is eventually contained in $\mathcal U$.  Lemma \ref{sequence family} implies that there exists a subsequence converging to $f$ within some family $\hat f$, so that subsequence is eventually contained in $\mathcal U(\hat f)$, and therefore eventually contained in $\mathcal U$. This implies that every sequence of curves converging to a curve in $\mathcal U$ is eventually contained in $\mathcal U$. 

\stop

 \begin{remark} Lemma \ref{open pullback} does not immediately imply that $\C\infty1$ convergence of a sequence $f_{i}$ to $f$ is equivalent to $ f_{i}$ eventually being in every open neighborhood of $f$ --- without further study, it is not obvious that there are enough open substacks.  For $f$ any stable curve, the equivalence of these two notions of convergence  follows from Proposition \ref{smooth model family}; see Lemma \ref{topological sequential convergence}.
 \end{remark}

\subsection{Stable curves, $\Msw$}

\

\begin{defn}\label{stable} Call a $\C\infty1$ curve $f\co\ex C(f)\longrightarrow \hat{\ex B}$ {\bf stable} if it has only a finite number of automorphisms and its smooth part\footnote{See Definition 3.14 of \cite{iec}    the smooth part functor $\totl\cdot$. The smooth part of $f\co \ex C(f)\longrightarrow \ex B$ is a nodal curve $\totl{f}\co \totl{\ex C(f)}\longrightarrow \totl{\ex B}$ in the possibly singular space $\totl{\ex B}$. Automorphisms of $\totl{f}$ that do not act compatibly on the tangent space at different sides of a node do not come from automorphisms of $f$.} $\totl f$ has only a finite number of automorphisms.

Let $\Msw(\hat {\ex B})$ be the substack of $\Ms(\hat{\ex B})$ consisting of families of  stable curves.
\end{defn}

\begin{remark}We only need the condition that $f$ has a finite number of automorphisms if $\ex C(f)$ is the exploded manifold $\ex T$ from Example 3.4 of \cite{iec}, which has smooth part a point. Otherwise, every nontrivial automorphism of $f$ is also a nontrivial automorphism of $\totl f$.  There are a number of possible candidates for a definition of a stable curve $f$. A weaker definition, agreeing with the above on holomorphic curves,  is to require that $f$ has only a finite number of automorphisms and that $f$ is not a nontrivial refinement\footnote{Refinements are a special kind of blowup, described in Section 10 of \cite{iec}.} of another curve.  A  much weaker definition is to just require that $f$ has a finite number of automorphisms; the stack of curves satisfying this much weaker definition of stability is not sufficiently well behaved for us. 

For example, the stack $\Msw(pt)$ of stable curves mapping to a point is an exploded orbifold using the two stronger definitions, but not the much weaker definition of `stable'. 
\end{remark}

We shall use the notation $\Mod\subset \Msw$ for the stack of stable holomorphic curves. 

\

\begin{remark} The ambient moduli stack of stable curves $\Msw$ is a well-behaved open substack of $\Ms$. In particular, $\Msw$ has topology pulled back from a Hausdorff topological space, whereas the same is not true for $\Ms$.  Each non-constant holomorphic curve $f$ in $\Ms(\hat {\ex B})$ has a stabilization $f^{st}$ in $\Msw(\hat {\ex B})$, and a degree-1 genus-and-end\footnote{Recall from Section 8 of  \cite{iec} that an `end' of an exploded curve is a stratum with tropical part a half infinite line, corresponding to a  `marked point'  or `puncture'  on the smooth part of the curve.} preserving map $\psi\co \ex C(f)\longrightarrow \ex C(\hat f)^{st}$ such that $f=\psi\circ f^{st}$; see Lemma \ref{curve stab}. For curves mapping to a point, this construction works in  families, and defines a stablization map $\Ms(pt)\longrightarrow\Msw(pt)$; see Lemma \ref{ev0}. If the domain $\ex C(f^{st})$ is not the exploded manifold $\ex T$, then within $\Ms(\hat {\ex B})$, there exists a family of curves $\hat f$ containing $f$ and $f^{st}$ such that the image of $f$ in $\totl{\ex F(\hat f)}$ is in the closure of  the image of curves isomorphic to $f^{st}$. It follows that every open substack containing $f$ also contains $f^{st}$. The set of stable curves within a family $\hat f$ in $\Ms$ is always open, so $\Msw$ is an open substack containing $f^{st}$ but not $f$ if $f$ is unstable. Accordingly, $\totl{\Ms_{top}}$ is not Hausdorff. In contrast, we prove in Lemma \ref{hausdorff} that $\totl{\Msw_{top}}$ is Hausdorff.  

\end{remark}
 \subsection{Decorated curves, $\dmsw$}

\

It is sometimes desirable to consider a stack of stable curves with some extra structure. For example, we may wish to label ends of curves  or consider extra structure on the tropical part of curves, as is necessary for proving the tropical gluing formulae from \cite{gfgw}. 

\begin{defn}\label{decorated} Say that a map $\fun\pi\co\dmsw\longrightarrow \Msw$ of stacks (over the category of $\C\infty1$ exploded manifolds) is an {\bf ambient  moduli stack of decorated curves} if it obeys the following:
\begin{itemize}
\item  Given any family $\hat f$ in $\Msw$, there exists a pulled-back family of decorated curves $\fun\pi^{*}\hat f$ in $\dmsw$ along with a map $\fun\pi(\fun\pi^{*}\hat f)\longrightarrow \hat f$ so that the following holds: Given any $\hat g$ in $\dmsw$ and map $\psi\co \fun\pi(\hat g)\longrightarrow \hat f$ in $\Msw$, there exists a unique map $\fun\pi^*\psi\co \hat g\longrightarrow \fun\pi^{*}\hat f$ in $\dmsw$ making the following diagram commute:
\[\begin{tikzcd}[row sep= small]\hat g \rar[dotted]{ \fun\pi^*\psi}\dar{\fun\pi}& \fun\pi^{*}\hat f\dar{\fun\pi}
\\ \fun\pi(\hat g)\ar{dr}{\psi}\rar[dotted]&\fun\pi(\fun\pi^{*}\hat f)\dar 
\\ & \hat f\end{tikzcd}\]  
\item The derivative of the map 
\[\ex F(\fun\pi^{*}\hat f)\longrightarrow \ex F(\hat f)\]
at any point in $\ex F(\fun\pi^{*}f)$ is bijective. Moreover, this map is proper in the sense that the inverse image of any compact  closed subset of $\ex F(\hat f)$ is a compact subset of $\ex F(\fun\pi^{*}\hat f)$.
 
\end{itemize}
\end{defn}

Using the general language of Definition \ref{representable map}, $\fun\pi\co \dmsw\longrightarrow \Msw$ is a proper representable  submersion which is also an immersion. 

From now on, we  use  $\dmsw$ to refer to an ambient moduli stack of decorated curves satisfying the above definition. We shall also use the notation $\dmod\subset \dmsw$ for the stack of stable holomorphic decorated curves.   Of course, $\Msw$ itself is an ambient moduli stack of (trivially) decorated curves, so any theorem stated for $\dmsw$ (or $\dmod$) applies also to $\Msw$ (or $\Mod$).

 In the category of smooth manifolds, the conditions on the map $\ex F(\fun\pi^{*}\hat f)\longrightarrow \ex F(\hat f)$ make it into a finite cover. This is not quite the case here, because $\fun\pi$ need not be etale; see Remark \ref{submersion etale}. Moreover,  recall (Definition 3.15 of \cite{iec}), that the correct analogue of `compact' or `proper' in the category of exploded manifolds is `complete'. An example of a map which is proper but not complete is the inclusion $\et 1{(0,1]}$ into $\et 1{[0,1]}$; see examples 3.8 and 3.12 in \cite{iec}.

The universal property defining $\fun\pi^{*}\hat f$ implies that  any morphism $\hat g\longrightarrow \hat f$ in $\Msw$ lifts canonically to $\fun\pi^{*}\hat g\longrightarrow \fun\pi^{*}\hat f$ in $\dmsw$. In particular, the group of automorphisms of $\hat f$ act on $\fun\pi^{*}\hat f$. Another important consequence of this universal property is that if $\hat f/G$ represents a substack of $\Msw$ (in the sense of Definition \ref{represented by family}), then $\fun\pi^{*}\hat f/G$ represents a substack of $\dmsw$.

\

\begin{remark}\label{decoration base change} Given any base change\footnote{Lemma  10.4 of \cite{iec} implies that we can change the base of families using fiber products as usual.} of our targets, 
\[\begin{tikzcd}\hat{\ex B}'\rar\dar&\hat {\ex B}\dar
 \\ \ex B_0'\rar& \ex B_0\end{tikzcd}\]
if we have $\dmsw(\hat {\ex B})$ defined, we can define $\dmsw(\hat{\ex B'})$ as the fiber product of $\dmsw(\hat {\ex B})$ with $\Msw(\hat{\ex B'})$ over $\Msw(\hat{\ex B})$. So,  a decorated family of curves in $\hat{\ex B'}$ is a family of curves in $\hat {\ex B'}$ with a choice of decoration (when this family of curves is considered as a family in $\hat {\ex B}$).  

Be warned that this definition of $\dmsw(\hat{\ex B'})$ sometimes depends on the particular choice of map $\hat{\ex B}'\longrightarrow \hat{\ex B}$.
\end{remark}

\

%
%

\subsection{Tangent space of $\Msw$ and $\dmsw$ }\label{tangent section}

\

\

In this section, we  define the tangent space of the open substack $\Msw\subset\Ms$ of stable curves. The discussion also applies  to any ambient moduli stack of decorated curves $\dmsw$. In fact, the tangent space will never depend on the decorations because the lifting property of Definition \ref{decorated} implies that a deformation (parametrized by $\mathbb R$) of a decorated family $\hat f$ of curves in $\dmsw$ is equivalent to a deformation of the undecorated image of $\hat f$ in $\Msw$. We  first discuss the tangent space to a single curve $f$ with codomain $\ex B$, and then define the relative tangent space in the case of a family of codomains $\hat {\ex B}\longrightarrow \ex B_0$ and the (relative) tangent sheaf in the case of a family of curves.

\subsubsection{$T_{f}\Msw(\ex B)$}

\

In this section, we give an explicit construction of the tangent space $T_{f}\dmsw(\ex B)$  as a vector-space, with a natural action of the automorphism group of $f$. The tangent space $T_f\Msw(\ex B)$ is a special case. Lemmas \ref{tangent1} and \ref{tangent2} will  imply that this $T_{f}\dmsw(\ex B)$ is equivalent to the following  natural definition:

\begin{defn}\label{tangent as set} Suppose that $\hat f_{i}$ are families  parametrized by $\mathbb R$ in a stack $\mathcal X$ over the category of exploded manifolds, with a given inclusion  $f\longrightarrow \hat f_{i}$ over $0\in\mathbb R$. Say that $\hat f_{1}$ is {\bf tangent to} $\hat f_{2}$ at $f$ if there exists a family $\hat g$ (parametrised by some $\mathbb R^n$) and a commutative diagram
\[\begin{tikzcd}[row sep=tiny] & \hat f_{1}\ar{dr} 
\\ f\ar{ur}\ar{dr}&&\hat g
\\ &\hat f_{2}\ar{ur}\end{tikzcd}\] 
such that the corresponding maps $\ex F(\hat f_{i})=\mathbb R\longrightarrow \ex F(\hat g)$ have the same derivative at $0$.

Define the deformations of $f$,  $\mathcal {DEF}_{f}$,  to be the set\footnote{Assuming that $\mathcal X$ is small enough that there is a set of such deformations of $f$ --- any stack of interest to us is equivalent to such a small stack.} of families $\hat f$ with $\ex F(\hat f)=\mathbb R$, and with a chosen isomorphism of $ f$ with $\hat f\rvert_0$. Then define $T_f\mathcal X$ to be the quotient of $\mathcal {DEF}_{f}$  by the equivalence relation generated by being tangent at $f$.

Given any map of stacks $\fun\Phi\co\mathcal X\longrightarrow \mathcal X'$ and individual object $f\in \mathcal X$, there is an induced derivative $T_f\fun\Phi\co T_f\mathcal X\longrightarrow T_{\fun\Phi(f)}\mathcal X'$ sending the equivalence class of $\iota\co f\rightarrow \hat f$ to the equivalence class of $\fun\Phi(\iota)\co \fun\Phi(f)\longrightarrow \fun\Phi(\hat f)$.

Similarly, given any isomorphism $\psi\co f\longrightarrow f'$, there is a  bijection $T_{f}\psi\co T_{f}\mathcal X\longrightarrow T_{f'}\mathcal X$ sending the equivalence class of $\iota\co f\rightarrow \hat f$ to the equivalence class of $\iota\circ \psi^{-1}$.
\end{defn}

In the case of a stack represented by an exploded manifold, this definition coincides with the usual tangent space. Moreover, it works for exploded orbifolds; see Definition \ref{orbifold}.

\begin{example}\label{orbifold tangent space} Let $\fun\pi\co \ex U\longrightarrow \ex U/G$ be the quotient map from Example \ref{quotient map}. Then $T_{x}\fun\pi\co T_{x}\ex U\longrightarrow T_{\fun\pi(x)}(\ex U/G)$ is a bijection. In $\ex U$, deformations of $x$ are just  smooth maps $\gamma\co \mathbb R\longrightarrow \ex U$ with $\gamma(0)=x$. In $\ex U/G$, these are sent to the families  with domain the trivial $G$--bundle over $\mathbb R$  and  $\fun\pi(\gamma)(t,g)=g*\gamma(t)$. Similarly $\fun\pi(x)$ is the map $G\longrightarrow U$ given by $g\mapsto g*x$. Any deformation of $\fun\pi(x)$ is therefore canonically trivialized, and is $\fun\pi(\gamma)$ for a unique deformation $\gamma$ of $x$. Similarly deformations of $x$ parametrized by $\mathbb R^{n}$ are in bijection with deformations of $\fun\pi(x)$ parametrized by $\mathbb R^{n}$, so two deformations are tangent in $\ex U$ if and only if they are tangent in $\ex U/G$. 

More generally, given a surjective etale map $\fun\Phi\co \ex U\longrightarrow \mathcal X$ (so the exploded orbifold is the quotient stack of the groupoid $\ex U\times_{\mathcal X}\ex U\rightrightarrows \ex U$) the derivative $T_{x}\fun\Phi\co T_{x}\ex U\longrightarrow T_{\fun\Phi(x)}\mathcal X$ is a bijection. This is in agreement with defining $T\mathcal X$ as the quotient stack of the groupoid $T\lrb{\ex U\times_{\mathcal X}\ex U}\rightrightarrows T\ex U$.

\end{example}

\

Let us now give an explicit description for $T_{f}\dmsw(\ex B)$,  equipping the set from Definition \ref{tangent as set} with the structure of a vector-space.
 Given a curve $f\co\ex C\longrightarrow \ex B$,   the tangent space $T_{f}\dmsw(\ex B)$ of $\dmsw(\ex B)$ at $f$ is defined using the following short exact sequence,  discussed below.
  \begin{equation} \label{hles}0\longrightarrow \Gamma(T\ex C)\xrightarrow{h} \Gamma^{0,1}(T\ex C\otimes T^*\ex C)\times \Gamma( f^{*} T\ex B)\longrightarrow T_{f}\dmsw\longrightarrow 0\end{equation}

 Let $\Gamma(f^{*}T{\ex B})$ denote the space of $\C\infty1$ sections of $f^{*}T{\ex B}$. The action of the almost complex structure $J$ on sections of $f^{*}T\ex B$ makes $\Gamma(f^{*}T{\ex B})$  a complex vector-space. Think of $\Gamma(f^{*}T\ex B)$ as the infinitesimal variations of a map from a fixed domain\footnote{Note that for now, we are using $\ex C$ as shorthand for $\ex C(f)$, so $\ex C$ means an exploded manifold, not a functor.} $\ex C$.

Let $\Gamma^{0,1}(T\ex C\otimes T^*\ex C)$ denote the space of  $\C\infty1$, $j$-anti-linear sections $\alpha$ of $T\ex C\otimes T^{*}\ex C$ that vanish on edges\footnote{Recall, from Definition 8.3 \cite{iec}, that an edge of a curve $\ex C$  is a stratum isomorphic to $\et 1{(a,b)}$. The tropical part of such a stratum is an edge of the graph $\totb{\ex C}$. Each edge corresponds to a node or marked point in the nodal curve that is the smooth part of $\ex C$. For the stratified structure of an exploded manifold and the tropical part functor, see Section 4 of \cite{iec}. } of $\ex C$.
This $\Gamma^{0,1}(T\ex C\otimes T^*\ex C)$ is the space of infinitesimal variations of almost-complex structure on a fixed domain $\ex C$. The action of $j$ on the left makes $\Gamma^{0,1}(T\ex C\otimes T^*\ex C)$  a complex vector-space.

To obtain $T_{f}\dmsw$, we quotient $\Gamma(T\ex C)\times\Gamma^{0,1}(T\ex C\otimes T^*\ex C)$ by the reparametrization action: in particular, we quotient by the image of the map $h$ below.
\begin{equation}\label{hdef}\begin{array}{ccc}h\co\Gamma(T\ex C)&\longrightarrow &\Gamma^{0,1}(T\ex C\otimes T^*\ex C)\times \Gamma( f^{*} T\ex B)
\\ v&\mapsto & (L_{v}j,df(v))\end{array}\end{equation}
 If $f$ is in $\dmsw$, then $f$ has no infinitesimal automorphisms and $h$ is injective. We define $T_{f}\dmsw$ to be the quotient of $\Gamma^{0,1}(T\ex C\otimes T^*\ex C)\times\Gamma(f^{*}T\ex B)$ by the image of $h$.  We will show, in  Lemmas \ref{tangent1} and \ref{tangent2}, that this $T_f\dmsw$ is the tangent space at $f$  in the sense of  Definition \ref{tangent as set}.

 The following lemma shows that change of $j$ under the flow of a vector-field $v$,  $L_{v}j$,  can be regarded as $2j\circ \dbar v$.

\begin{lemma}\label{Lvdbar}  Let $\nabla$ be any holomorphic connection on $T\ex C$. Then 
\[L_{v}j=j\circ \nabla v-  \nabla v\circ j\]
and, in local holomorphic coordinates
\[L_{v}j=2j\circ \dbar v\ .\]
\end{lemma}
\pf

If $e$ is any holomorphic vector-field, then $j\circ\nabla e-(\nabla e)\circ j=0$, so $j\circ \nabla he-(\nabla he)\circ j=j\circ 2(\dbar h) e$. Therefore, in holomorphic coordinates, where $v$ can be considered as a complex function,  we have
\[j\circ \nabla v-\nabla jv=j\circ  2\dbar v\ .\]
Now we can calculate in holomorphic  coordinates, where $j=\partial_{y}\otimes dx -\partial_{x}\otimes dy$. Write $v=v_{1}\partial_{x}+v_{2}\partial_{y}$.  Then 
\[\begin{split}L_{v}j&= -(\partial_{y}v_{1}\partial_{x}+\partial_{y}v_{2}\partial_{y}))\otimes dx +\partial_{y}\otimes (\partial_{x}v_{1}dx+\partial_{y}v_{1}dy) 
\\&\ \  +(\partial_{x}v_{1}\partial_{x}+\partial_{x}v_{2}\partial_{y})\otimes dy-\partial_{x}\otimes (\partial_{x}v_{2}dx+\partial_{y}v_{2}dy)
\\&=(\partial_{y}v_{1}+\partial_{x} v_{2})\lrb{-\partial_{x}\otimes dx+\partial_{y}\otimes dy}
\\ &\ \ +(\partial_{x}v_{1}-\partial_{y}v_{2})\lrb{\partial_{y}\otimes dx+\partial_{x}\otimes dy}
\end{split}\]
On the other hand, 
\[\begin{split}j\circ 2\dbar v&=j\circ dv-dv\circ j
\\&= \partial_{x}v_{1}\partial_{y}\otimes dx+\partial_{y}v_{1}\partial_{y}\otimes dy-\partial_{x}v_{2}\partial_{x}\otimes dx-\partial_{y}v_{2}\partial_{x}\otimes dy
\\ & \ \ +\partial_{x}v_{1}\partial_{x}\otimes dy-\partial_{y}v_{1}\partial_{x}\otimes dx+\partial_{x}v_{2}\partial_{y}\otimes dy-\partial_{y}v_{2}\partial_{y}\otimes dx 
\\&=L_{v}j\end{split}\]

\stop

Lemma \ref{Lvdbar} implies that $L_{jv}j=j\circ L_{v}j$, so the corresponding map 
\[\begin{array}{ccc}\Gamma(T\ex C)&\longrightarrow &\Gamma^{0,1}(T\ex C\otimes T^*\ex C)
\\ v&\mapsto & L_{v}j\end{array}\]
is complex linear, and therefore the map  $h$, from (\ref{hles}), (\ref{hdef}), is complex linear if and only if $f$ is holomorphic. Accordingly,   $T_{f}\dmsw$ is a complex vector-space whenever $f$ is holomorphic, and is otherwise a real vector-space.
 
 \

\subsubsection{$T_{f}\dmsw(\hat {\ex B})$ and $T_{f}\dmsw(\hat {\ex B})\ov{\ex B_{0}}$}

\

In the case of a family of targets $\hat{\ex B}\longrightarrow \ex B_0$, define the relative tangent space $T_{f}\dmsw(\hat{\ex B})\ov{\ex B_0}$ as $T_{f}\dmsw(\ex B)$, where $\ex B$ is the member of the family $\hat{\ex B}$ containing the image of $f$. The following is the defining short exact sequence: 
 \begin{equation}\label{relative tangent sequence}0\longrightarrow \Gamma(T\ex C)\xrightarrow{h}\Gamma^{0,1}(T\ex C\otimes T^*\ex C)\times \Gamma( f^{*} T_{vert}\hat{\ex B})\longrightarrow T_{f}\dmsw(\hat{\ex B})\ov{\ex B_0}\longrightarrow 0\end{equation}

\

We define the tangent space $T_{f}\dmsw(\hat{\ex B})$ similarly to the case of a single target $\ex B$, except we use the notation $\Gamma_{\ex B_{0}}(f^{*}T\hat {\ex B})$ to denote  sections of $ f^{*}T\hat{\ex B}$ that become constant  when composed with the derivative of the map $\hat{\ex B}\longrightarrow \ex B_0$. We use $\Gamma_{\ex B_{0}}(f^{*}T\hat {\ex B})$ instead of $\Gamma(f^{*}T\hat {\ex B})$ because we are interested in infinitesimal variations of $f$ as a curve  contained in a fiber of $\hat{\ex B}\longrightarrow \ex B_{0}$ instead of variations of $f$ as a map to $\hat{\ex B}$. So, the following is the defining exact sequence for $T_{f}\dmsw(\hat{\ex B})$.
\[0\longrightarrow  \Gamma(T\ex C)\xrightarrow{h}\Gamma^{0,1}(T\ex C\otimes T^*\ex C)\times \Gamma_{\ex B_0}( f^{*} T\hat{\ex B})\longrightarrow T_{f}\dmsw(\hat{\ex B})\longrightarrow 0\]

\

\subsubsection{Derivatives}\label{derivatives section}

\

 Given any family of curves $\hat f$ in $\dmsw$ containing $f$ and any vector $v$ in $T_{f}\ex F(\hat f)$, define an element, $[v]$ of $T_{f}\dmsw$ by differentiating $\hat f$ in the direction of $v$ as follows: Extend $v$ to a vector-field on $\ex F(\hat f)$, and lift it to a $\C\infty1$ section $v'$ of $T\ex C(\hat f)$, then define
 \[[v]:=[L_{v'}j\rvert_{\ex C(f)},df(v')]\]
 where the righthand side indicates the image of $(L_{v'}j\rvert_{\ex C(f)},df(v'))$ in $T_{f}\dmsw$ under the quotient by the image of $h$.
 
 We must verify that $[v]$ is well defined.  
  The fiberwise almost-complex structure $j$ is a section of $T^{*}_{vert}{\ex C(\hat f)}\otimes T_{vert} {\ex C(\hat f)}$, and  the Lie derivative of $j$ with respect to any lifted vector-field is again a section of $T^{*}_{vert}{\ex C(\hat f)}\otimes T_{vert} {\ex C(\hat f)}$, because the flow of lifted vector-fields respects the fibers of $\ex C(\hat f)\longrightarrow \ex F(\hat f)$. The Lie derivative $L_{v'}j\rvert_{\ex C(f)}$ does not depend on the choice of extension of $v$, because the flow of any vector-field vanishing on a fiber preserves that fiber, so $L_{v'}j\rvert_{\ex C(f)}$ does not depend on a choice of extension. Moreover, Lemma \ref{Lvdbar} implies that $L_{v'}j\rvert_{\ex C(f)}$ is a section of $\Gamma^{0,1}(T^{*}\ex C(f)\otimes T\ex C(f))$, because $v'$ is locally equal to a $j$--preserving vector-field plus a section of $T\ex C(f)$. So, $[v]:=[L_{v'}j\rvert_{\ex C(f)},df(v')]$ really is a element of $T_{f}\dmsw$. It remains to check that $[v]$ does not depend on our choice of lift $v'$: any other lift $v''$ of $v$ will differ on $\ex C(f)$ from $v'$ by a  section of $T\ex C(f)$, so $(L_{v''}j\rvert_{\ex C(f)},df(v''))$ will differ from $(L_{v'}j,df(v'))$ by a vector in the image of $h$. Therefore, $[v]:=[L_{v'}j\rvert_{\ex C( f)},df(v')]$ gives a well-defined element of $T_{f}\dmsw$.

Therefore for any  curve $f$ in $\hat f$, we have a well-defined linear map \[T_f\hat f\co T_{f}\ex F(\hat f)\longrightarrow T_{f}\dmsw\ \] sending $v$ to $[v]$. Moreover, this map restricts to a well-defined linear map \begin{equation}\label{derivative def}T_f\hat f\co T_{f}\ex F(\hat f)\ov{\ex B_0}\longrightarrow T_{f}\dmsw(\hat{\ex B})\ov{\ex B_0}\ .\end{equation} 

\

 Lemmas \ref{tangent1} and \ref{tangent2} below imply that, as a set,  $T_{f}\dmsw$ is the quotient of the set of families of curves parametrized by $\mathbb R$ containing $f$ at $0$ by an equivalence relation generated by declaring two such families equivalent if they are tangent at $0$ within a two-dimensional family of curves. It follows that given any $\C\infty1$ map \[\fun\Phi\co\dmsw\longrightarrow \ex A\] there is an induced linear tangent map  \[T_{f}\fun\Phi\co T_{f}\dmsw\longrightarrow T_{\fun\Phi(f)}\ex A\]
so that,  if $\fun\Phi(\hat f)$ is the induced map $\ex F(\hat f)\longrightarrow \ex A$,

 \[T_{f}\fun\Phi(\hat f)=T_{f}\fun\Phi \circ T_f\hat f\ .\]

\begin{lemma}\label{tangent1}Given any curve $f$ in $\dmsw(\hat{\ex B})$ and vector in $T_{f}\dmsw(\hat {\ex B})$, there exists a family of curves $f_t$ parametrized by $\mathbb R$ with that vector in the image of $T_ff_t$. 

\end{lemma}
\pf

We must construct a family with a given derivative. The definition of $T_{f}\dmsw$ does not depend on the particular decoration chosen. So, we shall prove this lemma for $T_{f}\Msw$.  Then Definition \ref{decorated} implies that to complete the proof for any $\dmsw$, we can construct the required family of curves in $\Msw$ then lift them to $\dmsw$.

 To avoid any issues with the precise nature of almost-complex structures at edges\footnote{Definition 8.3 of \cite{iec}.} of a curve, we first reduce to the case that the variation in almost-complex structure on the domain $\ex C$  is described by  a section of $\Gamma^{0,1}(T\ex C\otimes T^{*}\ex C)$
 vanishing in a neighborhood of all edges of $\ex C$. Consider a neighborhood of an edge of the domain $\ex C$ of $f$. A small such neighborhood  is isomorphic  an open subset of $\et 1{[0,l]}$ or $\et 1{[0,\infty)}$ with its standard complex structure.  A vector-field  $v$  in standard coordinates is then  the real and imaginary parts of $\tilde z\frac \partial{\partial\tilde  z}$ times $\C\infty1$ functions. As noted in Lemma \ref{Lvdbar}, in these coordinates $L_{v}j$  is $2j$ times the standard $\dbar$ operator. Therefore, Theorem \ref{f replacement}   implies that given any section in $\Gamma^{0,1}(T\ex C\otimes T^{*}\ex C)$, there exists a $\C\infty1$ vector-field $v$ so that $L_{v}j$ is equal to the given section in a neighborhood of each edge. We can therefore represent any given vector in $T_{f}\dmsw(\hat {\ex B})$ using a section of $\Gamma^{0,1}(T\ex C\otimes T^{*}\ex C)$ vanishing on a neighborhood of edges of $\ex C$.

Given any section  $\alpha$ in $\Gamma^{0,1}(T\ex C\otimes T^{*}\ex C)$ vanishing in a neighborhood of edges of $\ex C$, we can construct  a family of almost-complex structures $j_{t}$ on $\ex C$ so that $j_{t}$ is the original almost-complex structure when $t=0$ and near edges, and $\frac{\partial}{\partial t}j_{t}=\alpha$ at $t=0$. Also, we can construct a $\C\infty1$ family of maps $f_{t}$ from $\ex C$ to $\hat{\ex B}\longrightarrow \ex B_0$ so that $\frac\partial{\partial t}f_{t}$ at $t=0$ is any given $\C\infty1$ section of $f^{*}T\hat{\ex B}$ projecting to a constant map to $T\ex B_0$. It follows that given any vector in $T_{f}\Msw$, there exists a family of curves parametrized by $\mathbb R$ with the given derivative. 

\stop

\begin{lemma}\label{tangent2}
Suppose that $\hat f_{1}$ and $\hat f_{2}$ are two one-dimensional families in $\dmsw$, containing $f$, such that the image of $T_f\hat f_i$ within $T_{f}\dmsw$ coincides.  Then there exists another one-dimensional family $\hat f_{0}$ containing $f$ and two   $2$--dimensional families $\hat g_{i}$ with given maps  $\hat f_{0}\longrightarrow \hat g_{i}$ and  $\hat f_{i}\longrightarrow \hat g_{i}$ so that the maps 
\[\ex F(\hat f_{0})\longrightarrow \ex F(\hat g_{i})\]
\[\ex F(\hat f_{i})\longrightarrow \ex F(\hat g_{i})\]
are tangent at $f$.
\end{lemma}

\pf
Again, the lifting property of Definition \ref{decorated} implies that it suffices to prove this lemma for $\Msw$.

The reason that an extra family $\hat f_{0}$ is needed, and the heart of the technical problem to be overcome, is the following observation: Given a smooth section of an infinite-dimensional vector-bundle over the real line, there may not exist a finite-dimensional sub-vector-bundle containing the given smooth section. This problem disappears if we require some kind of non-degeneracy at zeros of the section.    

Without losing generality, assume our $1$--dimensional families $\hat f_i$ are parametrised by $\mathbb R$,  and assume that that $f$ is the curve over $0\in \mathbb R=\ex F(\hat f_i)$.  Let $t$ indicate the coordinate parametrizing $\ex F(\hat f_{i})$. The following claim is a version of Hadamard's lemma:
\begin{claim}\label{Hadamard} If a $\C\infty1$ section $\nu$ of a vector-bundle over $\ex C(\hat f_{i})$ vanishes at $t=0$, then  
$\nu=t\nu'$ where $\nu'$ is also a $\C\infty1$ section. If $\nu$ vanishes on all edges of curves in $\ex C(\hat f_{i})$, then $\nu'$ does too. \end{claim}

To prove Claim \ref{Hadamard}, choose a connection $\nabla$ on our vector-bundle, and a lift, $v$, of $\frac \partial{ \partial t}$ to a $\C\infty1$ vector-field on $\ex C(\hat f_{i})$. The flow of $v$ identifies $\ex C(\hat f_{i})$ with $\ex C(f)\times \mathbb R$ and together with $\nabla$ allows us to trivialize the vector-bundle in the $\mathbb R$ direction. Then we can write
\[\nu(z,t)=\int_{0}^{1}(\nabla_{tv}\nu)(z,ts)ds=t\int_{0}^{1}(\nabla_{v}\nu)(z,ts)ds\ .\]
 Note that $\nabla_{v}\nu$ is a $\C\infty1$ section of our vector-bundle vanishing on edges of curves in $\ex C(\hat f_{i})$ if $\nu$ does, therefore
 \[\nu'(z,t):=\int_{0}^{1}(\nabla_{v}\nu)(z,ts)ds\]
is a $\C\infty1$ section of our vector-bundle vanishing on edges of curves in $\ex C(\hat f_{i})$ if $\nu$ does. This completes the proof of Claim \ref{Hadamard}.

 Claim \ref{Hadamard} implies that, if a section of a vector-bundle over $\ex C(\hat f_{i})$ vanishes at $t=0$ to order $n-1$, but has nonvanishing $n$th derivative, then it is equal to $t^{n}$ times a $\C\infty1$ section  not vanishing at $0$.

If the domain of $f$ is not stable, then we can choose some injective codimension--$2$ map $\iota\co S\longrightarrow \hat{\ex B}$  such  that $\ex C(f)$ with the extra marked points in $f^{-1}\iota(S)$ is stable, $f$ is transverse to $\iota$,    and all intersection points  are the image of some relatively compact open subset $U\subset S$. Then a neighborhood of $f$ in $\hat f_{i}$ will remain transverse to $\iota\co  U\longrightarrow \hat{\ex B}$; see \cite[Section 9]{iec}, noting that as $\ex F(\hat f_i)=\mathbb R$, there are no departures from usual transversality theory.  For simplicity,  assume that $\hat f_{i}$ remains everywhere transverse to $\iota$ --- as it suffices to prove our lemma for a neighborhood of $f$ in $\hat f_{i}$,  we shall repeatedly assume $\hat f_{i}$ is small enough as needed.

 As explained in  Section \ref{ev0 section},  there is a unique $\C\infty1$ map $s_{i}\co\mathbb R\longrightarrow \M$ so that the pullback of the universal curve over $\M$ is $\ex C(\hat f_{i})$ with extra edges at each of the points in $\hat f_{i}^{-1}(\iota (U))$. Moreover, as proved in Section \ref{ev0 section},  we can locally represent $\M$ as the quotient of a given family of stable curves by a finite group of automorphisms.   These two maps $s_{i}\co\mathbb R\longrightarrow \M$ are tangent at $0$. Chose another map $s_{0}\co\mathbb R\longrightarrow \M$  tangent to both $s_{i}$ at $0$, but not equal to either of them to second order. Then Claim \ref{Hadamard} implies that there exist two smooth maps $\tilde s_{i}\co\mathbb R^{2}\longrightarrow\M $ so that $s_{i}(t)$ is equal to $\tilde s_{i}(t,0)$, and $s_{0}(t)=\tilde s_{i}(t,t^{2})$. 

Define the domain $\ex C(\hat g_{i})$ of $\hat g_{i}$ to be the pullback of the universal curve by $\tilde s_{i}$, with the extra ends (corresponing to $\hat f_{i}^{-1}(S)$) removed. This removal operation replaces a neighborhood of an end with its smooth part,  replacing an open subset of $\et 11$ with the corresponding open subset  $\mathbb C=\totl{\et 11}$. Choose a Riemannian metric on $\M$. The topology induced by such a metric is significantly stronger than the usual exploded manifold topology. Small enough open balls in this metric are equivalent to the quotient of $\mathbb R^n$ by a finite group of diffeomorphisms, and we can always trivialise the universal curve over such open balls. 
Around the image of $f$ (using the strong topology induced by the metric),  identify all the domains of curves in the the universal curve over $\M$, and make the identifications holomorphic in a neighborhood of edges. Then, restricted to a neighborhood of $0\in\mathbb R^{2}=\ex F(\hat g_i)$, we can regard $\ex C(\hat g_{i})$, and hence $\ex C(\hat f_{i})$, to be given by a family of almost-complex structures on a fixed domain, and regard $\ex C(\hat f_{i})$ to be identical, but have different almost-complex structures for $i=0,1,2$. With these identifications, the fact that $\hat f_{i}$ are tangent for $i=1,2$ implies that the maps  $\hat f_{i}$ are equal to first order restricted to $\ex C(f)$. We can therefore choose a map $\hat f_{0}$ sending $f_{i}^{-1}(\iota(U))$ to $\iota(U)$,  equal $\hat f_{i}$ to first order at $\ex C(f)$, but not equal to either of them to second order at $\ex C(f)$. Claim \ref{Hadamard} then implies that there exist families $\hat g_{i}\co\ex C(\hat g_{i})\longrightarrow \ex B$ so $\hat f_{i}$ is the restriction of $\hat g_{i}$ to $\ex C(\hat f_{i})$, and $\hat f_{0}$ is the restriction of $\hat g_{i}$ to $\ex C(\hat f_{0})$.

  \stop
  
  \
  
Lemmas \ref{tangent1} and \ref{tangent2} allow us to differentiate any $\C\infty1$ map $\fun\Phi \co\dmsw\longrightarrow \ex X$ to obtain a linear map
 
 \[T_{f}\fun\Phi\co T_{f}\dmsw\longrightarrow T_{\fun\Phi(f)}\ex X\]
 defined so that for all $\C\infty1$ families of curves $\hat f$ in $\dmsw$ containing $f$, the diagram
 \[\begin{tikzcd}T_{f}\dmsw\rar{T_{f}\fun\Phi}& T_{\fun\Phi(f)}\ex X
 \\T_{f}\ex F(\hat f)\uar{T_f\hat f} \ar{ur}[swap]{T_{f}\fun\Phi(\hat f)}\end{tikzcd}\]commutes, where $\fun\Phi(\hat f)\co\ex F(\hat f)\longrightarrow \ex X$ is the map induced by $\fun\Phi$. Lemma \ref{tangent1} implies that $T_{f}\fun\Phi$ is unique, and Lemma \ref{tangent2} implies that it is well defined.
  
  \subsubsection{$\mathbb R$--nil vectors}
  
  \

  Recall, \cite[Definition 6.8]{iec}, that the integral vectors of an exploded manifold are the vectors, $v$, such that $\tilde z^{-1}v\tilde z$ is always an integer for any locally defined exploded function $\tilde z\co U\longrightarrow \mathbb C^*\e{\mathbb R}$. Such a vector always acts as the zero derivation on smooth or  $\C\infty1$ functions, and is an example of a $\mathbb R$--nil vector. 

\begin{defn}\label{R-nil}A $\mathbb R$--nil vector $v$ on an exploded manifold is a vector which acts as a zero derivation on any $\C\infty1$, $\mathbb R$--valued function. \end{defn}

There is a canonical complex structure on the $\mathbb R$--nil vectors at a point so that $(Jv)(\tilde z)=i(v\tilde z)$, and the $\mathbb R$--nil vectors at a point are always the complex-linear span of the integral vectors. Clearly, derivatives always send $\mathbb R$--nil vectors to $\mathbb R$--nil vectors, and such derivative maps are always complex with respect to the canonical complex structure on $\mathbb R$--nil vectors. 

The bundle of $\mathbb R$--nil vectors on a stratum of an exploded manifold have a canonical flat connection:  the connection  preserving the canonical complex structure and the lattice of integral vectors. So, a constant $\mathbb R$--nil vector-field  is some sum of complex numbers times integral vector-fields.
 
 \

 There is a similar notion of integral and $\mathbb R$--nil vectors on $T_{f}\dmsw$ and $T_{f}\dmsw\ov{\ex B_0}$. 
 
 \begin{defn}A vector $v$ in $T_{f}\dmsw$ is $\mathbb R$--nil (or integral) if there exists a family of curves $\hat f$ containing $f$ and an $\mathbb R$--nil (or integral) vector $v'\in T_f\ex F(\hat f)$ such that $v=T_f\hat f(v')$. \end{defn}
 
 As the canonical complex structure on $\mathbb R$--nil vectors is compatible with all exploded maps, it follows that there is a canonical complex structure on the $\mathbb R$--nil vectors within $T_{f}\dmsw$.

 \subsubsection{$T_{\hat f}\dmsw$}
 
  \

For a family of curves $\hat f$, it is not clear that the vector-spaces $T_{f}\dmsw$ for all $f$ in $\hat f$ fit together to form a vector-bundle over $\ex F(\hat f)$. On the other hand, there is a natural tangent sheaf, $T_{\hat f}\dmsw$ over the exploded manifold $\ex F(\hat f)$,  encoding first order  deformations of $\hat f$ parametrized by $\ex F(\hat f)$. This tangent sheaf is defined by the short exact sequence
\begin{equation}\label{Tsheaf}\begin{tikzcd} \Gamma(T_{vert}\ex C(\hat f))\rar[hook]{h}& \Gamma^{0,1}(T_{vert}\ex C(\hat f)\otimes T_{vert}^*\ex C(\hat f))\times \Gamma_{\ex B_0}( \hat f^{*} T \hat{\ex B})\dar
\\ & T_{\hat f}\dmsw(\hat{\ex B})\end{tikzcd}\end{equation}
and defined in the relative case by the following short exact sequence:
\begin{equation}\label{rTsheaf}\begin{tikzcd}\Gamma(T_{vert}\ex C(\hat f))\rar[hook]{h}& \Gamma^{0,1}(T_{vert}\ex C(\hat f)\otimes T_{vert}^*\ex C(\hat f))\times \Gamma( \hat f^{*} T_{vert}\hat{\ex B})\dar\\& T_{\hat f}\dmsw(\hat{\ex B})\ov{\ex B_0}\end{tikzcd}\end{equation}

Again, $\Gamma^{0,1}(T_{vert}\ex C(\hat f)\otimes T_{vert}^*\ex C(\hat f))$ indicates $\C\infty1$ sections of $(T_{vert}\ex C(\hat f)\otimes T_{vert}^*\ex C(\hat f))$ that vanish on edges of curves and that anti-commute with $j$. These sections  represent infinitesimal variations of complex structure on the domain. Also, $\Gamma_{\ex B_{0}}( \hat f^{*} T\hat{\ex B})$ indicates $\C\infty1$ sections of $\hat f^{*}T\hat{\ex B}$ lifting sections  of $T\ex B_0$ over $\ex F(\hat f)$ using the following diagram.
\[\begin{tikzcd}\ex C(\hat f)\dar\rar{\hat f} &\hat{\ex B}\dar&\lar T\hat{\ex B}\dar
\\ \ex F(\hat f)\rar& \ex B_{0}&\lar T\ex B_{0}\end{tikzcd}\]

 Sections of $T_{vert}\ex C(\hat f)$ represent the infinitesimal reparametrizations of $\ex C(\hat f)$ fixing the parametrization of $\ex F(\hat f)$.  The action of reparametrization is again given by the map $h$:
  \[h(v):=(L_{v}j,d\hat f(v))\]

With the above definition, $T_{\hat f}\dmsw$ is a sheaf of $\C\infty1(\ex F(\hat f))$--modules,  because  given any $\lambda\in \C\infty1(\ex F(\hat f))$, $h(\lambda v+w)=\lambda h(v)+h(w)$. Note that if $f$ is a curve in $\hat f$, we can restrict any section of $T_{\hat f}\dmsw$ to give a tangent vector in $T_{f}\dmsw$, and more generally, given any map $\hat g\longrightarrow \hat f$, there is an induced pullback map $T_{\hat f}\dmsw\longrightarrow T_{\hat g}\dmsw$ making the assignment $\hat f\mapsto T_{\hat f}\dmsw$ into a contravariant functor and a sheaf over $\dmsw$; see Definition \ref{sheaf}. Accordingly,  the sheaf $T_{\hat f}\dmsw$ should be regarded as the pullback of the tangent sheaf of $\dmsw$. When the family $\hat f$ has bubbling or node formation behavior, it is not clear that $T_{\hat f}\dmsw$ represents sections of a vector-bundle over $\ex F(\hat f)$, however we will use a notion of a finite-dimensional sub-vector-bundle. In particular, when subsheaf of $T_{\hat f}\dmsw$  is locally free and finitely generated over $\C\infty1(\ex F)$, we regard  it  a sheaf of sections of a finite-dimensional sub-vector-bundle of $T_{\hat f}\dmsw$.

\

Consider a $1$-dimensional deformation $\hat f_{t}$ of $\hat f$; that is, a family $\hat h$ such that  $\ex F(\hat h)=\mathbb R\times \ex F(\hat f)$, and the pullback of $\hat h$ to ${t}\times \mathbb R$ is $\hat f_t$, with $\hat f_0=\hat f$. We can take the derivative of this deformation with respect to $t$ to obtain a section of $T_{\hat f}\dmsw$ as follows.
Choose a lift, $\tilde v$, of  the vector field $\frac\partial{\partial t}$ on $\ex F(\hat h):=\mathbb R\times \ex F(\hat f)$ to a section of  $T\ex C(\hat h)$. Then $(L_{\tilde v}j,d\hat f_{t}(\tilde v))$ restricted to $\ex C(\hat f)\subset \ex C(\hat h)$ is contained in $\Gamma^{0,1}(T_{vert}\ex C(\hat f)\otimes T_{vert}^*\ex C(\hat f))\times \Gamma_{\ex B_0}( \hat f^{*} T \hat{\ex B})$. The choice of lift, $\tilde v$, is determined up to a choice of vertical vector-field , so $(L_{\tilde v}j,d\hat f_{t}(\tilde v))$ projects to a section of $T_{\hat f}\dmsw$ that is well defined independent of our choice of lift of $\frac\partial{\partial t}$. The following lemma shows that all sections of $T_{\hat f}\dmsw$ can be constructed in this way.

\begin{lemma}\label{tangent3} Given any section $v$ of $T_{\hat f}\dmsw$, there is a one-dimensional deformation $\hat f_{t}$ of $\hat f$  with derivative at $t=0$ equal to $v$.
\end{lemma}

\pf

The lifting property of Definition \ref{decorated} implies that, without loss of generality, we can prove this lemma for $\Msw$.

Claim \ref{Lvdbar} tells us that, in holomorphic coordinates,  $L_{w}j$ is $2j$ times the standard $\dbar$ operator applied to $w$. Around any curve $f$ in $\hat f$, Theorem \ref{f replacement}  implies that we can construct a vertical vector-field  $w$ so that $L_{w}j$ equals a given section of $\Gamma^{0,1}(T_{vert}\ex C(\hat f)\otimes T_{vert}^*\ex C(\hat f))$ on a neighborhood in $\ex C(\hat f)$ of all the edges of $\ex C(f)$. We can patch together such vector-fields for different $f$, using a fiberwise constant partition of unity, to obtain  a globally defined section $w$ of $T_{vert}\ex C(\hat f)$ so that $L_{w}j$ is equal to the given section of $\Gamma^{0,1}(T_{vert}\ex C(\hat f)\otimes T_{vert}^*\ex C(\hat f))$ on a neighborhood of the edges of all the curves in $\ex C(\hat f)$.
Therefore, we can reduce to the case that our section $v$ of $T_{\hat f}\Msw$ is equal to $(\theta,r)$ where $\theta$ is a section of $\Gamma^{0,1}(T_{vert}\ex C(\hat f)\otimes T_{vert}^*\ex C(\hat f))$ vanishing on a neighborhood of all edges of curves in $\ex C(\hat f)$, and $r$ is a $\C\infty1$ section of $\hat f^{*}T_{vert}\hat{\ex B}$. 

Extend $j$ to a family $j_{t}$ of almost-complex structures on $\ex C(\hat f)$ with derivative at $0$ equal to $\theta$, and extend $\hat f$ to a  family of maps $\hat f_{t}$ with derivative at $0$ equal to $r$. The family of maps $f_{t}$ with domain $(\ex C(\hat f), j_{t})$ is the required deformation of $\hat f$ with derivative at $0$ equal to $v$.

\stop
 
\

\subsection{$D\dbar$}\label{Ddbar section}

\

\

In this section, we construct the linearization,  $D\dbar$, of the $\dbar$ operator at holomorphic curves $f$. We restrict to holomorphic curves because, to define a linearization of the $\dbar$ operator at a non-holomorphic curve $f$, a connection is required on the sheaf $\Y$ over the ambient moduli stack $\dmsw$ (recall that $\Y$ is the pullback of the sheaf from Definition \ref{Ydef} under the map $\dmsw\longrightarrow \Ms$). It is not clear to me that such a connection always exists globally.

 \
 
 For a $\C\infty1$ family of curves $\hat f$ containing a homomorphic curve $f$, let us define a linear map
  \[D\dbar\co T_{f}\ex F(\hat f)\longrightarrow \Y(f)\ .\]
 First, regard $\dbar\hat f$ as a section of the vector-bundle $Y(\hat f)$ over $\ex C(\hat f)$, and choose a connection $\nabla$ on $Y(\hat f)$. Let $v$ be a vector in $T_f\ex F(\hat f)$. To define $D\dbar(v)$, lift $v$ to a section $v'$ of $T\ex C(\hat f)$ over $\ex C(f)$. Then $D\dbar(v):=\nabla_{v'}\dbar \hat f$. Because $\dbar \hat f$ vanishes on $\ex C(f)$, $\nabla_{v'}\dbar \hat f$ does not depend on the choice of $\nabla$ or the lift $v'$ of $v$. Moreover,  because $\dbar\hat f$ vanishes on edges of $\ex C(f)$, $\nabla_{v'}\dbar\hat f$ vanishes on all edges of curves in $\ex C(\hat f)$,  so it is a section in $\Y(f)$. 

 Our $D\dbar$ is natural in the following sense: given any commutative diagram,
 \[\begin{tikzcd}f\rar \ar{dr}&\hat f\dar
 \\ & \hat h\end{tikzcd}\]
 the following diagram commutes.
\[ \begin{tikzcd}\dar T_{f}\ex F(\hat f)\rar{D\dbar}&\Y(f)
\\  T_{f}\ex F(\hat h)\ar{ur}{D\dbar}\end{tikzcd}\]
Therefore, 
Lemmas \ref{tangent1} and \ref{tangent2} imply that the maps $D\dbar$  all factor through a fixed  linear map from $T_{f}\dmsw$. 
 \[T_{f}\ex F(\hat f)\longrightarrow T_{f}\dmsw\xrightarrow{D\dbar}\Y( f)\]
An elementary extension of the calculation of the linearized $\dbar$ operator in \cite[Lemma 1.7]{reg} gives  a formula for $D\dbar$ in terms of $\Gamma^{0,1}(T_{vert}\ex C( f)\otimes T_{vert}^*\ex C( f))\times \Gamma_{\ex B_0}(  f^{*} T \hat{\ex B})$.

\

Recall that both $\Y(f)$ and $T_{f}\dmsw\ov{\ex B_{0}}\subset T_{f}\dmsw$ are complex at holomorphic curves $f$. The restriction of $D\dbar$ to $T_{f}\dmsw\ov{\ex B_{0}}$ is not usually $\mathbb C$--linear  unless $J$ is integrable. We shall see that $D\dbar$ has a finite-dimensional kernel and cokernel throughout the linear homotopy of $D\dbar$ to its complex-linear part. This allows us to canonically orient the kernel of $D\dbar$ relative to the cokernel, and  construct orientations and a kind of almost-complex structure on embedded Kuranishi structures; see Section \ref{rcs section}.

\

\subsection{Embedded Kuranishi structures}\label{K section}

\

If $D\dbar$ is surjective at a holomorphic curve $f$, then we shall show that the moduli stack of holomorphic curves close to $f$ is represented by the quotient of some  $\C\infty1$ family of curves $\hat f$ by a finite group of automorphisms; see Definition \ref{represented by family}. If, however, $D\dbar$ is not surjective at $f$, we need to choose a nice subsheaf $V$ of $\Y$ in a neighborhood of $f$  so that the moduli stack of curves with $\dbar$ in $V$ is well behaved.  Below is a standard definition of a subsheaf $V$ of $\Y$. 

\begin{defn} Let $\mathcal U$ be a substack of $\dmsw$. A  subsheaf $V$ of $\Y$ on $\mathcal U$ associates, to each family $\hat f$ in $\mathcal U$, a subsheaf $V(\hat f)\subset \Y(\hat f)$ such that given any map $\psi\co \hat f\longrightarrow \hat g$,  $V(\hat f)$ is the sheaf of $\C\infty1(\ex F(\hat f))$--modules  generated by the image of  $V(\hat g)$ under  $\Y(\psi)\co \Y(\hat g)\longrightarrow \Y(\hat f)$.
\end{defn}

We want our subsheaves $V$ of $\Y$ to be pulled back from nice geometrically defined sheaves. For example, around a curve $f$ that  is embedded and has a stable domain, we might pull back $V$ from a sheaf defined over $\bar M_{g,n}\times \hat{\ex B}$. The following defines what we mean by pullback.

\begin{defn} Given a family $\hat{\ex A}\longrightarrow \ex X$ of exploded manifolds with a fiberwise almost-complex structure, define
\[\Gamma^{0,1}(T_{vert}^{*}\hat {\ex A}\otimes T_{vert}\hat{\ex B})\]
to be the sheaf of $\C\infty1(\ex X)$--modules on $\ex X$ with global sections consisting of anti-holomorphic, $\C\infty1$ sections of  $T_{vert}^{*}\hat {\ex A}\otimes T_{vert}\hat{\ex B}$ over $\hat{\ex A}\times \hat{\ex B}$  vanishing on integral vectors within $T_{vert}\hat{\ex A}$, (and with sections over $U\subset \ex X$ consisting of the same thing with $\hat {\ex A}$ replaced by the inverse image of $U$ in $\hat{\ex A}$).

\end{defn}
Given a family of curves $\hat f$ in $\hat{\ex B}$ and a fiberwise-holomorphic map 
\[\begin{tikzcd}\ex C(\hat f)\rar{\psi} \dar &\hat {\ex A}\dar
\\\ex F(\hat f)\rar&\ex X\end{tikzcd}\]
sections $\nu$ of $T^{*}_{vert}\hat{\ex A}$  pull back to sections $\psi^*\nu$ of $T^{*}_{vert}\ex C(\hat f)$ by composing with the derivative of $\psi$. As $\psi$ is fiberwise holomorphic, the pullback map, $\psi^{*}$,  is complex linear.  Any section  vanishing on integral vectors of $T_{vert}\hat{\ex A}$ pulls back to a section  also vanishing on integral vectors of $T_{vert}\ex C(\hat f)$.  

 A section $v$ of $T_{vert}\hat{\ex B}$ also pulls back to a section $\hat f^*v$ of $\hat f^{*}T_{vert}\hat{\ex B}$. This  map $\hat f^*$ is also complex, therefore there is an induced complex map: 
\[\psi^*\otimes \hat f^*\co \Gamma^{0,1}(T_{vert}^{*}\hat {\ex A}\otimes T_{vert}\hat{\ex B})\longrightarrow \Y(\hat f):=\Gamma^{0,1}(T_{vert}^{*}\ex C(\hat f)\otimes f^{*}T_{vert}\hat{\ex B}) \]

\begin{defn} \label{Vpullback}
We say that the {\bf pullback} of a section or subsheaf of $\Gamma^{0,1}(T_{vert}^{*}\hat {\ex A}\otimes T_{vert}\hat{\ex B})$ to  $\Y(\hat f)$ is the image of the section or subsheaf under the map $\psi^*\otimes \hat f^*$.
\end{defn}

The following defines what we mean by a `nice' subsheaf $V$ of $\Y$.

\begin{defn}\label{simply generated}
 
 We say that a subsheaf $V$ of $\Y$ on $\mathcal U$ is {\bf simply generated} if there exists a fiberwise-holomorphic map
\[\begin{tikzcd}\mathcal U^{+1}\rar{\fun\Phi^{+1}} \dar &\hat {\ex A}/G\dar
\\ \mathcal U\rar{\fun\Phi} &{\ex X}/G\end{tikzcd}\]
 with an effective $G$--action in the sense of Definition \ref{fhm},
 and  sections $v_{1},\dotsc,v_{n}$ of $\Gamma^{0,1}(T_{vert}^{*}\hat {\ex A}\otimes T_{vert}\hat{\ex B})$ such that the following holds:\begin{enumerate}\item For any family of curves $\hat f$ in $\mathcal U$,  the pulled back sections  \[((\fun\Phi(\hat f)^{+1})^*\otimes(\hat f\bd G)^*)v_{i} \] of $\Y(\hat f\bd G)$ are linearly independent, and generate $V(\hat f\bd G)$ as a sheaf of $\C\infty1(\ex F(\hat f\bd G))$--modules. (As in Definition \ref{fhm},  $\hat f\bd G$ indicates  the $G$--fold cover of $\hat f$ with a $G$--equivariant fiberwise-holomorphic map $\fun\Phi^{+1}(\hat f)\co \ex C(\hat f\bd G)\longrightarrow \ex A$.)
 
\item The sheaf of $\C\infty1(\ex X)$--modules over $\ex X$ generated by $v_{1},\dotsc,v_{n}$ is $G$--invariant.
\end{enumerate}
\end{defn}

Theorem \ref{V moduli stack}  states that if $V$ is a simply-generated subsheaf of $\Y$ transverse to $D\dbar$ at a holomorphic curve $f$, then there exists an open neighborhood $\mathcal O\subset \dmsw$ of $f$ and a $\C\infty1$ family of curves $\hat f$ with automorphism group $G$ such that $\hat f/G$ represents the substack, $\dbar^{-1}(V)\subset \mathcal O$, consisting  of curves $h$  with the property that $\dbar h\in V(h)$; see definitions \ref{family quotient} and \ref{represented by family}. The substack of holomorphic curves within  $\hat f/G$  therefore represents the  moduli stack of holomorphic  curves in $\mathcal O$. So,  the moduli stack of holomorphic curves in $\mathcal O$ is represented by $\{\dbar \hat f^{-1}(0)\}/G\subset \hat f/G$.

\begin{defn}Say that a subsheaf $V$ of $\Y$ is {\bf complex} if $V(f)$ is a complex-linear subspace of $\Y(f)$
 for all curves  $ f$.  \end{defn}

Recall that if $f$ is holomorphic,  $T_{f}\dmsw\ov{\ex B_0}$ is a complex vector-space, and there is a well-defined linearization of the $\dbar$ operator
\[D\dbar\co T_{f}\dmsw\ov{\ex B_0}\longrightarrow \Y(f)\]
constructed in Section \ref{Ddbar section}. This map $D\dbar$ is not necessarily $\mathbb C$--linear, however we shall show that there is a homotopy  \[(1-t)D\dbar+tD\dbar^{\mathbb C}\] from $D\dbar$ to its complex-linear part, transverse to a finite-dimensional $\mathbb C$--linear subspace $V$ of $\Y(f)$ for all $t\in[0,1]$, so 
\[K_t(f):=\lrb{(1-t)D\dbar+tD\dbar^{\mathbb C}}^{-1}(V)\]
is a family of finite-dimensional vector-spaces in $T_{f}\dmsw\ov{\ex B_0}$. 
  
\begin{defn}\label{strongly transverse}Say that a subsheaf $V$ of $\Y$ is {\bf strongly transverse} to $D\dbar$ at a holomorphic curve $f$ if 
\[(1-t)D\dbar+tD\dbar^{\mathbb C}\co T_{f}\dmsw\ov{\ex B_0}\longrightarrow  \Y(f)\]
is transverse to $V$ for all $t\in[0,1]$.
\end{defn}

Combining notions from definitions \ref{open substack}, \ref{decorated}, \ref{simply generated}, \ref{Ydef},   \ref{strongly transverse} and \ref{represented by family}, we make the following important definition of Kuranishi charts. 
\begin{defn}\label{K chart}A  {\bf Kuranishi  chart} $(\mathcal U,V,\hat f/G)$ on $\dmsw(\hat{\ex B})$ is 
\begin{itemize} \item an open substack \[\mathcal U\subset \dmsw(\hat{\ex B})\] 
\item a simply-generated complex  subsheaf $V$ of $\Y$, defined on $\mathcal U$,  and strongly transverse to $\dbar$ at all holomorphic curves in $\mathcal U$,
\item \label{k5} and a $\C\infty1$ family of curves $\hat f$ in $\mathcal U$ with automorphism group $G$ such that 
\begin{enumerate}\item $\hat f/G$ represents the  substack $\dbar^{-1}V\subset \mathcal U$, 
\item and so that the map  $\ex F(\hat f)\longrightarrow \ex B_0$ is a submersion.
\end{enumerate}
\end{itemize}

\end{defn}

 The exploded manifold $\ex F(\hat f)$ with the vector-bundle $V(\hat f)$, the section of this vector-bundle given by $\dbar$,  and the action of $G$ is our version of a Kuranishi chart defined by Fukaya and Ono in \cite{FO}. Apart from being embedded in an ambient moduli stack, the main difference is that we use exploded manifold charts, so when defining a virtual fundamental class we use the transversality and intersection theory for exploded manifolds; see \cite[Section 9]{iec}.
 
 Our final condition, that $\ex F(\hat f)\longrightarrow \ex B_0$ is a submersion, is included to ensure that any base change of our family of targets 
 \[\begin{tikzcd}\hat{\ex B}'\rar\dar&\hat {\ex B}\dar
\\ \ex B_0'\rar& \ex B_0\end{tikzcd}\]
pulls back  Kuranishi charts on $\Msw(\hat{\ex B})$ to Kuranishi charts on $\Msw(\hat{\ex B}')$; similarly, Kuranishi charts on $\dmsw(\hat{\ex B})$ pull back to Kuranishi charts on $\dmsw(\hat{\ex B}')$, but see Remark \ref{decoration base change}.

We include the condition of strong transversality to $V$ and require $V$ be complex  to ensure that the information of the stable almost-complex structure, defined by Fukaya and Ono in \cite{FOinteger}, is reflected in our  Kuranishi chart. This information can be used to orient a Kuranishi chart (as in  \cite[Section 3.1]{vfc}), and  to define integer counts of holomorphic curves using Fukaya and Ono's method from \cite{FOinteger}; see \cite{icc} for details, and   \cite{kuranishihomology,AMS2021, Bai2022} for related work.  Such structures and invariants are compatible with evaluation maps that are holomorphic submersions, defined below.

\begin{defn}\label{holomorphic submersion} A {\bf submersion} $\fun\Phi$  to a family $\ex X\longrightarrow \ex X_0$ of exploded manifolds or orbifolds is a commutative diagram of $\C\infty1$ maps
\[\begin{tikzcd}\dmsw(\hat {\ex B})\dar \rar{\fun\Phi}&\ex X\dar
\\ \ex B_0 \rar & \ex X_{0} \end{tikzcd}\]
such that $T_{f}\fun\Phi\co T_{f}\dmsw\longrightarrow T_{\fun\Phi(f)}\ex X$ is surjective for all $f$, where we use Definition \ref{tangent as set} and Example \ref{orbifold tangent space} in the orbifold case.

Say that a submersion $\fun\Phi\co \dmsw\longrightarrow \ex X$ is {\bf holomorphic} if $\ex X\longrightarrow \ex X_{0}$ has a fiberwise almost-complex structure, and for each holomorphic curve $f$, the map
\[T_{f}\fun\Phi\co T_{f}\dmsw\ov{\ex B_0}\longrightarrow T_{\fun\Phi(f)}\ex X\ov{\ex X_{0}}\]
is complex.
\end{defn}

 Examples of such holomorphic submersions include the maps $\fun{ev}^{+n}$ defined in Section \ref{+n section}, and  the usual evaluation map from Gromov--Witten theory, evaluating curves at marked points.

\begin{defn}\label{submersive def} Given a submersion \[\fun\Phi\co \dmsw\longrightarrow \ex X\] where $\ex X$ is a finite dimensional exploded manifold or orbifold, say that a Kuranishi chart $(\mathcal U,V,\hat f/G)$ on $\dmsw$ is $\fun\Phi$--submersive if the induced map 
\[\fun\Phi\co \ex F(\hat f)\longrightarrow \ex X\]
is a submersion, and if for all holomorphic curves $f$ in $\hat f$, $D\dbar$ restricted to $\ker T_{f}\fun\Phi$ is strongly transverse to $V(f)$ in the sense that  
\[((1-t)D\dbar+tD\dbar^{\mathbb C})(T_{f}\dmsw\ov{\ex B_0} \cap\ker T_{f}\fun\Phi)\]
is transverse to $V(f)$ for all $t$ in $[0,1]$.
\end{defn}

\begin{defn}\label{K compatible}Two Kuranishi charts $(\mathcal U_{1},V_{1},\hat f_{1}/G_{1})$,  $(\mathcal U_{2},V_{2},\hat f_{2}/G_{2})$ are 
{\bf compatible} if restricted to $\mathcal U_{1}\cap \mathcal U_{2}$, either $V_{1}$ is a subsheaf of $V_{2}$ or $V_{2}$ is a subsheaf of $V_{1}$.
\end{defn}

In the case that $V_{1}$ is a subsheaf  of $V_{2}$ , $\hat f_1\rvert_{\mathcal U_1\cap\mathcal U_2}$ is a family in the stack represented by $\hat f_2/G_2$ (see definitions \ref{family quotient} and \ref{represented by family}) so there is a unique (up to unique $2$-isomorphism) $\C\infty1$ transition map 
\[\hat f_{1}\rvert_{\mathcal U_{1}\cap \mathcal U_{2}}\longrightarrow \hat f_{2}/G_{2}\ .\]
 in the sense of Example \ref{family map}. In particular, this transition map can be expressed in the form 
 \[\begin{tikzcd}\hat f_{1}\times_{\dmsw}\hat f_2 \rar{\psi}\dar{\pi} &\hat f_2
 \\ \hat f_1\rvert_{\mathcal U_1\cap\mathcal U_2} \end{tikzcd}\]
where $\pi$ is a $G_1$--equivarant  $G_2$--bundle internal to $\dmsw$, and $\psi$ is a $G_1$--invariant and  $G_2$--equivariant  map.

\begin{defn}\label{K extend def}A Kuranishi chart $(\mathcal U,V,\hat f/G)$ is 
{\bf extendible} if it has an extension, which is a Kuranishi chart $(\mathcal U\x,V\x,\hat f\x/G)$ satisfying the following conditions.
\begin{itemize}\item There exists a continuous map 
\[\rho\co \mathcal U\x\longrightarrow (0,1]\]
so that \begin{itemize}
\item \[\mathcal U=\rho^{-1}((1/2,1])\ ;\]
\item for any $t>0$, all holomorphic curves in the closure of $ \rho^{-1}((t,1])$ within $\dmsw$ are contained within $\rho^{-1}([t,1])$;
\item for any $t>0$, the closure within $\dmsw$ of  the subset of $\hat f\x$ where $\rho>t$ is contained within $\hat f\x$; \footnote{In particular,  the closure of $\hat f$ within $\hat f\x$ is closed within $\dmsw$. This means wild behavior at  the boundary of $\hat f$ can be excluded by requiring  extensions to $\hat f\x$. }
\end{itemize}
\item $V$ is the restriction of $V\x $ to $\mathcal U$, and $\hat f$ is the restriction of $\hat f\x$ to $\mathcal U$.
\end{itemize}
\end{defn}

Note that extendible Kuranishi charts pull back to extendible Kuranishi charts. 
The condition of extendability is designed to prevent pathological behavior from occurring at the boundary of Kuranishi charts. We shall have cause\footnote{See for example  \cite[Proposition 2.3]{vfc}.} to repeatedly shrink the size of extensions during inductive constructions.  Restricting $(\mathcal U\x, V\x, \hat f\x/G)$  to  $\rho>t$ for any $t\in (0, 1/2)$ gives an extension of $(\mathcal U,V,\hat f/G)$.

As we shall have no reason to distinguish between $V$ and $V\x$, we sometimes use the notation $V$ to refer to $V\x$.

\begin{defn}\label{K compatible collection} A collection of extendible Kuranishi charts is {\bf locally finite} if there exists an extension $(\mathcal U\x ,V\x ,\hat f\x /G)$ of each Kuranishi chart $(\mathcal U,V,\hat f/G)$ such that each holomorphic curve has a neighborhood intersecting only finitely many of the  $\mathcal U\x $,  and each  $\mathcal U\x$ intersects only finitely many of the other $\mathcal U\x$.

 The collection is {\bf compatible}  if every pair of extended Kuranishi charts is compatible.

 The collection is said to {\bf cover} any substack  of  $\dmsw$  covered by  $\{\mathcal U\}$.
\end{defn}

Now, combining definitions \ref{K chart}, \ref{K compatible}, \ref{K extend def}, and \ref{K compatible collection}, we arrive at the definition of an embedded Kuranishi structure.

\begin{defn}\label{K def}
An {\bf embedded Kuranishi structure} on the moduli stack of stable holomorphic curves $\dmod\subset\dmsw$ is a countable, locally finite,  compatible collection of extendible Kuranishi charts $\{(\mathcal U_{i},V_{i},\hat f_{i}/G_{i})\}$ covering the holomorphic curves $\dmod\subset\dmsw$.

\end{defn}

\

In many cases, the moduli stack of stable holomorphic curves,   $\Mod$, has the following compactness property: The map $\Mod\longrightarrow\ex B_{0}$ is proper when $\Mod\subset\Msw$ is restricted to any connected component of $\Msw$; see \cite{cem} for details. Definition \ref{decorated} implies that in this case, the same compactness holds for the moduli stack of decorated stable holomorphic curves $\dmod\subset\dmsw$. When $\dmod$ satisfies this compactness property, Theorem \ref{K existence}  states that there exists an embedded Kuranishi structure on $\dmod$, and that, given any submersion $\fun\Phi\co \dmsw\longrightarrow \ex X$,  there exists a $\fun\Phi$--submersive embedded Kuranishi structure on $\dmod$.

Given any base change of our family of targets,
 \[\begin{tikzcd}\hat{\ex B}'\rar \dar&\hat {\ex B}\dar
 \\ \ex B_0'\rar& \ex B_0\end{tikzcd}\]
the pullback of any embedded Kuranishi structure on $\Mod(\hat {\ex B})\subset \Msw(\hat{\ex B})$ is an embedded Kuranishi structure on $\Mod(\hat{\ex B}')$.  Corollary \ref{K homotopy}  proves that any two embedded Kuranishi structures on $\dmod(\hat{\ex B} )$ are homotopic. That is,  there exists an embedded Kuranishi structure on $\dmod(\hat {\ex B}\times \mathbb R)$  pulling back to give each of the original embedded Kuranishi structures under the inclusions of $\hat {\ex B}$ over $0$ and $1$ in $\mathbb R$.

\

%
%
 
 Given an embedded Kuranishi structure $\{(\mathcal U_{i},V_{i},\hat f_{i}/G_{i})\}$,  we shall construct, in Section \ref{rcs section}, a canonical homotopy class of complex structure on $T_{f}\ex F(\hat f _{i})\ov{\ex B_{0}}$ for all holomorphic curves $f$ in $\hat f_{i}$ so that 
 
 \begin{itemize}\item all $\mathbb R$--nil vectors at holomorphic curves  are given the canonical complex structure;
 \item the action of $G_i$  preserves the complex structure:  if an element of $G_i$ sends  $f$ to $f' $, the corresponding map $T_{f}\ex F(\hat f _i)\ov{\ex B_0}\longrightarrow T_{f' }\ex F(\hat f  _i)\ov{\ex B_0}$ is complex;
\item if $f$ is a  holomorphic curve in $\hat f  _i$ and $\hat f  _j$, and  $V_i(f)\subset V_j(f)$, then the following is a short exact sequence of complex linear maps
\[0\longrightarrow T_f\ex F(\hat f  _i)\ov{\ex B_0}\longrightarrow T_f\ex F(\hat f  _j)\xrightarrow{D\pi_{V_i}\dbar}  V_j/V_i(f)\longrightarrow 0\] 
\item there exists a complex structure on $T\ex F(\hat f  _i)\ov{\ex B_0}$, defined on a neighborhood of the holomorphic curves, and  restricting to the given complex structure at holomorphic curves. 
\end{itemize}
 We shall also show in Section \ref{rcs section} that given a holomorphic submersion (from Definition \ref{holomorphic submersion}), 
 \[\begin{tikzcd}\dmsw\dar \rar{\fun\Phi}& \ex X\dar
 \\ \ex B_0\rar&\ex X_0\end{tikzcd}\]
we can construct this complex structure so that 
 \[T_f\fun\Phi\co T_f\ex F(\hat f_i)\ov{\ex B_0}\longrightarrow T_{\fun\Phi(f)}\ex X\ov{\ex X_0}\]
 is complex for all $f$ in $\hat f_i$.

\section{Trivializations and pre-obstruction models}

\label{tpo section}

This section summarizes the results of \cite{reg}  necessary for this paper. The notion, from  \cite[Definition 1.4]{reg}, of a trivialization $(\mathcal F,\Phi)$ associated to a family of curves
\[\begin{tikzcd}\ex C(\hat f)\dar \rar{\hat f} & \hat{\ex B}\dar
\\ \ex F(\hat f)\rar & \ex B_0 \end{tikzcd}\]
 allows us to identify sections of $\hat f^{*}T_{vert}\hat{\ex B}$ with families of curves parametrized by $\ex C(\hat f)$, and to identify $\dbar$ of such a family of curves with a section of $\Y(\hat f)$.

\begin{defn}\label{trivialization def} Given a $\C\infty1$ family  of curves $\hat f$, 
 a choice of trivialization $(\mathcal F,\Phi)$ is the following information:
 \begin{enumerate}
 \item
 A $\C\infty1$ map $\mathcal F$ such that the following diagram of $\C\infty1$ maps commutes
 \[\begin{tikzcd} \hat f^{*}{T_{vert}\hat{\ex B}}\rar{\mathcal F}\dar &\hat{\ex B} \dar
\\ \ex  F(\hat f)\rar &\ex B_{0}
 \end{tikzcd}\]
 and such that \begin{enumerate}\item\label{triv1} $\mathcal F$ restricted to the zero section equals  $\hat f$, 
 \item $T\mathcal F$ restricted to the canonical inclusion  $\hat f^{*}{T_{vert}\hat{\ex B}}\subset T\hat f^{*}{T_{vert}\hat{\ex B}}$ over the zero section is the identity,
 \item\label{injective derivative}  $T\mathcal F$ restricted to the vertical tangent space at any point of $\hat f^{*}T_{vert}\hat{\ex B}$ is injective.
 \end{enumerate} 
 \item A $\C\infty1$ vector-bundle map, $\Phi$ 
 \[\begin{tikzcd}\mathcal F^{*}{T_{vert}\hat{\ex B}}\rar{\Phi}\dar &\hat f^{*}{T_{vert}\hat{\ex B}}\dar{\pi}
 \\ \hat f^{*}{T_{vert}\hat{\ex B}}\rar{\pi} &\ex C(\hat f)\end{tikzcd}\]
 which is a $\mathbb C$--linear isomorphism on each fiber, and which is the identity when the vector-bundle  $\mathcal F^{*}{T_{vert}\hat{\ex B}}\longrightarrow \hat f^{*}{T_{vert}\hat{\ex B}}$ is restricted to the zero section of $\hat f^{*}T_{vert}\hat{\ex B}$.
\end{enumerate}
 
 A trivialization allows us to define $\dbar$ of a section \[\nu \co \ex C(\hat f)\longrightarrow\hat f^{*}{T_{vert}\hat{\ex B}}\] as follows: $\mathcal F\circ \nu$ is a map $\ex C(\hat f)\longrightarrow \hat {\ex B}$, so $\dbar (\mathcal F\circ \nu)$ is a section of 
 \[\Y(\mathcal F\circ\nu)=\Gamma^{0,1}\lrb{{T_{vert}^{*}\ex C(\hat f) }\otimes {(\mathcal F\circ \nu)^{*}T_{vert}\hat{\ex B}}}\ .\] Applying the map $\Phi$ to the second component of this tensor product  identifies  $\Y(\mathcal F\circ \nu)$ with $\Y(\hat f)$, and identifies $\dbar(\mathcal F\circ \nu)$ with a section of $\Y(\hat f)$. Define $\dbar \nu$ to be this section of $\Y(\hat f)$.

\end{defn}

For example, we can construct a trivialization by extending $\hat f$ to a map $\mathcal F$ satisfying the above conditions --- for instance by choosing a smooth connection on $T_{vert}\hat {\ex B}$ and reparametrising the exponential map on a neighborhood of the zero section in  $f^{*}T_{vert}\hat {\ex B}$ --- and defining $\Phi$ using parallel transport along a linear path to the zero section using a smooth $J$--preserving connection on $T_{vert}\hat{\ex B}$. 

Given a choice of  trivialization for the family $\hat f$ and a $\C\infty1$ section $\nu$ of $\hat f^{*}T_{vert}\hat{\ex B}$, there is an induced choice of trivialization for the family $\mathcal F(\nu)$.

\begin{defn}\label{pre-obstruction model}
A  $\C\infty1$  pre-obstruction model $(\hat f,V,\mathcal F,\Phi,\{s_{i}\})$, is given by
\begin{enumerate}
\item
 a $\C\infty1$ family of curves $\hat f$;
\item   a subsheaf $V$ of $\Y(\hat f)$ that  is   locally free, and finitely generated as a sheaf of $\C\infty1(\ex F((\hat f))$--modules;

\item a choice of trivialization $(\mathcal F,\Phi)$ for $\hat f$; and
 
 \item a finite collection $\{s_{i}\}$ of extra marked points on $\ex C(\hat f)$ corresponding to $\C\infty1$ sections \[s_{i}\co \ex F(\hat f)\longrightarrow \ex C(\hat f)\] so that, restricted to any curve $\ex C$ in $\ex C(\hat f)$, these marked points are distinct and contained inside the smooth components\footnote{Note that the extra marked points on $\ex C$ are not the same as ends of $\ex C$, so we do not modify the exploded manifold structure of $\ex C$ at these extra marked points.  See  \cite[Definition 8.3]{iec} for smooth component of $\ex C$.} of $\ex C$.
 
 \end{enumerate}
\end{defn}

 We shall usually use the notation $(\hat f,V)$ for a pre-obstruction bundle. We will refer to $V\subset\Y(\hat f)$ as a finite-rank sub-vector-bundle, even though $\Y(\hat f)$ might not really be a vector bundle over $\ex F(\hat f)$.

\begin{defn}Given any family of curves, $\hat f$, and a collection, $\{s_{i}\}$, of extra marked points on $\ex C(\hat f)$,  let  $X^{\infty,\underline1}(\hat f)$ indicate the space of $\C\infty1$ sections of $ f^{*}T_{vert}\hat{\ex B}$  vanishing on the extra marked points $\{s_{i}\}$ on ${\ex C}(\hat f)$.
\end{defn}

Note that both  $X^{\infty,\underline 1}(\hat f)$ and $\Y(\hat f)$ are complex vector-spaces because they consist of sections of complex vector-bundles.

We can restrict any pre-obstruction bundle $(\hat f,V)$ to a single curve $f$ in $\hat f$. The restriction of $V$ to this curve $f$ is a finite-dimensional linear subspace  $V(f)\subset \Y(f)$. 

 Let
 $D\dbar(f)\co X^{\infty,\underline 1}(f)\longrightarrow \Y(f)$ indicate the  derivative of $\dbar$ at $0\in X^{\infty,\underline 1}(f)$.   We are most interested in pre-obstruction bundles $(\hat f,V)$ containing  holomorphic curves $f$ such that $D\dbar(f)$ is injective and has image complementary to $V(f)$. If there are enough extra marked points to stabilise $\ex C(f)$,  then  $X^{\infty,\underline 1}(f)$ is a linear subspace of $T_{f}\Msw\ov{\ex B_{0}}\subset T_{f}\Msw$. Moreover,  if $f$ is also holomorphic, $X^{\infty,\underline 1}$ is a $\mathbb C$--linear subspace of $T_{f}\Msw\ov{\ex B_{0}}$, and $D\dbar(f)$ is the restriction of $D\dbar\co T_{f}\Msw\longrightarrow \Y(f)$ to this subspace.

To describe the importance of pre-obstruction models, we need the notion of a simple perturbation below. Recall, from Definition \ref{Ydef}, that the sheaf of $\C\infty 1(\ex F(\hat f))$--modules $\mathcal Y(\hat f)$ describes the sections of a vector bundle $Y(\hat f)$ over $\ex C(\hat f)$.

\begin{defn}\label{simple perturbation}
Given a trivialization for  $\hat f$, a {\bf simple perturbation} of $\dbar$ is a map \[\dbar'\co X^{\infty,\underline 1}(\hat f)\longrightarrow \Y(\hat f)\]
such that

\[\dbar'\nu=\dbar(\nu)+\Psi(\nu)\]

where $\Psi$ is a (usually nonlinear)  $\C\infty1$ map 
\[\begin{tikzcd}[column sep=tiny]\hat f^{*}{T_{vert}\hat {\ex B}}\ar{rr}{\Psi}\ar{dr}& &Y(\hat f)\ar{dl}
\\ &\ex C(\hat f) \end{tikzcd}\]
that agrees with  zero-section when restricted to edges of curves in $\ex C(\hat f)$.

\

\end{defn}

\begin{example}[Construction of a simple perturbation]\label{simple construction}\end{example}
 Let  $\theta$ be a section of $\Gamma^{0,1}\lrb{{T_{vert}^{*}{\ex C}(\hat f) }\otimes {T_{vert}\hat{\ex B}}}$ and suppose that $\hat f$ comes with a trivialization $(\mathcal F,\Phi)$. A section $\nu$ of $\hat f^{*}T_{vert}\hat{\ex B}$ defines a map \[(\id,\mathcal F(\nu))\co \ex C(\hat f)\longrightarrow \ex C(\hat f)\times \hat {\ex B}\] Pulling back the section $\theta$ over $(\id, \mathcal F(\nu))$ gives a section of  $\Y(\mathcal F(\nu))$, which we can identify as a section of $\Y(\hat f)$ using the map $\Phi$ from our trivialization. Therefore, we get a modification $\dbar'$ of the usual $\dbar$ equation on sections $\hat f^{*}T_{vert}\hat{\ex B}$ given by the trivialization 
\[\dbar'\nu:=\dbar\nu-\Phi\lrb{(\id,\mathcal F(\nu))^{*}\theta}\ .\]
As $\mathcal F$, $\theta$ and $\Phi$ are $\C\infty1$ maps, we can define a $\C\infty1$ map 
\[\psi\co \hat f^*T_{vert}\hat {\ex B}\longrightarrow Y(\hat f)\] by $\psi(\nu):=-\Phi\lrb{(\id,\mathcal F(\nu))^{*}\theta}$ so $\dbar'$ is a simple perturbation of $\dbar$.

\

The following theorem is the main result of \cite{reg}. It describes the solutions to $\dbar'\subset V$ in the case that $D\dbar'(f)$ is injective, and has image complementary to  $V(f)$ in the sense that it is transverse to $V(f)$, and intersects $V(f)$ only at $0$. To achieve such a situation, we generally add extra marked points to our family to ensure that $D\dbar'$ is injective, then choose $V$ to be complementary to  the image of $D\dbar'$.

\begin{thm}\label{regularity theorem}
Suppose that $(\hat f,V)$ is a $\C\infty1$ pre-obstruction model containing the curve $f$, so that $\dbar' f\in V(f)$, and 
\[D\dbar'(f)\co X^{\infty,\underline 1}(f)\longrightarrow \Y(f)\]
is injective and has image complementary to $V(f)$.

Then there exists  a neighborhood $\hat f'$ of $f$ within $\hat f$ and a solution $\nu\in X^{\infty,\underline 1}(\hat f')$ to the equation 
\[\dbar'\nu=0\mod V\ .\] 
 
Moreover, there exists a neighborhood $O$ of $0\in X^{\infty,\underline 1}(\hat f')$ so that $\nu$ is the unique solution to $\dbar'\nu=0\mod V$ within $O$. Moreover,  given any curve $f'\in \hat f'$, let  $\nu(f')$ and $O(f')$ be the restriction of  $\nu$ and $O$   to $X^{\infty,\underline 1}(f')$; then $\nu(f')$ is the unique solution to the equation $\dbar' \nu(f')=0\mod V(f')$ within $O(f')$.

\end{thm}

We also need the following corollary of Theorem \ref{regularity theorem}.
\begin{cor}\label{trt} Let $(\hat f,V)$ be  a $\C\infty1$ pre-obstruction model containing $f$, and let $\dbar'$ be a simple perturbation of $\dbar$ such that $\dbar' \hat f$ is tangent to $V$ at $f$, and $D\dbar'(f)$ is complementary to $V$. Then $\nu$ is tangent to $0$ at $f$, where  $\nu$ is the unique section of of $X^{\infty,\underline 1}(\hat f)$ such that $\dbar'(\nu)$ is a section of $V$. 

\end{cor}

\pf
We can restrict $(\hat f,V)$ to a pre-obstruction model $(\hat f',V)$ where $\hat f'$ is parametrized by $\mathbb R$ and $f$ is the curve over $0$.  We can also restrict $\dbar'$ to $(\hat f',V)$.
The uniqueness part of Theorem \ref{regularity theorem} implies that $\nu$ pulls back to the unique section $\nu'$ with $\dbar'\nu'\subset V$. It therefore suffices to prove that $\nu'$ is tangent to the zero section at $0$. 

Claim \ref{Hadamard} implies that $\dbar' \hat f'$ is equal to a section of $V$ plus $t^{2}\theta$ where $\theta$ is a section of $\Y(\hat f')$, and $t$ is the coordinate on $\mathbb R$. Consider the family of curves $\hat f'\times \mathbb R$; that is the pullback of $\hat f'$ under the projection $\ex F(\hat f')\times\mathbb R\longrightarrow\ex F(\hat f')$. We can pull back our pre-obstruction model, $\dbar'$ and $\theta$ to $\hat f'\times \mathbb R$. Theorem \ref{regularity theorem} then implies that in a neighborhood of $(f,0)$, there is a unique $\C\infty1$ solution $\psi$ to the equation \[(\dbar'\psi-x\theta)\in V\] where $x$ is the coordinate on the extra  $\mathbb R$ factor of $\hat f'\times\mathbb R$. The uniqueness part of Theorem \ref{regularity theorem} implies that 
 \[\psi(t,t^{2})=0\] 
and
\[\psi(t,0)=\nu'(t)\]
therefore $\nu'$ is tangent to the zero section at $t=0$, as required.

\stop

\

The following theorem, Theorem 1.2 of  \cite{reg}, implies that for any simple perturbation $\dbar'$ of $\dbar$,  $D\dbar'(f)$  is a Fredholm operator, and that $D\dbar'$ over a family $\hat f$ can be locally approximated by a $\C\infty1$ map of finite-rank vector-bundles.
\begin{thm}\label{f replacement}Given any family of curves  $\hat f$ with a trivialization and a simple perturbation $\dbar'$ of $\dbar$,  the restriction of $D\dbar'$ to $f$,
\[D\dbar'(f)\co X^{\infty,\underline 1}(f)\longrightarrow \Y(f)\]
has closed image,  finite-dimensional kernel and cokernel, and index
 \[\dim \ker D\dbar'(f)-\dim\coker D\dbar'(f) =2c_{1}-2n(g-1)\] where $c_{1}$ is the integral of the first Chern class of $(T_{vert}\hat{\ex B}, J)$ over the curve $f$, $2n$ is the rank of the vector-bundle $T_{vert}{\hat{\ex B}}$ and $g$ is the genus of the domain of $f$.

 Given any finite-rank  sub-vector-bundle $V\subset \Y(\hat f)$ so that  $D\dbar'(f)$ is transverse to $V(f)\subset \Y(f)$, there  exists a neighborhood $\hat f'$ of $f\in \ex F$, where $D\dbar'$ surjects onto $\Y(\hat f')/V$ with kernel  a finite-rank   sub-vector-bundle  $K=D\dbar'^{-1}(V)$ of $X^{\infty,\underline 1}(\hat f')$. 

\end{thm}

\begin{remark}\label{homotopy application}
The formula for $D\dbar'$ derived in \cite[ equations ($\star'$) and (25)]{reg} implies that the set of maps
\[X^{\infty,\underline 1}(\hat f)\longrightarrow \Y(\hat f)\]
equal to some $D\dbar'$ for some simple perturbation $\dbar'$ of $\dbar$ is convex and contains the complex map \[\frac12(\nabla\cdot+J\circ\nabla \cdot \circ j)\co X^{\infty,\underline 1}(\hat f)\longrightarrow \Y(\hat f)\]
for any $\C\infty1$ $J$--preserving connection $\nabla$ on $T_{vert}\hat{\ex B  }$. In particular, we can apply Theorem \ref{f replacement} to a linear homotopy from $D\dbar'$ to the above complex map. We shall use this homotopy to orient the kernel relative to the cokernel of $D\dbar$.
\end{remark}

\section{Evaluation maps}

\subsection{The  map $ \fun{ev}^{0}$}\label{ev0 section}
 
\

In this section, we consider the exploded manifold version of Deligne--Mumford space. Each map to a point is automatically holomorphic, so the moduli stack $\Msw(pt)$, of stable $\C\infty1$ families of exploded curves mapping to a point, agrees with the moduli stack $\M$ of stable families of holomorphic curves. We will show that $\M$, is an exploded orbifold, and construct a fiberwise-holomorphic  map \[\fun{ev}^{0}\co \Ms(\hat{\ex B})\longrightarrow \M\ .\]

Recall, \cite[Definition 8.3]{iec}, that the ends of an exploded curve $f\co \ex C\longrightarrow \ex B$ are encoded in the structure of the exploded manifold $\ex C$, and that the smooth part $\totl{\ex C}$ is a nodal curve with special points corresponding to the ends of $\ex C$. Accordingly $\M$ is analogous to the Deligne--Mumford stack of stable nodal curves with unspecified genus, and an unspecified number of unlabelled marked points. 

\begin{remark} We also show that this moduli stack of stable exploded curves $\M$ is (represented by) the explosion of the usual Deligne--Mumford space considered as a complex orbifold with normal-crossing divisors corresponding to the boundary. More specifically, the usual Deligne--Mumford space $\bar{ M}$ is a stack over the category of complex manifolds with normal-crossing divisors. It is locally represented by $U/G$ where $U$ is some complex manifold with normal-crossing divisors, and with some effective action of a finite group $G$. The inverse image $U^{+1}$ of $U$ under the forgetful map
\[\bar{ M}^{+1}\longrightarrow \bar{ M}\]
from Deligne--Mumford space with one extra marked point, is a family of nodal curves 
\[U^{+1}\longrightarrow U\]
with automorphism group $G$. The explosion functor applied to this family gives a family of exploded curves 
\[\expl U^{+1}\longrightarrow \expl U\]
with a group $G$ of automorphisms and this family locally represents the moduli stack of stable exploded curves; see Definition \ref{represented by family}. 

Thinking of Deligne--Mumford space as an orbifold locally equivalent to (the stack of holomorphic maps into)  $U/G$, we shall prove that $\M$ is locally equivalent to (the stack of all $\C\infty1$ maps into)  $\expl U/G$, and is, in this sense,  the explosion of Deligne--Mumford space.

\end{remark}

Definition 4.1 and Theorem 3.1 of \cite{uts} introduce the concept of families of curves with universal tropical structure. The tropical part of each curve $f$ in $\hat{\ex B}$ is a tropical curve $\totb f$ in $\totb {\hat{\ex B}}$, and the tropical part of family of curves $\hat f$ containing $f$ is an extension of $\totb f$ to a family of tropical curves $\totb{\hat f}$ in $\totb{\hat {\ex B}}$. Roughly speaking $\hat f$ has universal tropical structure if $\totb{\hat f}$ is the universal extension of $\totb { f}$, admitting a unique map from any other extension of $\totb{\hat f}$. This property, having universal tropical structure, is related  in log geometry to a family of log curves being basic, \cite{GSlogGW}, or having  minimal log structure, \cite{acgw,kim}. We shall need the following two results, respectively Theorem 4.8 and Lemma 4.6 from \cite{uts}:
 
 \begin{thm}\label{G uts family} For any stable curve $f$ in $\Msw(\hat {\ex B})$ with domain not equal to $\ex T$, there exists a family of curves $\hat f$
  containing $f$ with universal tropical structure  so that \begin{enumerate}
 \item there is a  group $G$ of automorphisms of $\hat f$  acting freely and transitively  on the set of maps of $f$ into $\hat f$;
 \item there is only one stratum\footnote{The stratification on an exploded manifold is discussed in Section 4 of \cite{iec}.} $\ex F_{0}$ of $\ex F(\hat f)$  containing the image of a map $f\longrightarrow \hat f$, and the smooth part of this stratum, $\totl{\ex F_{0}}$ is a single point;
 \item the action of $G$ on $\totl{\ex C(\hat f)}$ restricted to the inverse image of $\totl{\ex F_{0}}$ is effective, so $G$ is a subgroup of the group of automorphisms of $\totl{f}$. \end{enumerate}\end{thm}
 
Note that as $G$ is a subgroup of the group of automorphisms of $\totl{f}$ and $f$ is stable, $G$ is a finite group.

\begin{lemma}\label{uts map cor} 
 Let $\hat f$ be a family of curves with universal tropical structure at $f$.  Let $\hat h$ be a family of curves  containing  a  curve  $f'$ with a degree $1$ holomorphic map 
 \[\psi\co \ex C( f')\longrightarrow \ex C(f) \] 
 so that $f'=f\circ\psi$. Then, after possibly restricting to a neighborhood of $f'$ in $\hat h$, there exists an extension of $\psi$ to a map 
\[\begin{tikzcd}\ex C(\hat h)\rar{\Psi}\dar &\ex C(\hat f)\dar
\\\ex F(\hat h)\rar &\ex F(\hat f)
 \end{tikzcd}\]
 so that, in a metric\footnote{See remark 6.6 of \cite{iec}.} on $\hat{\ex B}$, the distance between the maps $\hat h$ and $\hat f\circ \Psi$ is bounded.

\end{lemma}

The following lemma is key for constructing the evaluation map $\fun{ev}^0$ and a local representation $\hat f/G$ for our moduli stack, $\mathcal M(pt)$,  of $\C\infty 1$ families of stable  curves mapping to a point. 

\begin{lemma}\label{rmap}Suppose that $\hat f$ is a family of stable curves in $\mathcal M(pt)$ with universal tropical structure  at some curve $f$ in $\hat f$ so that
the map 
 \[T_f\hat f\co T_{f}\ex F(\hat f)\longrightarrow T_{f}\M:=T_{f}\Msw(pt)\]
 is bijective.
Then given any family of curves  $\hat h$ containing a curve $f'$ with the same genus as $f$ and with a degree $1$ holomorphic map 
\[\psi\co \ex C(f')\longrightarrow \ex C(f)\] there exists an open neighborhood  of $f'$ in $\hat h$ with an extension of the above map to a fiberwise-holomorphic  map 
\[\begin{tikzcd}\ex C(\hat h)\rar{\tilde \psi}\dar& \ex C(\hat f)\dar
\\\ex F(\hat h)\rar &\ex F(\hat f)\end{tikzcd}\]

\end{lemma}

\pf 

By restricting to a neighborhood of $f'$ in $\hat h$ if necessary, Lemma \ref{uts map cor} gives a (not necessarily holomorphic) map 
\[\begin{tikzcd}\ex C(\hat h)\rar{\Psi}\dar &\ex C(\hat f)\dar
\\\ex F(\hat h)\rar &\ex F(\hat f)\end{tikzcd}\]
equal to  the given holomorphic map $\psi$ when restricted to $\ex C(f')$. Our goal is to modify $\Psi$ to a holomorphic map $\tilde \psi$, with the help of Theorem \ref{regularity theorem}.

In order to effectively use Theorem \ref{regularity theorem}, we first need to extend $\Psi$.
 Pull back $T\ex F(\hat f)$ along $\Psi$ to give a vector-bundle $E$ over $\ex F(\hat h)$, then pull back $\hat h$ over $E\longrightarrow \ex F(\hat f)$ to give a family $\hat h'$ of curves:
\[\begin{tikzcd}\ex C(\hat h')\rar\dar &\ex C(\hat h)\dar
\\\ex F(\hat h')=E\rar &\ex F(\hat h)\end{tikzcd}\]
Now extend $\Psi$ to a map $\hat h'$ such that the following diagram commutes
\[\begin{tikzcd}\ex C(\hat h')\rar{\hat h'}\dar &\ex C(\hat f)\dar
\\\ex F(\hat h')\rar &\ex F(\hat f)\end{tikzcd}\]
and such that, restricted to the zero-section, we have our old map $\Psi$, and normal to the zero-section,  the derivative of the map $\ex F(\hat h')\longrightarrow \ex F(\hat f)$ is the identity. As the vertical tangent space to $E$ at the zero-section is  the pullback of $T\ex F(\hat f)$,  this last condition makes sense.

Note that calling the above map $\hat h'$ involves a sleight of hand, because before, our notation implied that $\hat h'$ was a family of curves mapping to a point instead of into $\ex C(\hat f)$. At the zero section, we have a canonical inclusion of $\ex C(f')$ into $\ex C(\hat h')$ ---  we shall further abuse notation by calling $f'$ the restriction of $\hat h'$ to $\ex C(f')\subset \ex C(\hat h')$, so $f'$ now means a map $f'\co \ex C(f')\longrightarrow \ex C(\hat f)$, that is,  our original degree one holomorphic map $\psi\co \ex C(f')\longrightarrow\ex C(\hat f)$.

Consider $\ex C(\hat f)\longrightarrow \ex F(\hat f)$ as a family of targets for using Theorem \ref{regularity theorem}, and  choose a trivialization $(\mathcal F,\Phi)$ for $\hat h'$ in the sense of Definition \ref{trivialization def} using a connection $\nabla$ on $T_{vert}\ex C(\hat f)$,  holomorphic restricted on $\ex C(f)$. Let us check the injectivity of  the linearization, $D\dbar(f')$, of the $\dbar$ operator at $f'$ in $\hat h'$ using this trivialization. The domain,  $X^{\infty,\underline 1}(f')$ of $D\dbar(f')$  consists of $\C\infty1$ sections of the complex vector-bundle $(f')^{*}T\ex C(f)$ and $D\dbar(f')$ is  $\frac 12(\nabla+j\circ \nabla\circ j)$. Any component of $\ex C(f')$ collapsing to a point in $\totl{\ex C(f)}$ is a sphere,  and, on these components,  $D\dbar(f')$ is  the usual $\dbar$ operator acting on complex-valued functions. Hence, on collapsed components, $D\dbar(f')$  has no cokernel and has kernel equal to the constant functions. It follows that any element of the kernel of $D\dbar(f)$ must be constant on all such collapsed components, and must therefore be the pullback of a holomorphic vector-field  from $T\ex C(f)$ because $\ex C(f')\longrightarrow \ex C(f)$ has degree $1$. As $\ex C(f)$ is stable, it has no nonzero holomorphic vector-fields, so $D\dbar(f)$ has trivial kernel.

What is the cokernel of $D\dbar(f)$? Restricted to the inverse image of any small enough open subset of $\ex C(f)$, $D\dbar( f')$ is surjective: this is because all the  collapsed strata of $\ex C(f')$ are spheres or extra edges, and each connected component of $\ex C(f')$ which is collapsed has zero genus.  It follows that the cokernel of $D\dbar(f')$ is equal to the cokernel of  $\frac 12(\nabla +j\circ \nabla\circ j)$ acting on the space of vector-fields on $\ex C(f)$. We can identify this cokernel with $T_{f}\Msw(pt)$ as follows: $T_{f}\Msw(pt)$ is the quotient of $\Y(f)$ by the image of $L_{(\cdot)}j$, which is also the image of $\frac 12(\nabla+j\circ\nabla\circ j)$, because   Lemma \ref{Lvdbar} tells us that
  \[\nabla v+j\circ \nabla(jv) =- j\circ L_{v}j \]
 and $\frac 12(\nabla +j\circ \nabla\circ j)$ and $L_{(\cdot)}j$ are both complex linear. So, the cokernel of $D\dbar(f)$ equals $T_{f}\Ms$. Moreover, by assumption  $T_f\hat f\co T_{f}\ex F(\hat f)\longrightarrow T_f\M$ is bijective.  
 
 We now have that the derivative of $\dbar \hat h'$ at $f'\in \hat h'$ in the fiber direction of $E$ is an isomorphism onto the cokernel of $D\dbar(f')$.  So,  choose a pre-obstruction model $(\hat h',V,\mathcal F,\Phi)$ with $V(f')$  equal to  the image of the derivative of $\dbar\hat h'$ in the $E$ direction at $f'$. Then, Theorem \ref{regularity theorem} implies that, in a neighborhood of $f'$ in $\hat h'$, there exists a $\C\infty1$ section $\nu$ of $(\hat h')^{*}T_{vert}\ex C(\hat f)$  so that $\dbar\nu$ is a section of $V$. As $f'$ in $\hat h$ started off holomorphic, and the derivative of $\dbar\hat h'$ in the $E$ direction at $f'$ is in $V(f)$,  Corollary \ref{trt} then implies that $\nu$ vanishes to first order in the $E$ direction at $f$. Therefore the derivative of $\dbar\nu$ in the $E$ direction at $f'$ is equal to the derivative of $\dbar \hat h'$  in the $E$ direction at $f'$, which is an isomorphism onto $V(f)$.

Therefore, on a neighborhood of $f'$, the intersection of $\dbar\nu$  with the zero section is a section $s$ of $E\longrightarrow \ex F(\hat h)$. Note that $\mathcal F(\nu)\co \ex C(\hat h')\longrightarrow \ex C(\hat f)$ is a family of curves mapping to $\ex C(\hat f)\longrightarrow \ex F(\hat f)$.  The pullback $s^*\mathcal F(\nu)$ is a family of holomorphic curves in the family of targets, $\ex C(\hat f)\longrightarrow \ex F(\hat f)$, and, because $s$ is a section,  the domain of $s^*\mathcal F(\nu)$ is $\ex C(\hat h)$. Thus, $\tilde \psi:=s^*\mathcal F(\nu)$ is the required $\C\infty1$ holomorphic map from a neighborhood of $f'$ in $\ex C(\hat h)$ to $\ex C(\hat f)$.

\stop

\begin{lemma}\label{rsc} Suppose that $\hat f$ is a family of stable curves in $\M$ with a finite group $G$ of automorphisms so that for every curve $f$ in $\hat f$

\begin{enumerate}
\item\label{rsc1} $G$ acts freely and transitively on the set of  maps $f\longrightarrow\hat f$, and each of these maps have different smooth part $\totl{f}\longrightarrow \totl{\hat f}$;
\item\label{rsc2}  $\hat f$ has universal tropical structure at $f$;
 \item \label{rsc3}and the map 
 \[T_f\hat f\co T_{f}\ex F(\hat f)\longrightarrow T_{f}\M:=T_{f}\Msw(pt)\]
 is bijective.
\end{enumerate}

Then an open substack of $\M$ is represented by\footnote{See Definition \ref{represented by family}.} $\hat f/G$ and hence represented by $\ex F(\hat f)/G$.
\end{lemma}

\pf

 Lemma \ref{rmap} implies that, if  a curve $h$ in $\hat h$ is isomorphic to a curve in $\hat f$, then there is a map from a neighborhood of $h$ in $\hat h$ into $\hat f$  extending the given isomorphism. Let $\mathcal U$ indicate the substack of $\M$ consisting of families of curves isomorphic to some curve in $\hat f$. The above extension property implies that the curves in $\hat h$ which are in $\mathcal U$ form an open subset of $\hat h$. So $\mathcal U$ is open in the sense of Definition \ref{open substack}.

Now let $\hat h$ be a family of curves in $\mathcal U$ with two maps $\phi_{i}\co \hat h\longrightarrow \hat f$. Let us show that, if $\hat h$ is connected, then $\phi_{1}$ is $\phi_{2}$ composed with the action of an element of $G$. If these two maps are different at a curve $h$ in $\hat h$, Condition (\ref{rsc1}) implies that the $\phi_{i}(h)$ must differ by the action of an element of $G$, and therefore have different smooth parts. It follows that the two maps $\phi_{i}$ differ on an open  subset of $\ex F(\hat h)$.\footnote{We had to check that $\phi_{i}$ coincide on a closed subset because in general, two maps from an exploded manifold can be equal on a subset which is neither open nor closed.} Similarly, $\phi_{1}$ is equal to $\phi_{2}$ composed with the action of a given element of $G$ on a closed subset of $\ex F(\hat h)$. As restricted to any curve, $\phi_{1}$ is always  equal to $\phi_{2}$ composed with some element of $G$, and $G$ is finite, it follows that the subset on which $\phi_{1}$ is equal to $\phi_{2}$ composed with a given element of $G$ is both open and closed. In particular, if $\hat h$ is a connected family of curves, then $\phi_{1}$ is equal to $\phi_{2}$ composed with some element of $G$.

Around every curve $h$ in a family $\hat h$ of curves in $\mathcal U$, Lemma \ref{rmap} gives us the existence of $\abs G$ maps of a neighborhood of $h$ in $\hat h$ to $\hat f$,  permuted by the action of $G$ on $\hat f$. (Throughout this paper, $\abs G$ indicates the order of the group $G$.) As argued above, these maps are unique up to the action of $G$, so  they patch together to form  a unique\footnote{A unique $G$--bundle  is one that is unique up to canonical isomorphism. When we have such local uniqueness, the cocycle condition for a descent datum is automatically satisfied, so we can patch together our locally defined $G$--bundles to a unique (up to canonical isomorphism) global $G$--bunde.} $G$--bundle  $\hat h\bd G\longrightarrow \hat h$ with a $G$--equivariant map into $\hat f$. In other words, there exists a unique (up to canonical $2$-isomorphism) map of $\hat h$ into $\hat f/G$, in the sense of Example \ref{family map}. It follows that $\mathcal U$ is equivalent to the stack  $\hat f/G$.

\stop

\

Given any stable curve $f$ in $\M$, we could construct  $\M$ near $f$ by using Lemma \ref{G uts family} to construct a family containing $f$  obeying the first two criteria of Lemma \ref{rsc}, then extending this family to also satisfy the last condition of Lemma \ref{rsc}. Instead, we will  apply the explosion functor\footnote{For the explosion functor, see Section 5 of \cite{iec}. } to the universal curve over Deligne--Mumford space in a neighborhood of $\totl f$.

Let $\pi\co U^{+1}\longrightarrow U$ be a family of nodal curves with marked points (in the category of complex manifolds, or schemes over $\mathbb C$) which, when quotiented by its group $G$ of automorphisms, locally represents the universal curve over Deline-Mumford space, as constructed by Deligne, Mumford and Knudsen in \cite{DeligneMumford,KnudsenDM}, or as constructed more geometrically by Robbin and Salamon in \cite{SalamonDM}, where $\pi\co U^{+1}\longrightarrow U$ is called a universal unfolding. 

In either case, the following holds:
\begin{enumerate}
\item\label{dmpG} $G$ acts freely and transitively on the set of inclusions of a given nodal curve into $U^{+1}\longrightarrow U$.
\item $U$ and $U^{+1}$ are (or may be considered as) complex manifolds and $\pi$ is a holomorphic map.
\item $U$ minus the set of smooth curves is a normal-crossing divisor $D$, and $\pi^{-1}D$ is also a normal-crossing divisor, as is the union, $D'$, of $\pi^{-1}D$ with the locus of all marked points.
\item \label{dmp}If  a curve $\pi^{-1}(p)$ has $n$ nodes, then around $p$ there exist holomorphic coordinates $(z_{1},\dotsc,z_{m})$ centered on $p=(0,\dotsc,0)$, such that $D$ is locally $z_{1}\dotsb z_{n}=0$. At the $i$th node of this curve, $\pi^{-1}(p)$, there are local holomorphic coordinates such that $\pi^{*}z_{j}$ are coordinate functions for $j\neq i$, and there are two extra coordinate functions $z^{\pm}_{i}$ so that 
\[\pi^{*}z_{i}=z^{+}_{i}z_{i}^{-}\ .\]
Away from nodes, all the $z_{i}$ pull back to be coordinate functions. And,  around a marked point in $\pi^{-1}(p)$, there are local coordinates $(z,\pi^{*}z_{1},\dotsc,\pi^{*}z_{m})$ so that $D'=\{z\prod_{i}\pi^{*}z_{i}=0\}$. 
\end{enumerate}

We can therefore apply the explosion functor (described in \cite[Section 5]{iec}) to $\pi$.
\[\expl \pi\co \expl U^{+1}\longrightarrow \expl U\] 

\begin{lemma}\label{explodeDM} $\expl\pi$ satisfies the conditions of Lemma \ref{rsc}. In particular,   $\expl U/G$ represents an open substack $\mathcal U$ of $\M$  consisting of curves isomorphic to some curve in $\expl \pi$.
\end{lemma}

\pf

The property (\ref{dmp}) of $\pi$ implies that  $\expl \pi$ is a family of curves. As the smooth part of $\expl \pi$ is equal to $\pi$, $\expl \pi$ is a family of stable curves, and the property (\ref{dmpG}) of $\pi$ implies condition (\ref{rsc1}) of Lemma \ref{rsc}.

 Property (\ref{dmp}) of $\pi$ implies that the tropical structure $\mathcal P(x)$ of $\expl U$ at any point $x\in U$ is equal to $[0,\infty)^{n}$ where $n$ is the number of nodes of the curve $ \pi^{-1}(x)$. Furthermore,  at the $i$th internal edge $e_{i}$ of $\expl \pi^{-1}(x)$, the tropical structure of $\expl U^{+1}$ is given by the fiber product  
 \[\begin{tikzcd} \mathcal P(e_{i})\rar \dar & {[0,\infty)}^{2}\dar{a+b}
\\ \mathcal P(x)\rar & {[0,\infty)} \end{tikzcd}\]
where the bottom arrow is projection onto the $i$th factor of $\mathcal P(x)=[0,\infty)^{n}$. Therefore, the tropical structure
  $\expl \pi$  restricted to any curve $f$ is the universal extension of the tropical structure of $f$; see  \cite[Remark 3.3]{uts}. So,   condition (\ref{rsc2}) of Lemma \ref{rsc} holds. 

All that remains is to show that condition (\ref{rsc3}) of Lemma \ref{rsc} holds. In particular, we must show that for any curve $f$ in $\expl \pi$, the map 
 \[T_f\expl\pi\co T_{f}\expl U\longrightarrow T_{f}\M\]
 is bijective. 
 
  Recall, from Section \ref{tangent section}, that $T_{f}\M$ is the cokernel of the map $L_{\cdot}j$, and Claim \ref{Lvdbar} implies that  $T_{f}\M$ is also  the cokernel of the $\dbar$ operator acting on $\C\infty1$ vector-fields on $\ex C(f)$. The complex dimension of $T_{f}\M$ is $3g-3+k$ --- this follows from Theorem \ref{f replacement}  and a calculation of the first Chern class of $T\ex C(f)$ as $2-2g-k$ where $k$ is the number of ends of $\ex C(f)$ and $g$ is the genus of $\ex C(f)$. As we include the image of marked point sections in the divisor of $U^{+1}$, the marked points in  $U^{+1}\longrightarrow U$ become ends of the curves in $\expl U^{+1}\longrightarrow \expl U$. In particular, the number $k$ of ends agrees with the number of marked point sections.   So, $T_{f}\M$ has the same dimension as Deligne--Mumford space and therefore $\expl U$. It follows that to check that  $T_{f}\expl \pi$ is bijective, we need only check that it is injective.
 
 Let $v$ be in the kernel of the map $T_f\expl \pi\co T_{f}\expl U\longrightarrow T_{f}\M$. The description of this map (\ref{derivative def}) from Section \ref{derivatives section} implies that there must be a lift $v'$ of $v$ to a section of $T\expl U^{+1}$ restricted to $\ex C(f)$ so that $L_{v'}j=0$. It remains to show that such a $v'$ must be $0$. 
 
 Use $\totl{v'}$ to indicate the image of $v'$ under the derivative of the smooth part map $\expl U^{+1}\longrightarrow U^{+1}$, and $\totl v$ to indicate the image of $v$ under the derivative of the map $\expl U\longrightarrow U$. We have that $\totl {v'}$ is a lift of $\totl v$ and that $L_{\totl {v'}}j=0$ where $j$ now indicates the fiberwise almost-complex structure on $U^{+1}\longrightarrow U$ induced by the complex structure on $U^{+1}$. As defined, it is not obvious that $\totl{v'}$ is smooth at nodes and marked points of $\totl{\ex C(f)}$, however, in the coordinates around nodes from property \ref{dmp}, such a vector-field  must be a constant vector-field plus a vector-field which is continuous and holomorphic away from the node or marked point, so it must be smooth. Moreover, $\totl {v'}$ is tangent to each divisor.  As $U^{+1}\longrightarrow U$ represents the moduli stack of stable marked nodal curves, it follows that $\totl{v'}$ and $\totl{v}$ must be $0$.  
 
 We now have that  $\totl{v'}=0$. If $\totl{\ex C(f)}$ has no nodes, then around $f$ the smooth part map $\expl U\longrightarrow U$ is an isomorphism, so $\totl{v}=0$ implies that $v=0$. Now suppose that $\totl{\ex C(f)}$ has $n$ nodes, and use the coordinates from property (\ref{dmp}) of Deligne--Mumford space. In particular, around the $i$th node we have coordinates including $z_{i}^{\pm}$ such that $\pi^{*}z_{i}=z_{i}^{+}z_{i}^{-}$. These coordinates correspond to coordinates $\tilde z_{i}^{\pm}$ on $\expl U^{+1}$ and a coordinate $\tilde z_{i}=\tilde z_{i}^{+}\tilde z_{i}^{-}$ from $\expl U$. As $ z_{i}^{\pm}$ are coordinates around the $i$th node, and  $\totl{v'\tilde z_{i}^{\pm}}=\totl{v'}z_{i}^{\pm}=0$, it follows that $v'\tilde z_{i}^{\pm}=0$, and therefore $v\tilde z_{i}=0$. As the divisor $D$ in $U$ is defined by the product of these $z_{i}$ corresponding to each node, it follows that $v$ must be equal to $0$. We have now shown that the map $T_{f}\expl \pi\co U\longrightarrow T_{f}\mathcal M(pt)$ is injective, and hence an isomorphism. And, our proof is complete.
 
 \stop
 
 \begin{cor} $\M$ is a $\C\infty1$ exploded orbifold\footnote{See Definition \ref{orbifold}.}, and the explosion of Deligne--Mumford space.
 \end{cor}
 
 \pf We have shown that $\M$ is locally represented by $\hat f/G$, where $\hat f$ is obtained by exploding a holomorphic family of stable nodal curves $U^{+1}\longrightarrow U$, whose quotient by $G$ represents an open substack of Deligne--Mumford space (with the analytic topology). As Deligne-Mumford space has a second countable and Hausdorff topology, it follows that $\totl{\M_{top}}$ is second countable and Hausdorff, so $\M$ is an exploded orbifold. 
 
%
 \stop
 
 \begin{lemma}\label{curve stab} Suppose that $f$ is a (possibly unstable) holomorphic curve in $\Ms(\ex B)$. Then exactly one of the following two options holds:
 \begin{enumerate}\item \label{nostab} $f$ is constant, and has genus zero with less than three ends, or has genus 1 and no ends; or 
 \item  there exists a unique\footnote{Unique up to unique isomorphism.} stable curve $f^{st}$ in $\Msw(\ex B)$ and a degree $1$  holomorphic map
 \[\psi\co \ex C(f)\longrightarrow \ex C(f^{st})\] 
 such that
 \[f=f^{st}\circ \psi\]
 and such that $f^{st}$ has the same genus, and the same number of ends as $f$. 
\end{enumerate}
 \end{lemma}
 \pf

 For this proof, call a holomorphic map $\psi\co \ex C(f)\longrightarrow \ex C(f')$ a \emph{partial stabilisation} of $f$ if it has degree $1$, preserves the number of ends and genus, and satisfies
 \[f=f'\circ \psi\ . \]
 Call this map a \emph{stabilisation}  if  $f'$ is stable. We need to construct a stabilisation of $f$, and prove that it is unique up to unique isomorphism. Note that any stabilisation of a stable curve must be an isomorphism.

 The idea is to `remove' all unstable strata using a series of partial stabilisations. Recall that $\totl f\co \totl{\ex C(f)}\longrightarrow \totl{\ex B}$ is a map from a nodal Riemann surface with marked points at the image of ends of $\ex C(f)$, and each stratum in $\ex C(f)$ is the inverse image of a node, marked point, or component of the complement of all nodes and marked points in $\totl{\ex C(f)}$. We refer to these last strata of $\ex C(f)$ as smooth strata, because they are smooth manifolds.   An edge of $\ex C(f)$ refers to the inverse image of a node or marked point, because such strata correspond to edges of the tropical curve $\totb{\ex C(f)}$.  Each unstable stratum  consists of an unstable punctured Riemann surface on which $\totl f$ is constant. The set of points in $\ex B$ sent to a single point in $\totl{\ex B}$ and $\totb{\ex B}$ is always $(\mathbb C^*)^k$ for some $k$. All holomorphic maps from a connected unstable Riemann surface to $(\mathbb C^*)^k$ are either constant, or the Riemann surface is isomorphic to $\mathbb C^*$ and the holomorphic map is a monomial map $z\mapsto (c_1z^{\alpha_1},\dotsc,c_kz^{\alpha_k})$. So, $f$ on such unstable components is either constant, or a monomial map.

\begin{itemize}

\item
If a stratum $S$ of $\ex C(f)$ is a sphere attached to exactly one edge, and $\totl f$ is constant on $S$, we have that $f$ is also constant on $S$. There are then two cases:  either $f$ is constant and has genus 0 and exactly one end, or there are  holomorphic coordinates on a neighborhood of the edge modelled on an open subset of $\et 1{[0,l]}$ with coordinate $\tilde z$ so that $\totl{\tilde z}$ gives local coordinates on the smooth statum $S'$ of $\ex C(f)$ attached to the other end of the edge, and $\totl{\e l\tilde z^{-1}}$ gives local coordinates on $S$ . In the first case, $f$ satisfies (\ref{nostab}) and has no stabilisation, but in the second case, we can perform a partial stabilisation as follows. Replace this coordinate chart with the corresponding open subset  of $\mathbb C$ with coordinate $z=\totl{\tilde z}$. There is an obvious partial stabilisation from our old curve to this new one $f'$, given in this coordinate chart by $\tilde z\mapsto\totl{\tilde z}$, and sending our unstable sphere $S$ and the edge attached to it to  $\{z=0\}$.  (This map  is the identity everywhere else.) 

As stable curves have no constant once-punctured-sphere strata, the resulting partial stabilisation $\ex C(f)\longrightarrow \ex C(f')$ has the property that any stabilisation of $f$ factors uniquely through it.

\

\item
If a stratum of $\ex C(f)$ is a sphere attached to exactly two edges, and $\totl f$ is constant on this stratum,  there exists a holomorphic identification of a neighborhood of this stratum with a refinement of an open subset of $\et 1{[0,l]}$, $\ex T$ or $\et 1{[0,\infty)}$, and $f$ factors through this refinement map because it is a monomial map on this stratum. Replace this open set with the corresponding open subset of $\et 1{[0,l]}$, $\ex T$ or $\et 1{[0,\infty)}$. The partial stabilisation from the old exploded curve to the new curve $f'$ is this refinement map. (Refer to Section 10 of \cite{iec} for the definition of refinements.)

Again, the resulting partial stabilisation $\ex C(f)\longrightarrow \ex C(f')$ has the property that any stabilisation of $f$ must factor uniquely through it.
\end{itemize}

Each of the above types of partial stabilisations reduces the number of smooth strata by one, so after applying maps of the above type a finite number of times, we either obtain a constant curve satisfying (\ref{nostab}), or obtain a connected exploded curve $\ex C(f^{st})$ with no smooth one-or-two-puctured-sphere strata on which $\totl f$ is constant. Our curve $f^{st}$ is therefore stable because either $\ex C(f^{st})=\ex T$ and $f^{st}$ is non-constant, or  $\totl f$ has a finite number of automorphisms. Moreover, the stabilization  $\psi\co \ex C(f)\longrightarrow \ex C(f^{st})$ has the property that any stabilisation $\psi'\co \ex C(f)\longrightarrow \ex C(f')$ of $f$ factors uniquely through it. As the resulting map $\ex C(f^{st})\longrightarrow \ex C(f')$ is a stabilisation of a stable curve, it must be an isomorphism. So our stabilisation is unique up to unique isomorphism, as required. 

\stop
 
 In general, the stabilisation does not work in families, but given a family $\hat f$ of curves in $\Ms(\hat {\ex B})$, we can stabilise the domain of $\hat f$ to obtain a family in $\M$. 
 
 \begin{lemma} \label{ev0}Suppose that $\hat f$ in $\Ms(pt)$ is a connected $\C\infty1$ family of curves  for which $2g+n\geq3$, where $g$ is the genus and $n$  the number of ends of curves in $\hat f$.   Then there exists a unique\footnote{Unique up to canonical isomorphism.} stabilization, $\hat f^{st}$, of $\hat f$, which is  a $\C\infty1$ family of stable curves $\hat f^{st}$ in $\M$ with
 a degree $1$, fiberwise-holomorphic, $\C\infty1$ map
 \[\begin{tikzcd}\ex C(\hat f)\rar\dar &\dar\ex C(\hat f^{st})
 \\\ex F(\hat f)\rar{\id} &\ex F(\hat f^{st})=\ex F(\hat f)\end{tikzcd}\]
 preserving the genus and number of ends of fibers. There is then a fiberwise holomorphic stabilisation map
 \[\begin{tikzcd}\Ms(pt)^{+1}\rar{\fun{st}^{+1}}\dar & \M^{+1}\dar
 \\ \Ms(pt)\rar{\fun{st}} &\M \end{tikzcd}\]
so that $\fun{st}(\hat f)=\hat f^{st}$.

After composing with the forgetful map $\Ms(\hat {\ex B})\longrightarrow \Ms(pt)$,  the above defines a fiberwise-holomorphic\footnote{See Definition \ref{curve fhm}.} map 
\[\begin{tikzcd}\Ms(\hat{\ex B})^{+1}\ar[bend left]{rr}{(\fun{ev}^0)^{+1}}\dar\rar & \Ms(pt)^{+1}\rar{\fun{st}^{+1}}\dar&\M^{+1}\dar
\\ \Ms(\hat{\ex B})\rar \ar[bend right]{rr}{ \fun{ev}^{0}}&\Ms(pt)\rar{\fun{st}}&\M\end{tikzcd}\]
  \end{lemma}
on the components of $\Ms(\hat{\ex B})$ for which $2g+n\geq3$.

\pf

Lemma \ref{curve stab} applied to a curve $f\in \Ms(pt)$ gives us a stabilisation $f^{st}\in \M$ of this curve. 
Let us extend our stabilization of $f$ to a family $\hat f$ of curves containing $f$. Let $\hat g$ be a family of stable curves containing $f^{st}$, and satisfying the requirements of Lemma \ref{rsc} --- such a family of curves exists, as proved by Lemma \ref{explodeDM}. By restricting $\hat f$ to a smaller neighborhood of $f$ if necessary, Lemma \ref{rmap} then gives us a fiberwise holomorphic map 
\[\begin{tikzcd}\ex C(\hat f)\rar{\phi^{+1}}\dar &\ex C(\hat g)\dar
\\ \ex F(\hat f)\rar{\phi} &\ex F(\hat g)\end{tikzcd}\]
  extending the given stabilisation $\ex C(f)\longrightarrow \ex C(f^{st})\subset \ex C(\hat g)$. So long as $\hat f$ is connected, $\phi^{+1}$ is fiberwise degree 1 and preserves genus and number of ends. We can therefore pull back $\hat g$ over the map $\phi$ to obtain a stabilization $\phi^*{\hat g}$ of $\hat f$:
\[\begin{tikzcd}\ex C(\hat f)\dar\rar{\psi}&\ex C(\hat f^{st})=\ex C(\phi^*\hat g)\dar
\\ \ex F(\hat f)\rar{\id}&\ex F(\hat f)\end{tikzcd}\]

We must now verify that the assignment $\hat f\mapsto \hat f^{st}$ is a functor, and defines a fiberwise holomorphic map.
Suppose that $\hat h$ is a family with a stabilization $\hat h^{st}$ and a map $\alpha\co \hat h\longrightarrow \hat f$. We must construct a canonical map $\alpha^{st}\co \hat h^{st}\longrightarrow \hat f^{st}$.  We have assumed that there is a group $G$ of automorphisms of $\hat g$ so that given any curve $h$ in $\hat h$, $G$ acts freely and transitively on the set of maps $h^{st}\longrightarrow \hat g$, and each of these maps has a different smooth part. As every fiberwise-degree-1 holomorphic map $\ex C(h)\longrightarrow \ex C(\hat g)$ factors uniquely through a map  $h^{st}\longrightarrow \hat g$, it follows that $G$ acts freely and transitively on the set of fiberwise-degree-$1$ holomorphic maps $\ex C(h)\longrightarrow \ex C(\hat g)$, and each of these maps has a different smooth part. As in the proof of Lemma \ref{rsc}, Lemma \ref{rmap} then implies that, if $\hat h$ is small enough, there are exactly $\abs G$ degree 1 fiberwise-holomorphic maps $\ex C(\hat h)\longrightarrow \ex C(\hat g)$, permuted by the action of $G$, and corresponding to $\abs G$ maps $\hat h^{st}\longrightarrow \hat g$, also unique up to the action of $G$. There must therefore be exactly one of these maps $\alpha_1$ such that the following diagram commutes 
\[\begin{tikzcd}\ex C(\hat h)\rar{\ex C(\alpha)}\dar&\ex C(\hat f)\dar{\phi^{+1}}
\\ \ex C(\hat h^{st})\rar{\ex C(\alpha_1)}& \ex C(\hat g) \end{tikzcd}\]
Then, as $\hat f^{st}=\phi^*\hat g$ and $\ex F(\alpha)=\ex F(\alpha_1)$,  the map $\alpha_1\co \hat h^{st}\longrightarrow \hat g$ factors uniquely as a compostion  $\hat h^{st}\xrightarrow{\alpha^{st}} \hat f^{st}\longrightarrow \hat g$ of  maps such that the following diagram commutes
\[\begin{tikzcd}\ex C(\hat h)\dar\rar{\ex C(\alpha)}&\ex C(\hat f)\dar
\\\ex C(\hat h^{st})\rar{\ex C(\alpha^{st})}&\ex C(\hat f^{st})\end{tikzcd}\]
The unique-factorization property of stabilizations of individual curves implies that the map $\alpha^{st}\co\hat h^{st}\longrightarrow \hat f^{st}$ is the unique\footnote{Unique here means that $\hat h^{st}$ is defined up to canonical isomorphism, and with such a choice of $\hat h^{st}$, the map $\alpha^{st}$ to $\hat f^{st}$ is unique.} map so that the above diagram commutes.  In this argument, we assumed that $\hat h$ was `small enough', so we have only constructed the $\alpha^{st}$ locally, however, the uniqueness of this map implies that all such local constructions patch together\footnote{Our uniqueness ensures that every collection of local choices for $\hat h^{st}$ give a descent datum for $\hat h^{st}\longrightarrow \hat f^{st}$ in the stack $\M/\hat f^{st}$ because cocycle conditions are automatically satisfied.} to a globally defined map $\alpha^{st}$.

To summarize, we have shown that every `small enough' family $\hat f$ has a stabilization $\hat f^{st}$, and that if $\hat h$ also has a stabilization $\hat h^{st}$ and there is a map $\alpha\co\hat h\longrightarrow \hat f$, then there is a unique map $\alpha^{st}\co \hat h^{st}\longrightarrow \hat f^{st}$ so that the above diagram  commutes. This uniqueness of locally defined stabilizations implies that locally defined stabilizations glue together. Therefore,  every family of curves $\hat f$ satisfying the requirements of our lemma has a stabilization $\hat f^{st}$ and  given any map $\alpha\co \hat h\longrightarrow \hat f$, there exists a unique map $\alpha^{st}\co\hat h^{st}\longrightarrow \hat f^{st}$ so that the above diagram commutes. This uniqueness implies that $(\alpha_1\circ\alpha_2)^{st}=\alpha_1^{st}\circ\alpha_2^{st}$, so we have defined a map of stacks $\mathbf{st}\co \Ms(pt)\longrightarrow \M$.

 Moreover, the fiberwise holomorphic maps $\psi\co \ex C(\hat f)\longrightarrow \ex C(\hat f^{st})$ lift this map to a fiberwise holomorphic map in the sense of Definition \ref{curve fhm}. In particular, let $s$ be a section of $\ex C(\hat f)\longrightarrow \ex F(\hat f)$. We  can define $\mathbf {st}^{+1}(\hat f,s)$ to be $(\hat f^{st},\psi\circ s)$. There is an obvious forgetful map $\fun{Dom}\co \Ms(\hat{\ex B})\longrightarrow \Ms(pt)$ which replaces a family $\hat h$ of curves in $\hat{\ex B}$ with its domain $\ex C(\hat h)\longrightarrow \ex F(\hat h)$. This forgetful map extends to a fiberwise holomorphic map, $(\fun{Dom},\fun{Dom}^{+1})$ and composing this with the fiberwise holomorphic map $(\mathbf{st},\mathbf{st}^{+1})$, we obtain the required fiberwise holomorphic map $(\fun{ev}^0,(\fun{ev}^0)^{+1})$.
\[\begin{tikzcd}\Ms(\hat{\ex B})^{+1}\ar[bend left]{rr}{(\fun{ev}^0)^{+1}}\dar\rar{\fun{Dom}^{+1}} & \Ms(pt)^{+1}\rar{\mathbf{st}^{+1}}\dar&\M^{+1}\dar
\\ \Ms(\hat{\ex B})\rar{\fun{Dom}} \ar[bend right]{rr}{ \fun{ev}^{0}}&\Ms(pt)\rar{\mathbf {st}}&\M\end{tikzcd}\]
defined on the connected components of $\Ms(\hat{\ex B})$ for which $2g+n\geq3$.

\stop

\subsection{The evaluation maps $\fun{ev}^{+n}$ and adding extra marked points to families}\label{+n section}

\
 
 In what follows, we define an `evaluation map' for a family of curves using a functorial construction of a family of curves $\hat f^{+n}$, with $n$ extra (labeled) ends, from a given family of curves $\hat f$. In Definition \ref{f^{+n}}, we list the properties required of $\hat f^{+n}$; we then prove existence and show that this family is unique (up to unique isomorphism) in Lemma \ref{+1 functorial}. 
 
\begin{defn}
Given a submersion\footnote{As usual, a submersion of exploded manifolds means a map with surjective derivative. See Section 6 of \cite{iec} for a discussion of tangent spaces of exploded manifolds, and see Section 9 for a discussion of fiber products.} $f\co \ex D\longrightarrow \ex E$ of exploded manifolds,  use the following notation for the fiber product of $\ex D$ over $\ex E$ with itself $n$ times:
\[\sfp {\ex D}{\ex E}n:=\ex D\fp ff \ex D\fp ff \dotsb \fp ff \ex D\]

\end{defn}

For a family of curves $\hat f$ in $\mathcal U$, we have already defined the family $\hat f^{+1}$ in $\mathcal U^{+1}$ in Definition \ref{basic+1def}. In particular, we constructed  $\hat f^{+1}$ as 
\[\hat f^{+1}:=(\pi_{\ex F(\hat f)}^*\hat f,\Delta)\ ,\] where $\pi_{\ex F(\hat f)}^*\hat f$ is some family  of curves in $\mathcal U$ with domain $\ex C(\hat f)\times_{\ex F(\hat f)}\ex C(\hat f)\longrightarrow \ex C(\hat f)$ and  $\Delta$ is the natural diagonal section of this domain.  In what follows, we upgrade this definition by `exploding' the domain along the image of $\Delta$ to create a family of curves with one extra end. So, we replace the above definition of $\hat f^{+1}$.

\begin{defn}\label{f^{+n}}
Given a family of curves $\hat f$ in $\hat{\ex B}\longrightarrow \ex B_{0}$, define the family of curves $\hat f^{+1}$ to be a family of curves in $\hat{\ex B}\times_{\ex B_{0}}\hat{\ex B}$ with one extra end 
\[\begin{tikzcd} \ex C(\hat f^{+1})\dar{\pi_{\ex F(\hat f^{+1})}}\rar{\hat f^{+1}}&\hat{\ex B}\times_{\ex B_{0}}\hat{\ex B}\dar
\\\ex F(\hat f^{+1})=\ex C(\hat f)\rar{\hat f}&\hat{\ex B}
\end{tikzcd}\]
and satisfying the properties enumerated below.
\begin{enumerate}
\item The fiber of $\pi_{\ex F(\hat f^{+1})}\co \ex C(\hat f^{+1})\longrightarrow \ex F(\hat f^{+1})$ over a point $p\in\ex F(\hat f^{+1})=\ex C(\hat f)$ is equal to the fiber of $\pi_{\ex F(\hat f)}\co \ex C(\hat f)\longrightarrow \ex F(\hat f)$ containing $p$ with an extra end at the point $p$. In other words,  there is a canonical holomorphic map between these curves
\[\pi_{\ex F(\hat f^{+1})}^{-1}(p)\longrightarrow \pi_{\ex F(\hat f)}^{-1}\pi_{\ex F(\hat f)}p \]
 that is an isomorphism on the complement of $p$ and  its  inverse image, but which sends the extra end of $\pi_{\ex F(\hat f^{+1})}^{-1}(p)$ to $p$.
 
\item There exists a fiberwise-holomorphic, degree $1$ map 
\[\begin{tikzcd}\ex C(\hat f^{+1})\arrow{r}\dar &{\ex C}(\hat f)\times_{\ex F(f)}{\ex C}(\hat f)\dar
\\ \ex C(\hat f)\arrow{r}{\id}&\ex C(\hat f) \end{tikzcd}\] so that the following diagram commutes.
\[\begin{tikzcd}
{\ex C}(\hat f^{+1})\dar\ar{dr}{\hat f^{+1}}\ar[bend right=90]{dd}
\\ {\ex C}(\hat f)\times_{\ex F(\hat f)}{\ex C}(\hat f)\dar  \arrow{r}{\hat f\times \hat f} & \dar\hat{\ex B}\times_{\ex B_{0}}{\hat {\ex B}}
\\ {\ex C}(\hat f)\arrow{r}{\hat f}\dar&\hat{\ex B}\dar
\\ \ex F  \arrow{r} & \ex B_{0}
\end{tikzcd}\]

\end{enumerate}
 
 Define $\hat f^{+0}$ to be $\hat f$, and for positive integers $n$, define $\hat f^{+n}$  inductively using 
  \[\hat f^{+n}=\lrb{\hat f^{+(n-1)}}^{+1}\]
 so $\hat f^{+n}$ is a family of curves in  $\sfp{\hat{\ex B}}{\ex B_{0}}{n+1}\longrightarrow \sfp{\hat{\ex B}}{\ex B_{0}}{n}$.
 \[\begin{tikzcd}\ex C(\hat f^{+n})\arrow{r}{\hat f^{+n}}\dar&\sfp{\hat{\ex B}}{\ex B_{0}}{n+1}\dar
\\ \ex C(\hat f^{+(n-1)})\dar[dotted] \arrow{r}{\hat f^{+(n-1)}}&\dar[dotted]\sfp{\hat {\ex B}}{\ex B_{0}}{n}
\\ \ex C(\hat f)\rar{\hat f^{+0}}\dar &\hat{\ex B}\dar
\\\ex F(\hat f)\rar & \ex B_{0}
\end{tikzcd}\]

\end{defn}

In what follows, we construct
$\hat f^{+1}$ satisfying the above requirements;  then we  show that such a family is unique (up to unique isomorphism) and show that the construction is functorial.

Construct the total space of the domain, $\ex C(\hat f^{+1})$ by `exploding' the diagonal of $\ex C(\hat f)\times _{\ex F(\hat f)}\ex C(\hat f)$ as follows:  Consider the diagonal map $\Delta\co \ex C(\hat f)\longrightarrow \ex C(\hat f)\times _{\ex F(\hat f)}\ex C(\hat f)$. The image of the tropical part of this map, $\totb\Delta$, defines a subdivision of the tropical part of   $\ex C(\hat f)\times _{\ex F(\hat f)}\ex C(\hat f)$. As noted in Section 10 of \cite{iec}, any such subdivision  determines a unique refinement $ {\ex C}'\longrightarrow \ex C(\hat f)\times _{\ex F(\hat f)}\ex C(\hat f)$. Note that the diagonal map to  this refinement ${\ex C}'$ is still defined, 
\begin{equation}\label{+1diag}\begin{tikzcd}  & \ex C'\dar
\\\ex C(\hat f)\ar{ur}{\Delta'}\arrow{r}[swap]{\Delta}&\ex C(\hat f)\times _{\ex F(\hat f)}\ex C(\hat f)
\end{tikzcd}\end{equation}
and a  neighborhood of the image of $\Delta'$ in $\ex C'$ is isomorphic to a neighborhood of $0$ in a $\mathbb C$--bundle over $\ex C(\hat f)$.
 Now `explode' the image of the diagonal $\Delta'$ in $\ex C'$ to make $\ex C(\hat f^{+1})\longrightarrow \ex C'$ as follows: 
Choose coordinate charts on $\ex C'$ so that any coordinate chart intersecting the image of the diagonal is  some subset of  $\mathbb C\times U$ where $U$ is a coordinate chart on $\ex C(\hat f)$, the projection to $\ex C(\hat f)$ is the obvious projection to $U$, the complex structure on the fibers of this projection is  the standard complex structure on $\mathbb C$, and  the image of the diagonal is $0\times U$. Replace these charts with the corresponding subsets of $\et 1{[0,\infty)}\times U$, and leave coordinate charts that do not intersect the image of the diagonal unchanged. Any transition map between coordinate charts of the above type is of the form $(z,u)\mapsto (g(z,u)z,\phi(u))$ where $g(z,u)$ is $\mathbb C^{*}$--valued and fiberwise-holomorphic in $z$. In the corresponding exploded charts, the transition map is  $(\tilde z,u)\mapsto(g(\totl{\tilde z},u)\tilde z,\phi(u))$, and  transition maps between other charts remain unchanged. This defines the domain $\ex C(\hat f^{+1})$ of $\hat f^{+1}$. The map $\ex C(\hat f^{+1})\longrightarrow \ex C'$ is given in the above coordinate charts by $(\tilde z,u)\mapsto (\totl{\tilde z},u)$. Composing this map with the refinement map $\ex C'\longrightarrow \ex C(\hat f)\times _{\ex F(\hat f)}\ex C(\hat f)$ then gives a degree-one fiberwise-holomorphic map 
\[\begin{tikzcd}\ex C(\hat f^{+1})\rar\dar &\ex C(\hat f)\times _{\ex F(\hat f)}\ex C(\hat f)\dar
\\ \ex C(\hat f)\rar{\id}&\ex C(\hat f)\end{tikzcd}\] 
The map $\hat f^{+1}\co \ex C(\hat f^{+1})\longrightarrow\hat {\ex B}\times _{\ex B_{0}}\hat{\ex B}$ is given by the above constructed map $\ex C(\hat f^{+1})\longrightarrow \ex C(\hat f)\times _{\ex F(\hat f)}\ex C(\hat f)$ composed with the map 
\[\hat f\times \hat f\co \ex C(\hat f)\times _{\ex F(\hat f)}\ex C(\hat f)\longrightarrow \hat {\ex B}\times _{\ex B_{0}}\hat{\ex B}\] which is $\hat f$ in each component. All the above maps are smooth or $\C\infty 1$ if $\hat f$ is. This constructed family of curves $\hat f^{+1}$ obeys the requirements of Definition \ref{f^{+n}}.

\

The following lemma implies that $\hat f^{+n}$ is unique (up to unique isomorphism) and that the construction of $\hat f^{+n}$ is functorial.

\begin{lemma} \label{+1 functorial}Given a map of families $\psi\co \hat f\longrightarrow \hat g$ and families $\hat f^{+1}$ and $\hat g^{+1}$ satisfying the requirements of Definition \ref{f^{+n}} there is a unique induced map $\psi^{+1}\co \hat f^{+1}\longrightarrow \hat g^{+1}$ such that the diagram
\[\begin{tikzcd}\ex C(\hat f^{+1})\rar{\ex C(\psi^{+1})}\dar&\ex C(\hat g^{+1})\dar
\\ \dar \ex C(\hat f)\times_{\ex F(\hat f)}\ex C(\hat f) \rar & \ex C(\hat g)\times_{\ex F(\hat g)}\ex C(\hat g)\dar
\\\ex C(\hat f) \rar{\ex C(\psi)} &\ex C(\hat g)
\end{tikzcd}\]
 commutes.
 \end{lemma}

\pf

Both $\ex C(\hat g^{+1})$ and $\ex C(\hat g)\times _{\ex F(\hat g)}\ex C(\hat g)$ are families over $\ex C(\hat g)$, but the fiber in $\ex C(\hat g^{+1})$ over a point $p\in \ex C(\hat g)$ has one extra end. Away from this end, the map $\ex C(\hat g^{+1})\longrightarrow \ex C(\hat g)\times _{\ex F(\hat g)}\ex C(\hat g)$ is a fiberwise-holomorphic isomorphism. Therefore,  away from the extra end, the required map 
\[\begin{tikzcd}\ex C(\hat f^{+1})\rar{\ex C(\psi^{+1})}\dar&\ex C(\hat g^{+1})\dar
\\ \ex C(\hat f)\times_{\ex F(\hat f)}\ex C(\hat f) \rar & \ex C(\hat g)\times_{\ex F(\hat g)}\ex C(\hat g)
\end{tikzcd}\]
exists, is unique, and is fiberwise holomorphic.

On a neighborhood of the extra end in the fiber over $p\in \ex C(\hat g)$, there exists a fiberwise-holomorphic exploded coordinate function $\tilde z$ so that the extra end is at $\totl{\tilde z}=0$. The fiberwise-holomorphic function $\totl{\tilde z}$ is a fiberwise-holomorphic coordinate function on  $\ex C(\hat g)\times_{\ex F(\hat g)}\ex C(\hat g)$  vanishing on the image of the diagonal. Therefore, $\totl{\tilde z}$ pulls back to a fiberwise-holomorphic coordinate function on $\ex C(\hat f)\times_{\ex F(\hat f)}\ex C(\hat f)$  vanishing on the image of the diagonal. It follows that, if $\tilde z'$ is a locally defined fiberwise-holomorphic coordinate function on $\ex C(\hat f^{+1})$ so that the extra end is at $\totl{\tilde z'}=0$, then the pullback of $\totl{\tilde z}$ is  $h\totl{\tilde z'}$ where $h$ is some $\mathbb C^{*}$--valued fiberwise-holomorphic function. Therefore there locally exists a unique map
\[\begin{tikzcd}\ex C(\hat f^{+1})\rar{\ex C(\psi^{+1})}\dar&\ex C(\hat g^{+1})\dar
\\ \ex C(\hat f)\times_{\ex F(\hat f)}\ex C(\hat f) \rar & \ex C(\hat g)\times_{\ex F(\hat g)}\ex C(\hat g)
\end{tikzcd}\]
which is fiberwise holomorphic and pulls back $\tilde z$ to $h\tilde z'$. As our locally defined maps all satisfy the same uniqueness property, they glue together to give the required unique map. Restricted to each fiber, this map is a holomorphic isomorphism. The fact that $\hat f^{+1}$ factors through $\ex C(\hat f)\times_{\ex F(\hat f)}\ex C(\hat f)$ implies that our map $\ex C(\hat f^{+1})\longrightarrow \ex C(\hat g^{+1})$ corresponds to a unique map $\hat f^{+1}\longrightarrow \hat g^{+1}$.

\stop

\
 
 Given a moduli stack of curves $\mathcal X$, we can use this construction to upgrade Definition \ref{+1 stack} to describe $\mathcal X^{+1}$ as a moduli stack of curves; see Definition \ref{stack of curves}. Recall that a family $(\hat f,s)$ in $\mathcal X^{+1}$ is a family $\hat f$ in $\mathcal X$, and a section $s\co \ex F(\hat f)\longrightarrow \ex C(\hat f)$. Morphisms $(\hat f,s)\longrightarrow (\hat g,s')$ in $\mathcal X^{+1}$ are given by morphisms $\psi\co \hat f\longrightarrow\hat g$ compatible with the sections.

\begin{defn}\label{+1 curve stack} If $\mathcal X$ is a moduli stack of curves in $\hat {\ex B}\longrightarrow \ex B_0$, define the moduli stack $\mathcal X^{+1}$ of curves in $\sfp{\hat {\ex B}}{\ex B_0}{2}\longrightarrow \hat {\ex B}$ to be the stack $\mathcal X^{+1}$ together with the map 
\[\fun\Phi\co \mathcal X^{+1}\longrightarrow  \Ms(\sfp{\hat {\ex B}}{\ex B_0}{2}) \]
defined on families by
 \[\fun\Phi(\hat f,s)=s^*\hat f^{+1} \ .\]
To define $\fun\Phi$ on morphisms,  given $\psi\co \hat f\longrightarrow \hat g$ inducing a morphism $(\hat f,s)\longrightarrow (\hat g,s') $ in $\mathcal X^{+1}$, define 
$\fun\Phi(\psi)\co s^*\hat f\longrightarrow (s')^*\hat g$ to be the unique morphism such that the following diagram commutes.
\[\begin{tikzcd} \ex C(s^*\hat f)\dar \rar{\fun\Phi(\psi)} & \ex C((s')^*\hat g)\dar 
\\ \ex C(\hat f^{+1})\rar{\psi^{+1}} &\ex C(\hat f^{+1}) \end{tikzcd}\]

For $n$ a positive integer, define the moduli stack of curves $\mathcal X^{+n}$ in $\sfp{\hat {\ex B}}{\ex B_0}{n+1}\longrightarrow \sfp{\hat {\ex B}}{\ex B_0}{n}$  to be $(\mathcal X^{+(n-1)})^{+1}$.

\end{defn}

Note that $\Msw(\hat {\ex B})^{+n}$ is a moduli stack of decorated curves in $\sfp{\hat {\ex B}}{\ex B_0}{n+1}\longrightarrow \sfp{\hat {\ex B}}{\ex B_0}{n}$, because the domain of a curve in $\Msw(\hat {\ex B})^{+n}$ has $n$ distinguished ends (as well as an unspecified number of extra, unlabelled ends).
 
 The stabilisation map $\fun{st}\co \Ms(pt)\longrightarrow \M$ from Lemma \ref{ev0} extends to a fiberwise holomorphic tower of maps
 \[\begin{tikzcd}\Ms(pt)^{+n}\ar{rr}{\fun{st}^{+n}}\dar && \M^{+n} \dar
 \\ \Ms(pt)^{+(n-1)}\ar{rr}{\fun{st}^{+(n-1)}} && \M^{+(n-1)} \end{tikzcd}\]
defined so that, given a family $\hat h$ in $\Ms(pt)^{+n}$,  $\fun{st}^{+n}(\hat h)$ is $\hat h^{st}$, along with the marking distinguishing the $n$ special ends of the family $\hat h$ so that the resulting stabilised family is in $\M^{+n}$. As with $\fun{st}$, this map is only defined on components of $\Ms(pt)^{+n}$ describing curves with sufficient topology for the stabilisation to be defined, however it is always defined when $n\geq 0$.

Similarly, taking the domain of a family in $\Msw(\hat {\ex B})^{+n}$, but remembering the $n$ special ends, defines a map of stacks 
\[\fun{Dom}^{+n}\co \Msw(\hat {\ex B})^{+n}\longrightarrow \Ms(pt)^{+n}\]
Composing these two maps we define a map of stacks $(\fun{ev}^0)^{+n}\co \Msw(\hat {\ex B})^{+n}\longrightarrow \M^{+n} $.
\[\begin{tikzcd} \Msw(\hat {\ex B})^{+n}\ar[bend left]{rr}{(\fun{ev}^0)^{+n}} \rar{\fun{Dom}^{+n}} & \Ms(pt)^{+n}\rar{\fun{st}^{+n}} & \M^{+n}\end{tikzcd}\]

There is also a natural map of stacks
\[\fun{Base}\co \Msw(\hat {\ex B})^{+n}\longrightarrow  \sfp{\hat {\ex B}}{\ex B_0}{n}\]
sending a family 
\[\begin{tikzcd}\ex C(\hat h)\dar\rar{\hat h} & \sfp{\hat {\ex B}}{\ex B_0}{n+1}\dar
\\\ex F(\hat h)\rar & \sfp{\hat {\ex B}}{\ex B_0}{n} \end{tikzcd}\]
to the map $\ex F(\hat h)\longrightarrow \sfp{\hat {\ex B}}{\ex B_0}{n}$ in the bottom row of the above diagram. Combining $(\fun{ev}^0)^{+n}$ and $\fun{Base}$, we get a natural map of stacks
\[\fun{ev}^{+n}:=((\fun{ev}^0)^{+n},\fun{Base})\co \Msw(\hat {\ex B})^{+n}\longrightarrow \M^{+n}\times \sfp{\hat {\ex B}}{\ex B_0}{n}\ \]
defined whenever $n\geq 3$, or, in the case that $n<3$, defined on components of $\Msw(\hat {\ex B})^{+n}$ describing curves with sufficient topology for the stabilisation to exist. Applying $\fun{ev}^{+n}$ to $\hat f^{+n}$, we get the evaluation map
\begin{equation}\label{ev def} \fun{ev}^{+n}_{\hat f}:=((\fun{ev}^{0})^{+n}(\hat f^{+n}),\hat f^{+(n-1)})\co \ex F(\hat f^{+n})\longrightarrow \M^{+n}\times\sfp{ \hat{\ex B}}{\ex B_{0}}n \ .\end{equation}

\section{Core families}
\label{core section}

Core families are an essential ingredient in our  local description of the ambient moduli stack $\Msw(\hat {\ex B})$ of stable, not-necessarily-holomorphic  $\C\infty1$ curves. Some such notion is necessary, because the `space' of stable curves in $\hat{\ex B}\longrightarrow \ex B_{0}$ of a given regularity is not locally modelled on (an orbifold version of) a Banach space --- this is because the domains of curves  are not fixed, and because of  bubble and node formation.  The ambient `space' of stable curves could be described as a `orbifold' by using an exploded manifold adaption   of the theory of polyfolds developed by Hofer, Wysocki and Zehnder in  \cite{polyfold0,polyfold1,polyfold2,polyfoldgw,polyfoldsc,polyfoldint}; given such a description, families in our ambient moduli stack would then be $\C\infty 1$ maps of finite dimensional exploded manifolds into this ambient polyfold.    Such an adaption of the theory of polyfolds to the exploded setting is a worthwhile direction for further research, not explored in this paper.

Given an open substack $\mathcal O$ of our (infinite dimensional) ambient moduli stack, a core family $((\fun\Phi,\fun\Phi^{+1}), \hat f/G, \{s_i\})$  for $\mathcal O$ provides a way of capturing the interesting structure of $\mathcal O$ with  a (finite dimensional) family $\hat f$ in $\mathcal O$.  In particular, there is a map of stacks $\fun\Phi\co O\longrightarrow \ex F(\hat f)/G$ that plays the role of a retraction of the infinite dimensional stack $\mathcal O$ onto a simple finite-dimensional substack $\ex F(\hat f)/G\equiv \hat f/G\subset \mathcal O$, and this map lifts to a fiberwise holomorphic map $(\fun\Phi,\fun\Phi^{+1})$. Necessarily, this retraction only respects the domains of families in $\mathcal O$, and not the corresponding maps into $\hat{\ex B}$ so that is is not in general true that $\hat h=\hat f/G\circ \fun\Phi_{\hat h}^{+1}$. (Here we use $f/G$ to indicate the map $f/G\co \ex C(\hat f)/G\longrightarrow \hat{\ex B}$.) Nevertheless it is true that $\hat h$ is a deformation of $\hat f/G\circ\fun\Phi_{\hat h}^{+1}$, and these two maps agree on the inverse image of the `marked point' sections $s_i\co \ex F(\hat f)\longrightarrow \ex C(\hat f)$.  Moreover, the fiberwise holomorphic map $(\fun\Phi,\fun\Phi^{+1})$ is determined by these sections $\{s_i\}$, so that $(\fun\Phi,\fun\Phi^{+1})$ is constant on one-parameter deformations $f_t$ of  a curve $f$ in $\hat f$, so long as this   one-parameter deformation of $f$  does not change the domain of $f$, or modify the map $f$ on the image of the marked point sections. 
  
\begin{defn}\label{core family}
A core family   $((\fun\Phi,\fun\Phi^{+1}), \hat f/G,\{s_{i}\})$ for an open substack $\mathcal O$ of $\Ms(\hat {\ex B})$ or $\dmsw(\hat{\ex B})$ is the following data:
\begin{itemize}
 \item a  $\C\infty 1$ family, $\hat f$, of stable curves in $\mathcal O$ with a group, $G$, of automorphisms,
 \[\begin{tikzcd}\ex C(\hat f)\rar{\hat f}\dar&\dar\hat{\ex B}
\\ \ex F(\hat f)\rar &\ex B_{0}
 \end{tikzcd}\]
 \item a finite, nonempty $G$--invariant collection of $\C\infty 1$ `marked point' sections $s_{i}\co \ex F(\hat f)\longrightarrow{\ex C}(\hat f)$ which do not intersect each other, and which do not intersect the edges of  curves in $ \ex C(\hat f)$,
 \item and a fiberwise-holomorphic map comprised of a strictly commuting diagram of stacks 
 \[\begin{tikzcd}\mathcal O^{+1}\dar \rar{\fun\Phi^{+1}}& \ex C(\hat f)/G\dar
 \\\mathcal O\rar{\fun\Phi} &\ex F(\hat f)/G \end{tikzcd}\]
 with an effective $G$--action in the sense of Definition \ref{fhm},

  \end{itemize}
   such that the following holds:
  \begin{enumerate}
 
 \item\label{cc1} The fiberwise holomorphic map $(\fun\Phi,\fun\Phi^{+1})$ is a retraction onto $\hat f/G\subset \mathcal O$ in the following sense. 
 Using the notation of Definition \ref{fhm},
 \[\fun\Phi_{\hat f}^{+1}\co \ex C(\hat f)\longrightarrow \ex C(\hat f)/G\] 
 is (canonically $2$-isomorphic to) the quotient map 
 \[\fun\pi\co \ex C(\hat f)\longrightarrow \ex C(\hat f)/G\]
  from Example \ref{quotient map}. So, the following is a fiber product diagram.
 \[ \begin{tikzcd} \ex C(\hat f\bd G)\rar{\fun\Phi^{+1}(\hat f)}\dar &\ex C(\hat f)\dar{\fun\pi}
\\ \ex C(\hat f)\rar{\fun\pi}& \ex C(\hat f)/G \end{tikzcd}\] 

\item\label{cc2} 

For any family  of curves $\hat h$ in $\mathcal O$, there exists a $1$ parameter family of maps 
\[\hat h_t\co \ex C(\hat h)\longrightarrow \hat{\ex B}\]
corresponding to a  $\C\infty1$ family of curves 
parametrized by $\mathbb R\times \ex F(\hat h)$, with domain $\mathbb R\times \ex C(\hat h)$, such that 
\begin{itemize}\item $\hat h_0$ is $\hat h$,
\item $\hat h_1$ is\footnote{As this is a composition of maps of stacks, we can only say that $h_{1}$ represents the composition, however this uniquely specifies the map $h_{1}$.} the composition
\[\begin{tikzcd}\ar[bend left]{rr}{\hat h_1} \ex C(  \hat h)\rar[swap]{\fun\Phi^{+1}_{ \hat h}}&\ex C(\hat f)/G\rar[swap]{\hat f/G}& \hat{\ex B}\end{tikzcd} \ .\]

\item and $\hat h_t$ is independent of $t$ on the special points $(\fun\Phi_{\hat h}^{+1})^{-1}s_i(\ex F(\hat f))$. 
\end{itemize}
\item\label{cc3} 
Conversely, given any curve $f$ in $\hat f$ and a one--parameter family of maps $f_t\co \ex C(f)\longrightarrow \hat{\ex B}$ independent of $t$ on the image of $s_i$, and with $f_0=f$, consider the open subset of $\mathbb R$ comprised of $t$ such that $f_t$ is in $\mathcal O$.  For all such $t$ 
\[\fun\Phi_{f_t}^{+1}\co \ex C(f)\longrightarrow  \ex C(\hat f)/G\]
is independent of $t$.


%

\end{enumerate}
\end{defn}

 Proposition \ref{smooth model family}, constructs a core family  containing any given stable holomorphic curve with at least one smooth component. Note that conditions \ref{cc1} and \ref{cc3}  determine  $(\fun\Phi,\fun\Phi^{+1})$  (up to unique  $2$-isomorphism) from $\hat f/G$ and $\{s_{i}\}$, so we can specify a core family using only $(\hat f/G,\{s_{i}\})$. So, where convenient, we will refer to a core family $((\fun\Phi,\fun\Phi^{+1}),\hat f/G,\{s_i\})$ using only $(\hat f/G,\{s_i\})$.

Condition (\ref{cc1}) implies that $\hat f/G$ is a substack of $\mathcal O$. Any map $\psi\co\hat h\longrightarrow \hat f$ induces a commutative diagram
\[\begin{tikzcd}\ex C(\hat h\bd G)\rar{\ex C(\psi\bd G)} \ar{dr}[swap]{\fun\Phi^{+1}(\hat h)}& \ex C(\hat f\bd G)\dar{\fun\Phi^{+1}(\hat f)}
\\ & \ex C(\hat f)\end{tikzcd}\]
and as a consequence, if such a $\psi$ exists, there are precisely $|G|$ such maps, permuted by the action of $G$. We also can conclude from the above commutative diagram that $\hat f\circ \fun\Phi^{+1}(\hat h)=\hat h\bd G$, because $\hat f\circ \fun\Phi^{+1}(\hat f)=\hat f\bd G$. 
 As spelled out in Definition \ref{family quotient stack}, a family in $\hat f/G$ consists of a $G$--bundle $\hat h\bd G\longrightarrow \hat h$ internal to $\mathcal O$ and a $G$--equivariant map $\hat h\bd G\longrightarrow \hat f$. The above argument implies that given $\hat h$, such a family in $\hat f/G$ is unique (up to unique isomorphism). We can use this fact to verify that the map $\fun\Phi_{\mathcal X}\co \hat f/G\longrightarrow \mathcal O$ is an isomorphism onto a substack. In particular, applying $\fun \Phi_{\mathcal X}$ to such a family $\hat h\bd G\longrightarrow \hat f$ returns the family $\hat h$ in $\mathcal O$. Then $\fun\Phi(\hat h)$ is a family in $\ex F(\hat f)/G\equiv\hat f/G$ but, as argued above,  there exists a unique isomorphism between this family and our original family in $\hat f/G$. This implies that $\fun\Phi_{\mathcal X}\co \hat f/G\longrightarrow O$ is an inclusion onto a substack, with inverse $\fun\Phi$ followed by the equivalence $\ex F(\hat f)/G\equiv \hat f/G$. 

Unpacking the definition of a fiberwise holomorphic map $(\fun\Phi,\fun\Phi^{+1})$, for each family of curves $\hat h\in\mathcal O$,  we have a canonical map of stacks $\fun\Phi^{+1}_{\hat h}\co \ex C(\hat h)\longrightarrow \ex C(\hat f)/G$, from (\ref{phi_f}), which amounts to a $G$--bundle $\ex C(\hat h\bd G)$ over $\ex C(\hat h)$ together with a $G$--equivariant, fiberwise-holomorphic map of exploded manifolds
\[\fun\Phi^{+1}(\hat h)\co \ex C(\hat h\bd G)\longrightarrow \ex C(\hat f)\]
Note that we do not require the map $\fun\Phi^{+1}(\hat h)$ to be compatible with the maps $\hat f$ or $\hat h$, and we do not require that this map be an isomorphism on fibres, although both these things are true in the special case that $\hat h$ is contained in the substack $\hat f/G$. In general, there is an open substack of $\mathcal O$ comprised of families $\hat h$ such that $\fun\Phi^{+1}(\hat h)$ is an isomorphism on fibres. 
 
 \begin{remark}\label{decorated core family} Let $\fun\pi\co \dmsw\longrightarrow \Msw$ be a decorated moduli stack of curves. 
The lifting property of Definition \ref{decorated} implies that a deformation of a  family of decorated curves is uniquely defined by deforming the underlying undecorated family. Let  $\hat f/G$ be a core family for $\mathcal O\subset \Msw$, and suppose that $\mathcal O$ is small enough that for each family of curves $\hat h$ in $\mathcal O$,  the map $\fun\Phi^{+1}(\hat h)\co \ex C(\hat h\bd G)\longrightarrow \ex C(\hat f)$ is an  isomorphism on each fibre.  Consider the pullback  $\fun\pi^*\hat f$ of $\hat f$ to $\dmsw$. Then, the lifting property of Definition \ref{decorated} implies that $\fun\pi^*\hat f/G$ is a core family for the inverse image of $\mathcal O$ in $\dmsw$.  In particular, we can define
a fiberwise holomorphic map $\fun \pi^*\fun\Phi^{+1}\co \fun\pi^*\mathcal O^{+1}\longrightarrow \ex C(\fun\pi^*\hat f)/G$ as follows. Given a family $\hat h$ in $\fun\pi^{-1}\mathcal O$, let $(\fun\pi\hat h)_t$ be the $\mathbb R$--deformation of $\hat h$ from Definition \ref{core family} part \ref{cc2}, and let $\hat h_t$ be the unique lift of this to a $\mathbb R$--deformation of $\hat h$ such that $\fun\pi(\hat h_t)=(\fun\pi\hat h)_t$. Then $\fun\Phi^{+1}_{\fun\pi\hat h_1}\co \ex C(\fun\pi\hat h)\longrightarrow \ex C(\hat f)/G$ is compatible with the maps $\hat h_1$ and $\hat f/G$, and is an isomorphism on each fibre,  so defines a map $\phi\co \fun\pi\hat h_1\longrightarrow \hat f/G$ and therefore a map $\fun\pi^*\phi\co \hat h_1\longrightarrow \fun\pi^*\hat f/G$. We can then define $\fun\pi^*\fun\Phi^{+1}_{\hat h}$ to be $\ex C(\fun\pi^*\phi)\co \ex C(\hat h)\longrightarrow \ex C(\fun\pi^*\hat f)/G$. Condition \ref{cc3} of Definition \ref{core family} ensures that $\fun\pi^*\fun\Phi^{+1}_{
\hat h}$ is independent of choice of deformation $(\fun\pi\hat h)_t$, and that $\fun\pi^*\fun\Phi$ satisfies Definition \ref{fhm}. Because $\fun\Phi^{+1}_{\hat f}$ represents the functor $\ex C$ applied to the quotient map $ \fun q\co \hat f\longrightarrow \hat f/G$, we have that $\fun\pi^*\fun\Phi^{+1}_{\fun\pi^*\hat f}$ represents $\ex C$ applied to  $\fun\pi^*\fun q\co \fun\pi^*\hat f\longrightarrow \fun\pi^*\hat f/G$, and  Definition \ref{decorated} implies that $\fun\pi^*\fun q$ is also a quotient map.

\end{remark}

\begin{remark} \label{deformation remark} Suppose we have a fiberwise holomorphic map  $((\fun\Phi,\fun\Phi^{+1}),\hat f/G,\{s_i\})$ for an open substack $\mathcal U$ satisfying conditions \ref{cc1} and \ref{cc3} of Definition \ref{core family}, and satisfying the following weaker version of Condition \ref{cc2} for all $\hat h\in \mathcal U$: 
\[\hat h\rvert_{(\fun\Phi^{+1}_{\hat h})^{-1}S}=\hat f/G\circ \fun\Phi^{+1}_{\hat h}\rvert_{(\fun\Phi^{+1}_{\hat h})^{-1}S}\]
where  $S$ denotes the image of the marked point sections $\{s_i\}$ in $\ex C(\hat f)/G$, and
\begin{equation}\label{tropical cond}\totb{\hat h}=\totb{\hat f}/G\circ \totb{\fun\Phi^{+1}_{\hat h}}\ ,\end{equation}
 Here, the notation $\totb{\hat h}$ indicates the tropical part of a map $\hat h$; see \cite[Section 4]{iec}, and note that deforming a map does not change its tropical part.  If this weaker condition holds, then there exists an open neighborhood $\mathcal O$ of $\hat f/G$ on which Condition \ref{cc2} holds.  In particular  we can choose a $G$--equivariant  trivialization $(\mathcal F,\Psi)$ for $\hat f$ as in Definition \ref{trivialization def}, and there exists a neighborhood $\mathcal O$ of $\hat f$ with each family  $\hat h$ in $\mathcal O$ uniquely expressible as 
\[\hat h=(\mathcal F/G)\circ (\nu_{\hat h})\]
for a unique  lift $\nu_{\hat h}$ of $\fun\Phi_{\hat h}$ as follows 
\[\begin{tikzcd}& \hat f^*T_{vert}\hat{\ex B}/G\dar
\\ \ex C(\hat h)\ar{ur}{\nu_{\hat h}}\rar{\fun\Phi^{+1}_{\hat h}}& \ex C(\hat f)/G\end{tikzcd}\]
such that $\nu_{\hat h}$ vanishes at the image of the sections $s_i$. Then, we have a canonical deformation $\hat h_t:=(\mathcal F/G)\circ ((1-t)\nu_{\hat h})$ satisfying condition \ref{cc2} of Definition \ref{core family}. Moreover, we can choose this neighborhood such that $\hat h_t$ remains in this neighborhood for $0\leq t\leq 1$.  Unpacking definitions, $\fun\Phi^{+1}_{\hat h}$ is a $G$--bundle $\hat h\bd G\longrightarrow \hat h$ with a $G$--equivariant map $\fun\Phi^{+1}(\hat h)\co \ex C(\hat h\bd G)\longrightarrow \ex C(\hat f)$ and $\nu_{\hat h}$ is a $G$--equivariant lift $\nu'_{\hat h}$ of this map such that the following diagram commutes
\[\begin{tikzcd}\hat{\ex B}& \lar[swap]{\mathcal F}\hat f^*T_{vert}\hat{\ex B}\dar
\\ \ex C(\hat h\bd G)\uar{\hat h\bd G}\ar{ur}{\nu'_{\hat h}}\rar{\fun\Phi^{+1}(\hat h)}& \ex C(\hat f)\end{tikzcd}\]
 and $\nu_{\hat h}'$ vanishes on the image of the sections $s_{i}$. 
 
 In particular, conditions \ref{triv1} and  \ref{injective derivative} of Definition \ref{trivialization def} imply that, given a metric on $\hat{\ex B}$     (as in Remark 6.6  of \cite{iec})  there exists some continuous positive function $m\co \ex C(\hat f)\longrightarrow (0,\infty)$ such that there exists a unique $\nu_{\hat h}'$ as above whenever $\dist (\hat h\bd G,\hat f\circ \fun\Phi^{+1}(\hat h))<\fun\Phi^{+1}(\hat h)\circ m$. Equation (\ref{tropical cond}) ensures that this distance is a continuous function, so this equation is satisfied on an open subset of $\ex C(\hat h\bd G)$. Define $\mathcal O$ as the substack consisting of curves $h$ satisfying the above such that $h_t$ remains in $\mathcal U$ for $t\in[0,1]$. As $\mathcal U$ is open, an open subfamily of $\hat h$  is in $\mathcal O$,  so $\mathcal O$ is open. Whenever $ h=\hat f\circ \fun\Phi^{+1}_{ h}$,  $ h$ is automatically in $\mathcal O$, so  $\mathcal O$ is an open neighborhood of $\hat f$, and also contains any refinement of a curve in $\hat f$ that is already in $\mathcal U$.

In the case that $\mathcal F$ is determined by exponentiation using a connection on $T_{vert}\hat{\ex B}$, we can think of $\nu_{\hat h}$ as a section of $(\hat f/G\circ \fun\Phi_{\hat h})^{*}T_{vert}{\hat {\ex B}}$ vanishing on  $\fun\Phi_{\hat h}^{-1}s_{i}$ such that 
\[\hat h=\text{Exp}_{\nu_{\hat h}}\circ \hat f/G\circ\fun\Phi_{\hat h} \ .\]

\end{remark}

\

\subsection{Construction of a core family}

\

The following theorem gives  sufficient criteria for when a given family with a collection of marked point sections is a core family:

\begin{thm}\label{core criteria}
Let  $\hat f$ be a family in $\Msw$   together with a group $G$ of automorphisms, and a finite nonempty  set of disjoint sections $s_{i}\co \ex F(\hat f)\longrightarrow \ex C(\hat f)$ not intersecting the edges of the curves in $\ex C(\hat f)$, such that the following conditions are satisfied:

 \begin{enumerate}
  \item \label{crit1} For all curves $f$ in $\hat f$,  the action of $G$ on the set of maps $f\longrightarrow \hat f$ is free and transitive.
 \item\label{crit2} For all curves $f$ in $\hat f$,  $\totl{\ex C(f)}$ with the extra\footnote{Recall, from \cite[Definition 8.3]{iec}, that $\totl{\ex C(f)}$ is a nodal Riemann surface with a special point at the image of each end of $\ex C(\hat f)$. In addition to these special points, we have the extra marked points from the sections $s_i$. } marked points from $\{s_{i}\}$ has no nontrivial automorphisms.
 \item \label{crit3}The action of $G$ preserves the set of sections $\{s_{i}\}$ --- so there is some  action of $G$ as a permutation group on the set of indices $\{i\}$ such that for all $g\in G$ and $s_{i}$,
 \[s_{i}\circ g=g\circ s_{g(i)}\]
 where the action of $g$ is on $\ex F(\hat f)$, $\ex C(\hat f)$ or the set of indices $\{i\}$ as appropriate. 
 \item\label{crit4} There exists a neighborhood $U$ of the image of the section
\[s\co \ex F(\hat f)\longrightarrow \ex F(\hat f^{+n})\]
defined by the $n$ sections $\{s_{i}\}$ such that
\[\fun{ev}^{+n}_{\hat f}\co \ex F(\hat f^{+n})\longrightarrow \M^{+n}\times\sfp{\hat {\ex B}}{\ex B_{0}}n\]
restricted to $U$ is an injective equi-dimensional submersion.
\item \label{crit5}The tropical structure of $\hat f$ is universal; see   \cite[Definition 4.1, Theorem 3.1]{uts}.
\item\label{crit6} For any curve $f$ in $\hat f$, there are exactly $\abs G$  points $x$ in $\ex F(f^{+n})$ such that   $\fun{ev}^{+n}_{f}(x)$ is in the closure of the image of $\fun{ev}^{+n}_{\hat f}\circ s$.
\end{enumerate}

Then $(\hat f/G,\{s_{i}\})$ is a core family for an open substack $\mathcal O$ of $\Ms$ containing  every refinement of any curve $f$ in $\hat f$.

\end{thm}

\pf 
We must construct $\fun\Phi$ and $\fun\Phi^{+1}$; an overview of the construction is as follows. Given a family of curves $\hat h$ close to $\hat f/G$, we want to construct  $\fun\Phi(\hat h)$ as a $G$--bundle $\fun\pi_{\hat h}\co \hat h\bd G\longrightarrow \hat h$ with a $G$--equivariant map $\fun\Phi(\hat h)\co \ex F(\hat h\bd G)\longrightarrow \ex F(\hat f)$. We construct  $\ex F(\hat h\bd G)$ as
\begin{equation}\label{hbdfp}\ex F(\hat h\bd G):= \ex F(\hat h^{+n})\fp{\fun{ev}^{+n}_{\hat h}}{\fun{ev}^{+n}_{\hat f}\circ s}\ex F(\hat f)\ .\end{equation}
As a fiber product, $\ex F(\hat h\bd G)$ has natural maps $\fun\Phi(\hat h)$ and $s_{\hat h}$ as follows.
\[\begin{tikzcd}\ex F(\hat h\bd G)\rar{\fun\Phi(\hat h)}\dar{s_{\hat h}}& \ex F(\hat f)\dar{\fun{ev}^{+n}_{\hat f}\circ s}
\\ \ex F(\hat h^{+n})\rar{\fun{ev}^{+n}_{\hat h}} & \M^{+n}\times\sfp{\hat {\ex B}}{\ex B_{0}}n \end{tikzcd}\]
 Criterion \ref{crit3} implies that the map 
 \[ \fun{ev}_{\hat f}^{+n}\circ s\co \ex F(\hat f)\longrightarrow \M^{+n}\times\sfp{\hat {\ex B}}{\ex B_{0}}n\] is $G$--equivariant for some action of $G$ permuting labels of extra marked points in $\M^{+n}\times\sfp{\hat {\ex B}}{\ex B_{0}}n$. The map $\fun{ev}^{+n}_{\hat h}$ is $G$ equivariant with this action of $G$ permuting labels of the extra marked points, so both $s_{\hat h}$ and $\fun\Phi(\hat h)$ are $G$--equivariant.   
 The map $\ex F(\fun\pi_{\hat h})\co \ex F(\hat h\bd G)\longrightarrow \ex F(\hat h)$ is defined as the composition of $s_{\hat h}$ with the natural projection forgetting the $n$ extra ends.
 \[\begin{tikzcd} \ar[bend left]{rr}{\ex F(\fun\pi_{\hat h})} \ex F(\hat h\bd G)\rar{s_{\hat h}} &\ex F(\hat h^{+n})\rar &\ex F(\hat h)\end{tikzcd} \]
This projection $\ex F(\hat h^{+n})\longrightarrow \ex F(\hat h)$ is  $G$--invariant so $\ex F(\fun\pi_{\hat h})$ is $G$--invariant. Below, we define an open substack $\mathcal O_2$ containing $\hat f$ and every refinement of a curve in $\hat f$ such that, for $\hat h$ in $\mathcal O_2$,  $\ex F(\hat h\bd G)$ is an exploded manifold, and  $\ex F(\fun\pi_{\hat h})$ is a  $G$--bundle --- and therefore defines a $G$--bundle $\fun\pi_{\hat h}\co \hat h\bd G\longrightarrow \hat h$ as the  pullback of  $\hat h$ over $\ex F(\fun\pi_{\hat h})$. (See Claims \ref{ebundle} and \ref{Oopen}.) 

We must verify that $\hat h\mapsto \fun\Phi(\hat h)$ defines a functor. Given any map $\psi\co\hat g\longrightarrow \hat h$, the naturally induced map $\psi^{+n}\co \hat g^{+n}\longrightarrow \hat h^{+n} $ is $G$--equivariant, and induces a $G$--equivariant map of fiber products $\ex F(\psi\bd G)\co \ex F(\hat g\bd G)\longrightarrow \ex F(\hat h\bd G)$ and therefore induces a $G$--equivariant map $\psi\bd G\co \hat g\bd G\longrightarrow \hat h\bd G$ of $G$--bundles. Moreover this map satisfies 
\begin{equation}\label{comp}\fun\Phi(\hat h)\circ\ex F(\psi\bd G)=\fun\Phi(\hat g)\  \text{ and }\ s_{\hat h}\circ \ex F(\psi\bd G)=s_{\hat g}\ .\end{equation}
 Because $(\psi_1\circ\psi_2)^{+n}=\psi_1^{+n}\circ \psi_2^{+n}$, it follows that $(\psi_1\circ \psi_2)\bd G=\psi_1\bd G\circ \psi\bd G_2$, so the above defines a map $\fun\Phi\co \mathcal O_2\longrightarrow \ex F(\hat f)/G$. 

Note that $s_{\hat h}$ induces a section of $(\hat h\bd G)^{+n}\longrightarrow \hat h\bd G$ and hence $n$ sections $s_{\hat h,i}$ of $\ex C(\hat h\bd G)\longrightarrow \ex F(\hat h\bd G)$. We lift $\fun\Phi$ to $\fun\Phi^{+1}\co\mathcal O_2^{+1}\longrightarrow \ex C(\hat f)/G$ by defining 
\[\begin{tikzcd} \ex C(\hat h\bd G)\rar{\fun\Phi^{+1}(\hat h)}\dar&  \dar\ex C(\hat f)
\\ \ex F(\hat h\bd G)\rar{\fun\Phi(\hat h)}& \ex F(\hat f)\end{tikzcd}\]
as the unique fiberwise holomorphic map which is fiberwise degree $1$, and pulls back the sections $s_i$ to $s_{\hat h,i}$, so
\begin{equation}\label{phi+1cond}\fun\Phi^{+1}(\hat h)\circ s_{\hat h,i}=s_i\circ \fun\Phi(\hat h)\ .\end{equation}
We show in Claim \ref{phi+1exist} that $\fun\Phi^{+1}(\hat h)$ exists and is $\C\infty1$. Condition \ref{crit2} ensures that it is unique, because there is at most one degree 1 holomorphic map from any curve with $n$ marked points to a curve in $\ex C(\hat f)$ with its marked points from $s_{i}$.  The defining property of this map ensures  that it is $G$--equivariant, and (\ref{comp}) implies that
\[\fun\Phi^{+1}(\hat h)\circ \ex C(\psi\bd G)=\fun\Phi^{+1}(\hat g) \]
so we have defined a fiberwise holomorphic map $\fun\Phi^{+1}\co \mathcal O_2^{+1}\longrightarrow \ex C(\hat f)/G$. This completes the overview of the construction of $(\fun\Phi,\fun\Phi^{+1})$. 

Let us verify Condition \ref{cc1} of Definition \ref{core family}. The map $s\co \ex F(\hat f)\longrightarrow \ex F(\hat f^{+n})$ induces a section, $s'$ of the $G$--bundle $\ex F(\hat f^{+n})\fp{\fun{ev}^{+n}_{\hat f}}{\fun{ev}^{+n}_{\hat f}\circ s}\ex F(\hat f)=\ex F(\hat f\bd G)\longrightarrow \ex F(\hat f)$ such that $\fun\Phi(\hat f)\circ s'$ is the identity. It follows that $\fun\Phi_{\hat f}\co \ex F(\hat f)\longrightarrow \ex F(\hat f)/G$ is equivalent to the quotient map from Example \ref{quotient map}. Similarly, this section $s'$ induces a section $\ex C(s')$ of the $G$--bundle $\ex C(\fun\pi_{\hat f})\co \ex C(\hat f\bd G)\longrightarrow \ex C(\hat f)$, and $\fun\Phi^{+1}(\hat f)\circ \ex C(s')\co \ex C(\hat f)\longrightarrow \ex C(\hat f)$ is the identity because $s_{\hat f}\circ s'=s$, so the identity satisfies (\ref{phi+1cond}). Therefore, $\fun\Phi^{+1}_{\hat f}\co \ex C(\hat f)\longrightarrow \ex C(\hat f/G)$ represents the quotient map from Example \ref{quotient map}, and Condition \ref{cc1} of Definition \ref{core family} holds.

Condition \ref{cc3} of Definition \ref{core family} also holds. For a curve $g$, the map $\fun\Phi^{+1}(g)$ is determined by the domain $\ex C(g)$ of $g$ and the $\abs{G}$ points in $(\fun{ev}^{+n}_g)^{-1}(\fun{ev}^{+n}_{\hat f}\circ s(\ex F(\hat f))$. It follows that for a family $f_{t}$ as in Condition \ref{cc3}, $\fun\Phi^{+1}_{f_{t}}$ is independent of $t$, so Condition \ref{cc3} is satisfied.

So,  $\fun\Phi$ satisfies conditions \ref{cc1} and \ref{cc3}. By construction, we have that $\hat f\circ \fun\Phi^{+1}(\hat h)$ agrees with $\hat h\bd G$ on the image of the sections $s_{\hat h,i}$. So, if $S\subset \ex C(\hat f)/G$ indicates the image of the sections $\{s_i\}$, we have that $\hat f\circ \fun\Phi^{+1}_{\hat h}\rvert_{(\Phi^{+1}_{\hat h})^{-1}S}=\hat h\rvert_{(\Phi^{+1}_{\hat h})^{-1}S}$.  Moreover, Construction \ref{O2def} and Claim \ref{phi+1exist} ensure that Equation (\ref{tropical cond}) holds, so  Remark \ref{deformation remark} implies that, by restricting to a smaller neighborhood $\mathcal O$ of $\hat f/G$ still containing every refinement of every curve in $\hat f$, $(\fun\Phi,\fun\Phi^{+1})$ will also  satisfy Condition \ref{cc2} of Definition \ref{core family}. We have now verified that $(\fun\Phi,\fun\Phi^{+1})$ satisfies all conditions of Definition \ref{core family}, but we still need to fill in details of the construction. 

Let us construct a substack $\mathcal O_1$ comprised of families $\hat h$ for which the fiber product $\ex F(\hat h\bd G)$ from (\ref{hbdfp}) is a $\C\infty1$ exploded manifold. 
Define $\mathcal O_{1}\subset \Ms$ to be the substack consisting of curves $h$ such that $\fun{ev}^{+n}_h $ intersects the closure of the image of $\fun{ev}^{+n}_{\hat f}\circ s$ exactly $\abs G$ times, and each of these intersections is a transverse intersection with the image of $\fun{ev}^{+n}_{\hat f}\circ s$. Conditions \ref{crit4} and \ref{crit6} ensure that any curve in $\hat f$ is in $\mathcal O_1$, and similarly any refinement of a curve in $\hat f$ is also in $\mathcal O_1$. 

For any family of curves $\hat h$ in $\mathcal O_{1}$ lemmas 9.7 and 9.8 of \cite{iec} imply that 
\[\ex F(\hat h\bd G):= \ex F(\hat h^{+n})\fp{\fun{ev}^{+n}_{\hat h}}{\fun{ev}^{+n}_{\hat f}\circ s}\ex F(\hat f)\]
is a $\C\infty1$ exploded manifold with a $G$ action discussed above, and moreover  the map $\ex F(\fun\pi_{\hat h})\co \ex F(\hat h\bd G)\longrightarrow \ex F(\hat h)$ is $G$ invariant, and a degree $\abs{G}$ equi-dimensional submersion. It is not clear, however, that this map is a $G$--bundle, because it might not be etale, and it is also not clear that $\mathcal O_{1}$ is open. (To prove these properties, we will make a further assumption on our substack, and use that $\hat f$ has universal tropical structure). Nevertheless, we can still define a family of curves $\hat h\bd G$ by pulling back $\hat h$ over the map $\ex F(\fun\pi_{\hat h})\co \ex F(\hat h\bd G)\longrightarrow \ex F(\hat h)$, and  $\fun\Phi(\hat h)\co \ex F(\hat h\bd G)\longrightarrow \ex F(\hat f)$ is still a well-defined $G$-equivariant map, and Equation (\ref{comp}) is still valid. 

Let us construct $\fun\Phi^{+1}(\hat h)\co\ex C(\hat h\bd G)\longrightarrow  \ex C(\hat f)$ for $\hat h$ a family in $\mathcal O_1$. There is a natural map 
\begin{equation}\tilde{\fun\Phi}^{+1}(\hat h)\co \ex C(s_{\hat h}^{*}\hat h^{+n})\longrightarrow \ex C(s^{*}\hat f^{+n})\end{equation}
induced by thinking of $\ex C(s^{*}\hat f^{+n})$ as a fiber product as follows.
Consider the map 
\begin{equation}\label{tilde ev}\fun{\tilde{ev}}^{+n}_{\hat f}\co \ex C(\hat f^{+n})\longrightarrow \M^{+(n+1)}\times\sfp{\hat{\ex B}}{\ex B_{0}}n\end{equation} constructed from the composition of  $\fun{ev}^{+(n+1)}_{\hat f}$ with the projection map
\[\begin{tikzcd}[row sep=tiny] && \M^{+(n+1)}\times\sfp{\hat {\ex B}}{\ex B_{0}}{n+1}\ar{dd}
\\\ex C(\hat f^{+n})\ar[end anchor= west]{urr}{\fun{ev}^{+(n+1)}_{\hat f}}\arrow[end anchor= west]{drr}[swap]{\fun{\tilde{ev}}^{+n}_{\hat f}}
\\&& \M^{+(n+1)}\times\sfp{\hat{\ex B}}{\ex B_{0}}n\end{tikzcd}\]
 forgetting the last factor of $\hat {\ex B}$ --- so, on the second component, $\fun{\tilde{ev}}^{+n}_{\hat f}$ is  the composition of the projection $\ex C(\hat f^{+n})\longrightarrow \ex F(\hat f^{+n})$ with the map $\hat f^{+(n-1)}$.
 Criteria \ref{crit2} and \ref{crit4} imply that $\fun{\tilde {ev}}^{+n}_{\hat f}$ is an injective equidimensional submersion when restricted to a neighborhood of $\ex C(s^{*}\hat f^{+n})\subset \ex C(\hat f^{+n})$, and that, on this neighborhood,  the following is a pullback diagram of families of curves. 
\[\begin{tikzcd}[column sep=large]\ex C(\hat f^{+n})\dar\arrow{r}{\fun{\tilde{ev}}^{+n}_{\hat f}}&\M^{+(n+1)}\times \sfp{\hat {\ex B}}{\ex B_{0}}{n}\dar
\\\ex F(\hat f^{+n})\arrow{r}{\fun{ev}^{+n}_{\hat f}}&\M^{+n}\times \sfp{\hat {\ex B}}{\ex B_{0}}{n}
\end{tikzcd}\]
Therefore, $\ex C(s^{*}\hat f^{+n})$ is  the following fiber product:
 \[\begin{tikzcd}[column sep=large]\ex C(s^{*}\hat f^{+n})\arrow{r}\dar&\M^{+(n+1)}\times \sfp{\hat {\ex B}}{\ex B_{0}}{n}\dar
\\\ex F(\hat f)\arrow{r}{\fun{ev}^{+n}_{\hat f}\circ s}&\M^{+n}\times \sfp{\hat {\ex B}}{\ex B_{0}}{n}
\end{tikzcd}\]
 Accordingly, for any $\hat h$ in $\mathcal O_1$ there is a unique $\C\infty1$ map $\tilde {\fun\Phi}^{+1}(\hat h)$ such that the following diagram commutes
 \[\begin{tikzcd}[column sep=large]\ex C(s_{\hat h}^{*}\hat h^{+n}) \arrow[bend left]{rr}{\tilde{\fun\Phi}^{+1}(\hat h)}\rar\dar & \ex C(\hat h^{+n})\arrow[bend left]{rr}{\fun{\tilde{ev}}^{+n}_{\hat h}}\dar &\ex C(s^{*}\hat f^{+n})\arrow{r}{\fun{\tilde{ev}}^{+n}_{\hat f}}\dar&\M^{+(n+1)}\times \sfp{\hat {\ex B}}{\ex B_{0}}{n}\dar
\\ \ex F(\hat h\bd G)\rar{s_{\hat h}}\ar[bend right]{rr}{\fun\Phi(\hat h)} &  \ex F(\hat h^{+n}) \arrow[bend right]{rr}{\fun{ev}^{+n}_{\hat h}}&\ex F(\hat f)\arrow{r}{\fun{ev}^{+n}_{\hat f}\circ s}&\M^{+n}\times \sfp{\hat {\ex B}}{\ex B_{0}}{n}
\end{tikzcd}\]
Moreover, on each fiber, $\tilde{\fun\Phi}^{+1}(\hat h)$ is holomorphic and degree $1$, because $\fun{\tilde{ev}}^{+n}_{\hat h}$ is fiberwise holomorphic and degree $1$.

The following claim tells us that, for curves in $\mathcal O_1$, we can forget the extra ends of $\ex C(s^{*}\hat f^{+n})$ and $\ex C(s_{\hat h}^{*}\hat h^{+n})$ to define $\fun\Phi^{+1}(\hat h)$.

\begin{claim}\label{phi+1exist} For $\hat h$ in $\mathcal O_{1}$, there exists a $\C\infty1$ map $\fun\Phi^{+1}(\hat h)$ such that the following diagram commutes. 
\[\begin{tikzcd}
\ex C(s_{\hat h}^{*}\hat h^{+n})\arrow{r}{\tilde{\fun\Phi}^{+1}(\hat h)}\dar&\ex C(s^{*}\hat f^{+n})\dar
\\\ex C(\hat h\bd G)\arrow{r}{\fun\Phi^{+1}(\hat h)}\dar& \ex C(\hat f)\dar
\\ \ex F(\hat h\bd G)\arrow{r}{\fun\Phi(\hat h)} &\ex F(\hat f)\end{tikzcd}\]
Moreover, $\fun\Phi^{+1}(\hat h)$ is fibewise holomorphic and degree $1$, and pulls back the sections $s_{i}$ to $s_{\hat h,i}$, so
\begin{equation}\label{ta} s_{i}\circ \fun\Phi (\hat h)=\fun\Phi^{+1}(\hat h) \circ s_{\hat h,i} \end{equation}
\end{claim}

Recall, from Section \ref{+n section}, that the domain of $\hat f^{+n}$ has $n$ extra ends. To prove Claim \ref{phi+1exist}, note that the pullback under $\tilde {\fun\Phi}^{+1}$ of the $n$ extra ends of $\ex C(s^{*}\hat f^{+n})$ are the $n$ extra ends of $\ex C(s_{\hat h}^{*}\hat h^{+n})$. Moreover, the image of all these extra ends in $\ex C(\hat f)$ do not intersect, and  do not intersect the edges of curves. It follows that  the same property holds for the image of the extra ends in $\ex C(\hat h\bd G)$: if one was on an edge of a curve $g$ in  $\ex C(\hat h\bd G)$, (\ref{ta}) would imply that that $g$ is constant on that edge, so it would be impossible for $ \fun{ev}^{+n}_g$ to be transverse to $\fun{ev}^{+n}_{\hat f}\circ s$.  In the construction of $\hat f^{+1}$, below Equation (\ref{+1diag}), there are two steps: a refinement of $\ex C(\hat f)\times_{\ex F(\hat f)}\ex C(\hat f)$ at the image of the diagonal, and the subsequent explosion of the diagonal.  Because our new ends stay away from each other and edges of curves in $\ex C(\hat f)$ or $\ex C(\hat h)$, the refinement step will not affect our domains, so we only need to undo the explosion step. Locally, around one of our extra ends we have coordinates on $\ex C(s_{\hat h\bd G}^{*}\hat h^{+n})$ and $\ex C(s^{*}\hat f^{+n})$ in the form of open subsets of $\ex T^{1}_{[0,\infty)}\times U$ and $\et 1{[0,\infty)}\times V$, where $U$ and $V$ are open subsets of $\ex F(\hat h\bd G)$ and $\ex F(\hat f)$ respectively, and $\fun\Phi^{+1}(\hat h)$ is in the form $(\tilde z,u)\mapsto (\tilde z b,\fun\Phi(\hat h)(u))$ with $b$ a $\mathbb C^{*}$--valued and fiberwise holomorphic function on our coordinate chart.  Noting that the smooth part of $\et 1{[0,\infty)}$ is $\mathbb C$, the image of our coordinates in $\ex C(\hat h\bd G)$ and $\ex C(\hat f)$ are the corresponding open subsets of $\totl{\ex T^{1}_{[0,\infty)}}\times U$ and $\totl{\et 1{[0,\infty)}}\times V$ respectively. The holomorphic function $b$ factors through our chart on $\ex C(\hat h\bd G)$, and   in these coordinates $\fun\Phi^{+1}(\hat h)$ is then $(\totl{\tilde z},u)\mapsto (\totl{\tilde z}b ,\fun\Phi(\hat h)(u))$. This is the unique map fitting into the commutative diagram above, so it follows that $\fun\Phi^{+1}(\hat h)$ is globally well defined. Of course, $\fun\Phi^{+1}(\hat h)$ is fiberwise holmorphic and degree 1 because $\tilde{\fun\Phi}^{+1}(\hat h)$ is, and by construction $s_{i}\circ \fun\Phi(\hat h)=\fun\Phi^{+1}(\hat h)\circ s_{\hat h,i}$, because the $i$th extra end of $s_{\hat h}^{*}\hat h^{+n}$ is sent to the $i$th extra end of $s^{*}\hat f^{+n}$. So, Claim \ref{phi+1exist} is valid.

\begin{construction}\label{O2def} Define a substack $\mathcal O_2\subset \mathcal O_1$ consisting of curves $h$ such that there exists a family of curves $ h_t$   parametrized by $\mathbb R$ with domains $\ex C( h_t)$ identified with $\ex C(h)$ and satisfying \begin{itemize}
\item $h_0=h$,
\item  $s_{h}^*h_1^{+n}=s^*\hat f^{+n}\circ \tilde {\fun\Phi}^{+1}(h)$.
\end{itemize}
\end{construction}
Note that this implies that $\totb{s_{h}^*h^{+n}}=\totb{\hat f^{+n}}\circ \totb{\tilde {\fun\Phi}^{+1}(h)}$ for $h$ in $\mathcal O_2$, because $\totb {s_h^*h^{+n}_t}$ is independent of $t$. We will need this fact to prove that $\ex F(\hat h\bd G)\longrightarrow \ex F(\hat h)$ is etale. As $\mathcal O_1$ contains every refinement of every curve in $\hat f$, $\mathcal O_2$ contains all these refinements too. We will show in Claim \ref{Oopen} that $\mathcal O_2$ is open, and show in Claim \ref{ebundle} that $\fun\pi_{\hat h}\co \hat h\bd G\longrightarrow \hat h$ is a $G$--bundle for $\hat h$ in $\mathcal O_2$.

\begin{claim}\label{12 3}Let $\hat h$ be a family of curves containing a curve $h$ in $\mathcal O_2$.  Recall, from Construction \ref{O2def},  the deformation $h_t$ of $h$. This  deformation $h_t$   extends to a deformation $\hat h_t$ of $\hat h$. In other words, $\hat h_t$ is a smooth family of maps $\hat h_t\co \ex C(\hat h)\longrightarrow \hat{\ex B}$ parametrized by $\mathbb R$, and the restriction of $\hat h_t$ to $\ex C(h)$ is $h_t$. \end{claim}

 To extend the deformation $h_t$, we use the following fact.

\begin{claim}\label{section extension}Let  $V$ be a vector-bundle over the total space of a family  $\pi\co \hat{\ex A}\longrightarrow \ex A_0$. Then any $\C\infty1$ section of $V$ over a fiber of $\pi$ extends to a $\C\infty1$ section of $V$ over $\hat{\ex A}$. 
\end{claim}
To prove Claim \ref{section extension},
the key case is extending from a fiber over the central stratum of a coordinate chart $\et mP$. As $\pi$ is a family in the sense of  \cite[Definition 10.1]{iec}, a section of $V$ over such a fiber extends uniquely to a section over $\pi^{-1}\et m{P^\circ}$ where $P^\circ$ is the interior of the polytope $P$. The section can then locally extended using the operator $(1-\Delta_{\totb\pi^{-1}(P^\circ)})$, where $\totb\pi^{-1} (P^\circ)$ is the set of strata in the inverse image of $P^\circ$, and $\Delta$ is the operator from  \cite[Definition 7.3]{iec}. Local extensions can then be patched together using a partition of unity; so we can conclude Claim \ref{section extension}.
%

To prove Claim \ref{12 3}, we can obtain $h_t$ using the flow of a $\C\infty1$ vector-field on $\ex C( h)\times\mathbb R\times \hat{\ex B}$ in the form $(0,\partial_ t,v(p,t))$ where for all $(p,t)\in\ex C( h)\times \mathbb R$, $v(p,t)$ is a section of $T_{vert}\hat {\ex B}$ with 
\[v(p,t)\rvert_{ h_t(p)}=\frac{\partial h_t(p)}{\partial t}\ . \]
Using Claim \ref{section extension}, we can extend this vector-field to a $\C\infty1$ vector-field on $\ex C(\hat h)\times\mathbb R\times \hat{\ex B}$ in the form $(0,\partial_t,\tilde v(p,t))$ with $\tilde v(p,t)$ a section of $T_{vert}\hat{\ex B}$. The flow of this vector-field then gives an extension of $h_t$ to a  deformation $\hat h_t$ of $\hat h$ satisfying 
\[\tilde v(p,t)\rvert_{ \hat h\bd G_t(p)}=\frac{\partial \hat h_t(p)}{\partial t}\ .\]
 This completes the proof of Claim \ref{12 3}.

\begin{claim}\label{etale neighborhood} Let $h$ be in $\mathcal O_2$, and let $\hat h$ be a family of curves containing $h$. Consider the following fiber product.  
\[\begin{tikzcd}\ex F(\hat h^{+n})\fp {\fun{ev}^{+n}_{\hat h}} {\fun{ev}^{+n}_{\hat f}}\ex F(\hat f^{+n})\rar{\pi_2}\dar{\pi_1} &\ex F(\hat f^{+n})\dar{\fun{ev}^{+n}_{\hat f}}
\\ \ex F(\hat h^{+n})\rar{\fun{ev}^{+n}_{\hat h}} & \M^{+n}\times\sfp{\hat {\ex B}}{\ex B_{0}}n   \end{tikzcd}\]
There exists an open neighborhood $O$ of 
\[\ex F(h\bd G)=\ex F( h^{+n})\fp{\fun{ev}^{+n}_ h}{\fun{ev}^{+n}_{\hat f}\circ s}\ex F(\hat  f)\subset \ex F(\hat h^{+n})\fp{\fun{ev}^{+n}_{\hat h}}{\fun{ev}^{+n}_{\hat f}}\ex F(\hat f^{+n})\] 
such that the restriction of $\pi_1$ to $O$ is an isomorphism onto an open subset. 
\end{claim}
Recall from Criterion \ref{crit4} that $\fun{ev}^{+n}_{\hat f}$ is an injective equi-dimenisional submersion on a neighborhood $U$ of $s(\hat f)$; so, $\pi_1\co\ex F(\hat h^{+n})\fp {\fun{ev}^{+n}_{\hat h}} {\fun{ev}^{+n}_{\hat f}}U\longrightarrow \ex F(\hat h^{+n})$ is also an injective equidimensional submersion. Remark \ref{submersion etale} tells us that this injective equidimensional submersion $\pi_1$ will be an isomorphism in a neighborhood of each point $x$ where its tropical structure $\mathcal P\pi_1(x)$ is an isomorphism.  From Claim \ref{12 3}, we have a family of curves $\hat h_1\co \ex C(\hat h)\longrightarrow \hat{\ex B}$ whose restriction to $\ex C(h)$ is the  curve $h_1$ from Construction \ref{O2def} with  $h_1\bd G=\hat f\circ \fun\Phi^{+1}(h)$. The $\abs G$ different points $x$ in $\ex F(h\bd G)$ correspond to  $\abs G$ different maps $\fun\Phi^{+1}(h)_x\co \ex C(h)\longrightarrow \ex C(\hat f)$ with $h_1=\hat f\circ \fun\Phi^{+1}(h)_x$. By assumption, $\hat f$ has universal tropical structure so, after replacing $\hat h$ by an open subfamily if necessary,  Lemma \ref{uts map cor} implies that each of these maps $\fun\Phi^{+1}(h)_x$ extend  to maps $\fun\Phi^{+1}_x\co \ex C(\hat h)\longrightarrow \ex C(\hat f)$ such that $\totb {\hat f}\circ \totb{\fun\Phi^{+1}_x}=\totb{\hat h_1}=\totb{\hat h}$. (This is the point where we needed the deformation from Construction \ref{O2def}.) Without loss of generality, we can assume that $\hat h$ is connected, and so $\fun\Phi^{+1}_x$ is fiberwise genus-preserving,  degree one, and sends ends to ends. It follows that the induced maps $\fun\Phi^{+n}_x\co \ex F(\hat h^{+n})\longrightarrow \ex F(\hat f^{+n})$ satisfy
\begin{equation}\label{lt}\totb{\fun{ev}^{+n}_{\hat f}}\circ \totb{\fun\Phi^{+n}_n}=\totb{\fun{ev}^{+n}_{\hat h}} \ .\end{equation} 
This map $\fun\Phi_x^{+n}$ induces a section $l_x$ of $\pi_1$ defined on a neighborhood of $\pi_1(x)$ with $l_x(\pi_1(x))=x$.
Equation (\ref{lt}) implies that $\mathcal P\pi_1(x)$ is inverse to $\mathcal Pl_x(\pi_1(x))$, so Remark \ref{submersion etale}  implies $\pi_1$ is etale in a neighborhood of $x$.  Claim \ref{etale neighborhood} follows.

\begin{claim}\label{ebundle}If $h$ is a curve in $\mathcal O_2$ and   $\hat h$ is a family of curves containing $h$, then $\ex F(\pi_{\hat h})\co \ex F(\hat h\bd G)\longrightarrow \ex F(\hat f)$ is etale in a neighborhood of $\ex F(h\bd G)\subset \ex F(\hat h\bd G)$. In particular, if $\hat h$ is  in $\mathcal O_2$, then $\ex F(\pi_{\hat h})\co \ex F(\hat h\bd G)\longrightarrow \ex F(\hat h)$ is etale, and hence a $G$--bundle.
\end{claim}

Claim \ref{ebundle} follows from Claim \ref{etale neighborhood}. Because $ h$ is in $\mathcal O_1$, we have that  $\fun{ev}^{+n}_h $ intersects $\fun{ev}^{+n}_{\hat f}\circ s(\ex F(\hat f))$ transversely exactly $\abs{G}$--times, so $\ex F(\hat h\bd G)\longrightarrow \ex F(\hat f)$  has derivative an isomorphism at $\ex F(h\bd G)$. As all the marked point sections defining $s\co \ex F(\hat f)\longrightarrow \ex F(\hat f^{+n})$ are non-intersecting and in smooth components of curves, it follows that the tropical structure of $s$ is an isomorphism $\mathcal Ps(x)$ at each point $x$. It follows the the tropical structure of $s_{\hat h}\co \ex F(\hat h\bd G)\longrightarrow \ex F(\hat h^{+n})$ at any point in $\ex F(h\bd G)$  coincides with the tropical structure of the map $\pi_2$, which is an isomorphism, by Claim \ref{etale neighborhood}. As with any family, the tropical structure of  the projection $\pi\co \ex F(\hat h^{+n})\longrightarrow \ex F(\hat h)$ at any point is either an isomorphism, or a projection from a higher dimensional polytope.  As $\ex F(\pi_{\hat h})$ has bijective derivative at each point in $\ex F(h\bd G)$, its tropical structure can't be a nontrivial projection. As $\ex F(\pi_{\hat h})=\pi\circ  s_{\hat h}$,   it follows that the tropical structure of $\ex F(\pi_{\hat h})$ is also an isomorphism. Remark \ref{submersion etale} then implies that $\ex F(\pi_{\hat h})$ is etale in a neighborhood of $\ex F(h\bd G)\subset \ex F(\hat h\bd G)$. If  $\hat h$ is in $\mathcal O_2$, it follows that $ \ex F(\pi_{\hat h})$ is etale. This completes the proof of Claim \ref{ebundle}.

\begin{claim}\label{Oopen} $\mathcal O_2$ is an open substack of $\Msw$.
\end{claim}

Let $h$ be in $\mathcal O_2$ and let $\hat h$ be a family of curves containing $h$. To prove Claim \ref{Oopen}, we must show that a neighborhood of $h$ within $\hat h$ is contained in $\mathcal O_2$. 
 Around any point in $\ex F( h\bd G)\subset\ex F(\hat h\bd G)$, Claim \ref{ebundle} gives that $\ex F(\hat h\bd G)\longrightarrow \ex F(\hat h)$ is etale. Moreover, for all curves $g$ in a neighborhood of $h$ in $\hat h$, $\fun{ev}^{+n}_g$ remains transverse to $\fun{ev}^{+n}_{\hat f}\circ s$ at $\abs G$ points. Moreover, on an open neighborhood, there are no further intersections with the closure of the image of $\fun{ev}^{+n}_{\hat f}\circ s$. If follows that a neighborhood of $h$ in $\hat h$ is in $\mathcal O_{1}$. 
 
  It remains to show that a neighborhood of $h$ in $\hat h$ is in $\mathcal O_{2}$. Without loss of generality, we can assume that $\hat h$ is in $\mathcal O_{1}$. As $\totb{s_{h}^*h^{+n}}=\totb{\hat f^{+n}}\circ \totb{\tilde {\fun\Phi}^{+1}(h)}$ and $\fun{ev}^{+n}_{\hat h}\circ s_{\hat h}=\fun{ev}^{+n}_{\hat f}\circ s\circ \fun\Phi(\hat h)$, it follows that, after restricting $\hat h$ to an open neighborhood of $h$ if necessary, 
  \begin{equation}\label{te}\totb{s_{\hat h}^*\hat h^{+n}}=\totb{\hat f^{+n}}\circ \totb{\tilde {\fun\Phi}^{+1}(\hat h)} \ .\end{equation}
This is because  a small enough open neighborhood of $\ex C(s_{h}^{*}h^{+n})\subset \ex C(s_{\hat h}^{*}\hat h^{+n})$ has tropical part consisting of the closure of the strata containing $\totb{\ex C(s_{h}^{*}h^{+n})}$. On these strata,  the two tropical maps on either side of equation (\ref{te}) are determined by their restriction to  $\totb{\ex C(s_{h}^{*}h^{+n})}$ and to the strata corresponding to the $n$ extra ends, where these maps agree.

   Consider the extension $\hat h_{t}$ of the deformation $h_{t}$ from Claim \ref{12 3}. The tropical part of $\hat h_{t}$ is independent of $t$, so we also have $\totb{s_{\hat h}^*\hat h_{1}^{+n}}=\totb{\hat f^{+n}}\circ \totb{\tilde {\fun\Phi}^{+1}(\hat h)}$. In particular, given any metric\footnote{See Remark 6.6 of \cite{iec}.} on $\hat {\ex B}$, the geodesic distance between  $s_{\hat h}^{*}\hat h_{1}^{+n}$ and $\hat f^{+n}\circ \tilde {\fun\Phi}^{+1}(\hat h)$ is a continuous function on the domain, $\ex C(\hat h\bd G)$. As $ s_{h}^{*}h_{1}^{+n}=\hat f\circ \tilde {\fun\Phi}^{+1}(h)$, it follows that this distance is small on a neighborhood of $\ex C(h\bd G)$, and in particular, for all curves $g$ in some neighborhood of $h$, our deformation $\hat h_{t}$ restricted to $\ex C(g)$ can be followed by a deformation to a curve $g_{1}$ with $ s_{g}^{*}g_{1}^{+n}=\hat f\circ \tilde {\fun\Phi}^{+1}(g)$. In particular, this neighborhood is in $\mathcal O_{2}$, and Claim \ref{Oopen} follows.

This completes the proof of Theorem \ref{core criteria}.

 \stop

\

The following proposition constructs a core family containing a given stable  curve with at least one smooth component (so its domain is not $\ex T$). Note that although this proposition constructs a core family on the ambient moduli stack of undecorated curves,   Remark \ref{decorated core family} implies that core families also exist around any  curve in the ambient moduli stack of  decorated stable curves,  $\dmsw$ (so long as this curve has at least $1$ smooth component).
 
\begin{prop}
\label{smooth model family}

Given a curve  $f $ in $\Msw$ with a domain not equal to $\ex T$, and a collection of marked points $\{p_{j}\}$ in the interior of  the smooth components of $\ex C(f)$, 
  there exists a $\C\infty1$ core family $(\hat f/G,\{s_{i}\})$,
   with $\hat f$ a family containing $f$ so that the restriction of $\{s_{i}\}$ to $f$ contains the given marked points $\{p_{j}\}$.

\end{prop}

\pf The proof proceeds by constructing $(\hat f/G,\{s_i\})$ satisfying the requirements of Theorem \ref{core criteria}.
Theorem \ref{G uts family} allows us to construct a family of curves $\hat f'$ containing $f$,  with universal tropical structure, and   with a finite group $G'$ of automorphisms  acting freely and transitively on the set of maps of $f$ into $\hat f'$. Moreover, Theorem \ref{G uts family} gives that  the smooth part\footnote{Definition 3.14 of \cite{iec}.} of the stratum of $\ex F(\hat f')$ containing $f$ consists of a single point, and $G'$ is a  subgroup of the group $G$ of automorphisms of the smooth part $\totl{f}$ of $f$. We can also assume that there is only one nodal curve in $\totl{\hat f'}$  isomorphic to $\totl{f}$.
 
  Let $\hat f_{0}$ be the quotient of $G\times \hat f'$ by the equivalence relation  $(g,\hat f')\simeq (gh^{-1},h*\hat f')$ for any $h\in G'$. In other words, $\hat f_{0}$ is $\abs {G/G'}$ disjoint copies of $\hat f'$. $G$ acts as a group of automorphisms  on $\hat f_{0}$ by multiplying the $G$ factor on the left; this $G$--action is free and transitive on the set of maps of $f$ into $\hat f_{0}$. Moreover,  this new family $\hat f_{0}$ also has universal tropical structure because any open subset of a family with universal tropical structure also has  universal tropical structure.

Choose an inclusion of $f$ into $\hat f_{0}$, then choose a $G$--invariant collection of $n$ non-intersecting sections $s_{i}$ of $\ex C(\hat f_{0})\longrightarrow \ex F(\hat f_0)$ such that the intersection of these sections $s_{i}$ with $\ex C(f)\subset \ex C(\hat f_{0})$ corresponds to a set of marked points $\{p_{i}\}$ satisfying  the following conditions:

\begin{enumerate}\item These marked points $\{p_{i}\}$ contain the set of marked points given in the statement of this proposition.
\item\label{ps0} Each $p_{i}$ is in a smooth component of $\ex C(f)$.
\item\label{psg} The action of $G$ as the automorphism group of $\totl f$ permutes the marked points $p_{i}$, and the action of $G$ on the set of sections $s_{i}$ is compatible in the following sense:  if $g$ as an automorphism of $\totl{\ex C(f)}$ sends $p_{i}$ to $p_{j}$, then the action of $g$ on $\ex F(\hat f_{0})$ followed by  $s_{j}$ is equal to $s_{i}$ followed by  the action of $g$ on $\ex C(\hat f_{0})$. Representing the various actions of $g$ simply as $g$, we can write this condition as
\[g*p_{i}=p_{j}\ \ \ \ 
\text{ implies that}\ \ \ 
s_{j}\circ g=g\circ s_{i}\ .\]
\item\label{ps1} $\ex C( f)$, together with the sub collection of points $p_{i}$ such that $d\totl f$ is injective at $p_{i}$, is stable.
\item\label{ps2} The nodal Riemann surface $\totl{\ex C( f)}$ with the extra marked points $\{p_{i}\}$ has no automorphisms.
\item\label{ps3} There is at least one marked point on each smooth component of $\ex C(f)$.
\end{enumerate}

Items \ref{ps0}, \ref{ps1}, \ref{ps2} and \ref{ps3} above remain true for the marked points obtained by intersecting the image of $\{s_{i}\}$ with $\ex C(f')$ for $f'$ in a neighborhood of $f$.

Following the notation of the proof of Theorem \ref{core criteria}, let $s\co \ex F(\hat f_{0})\longrightarrow \ex F(\hat f_{0}^{+n})$ be the map determined by the $n$ sections, $\{s_{i}\}$ ---  so the domain of the family of curves $s^{*}\hat f_{0}^{+n}$ is $\ex C(\hat f_{0})$ with extra ends at the images of $s_{i}$.

\begin{claim}\label{rigidify} If a curve $h$ in $f^{+n}$ has smooth part $\totl h$  isomorphic to the smooth part of a curve in $s^{*}\hat f_{0}^{+n}$, then $h$ is actually isomorphic to a curve in $s^{*}\hat f_{0}^{+n}$.
\end{claim}

To prove Claim \ref{rigidify}, forget the $n$ extra ends of $ h\in \hat f^{+n}$, and consider the corresponding isomorphism  of $\totl{f}$ with a curve in $\totl{\hat f_{0}}$.   As there was only one curve in $\totl{\hat f'}$ isomorphic to $\totl f$, this isomorphism must decompose as an automorphism $g_{1}$ of $\totl f$ followed by our chosen inclusion $\totl f\longrightarrow \totl{\hat f_{0}}$, followed by the action of some $g_{2}\in G$. Item \ref{psg} implies that the pullback of $\totl{s_{i}}$ under such an isomorphism is equal to the pullback of $\totl{s_{i}}$ under our chosen inclusion $\totl{f}\longrightarrow \totl{\hat f_{0}}$ followed by the action of $g_{2}g_{1}$. Therefore, the extra ends of  $h$ are located at the pullback of the sections $s_{i}$ via our chosen inclusion of $f$ into $\hat f_{0}$ followed by the action of an element of $G$. So, $h$ is isomorphic to a curve in $s^{*}\hat f_{0}^{+n}$, and Claim \ref{rigidify} is true.

\

Item \ref{ps2} implies that each of the $\abs{G}$ inclusions of $f$ into $\hat f_{0}$ corresponds to a different intersection of $\fun{ev}^{+n}_f$ with the image of $\fun{ev}^{+n}_{\hat f_{0}}\circ s$. However, some of these intersections may have the same smooth part. Claim \ref{rigidify} implies that no other curve in $f^{+n}$ has the same smooth part as a curve in $s^{*}\hat f_{0}^{+n}$. We can therefore add extra marked points satisfying the above properties until $\totl{\fun{ev}^{+n}_{f}}$ has precisely $\abs G$ intersections with the image of $\totl{\fun{ev}^{+n}_{\hat f_{0}}\circ s}$.

Let us verify that $\fun{ev}^{+n}_{\hat f_{0}}\circ s$ has injective derivative. As the smooth part of the strata of $\ex F(\hat f_{0})$ containing $f$ is $0$--dimensional,  $T_{f}\ex F(\hat f_{0})$ consists of only $\mathbb R$--nil vectors, so the derivative of $\fun{ev}^{+n}_{\hat f_{0}}\circ s$ applied to $T_{f}\ex F(\hat f_{0})$ is determined by integral vectors. Let us determine what happens to integral vectors by considering  the tropical structure of the map $\fun{ev}^{+n}_{\hat f_{0}}\circ  s$
at the curve $f$; \cite[Definition 4.4]{iec}. This is some integral-affine map $P_{u}\longrightarrow P'$, where the polytope $P_u$ is the tropical structure of $\ex F(\hat f_0)$ at $f$, and therefore parametrises a moduli space of tropical curves containing $\totb f$. The polytope $P'$ is part of the tropical structure of $ \M^{+n}\times\sfp{\hat {\ex B}}{\ex B_{0}}n $; it  records the image of $\hat f_{0}\circ s_{i}$ in $\totb{\hat {\ex B}}$, and because of item \ref{ps1} above, it also records the length of the internal edges of curves in $\totb{\hat f_{0}}$. 
Because $\hat f_{0}$ has universal tropical structure, Remark 3.3 of \cite{uts} implies that the map $P_{u}\longrightarrow P'$ is injective and sends integral vectors on $P_{u}$ to a full sublattice of the integral vectors on $P'$. In other words, $\fun{ev}^{+n}_{\hat f_{0}}\circ s$ sends integral vectors in its domain to a full sublattice of the integral vectors in its codomain. Therefore,  if $\hat f_{0}$ is chosen small enough,  $\fun{ev}^{+n}_{\hat f_{0}}\circ s$ is injective and has injective derivative.

Now let us prove that $\fun{ev}^{+n}_{\hat f_{0}}$ has injective deriative on a neighborhood of the image of $s$.
Item \ref{ps1},  above, ensures that, at the point  $s(f)\in\ex F(\hat f^{+n}_{0})$, the derivative of the smooth part of $\fun{ev}^{+n}_{f}$
is injective. As each of the $p_{i}$ are distinct and on smooth components of $\ex C(f)$, a neighborhood of $s(f)$ in $\ex F(f^{+n})$ is isomorphic to $\mathbb R^{2n}$. So, the derivative of $\fun{ev}^{+n}_{f}$ at this point is injective and has no nontrivial $\mathbb R$--nil vectors in its image. As the derivative of $\fun{ev}^{+n}_{\hat f_{0}}\circ s$ is injective and has only $\mathbb R$--nil vectors in its image, $\fun{ev}^{+n}(\hat f_0)$ has injective derivative at the point $s(\ex F(f))$. Therefore, so long as $\hat f_{0}$ is chosen small enough,   
$\fun{ev}^{+n}_{\hat f_{0}}$
 is injective and has injective derivative, when restricted to a neighborhood of the image of $s$.

To satisfy Criterion \ref{crit4} of Theorem \ref{core criteria}, we  now extend $(\hat f_{0}, \{s_{i}\})$ to $(\hat f,\{s_{i}\})$  so that $\fun{ev}^{+n}_{\hat f}$ has bijective derivative on a neighborhood of the image of $s$. The action of $G$ permutes the $n$ sections $\{s_{i}\}$ --- so there is an action of $G$ on $\hat f_{0}^{+n}$ lifting the action of $G$ on $\hat f_{0}$ and permuting the  extra end labels  so that $s\co \ex F(\hat f_{0})\longrightarrow\ex F(\hat f_{0}^{+n})$ is  $G$--equivariant. There is a corresponding  $G$--action  on $\M^{+n}\times \sfp{\hat{\ex B}}{\ex B_{0}}n$ so that \[\fun{ev}^{+n}_{\hat f_{0}}\co \ex F(\hat f_{0}^{+n})\longrightarrow \M\times \sfp{\hat{\ex B}}{\ex B_{0}}n\] is also $G$--equivariant. Choose a  $G$--invariant metric\footnote{See Remark 6.6 of \cite{iec} for a discussion of metrics on exploded manifolds. In the topology induced by a metric, an exploded manifold is a disjoint union of manifolds. Because metrics can be averaged, there is no obstruction to constructing a metric on an exploded orbifold using partitions of unity, or Proposition 2.3 of \cite{vfc}. In general, $U$ will not be an open neighborhood in the usual topology on $\M^{+n}\times \sfp{\hat{\ex B}}{\ex B_{0}}n$, as the topology induced by a metric on an exploded manifold is finer than the usual topology.} on $\M^{+n}\times \sfp{\hat{\ex B}}{\ex B_{0}}n$, and let $U$ be some small $G$--invariant tubular neighborhood  of the
image of $\fun{ev}^{+n}_{\hat f_{0}}$ restricted to a neighborhood of the section $s$.  Let $V$ be the restriction of this disk bundle $U$ to the image of $\fun{ev}^{+n}_{\hat f_{0}}\circ s$, (so $V$ has codimension $2n$ in $U$). As $\fun{ev}^{+n}_{\hat f_{0}}$ and $\fun{ev}^{+n}_{\hat f_{0}}\circ s$ are $G$--equivariant and the metric used to define our tubular neighborhood is $G$--invariant, $V$ is $G$--invariant. Define $\ex F(\hat f)$ to be $V$, and define $\ex C(s^{*}\hat f)$ by the following pullback diagram. 
\[\begin{tikzcd}\ex C(s^{*}\hat f)\arrow{r}\dar &\M^{+(n+1)}\times \sfp{\hat{\ex B}}{\ex B_{0}}n\dar
\\ \ex F(\hat f)\arrow[hook]{r} &\M^{+n}\times \sfp{\hat{\ex B}}{\ex B_{0}}n\end{tikzcd}\]
The action of $G$ on $\M^{+n}\times \sfp{\hat{\ex B}}{\ex B_{0}}n$ lifts to an action of $G$ on $\M^{+(n+1)}\times \sfp{\hat{\ex B}}{\ex B_{0}}n$  permuting the same end labels.
Therefore, there is an action of $G$ on the family of curves $\ex C(s^{*}\hat f)\longrightarrow \ex F(\hat f)$ making the above diagram $G$--equivariant. By removing the extra edges in $\ex C(s^{*}\hat f)$, and remembering their location with sections $s_{i}$, we get a $G$--invariant family of curves $\ex C(\hat f)\longrightarrow \ex F(\hat f)$ with a $G$--invariant set of sections $\{s_{i}\}$.

There is a $G$--equivariant inclusion of $\ex C(\hat f_{0})$ as a subfamily of $\ex C(\hat f)$ via the diagram
\[\begin{tikzcd}[column sep=huge]\ex C(\hat f_0)\arrow{r}\dar &\ex C(\hat f)\dar
\\ \ex F(\hat f_{0})\arrow{r}{\fun{ev}^{+n}_{\hat f_{0}}\circ s}&V=\ex F(\hat f)\end{tikzcd}\]
and, the restriction of the sections $s_{i}$ of $\ex C(\hat f)\longrightarrow \ex F(\hat f)$ are our original sections $s_{i}$. As $\ex F(\hat f)$ is  a disk bundle over $\ex F(\hat f_{0})$, $\ex C(\hat f)$ is  a disk bundle over $\ex C(\hat f_{0})$. We can therefore extend the map $\hat f_{0}$ to a map  $\hat f$
\[\begin{tikzcd}\ex C(\hat f)\arrow{r}{\hat f}\ar{d}&\hat{\ex B}\arrow{d}
\\ \ex F(\hat f)\arrow{r} &\ex B_{0}\end{tikzcd}\]
so that $\fun{ev}^{+n}_{\hat f}\circ s\co \ex F(\hat f)\longrightarrow \M\times \sfp{\hat {\ex B}}{\ex B_{0}}n$ is the identity inclusion of $V$. 
As this condition is $G$--equivariant and the original map $\hat f_{0}$ is $G$--invariant, we can construct our map $\hat f$ to be $G$--invariant.

As $\hat f$ is just the extension of $\hat f_{0}$ to a disk bundle,  Lemma 4.4 of \cite{uts} implies that $\hat f$ has universal tropical structure. Therefore $\hat f$ satisfies condition \ref{crit5} of Theorem \ref{core criteria}.

By construction,  $\fun{ev}^{+n}_{\hat f}\circ s$ is injective, and  the derivative of $\fun{ev}^{+n}_{\hat f}$ at $s(f)$ is an isomorphism. Therefore, by restricting $\hat f$ to a smaller $G$--invariant neighborhood of $f$ if necessary, $\fun{ev}^{+n}_{\hat f}$ is an injective equidimensional submersion in a neighborhood of the image of $s$. In other words, $\hat f$ satisfies condition \ref{crit4} of Theorem \ref{core criteria}.

Condition \ref{crit3}  of Theorem \ref{core criteria} is satisfied because $s$ is $G$--equivariant. Condition \ref{crit2} is satisfied because of item \ref{ps2} from the construction of $\{s_{i}\}$.

 We shall now verify Condition \ref{crit6} of Theorem \ref{core criteria}. We have already established that there are precisely $\abs G$ intersections of $\totl{\fun{ev}^{+n}_{f}}$ with $\totl{\fun{ev}^{+n}_{\hat f_{0}}\circ s}$ corresponding to the $\abs G$ maps of $f$ into $\hat f_{0}$.  The corresponding intersections of $\fun{ev}^{+n}_{f} $ with the image of $\fun{ev}^{+n}_{\hat f}\circ s$ are transverse (and 0-dimensional), and the same holds for the intersections of $\totl{\fun{ev}^{+n}_{f}}$ with the image of $\totl{\fun{ev}^{+n}_{\hat f}\circ s}$, because the smooth part of the derivative of $\fun{ev}^{+n}_{f}$ at the intersection points is still injective.  By restricting  $\hat f$ to a smaller $G$--equivariant neighborhood of $\hat f_{0}$ if necessary, we  therefore get  that there are precisely $\abs G$ intersections of $\totl{\fun{ev}^{+n}_{f}}$ with the closure of the image of $\totl{\fun{ev}^{+n}_{\hat f}\circ s}$, and that for any $f'$ sufficiently close to $f$ in $\hat f$, there are also  $\abs G$ transverse intersections  of $\fun{ev}^{+n}_{f'}$ with $\fun{ev}^{+n}_{\hat f}\circ s$. Furthermore, we can ensure that there are no further intersections of $\totl{\fun{ev}^{+n}_{f'}}$ with the closure of $\totl{\fun{ev}^{+n}_{\hat f}\circ s}$. Therefore, by further reducing the size of $\hat f$, we can ensure that for all $f'$ in $\hat f$,  $\fun{ev}^{+n}_{f'}$ intersects the closure of the image of $\fun{ev}^{+n}_{\hat f}\circ s$ exactly $\abs G$ times. In other words,  Condition \ref{crit6}  of Theorem \ref{core criteria} holds.

 To verify Condition \ref{crit1} of Theorem \ref{core criteria}, we must verify that for all $f'$ in $\hat f$, the action of $G$ on the set of maps $f'\longrightarrow f$ is free and transitive. As we have already shown that there are are precisely $\abs G$ intersections of $\fun{ev}^{+n}_{f'}$ with $\fun{ev}^{+n}_{\hat f}\circ s$, there are at most $\abs G$ maps $f'\longrightarrow f$, and it remains to verify that the action of $G$ on the set of these maps is free. This is easy, because the action of $G$ does not fix the image of $\totl{\ex C(f)}$ in $\totl{\ex C(\hat f)}$ under the inclusion $f\longrightarrow \hat f$, so the action of $G$ can not fix any curve in a $G$--equivariant neighborhood of $f$ within $\hat f$. Therefore, Condition \ref{crit1} of Theorem \ref{core criteria} will hold if we restrict $\hat f$ to a small enough $G$--invariant open neighborhood of $f$.

We have now verified that $(\hat f/G,\{s_{i}\})$ satisfies all the conditions of Theorem  \ref{core criteria}, so $(\hat f/G,\{s_{i}\})$ is a core family of curves for some open neighborhood $\mathcal O$ of $\hat f/G$.

\stop

\subsection{Topology of $\Msw$}\label{bump section}

\

\begin{lemma}\label{topological sequential convergence}Suppose that a sequence of curves $f_i$ in $\dmsw$ is eventually contained in every open substack of $\dmsw$ that contains $f$. Then $f_i$ converges in $\C\infty1$ to $f$ in the sense of  \cite[Definition 11.4]{iec}. \end{lemma}

\pf

 The case when $\ex C(f)=\ex T$ is easy, because curves with domain $\ex T$ form an easily analysable open substack of $\dmsw$; see Section \ref{T section}. Accordingly, we assume that $\ex C(f)\neq \ex T$, and construct a core family $(\hat f/G,\{s_i\})$ for an open substack $\mathcal O\subset \dmsw$ using Proposition \ref{smooth model family} and Remark \ref{decorated core family}.

By assumption $ f_i$  is eventually contained in $\mathcal O$. Moreover, because the inverse image of any open substack of  $\ex F(\hat f)/G$ under the map $\fun\Phi\co \mathcal O\longrightarrow \ex F(\hat f)/G$ is open, (see Example \ref{open inverse image}), the image of $\fun\Phi(f_i)$ in $\ex F(\hat f/G)$ converges to the $G$--orbit of $f\in \ex F(\hat f)$. So, by a judicious resolution of the $G$--fold ambiguity of $\fun\Phi(f_i)$ and $\fun\Phi^{+1}(f_i)$ from Definition \ref{core family}, we can pick fiberwise-holomorphic maps 
\[\begin{tikzcd} \ex C(f_i)\rar{\fun\Phi^{+1}(f_i)}\dar & \ex C(\hat f)\dar
\\ \ex F(f_i)\rar & \ex F(\hat f)\end{tikzcd}\] whose image converge to $\ex C(f)\subset \ex C(\hat f)$. Given any family $\hat h$ and fiberwise holomorphic map $x\co \ex C(\hat h)\longrightarrow \ex C(\hat f)$,  there is an open subset of $\ex F(\hat h)$ comprised of points $p\in\ex F(\hat f)$ such that $x$ restricted to the fiber over $p$  an isomorphism onto a fiber of $\ex C$. This is because the degree of these maps is constant, and a degree--1 holomorphic map $x$ between exploded curves is an isomorphism if its smooth part $\totl{x}$ does not contract any components.   Accordingly, there is an open subset of $\mathcal O$ such that the maps $\fun\Phi^{+1}(\hat h)$ are isomorphisms restricted to each fiber, so our maps $\fun\Phi^{+1}(f_i)$ are eventually isomorphisms onto fibers.

Choose a family of metrics on $\hat{\ex B}\longrightarrow \ex B_0$, so we exponentiate vectors.  Remark \ref{deformation remark} implies that there is an open substack of $\mathcal O$ comprised of curves $h$ for which $h=\text{Exp}_{\psi_h}\circ \hat f/G\circ \fun\Phi_h^{+1}$ where $\psi_h$ is a section of $(\hat f/G\circ \fun\Phi_h^{+1})^*T_{vert}\hat{\ex B}$ vanishing on the image of the sections $s_i$, and  $\psi_h$ is small in any appropriately continuous norm.  Moreover, we can choose a sequence of such norms so that convergence in these norms implies $\C\infty1$ convergence. As argued in the proof of  Lemma \ref{sequence family} we can then extend each $\psi_{f_i}$ to a section $\psi_i$ of $\hat f^*T_{vert} \hat{\ex B}$ such that $\psi_i$ converges in $\C\infty 1$ to $0$, and therefore $f_i$ converges to $f$ in the sense of Definition 11.4 of \cite{iec}.

\stop

We need the following technical lemma to prove that $\Msw$ has topology pulled back from a Hausdorff topological space, $\totl{ \Msw_{top}}$.

\begin{lemma}\label{double family} Let $\hat f_1$  and $\hat f_2$ be $\C\infty1$ families of curves in $\Ms$, and let $g_k\co\ex C(g_k)\longrightarrow \ex C(\hat f_1)\times \ex C(\hat f_2)$ be a sequence of stable holomorphic curves in the family of targets $\ex C(\hat f_1)\times \ex C(\hat f_2)\longrightarrow \ex F(\hat f_1)\times \ex F(\hat f_2)$ such that the sequence of points  $\ex F(g_k)$ in $\ex F(\hat f_1)\times \ex F(\hat f_2)$ converges. Suppose further that the projection of $g_k$ to each $\ex C(\hat f_i)$ is a degree 1 map, and that the genus and number of ends of $\ex C(g_k)$ is the same as the curves in $\hat f_i$.  Then there exist refinements\footnote{See definitions 10.5 and 10.9 of \cite{iec} for refinements of curves. The need to take refinements is because of an assumption of Theorem 6.1 in \cite{cem} that an exploded manifold be basic. It should be possible to remove this assumption using a family of curves with universal tropical structure, from Theorem 4.8 of \cite{uts}, then we would get that a subsequence of $f_i$ converges to a stable holomorphic curve.} $g'_{k}$ of $g_{k}$ and  a subsequence of $\{g_{k}'\}$ converging to a holomorphic curve in $\Ms(\ex C(\hat f_1)\times \ex C(\hat f_2))$. \end{lemma}

\pf
We use Theorem 6.1 from \cite{cem} on the holomorphic  curves $g_i$ in  the family of targets $\ex C(\hat f_1)\times \ex C(\hat f_2)\longrightarrow \ex F(\hat f_1)\times \ex F(\hat f_2)$.  One condition of Theorem 6.1 is that the targets be basic\footnote{Definition 4.6 of \cite{iec}.}, but a curve $\ex C(f)$ is basic if and only if its tropical part $\totb{\ex C(f)}$ has no edge with both ends attached to a single vertex. Assuming that $\ex C(\hat f_i)\longrightarrow \ex F(\hat f_i)$ are indeed basic, Theorem 6.1 of \cite{cem} implies that some subsequence of these curves converges in $\C\infty1$ to a stable holomorphic curve. We now deal with the case that $\ex C(\hat f_i)$ is not basic.

\begin{claim}\label{bc} There exists a basic family of curves $\hat f_{i}'$, an open subset $U_{1}\times U_{2}\subset \ex F(\hat f_{1})\times \ex F(\hat f_{2})$ containing the limit of $\ex F(g_{k})$ and complete finite degree maps $r_{i}$ and $r_{i}^{+1}$ such that the following diagram commutes
\[\begin{tikzcd}\ex C(\hat f'_{i})\dar \rar{r_{i}^{+1}}& \ex C(\hat f_{i}\rvert_{U_i})\subset \ex C(\hat f_i)\dar
\\ \ex F(\hat f_{i}')\rar{r_{i}} & U_i\subset \ex F(\hat f_{i})\end{tikzcd}\]
and such that the induced map $\ex C(\hat f_{i}')\longrightarrow \ex C(r_{i}^{*}\hat f_{i})$ is a refinement map and fiberwise holomorphic. 
\end{claim}

Before proving Claim \ref{bc}, let us verify that Lemma \ref{double family} follows. Pass to a subsequence so that $\ex F(g_{k})$ is contained $U_1\times U_2$. As $(r_{1},r_{2})$ is complete and finite degree, we can choose a sequence of points $p_{k}\in (r_{1},r_{2})^{-1}(\ex F(g_{k}))$ that converge within $\ex F(\hat f_{1}')\times \ex F(\hat f_{2}')$. Note that $\ex C(\hat f'_{1})\times \ex C(\hat f'_{2})\longrightarrow \ex C(r_{1}^{*}\hat f_{1})\times \ex C(r_{2}^{*}\hat f_{2})$ is a refinement map, so every stable holomorphic curve in $\ex C(r_{1}^{*}\hat f_{1})\times\ex C(r_{2}^{*}\hat f_{2})$ lifts uniquely to a stable holomorphic curve in $\ex C(\hat f'_{1})\times \ex C(\hat f'_{2})$ with a refined domain. (The domain of this lift is the fiber product with the refinement map.) In particular, we obtain a sequence of stable holomorphic curves $g''_{k}$ in $\ex C(\hat f'_{1})\times \ex C(\hat f'_{2})$ over the converging sequence of points $p_{k}$, such that $r^{+1}\co \ex C( g_{k}')\longrightarrow \ex C(g_{k})$ is a refinement map, and the curve $g_k'=r^{+1}\circ g''_{k}$ is a refinement of $g_{k}$. Theorem 6.1 of \cite{cem} gives that a  subsequence of $g_{k}''$ converges to some stable curve $g''$, so the corresponding subsequence of $ g'_{k}$ converges to $r^{+1}\circ g''$, which is a refinement of a stable holomorphic curve in $\ex C(\hat f_1)\times \ex C(\hat f_2)$.  

To complete the proof of Lemma \ref{double family}, it remains to prove Claim \ref{bc}.
If $\ex C(\hat f_{i})$ is not basic, choose  a coordinate chart $\mathbb R^{n}\times \et mP$ on $\ex F(\hat f_i)$ with the limit of the points $\ex F(g_{k})$ contained in the stratum with tropical part the interior of the $m$--dimensional  polytope $P$. For notational convenience, we can assume that this chart covers $\ex F(\hat f_{i})$. Local charts on $\ex C(\hat f_{i})\longrightarrow \mathbb R^{n}\times \et mP$ containing interior edges of curves in $\hat f_{i}$ are then isomorphic to maps $\mathbb R^{n}\times \et {m+1}Q\longrightarrow \mathbb R^{n}\times \et mP$, in the form $(x,\tilde z_{1},\dotsc \tilde z_{m},\tilde z_{m+1})\mapsto (z,\tilde z_{1},\dotsc,\tilde z_{m})$, where the polytope $Q$ is a subset of the inverse image of $P$ under the projection $\mathbb R^{m}\times \mathbb R\longrightarrow \mathbb R^{m}$, constructed by intersecting with two half spaces, creating two extra faces,  both projecting isomorphically to $P$ --- so these faces must be cut out by an equation in the form 
\begin{equation}\label{iap}\alpha_{1}\totb{\tilde z_{1}}+\dotsb+\alpha_{m}\totb{\tilde z_{m}}+\totb{\tilde z}_{m+1}=c\end{equation}
where $c\in \mathbb R$,  $\alpha_{1},\dotsb,\alpha_{m}$ are integers,  and $\alpha_{m+1}=1$. This condition corresponds to the deriviative of our family being surjective on integral vectors, a condition of Definition 10.1 of \cite{iec}.  Our exploded manifold  $\ex C(\hat f_{i})$ is non-basic if and only if there is some such stratum with both these faces glued to a single stratum in $\totb{\ex C(\hat f_{i})}$. We can deal with this problem by refining  $\ex C(\hat f_{i})$ in such a coordinate chart ---  subdividing $Q$ into two polytopes meeting at another face projecting isomorphically to $P$. However,  the equation for this new face must be in the form (\ref{iap}). If the top and bottom faces of $Q$ meet, there may be no such plane in the form (\ref{iap}) between the top and bottom faces.  If the function $h\co P\longrightarrow \mathbb R$ measuring the distance between the top and bottom faces has derivative divisible by $2$, the plane halfway between the top and bottom faces will be in the form (\ref{iap}), so we can bisect $Q$ in this plane, refining our coordinate and obtaining a new family of curves. Similarly, if $h$ is divisible by $k$, then we can subdivide $Q$ by a plane $1/k$-th the distance between the top and bottom planes.

We are left with the case that $h$ has primitive integral derivative. After a $\mathbb Z$--affine change of coordinates, we can assume that  $h=\totb{z_{1}}$. In this case we make a base change of our family using a complete degree 2 map $r\co\mathbb R^{n}\times \et m{P'}\longrightarrow \mathbb R^n\times\et mP$, in the form $r(x,\tilde z_{1},\dotsc,\tilde z_{m})=(x,\tilde z_{1}^{2},\tilde z_{2},\dotsc,\tilde z_{m})$, obtaining a new family $\hat f_{i,1}$ of curves such that the following is a pullback diagram
\[\begin{tikzcd}\ex C(r^{*}\hat f_{i})\dar\rar &\ex C(\hat f_{i})\dar
\\ \mathbb R^{n}\times \et m{P'}\rar{r}& \mathbb R^{n}\times \et mP\end{tikzcd}\]
Then we can refine $\ex C(\hat f_{i,1})$ by bisecting $Q\times_{P}P'$ as above. Similarly, if  there are other problematic edges, by using some complete finite degree map $r$ we can pull back our family, then refine the resulting family $r^{*}\hat f_{i}$ to obtain a basic family of curves $\hat f_{i}'$.  This completes the proof of Claim \ref{bc}, and therefore the proof of Lemma \ref{double family}.

\stop

Recall Example \ref{totlX}, where  we defined a notion of topologically indistinguishable curves, and a topological space $\totl{\mathcal X_{top}}$ with points equivalence classes of topologically indistinguishable curves in $\mathcal X$.

\begin{lemma}\label{hausdorff} Two curves $f_i$ in $\dmsw$ are topologically indistinguishable if and only if there is some family $\hat f$ with maps $f_i\longrightarrow \hat f$ whose image in $\totl{\ex F(\hat f)}$ is the same point. Moreover, $\totl{{ \dmsw}_{top}}$ is Hausdorff. 
\end{lemma}
\pf Again, we neglect the easy case of a curve with domain $\ex T$; see Section \ref{T section} for a description of the moduli stack in this case. Let $f_i$ be two topologically indistinguishable curves. We must construct a family $\hat f$ containing $\hat f_i$ over topologically indistinguishable points in $\ex F(\hat f_i)$. 
Choose a core family $(\mathcal O,\hat f/G,s_i)$ containing $ f_1$ such that the action of $G$ fixes $\totl{f_1}\subset \totl{\ex F(\hat f)}$.  For any family $\hat h$ in $\mathcal O$, there is a closed subset of $\ex F(\hat h)$ comprised of $h$ such that $\hat f\circ \fun\Phi^{+1}(\hat h)=h\bd G$, and the image of  $\totl{\fun\Phi(\hat h)}$ in $\totl{\ex F(\hat f)}$ is $\totl {f_1}\subset \totl{\ex F(\hat f)}$. Accordingly, there is a closed substack of $\mathcal O$ comprised of such curves. As this closed substack contains $f_1$, it must also contain the topologically indistinguishable curve  $f_2$. Moreover, $\fun\Phi^{+1}(f_2)$ defines $\abs G$ maps $f_2\longrightarrow \hat f$, all sent to the same point in $\totl {\ex F(\hat f)}$. So, we now have two inclusions $f_i\longrightarrow \hat f$ sent to the same point in $\totl{\ex F(\hat f)}$, and the first part of our lemma is proved.

Let us verify that $\totl{{\dmsw}_{top}}$ is Hausdorff. Suppose that $ f_1$ and $f_2$ are curves in $\totl{{\dmsw}_{top}}$ without disjoint open neighborhoods. It follows from  Lemma \ref{topological sequential convergence} that there is some sequence of curves $g_k$ converging in $\C\infty1$ to both $ f_1$ and $f_2$. We must prove that $f_1$ and $f_2$ are topologically indistinguishable. The proof proceeds by constructing a curve $h$ and an isomorphism of $h$ with  curves $f_i'$ topologically  indistinguishable from $f_i$, and contained in some core family $\hat f_i/G_i$ that contains $f_i$.

Choose core families $\hat f_i/G_i$ containing $f_i$ with corresponding fiberwise holomorphic maps $\fun\Phi_i\co \mathcal O_i\longrightarrow \hat f_i/G_i$. We will construct the domain of our curve $h$ first as the domain of a stable holomorphic curve in $\ex C(\hat f_1)\times \ex C(\hat f_2)$. By passing to a subsequence, we assume that our sequence $g_k$ is contained in both $\mathcal O_1$ and $\mathcal O_2$. Recall that $\fun\Phi_i(g_k)$ is a $G_i$--bundle $g_k\bd G\longrightarrow g_k$ with a $G_i$--equivariant holomorphic map $\fun\Phi_i^{+1}(g_k)\co \ex C(g_k\bd G)\longrightarrow \ex C(\hat f_i)$. Choose a metric on $\hat{\ex B}$, and note that there is an open substack of $\mathcal O_i$ comprised of curves $g$ such that the distance between the maps $\hat f_i\circ \fun\Phi_i(g\bd G)$ and $g\bd G$  is less than $\epsilon$. Accordingly, we can assume that the distance in our metric between $\hat f_i\circ \fun\Phi_i(g\bd G_k)$ and $g\bd G_k$ converges to $0$. After a judicious choice of sections of $g\bd G_k\longrightarrow g_k$, and passing to a subsequence, we then get maps $\Phi_{i,k}\co \ex C(g_k)\longrightarrow \ex C(\hat f_i)$ that are isomorphisms onto fibers of $\ex C(\hat f_i)\longrightarrow \ex C(\hat f_i)$, and whose image converge to $\ex C(f_i)\subset \ex C(\hat f_i)$ (and also converge to the domain of any curve topologically indistinguishable from $ f_i$). 

 It follows from Lemma \ref{double family} that the holomorphic curves $(\Phi_{1,k},\Phi_{2,k})\co \ex C(g_k)\longrightarrow \ex C(\hat f_1)\times \ex C(\hat f_2)$ have a subsequence of refinements that converge in $\C\infty1$ to a refinement $\phi'$ of some stable holomorphic curve $\phi\co\ex C(h)\longrightarrow \ex C(\hat f_1)\times \ex C(\hat f_2)$. The domain of our required curve $h$ is defined to be the domain of this stable holomorphic curve $\phi$, and the map $h\co \ex C(h)\longrightarrow \hat {\ex B}$ is defined using the composition.
 \[\begin{tikzcd}\ex C( h)\dar{h} \rar{\phi} & \ex C(\hat f_1)\times \ex C(\hat f_2)\dar{\pi_1}
 \\ \hat{\ex B}& \lar{\hat f_1}\ex C(\hat f_1)\end{tikzcd}\]
    As the distance  in our metric between the maps $\hat f_1\circ \Phi_{1,k}\co \ex C(g_k)\longrightarrow \hat{\ex B}$ and  $\hat f_2\circ \Phi_{2,k}\co \ex C(g_k)\longrightarrow \hat{\ex B}$ converges to $0$ (and the same is true of our refinements of these curves) it follows that we can choose our limit $\phi'$ such that
\[\hat f_1\circ \pi_1\circ\phi'=\hat f_2\circ \pi_2\circ\phi'\]
where  $\pi_i$ is the projection $\ex C(\hat f_1)\times \ex C( f_2)\longrightarrow \ex C(\hat f_i)$. 

We have that $\pi_i\circ \phi\co \ex C( h)\longrightarrow \ex C( \hat f_i)$ is a degree 1, genus-and-end-preserving holomorphic map onto the domain of some curve $ f_i'$ in $\hat f_i$ topologically indistinguishable from $f_i$, and that $h=\hat f_i\circ \pi_i\circ\phi$. It remains to show that $\pi_i\circ \phi$ is an isomorphism onto $\ex C(f_i')$. If this degree 1 holomorphic map is not an isomorphism, there must be some component of the nodal curve $\totl{\ex C( h)}$ sent to a single point in some $\totl{\ex C(\hat f_i)}$, and this component must be a sphere with one or two special points. Moreover, this component must be sent to a single point in $\totl{\ex B}$, so it must be collapsed to a single point in  $\totl{\ex C(\hat f_1)\times \ex C(\hat f_2)}$, because curves in both $\hat f_i$  are stable. This collapsed component contradicts the fact that $\phi\co \ex C(h)\longrightarrow \ex C(\hat f_1)\times \ex C(\hat f_2)$  is stable. Therefore, $\pi_i\circ\phi$ defines an isomorphism $h\longrightarrow f'_i$. As the curves $f_1'$ and $f_2'$ are isomorphic, and topologically indistinguishable from $f_1$ and $f_2$, it follows that $f_1$ and $f_2$ are topologically indistinguishable, and $\totl{{ \dmsw}_{top}}$ is Hausdorff, as required.

\stop

In \cite{cem} we prove in many situations that  $\Mod(\hat {\ex B})\longrightarrow \ex B_{0}$ is proper when $\Mod$ is restricted to connected components of $\Msw(\hat {\ex B})$. In \cite{cem}, this properness was stated as follows: every sequence of holomorphic curves in a connected component of $\Msw$ over a convergent sequence in $\ex B_{0}$ must have a  subsequence that converges in $\C\infty1$. 
The next lemma allows us to conclude that such sequential compactness implies compactness.

\begin{lemma}\label{metrizable strata} $\totl{\Msw_{top}}$ is a countable union of strata, and on each stratum there is a metric that has the following relation to the topology on $\totl{\Msw_{top}}$. Let $\totl{\mathcal X_{top}}\subset \totl{\Msw_{top}}$ be a sequentially compact subspace. Then, the subspace topology of  $\totl{\mathcal X_{top}}$ intersected with each stratum coincides with the metric topology. 
\end{lemma}
\pf

Give a stable curve  $f$, there is a canonical family $\tilde f$ of curves topologically indistinguishable from $f$, constructed as follows: Restrict any core family  $\hat f$ containing $f$ to the set of curves in $\hat f$ with the same image in $\totl{\ex F(\hat f)}$.  This family $\tilde f$ has the property that any connected family of curves $\hat h$ topologically indistingishable from $f$ admits exactly $G$ maps into $\tilde f$, where $G$ is the subgroup of the automorphism group of $\hat f$ preserving the image of $f$ in $\totl{\ex F(\hat f)}$. Any two such $f\longrightarrow \tilde f$ are canonically isomorphic. 
 
There is a natural stratification of $\Msw$ such that two curves are in the same stratum if they are isomorphic to curves in a single family $\hat f$, and both these curves are contained in the same stratum of $\ex F(\hat f)$. For this lemma, we use a slightly coarser stratification, defined as follows:   Stable curves $f_{1}$ and $f_{2}$  will be in the same stratum if there exists an invertible map $\phi\co \ex C(\tilde f_{1})\longrightarrow \ex C(\tilde f_{2})$ such that the following diagram commutes
\begin{equation}\label{msd}\begin{tikzcd}\ex C(\tilde f_{1})\dar\rar{\phi} &\ex C(\tilde f_{2}) \dar 
\\ \ex F(\tilde f_{1})\rar & \ex F(\tilde f_{2})\end{tikzcd} \end{equation}
and such that $\totb{\tilde f_{1}}=\totb{\phi\circ \tilde f_{2}}$. If  $ f_{1}$ and $f_{2}$ are topologically indistinguishable, then  $\phi$ can be chosen holomorphic such that  that the maps $f_{1}$ and $f_{2}\circ \phi$ coincide. The strata of $\totl{\ex C(\tilde f)}=\totl{\ex C(f)}$ corresponding to vertices of $\totb{\ex C(f)}$ are punctured Riemann surfaces. On any unstable stratum, we  assume that $\phi$ is holomorphic. This can always be achieved.  If, however, a stratum is stable, it has a canonical complete hyperbolic metric, and we can therefore measure the size of $\dbar\phi$ using the metric on the domain and codomain. Define $ \abs{\dbar\phi}^{\infty}$ to be the supremum of $\abs{\dbar\phi}$ over all these stable strata.  

Choose a metric on $\hat{\ex B}$, and let $\text{dist}(x,y)$ be the distance in this metric between any two points $x$ and $y$ with the same tropical part in $\totb{\hat{\ex B}}$. Then, for any two stable curves  $f_{1}$ and $f_{2}$ in the same stratum, define
\[d(f_{1},f_{2}):=\inf_{\phi}\lrb{\sup_{x\in \ex C(\tilde f_{1})}\text{dist}(\tilde f_{1}(x),\tilde f_{2}\circ \phi(x)) +\abs{\dbar\phi}^{\infty}+\abs{\dbar\phi^{-1}}^{\infty}} \]
where the infimum is over all invertible maps $\phi$  such that (\ref{msd}) commutes and such that $\phi$ is holomorphic on  unstable components. This $d(\cdot,\cdot)$ is symmetric, non-negative, and satisfies the triangle inequality. Moreover, if $f$ and $f'$ are topologically indistinguishable, then $d(f,f')=0$, so $d(\cdot,\cdot)$ induces a (possibly degenerate) metric on each stratum of $\totl{\Msw_{top}}$. It remains to check that there are only countably many strata, and that $d(\cdot,\cdot)$ induces the subspace topology on any sequentially compact subset. 

\begin{claim}\label{metrica} Let $\mathcal X$ be a sequentially compact substack of $\Msw$, and let $h_{i}$ be a sequence of  curves in $\mathcal X$. If $d(h_{i},f)$ converges to $0$, then $h_{i}$ converges in $\C\infty1$ to $f$.
\end{claim}

To prove Claim \ref{metrica}, it suffices to prove that every subsequence of $\{h_{i}\}$ has a subsequence converging to $f$. By assumption, each subsequence  of $\{h_{i}\}$  must have a subsequence converging to some stable  curve in $\C\infty1$, so  without losing generality, we can assume that $h_{i}$ converge to some stable  curve $h$, and using Lemma \ref{sequence family}, we can assume that $h_{i}$ converge to $h$ within a family of curves $\hat h$. It remains to show that $h$ is topologically equivalent to $f$. Choose a core family $\hat f/G$ containing $f$. Each neighborhood of $\tilde f$ in $\hat f$ contains curves with all complex structures close to $\tilde f$, so there exists a sequence of holomorphic maps $\phi_{i}'\co \ex C(h_{i})\longrightarrow \ex C(\hat f)$ which are holomorphic  isomorphisms onto fibers of $\ex C(\hat f)\longrightarrow \ex F(\hat f)$ converging in $\totl{\ex F(\hat f)}$ to the image of $f$ such that 
\[\text{dist}(h, \hat f\circ \phi_{i})\] 
converges to $0$. Arguing, using Lemma \ref{double family} as in the proof of Lemma \ref{hausdorff},  implies that $h$ is topologically equivalent to $f$, completing the proof of Claim \ref{metrica}.

Note that Claim \ref{metrica} also implies that $d(f,h)=0$ if and only if $f$ and $h$ are topologically indistinguishable, so $d(\cdot,\cdot)$ is nondegenerate, and defines a metric on our stratum of $\totl{\Msw_{top}}$. 

\begin{claim}\label{metricb} If a sequence of stable curves  $f_{i}$ are in the same stratum as $f$, and $f_{i}$ converges to $f$, then $\lim_{i\to\infty}d(f_{i},f)=0$. 
\end{claim}

To prove Claim \ref{metricb}, choose a core family $((\fun\Phi,\fun\Phi^{+1}),\hat f/G,\{s_i\})$ for an open neighborhood $\mathcal O$ of $f$. As $f_{i}$ converges to $f$, $f_{i}$ (and therefore $\tilde f_{i}$) is eventually contained in $\mathcal O$, and the image of ${\fun\Phi}_{ f_{i}}$ converges to the image of $f$ in $\totl{\ex F(\hat f)/G}$. Condition \ref{cc2} of Definition \ref{core family} implies that, after passing to a subsequence, we can assume that  that $\hat f/G\circ{\fun\Phi}_{ f_{i}}^{+1}$ is a $\mathbb R$--deformation of $ f_{i}$, and  therefore these two stable curves are in the same stratum of $\Msw$. Use $f_{i}'\in\hat f$ for some choice of curve in $\hat f$ isomorphic to this $\hat f/G\circ{\fun\Phi}_{ f_{i}}^{+1}$ so that the image of $f_{i}'$ converges to $f$ in $\hat f$.  By assumption, $f_{i}$ and $f_{i}'$ are in the same stratum of $\Msw$ as $f$, and therefore $\totb{\ex F(\tilde f)}$ is isomorphic to $\totb{\tilde f_{i}'}$, so the tropical part of the strata of $\ex F(\hat f)$ containing these $f$ and $f_{i}'$ must have the same dimension. As the image of $f_{i}'$ also converges to the image of $f$ in $\ex F(\hat f)$, this implies that $f_{i}'$ and $f$ must be in the same stratum of $\ex F(\hat f)$. On a neighborhood of $\ex C(f)$ within this stratum,  $\ex C(\hat f)\longrightarrow \ex F(\hat f)$ is isomorphic to $\mathbb R^n\times \ex C(\tilde f)\longrightarrow \mathbb R^n\times \ex F(\tilde f)$ (where the complex structure on fibers may vary smoothly with the $\mathbb R^n$ coordinates).  Projecting out the $\mathbb R^n$--direction, we obtain isomorphisms $\phi_i\co \ex C(\tilde f_i)\longrightarrow\ex C(\tilde f)$ which we can use to verify that $d(f_i,f)$ converges to $0$.

 Claims \ref{metrica} and \ref{metricb} along with lemmas \ref{topological sequential convergence} and \ref{open pullback} imply that $d(\cdot,\cdot)$ induces the subspace topology on any sequentially compact subset of a stratum of $\totl{\Msw_{top}}$. We must now check that $\Msw$ has only countably many strata.

%
%

\begin{claim}\label{stratac} Two stable curves $f$ and $h$ are in the same stratum of $\Msw$ if and only if there exists a strata-preserving homeomorphism $b\co \totl{\ex C(f)}\longrightarrow \totl{\ex C(h)}$ such that
\begin{itemize}
\item for each stratum $s$ of $\totl{\ex C(f)}$,  $\totl h$ sends $b(s)$ to the stratum of $\totl{\hat {\ex B}}$ containing $\totl{f}(s)$;
\item  The tropical structure of the curves $f$ and $h$ agree in the sense explained below. Note that   $0$-dimensional strata of $\totl{\ex C(f)}$  correspond to edges of $\totb{\ex C(f)}$. If $s$ corresponds to an internal edge, then a neighborhood of the stratum $s\subset \totl{\ex C(f)}$ is homeomorphic to two disks joined at the origin. An orientation of this edge  corresponds to a choice of one of these disks, so $b$ induces a bijection between oriented edges  $e$ of $\totb{\ex C(f)}$ and oriented edges $b(e)$ of $\totb{ \ex C(h)}$. We require that $\totb f(e)$ and $\totb h(b(e))$ are in  the same stratum of $\totb{\hat{\ex B}}$, and we require 
that the derivatives $\alpha_e$ and $\alpha_{b(e)}$ of $\totb f$ and $\totb h $ agree. So,  the last condition on our homeomorphism is that 
\[\alpha_e=\alpha_{b(e)}\]\end{itemize}
\end{claim}

Up to homeomorphism, there are only countably  many nodal curves with marked points. Moreover,  each  $\ex C(f)$ only has finitely many strata, and $\totb{\ex B}$ has at most countably many strata, and each stratum only has countably many integral vectors. So,   Claim \ref{stratac} implies that $\Msw$ has countably many strata, and to finish the proof of Lemma \ref{metrizable strata}, we need only prove Claim \ref{stratac}.

Let us identify $\totb {\ex F(\tilde f)}$ using Remark 3.3 of \cite{uts}, but simplifying using that $\hat{\ex B}$ is basic. For each vertex $v$ of $\totb{\ex C(f)}$, let the polytope $P_{v}$ be the closure of the stratum $P_{v}^{\circ}$ of $\totb{\ex B}$ containing $v$. There is a natural map $A_{v}\co\totb{\ex F(\tilde f)}\longrightarrow P^{\circ}_{v}$ given by evaluating $\totb{\tilde f}$ on the stratum of $\totb{\ex C(\tilde f)}$ containing $v$.  Moreover, for each internal edge $e$ of $\totb{\ex C(f)}$, there is a natural map $l_{e}\co \totb{\ex F(\tilde f)}\longrightarrow(0,\infty)$ that records the length of the edges of fibers of $\totb{\ex C(\tilde f)}\longrightarrow \totb{\ex F(\tilde f)}$ corresponding to $e$. Remark 3.3 of \cite{uts} tells us that $\prod_{v}A_{v}\prod_{e}l_{e}$ is a $\mathbb Z$--affine embedding of $\ex F(\tilde f)$ into $\prod_{v}P^{\circ }_v\times (0,\infty)^{k}$, where $\totb{f}$ has $k$ internal edges, and tells us that this embedding is onto a sub-polytope cut out by $\mathbb Z$--affine equations. In this case, the $\mathbb Z$--affine equations are also very easy to describe: Orient an internal edge of $\totb{\ex F(f)}$, so that it has a well-defined derivative $\alpha_{e}$, and is attached at its beginning to a vertex $v_{1}$ and at its end to a vertex $v_{2}$. Note that polytopes $P_{v_{i}}$ are faces of the polytope $P_{e}$ containing the image of $e$. Then, for any $x\in \totb{\ex F(\tilde f)}$, 
\[A_{v_{1}}(x)+l_{e}(x)\alpha_{e}=A_{v_{2}}(x) \]
and such equations for all internal edges $e$ are all the equations that cut out $\ex F(\tilde f)$ as a sub-polytope of $\prod_{v}P^{\circ }_v\times (0,\infty)^{k}$. As described in Remark 3.3 of \cite{uts}, $\totb{\ex C(\tilde f)}$ and $\totb{\tilde f}$ is also determined by this construction. In particular,  given points $A_{v}(x)\in P_{v}^{\circ}$ and lengths $l_{e}(x)$ satisfying these equations, there is a unique tropical curve that is a continuous deformation of $\totb f$ with these edge lengths and vertices at the chosen positions in $P_{v}^{\circ}$, and this curve is the fiber of $\totb{\tilde f}$ over $x\in \totb{\ex F(\tilde f)}$. From this construction, what we need is the following:  given a homeomorphism between $\totl{\ex C(f)}$ and $\totl{\ex C(h)}$ satisfying the conditions of Claim \ref{stratac}, we there is an isomorphism between $\totb{\tilde f}$ and $\totb{\tilde h}$, and in particular an isomorphism $\totb\phi\co \totb{\ex C(\tilde f)}\longrightarrow \totb{\ex C(\tilde h)}$ such that $\totb {\tilde h}\circ \totb\phi=\totb {\tilde f}$. 

By choosing a different curve in $\tilde h$ if necessary, we can now reduce to the case $f$ and $h$ have the same tropical part, so  $\totb \phi$ sends  $\totb{\ex C( f)}$ to $\totb{\ex C(h)}\subset \totb{\ex C(\tilde h)}$. The local model for edges of curves depends only on their length, so there exists a holomorphic isomorphism $\phi_0$ from a neighborhood of the edges of $\ex C(f)$ to a neighborhood of the edges of $\ex C(h)$ such that the tropical part of $\phi_0$ agrees with $\totb\phi$.   We can choose a smooth representative of our homeomorphism $b\co \totl{\ex C(f)}\longrightarrow \totl{\ex C(h)}$ that agrees with $\totl{\phi_0}$ on  a neighborhood of all special points in $\totl{C(f)}$. There is then a unique smooth isomorphism $\ex C(f)\longrightarrow \ex C(h)$ with smooth part $b$ and tropical part $\totb{\phi}$. Moreover, this isomorphism extends uniquely to an isomorphism $\phi\co \ex C(\tilde f)\longrightarrow \ex C(\tilde h)$ with tropical part $\totb\phi$. This $\phi$ is the isomorphism required to show that $f$ and $h$ are in the same stratum of $\Msw(\hat {\ex B})$, which concludes the proof of Claim \ref{stratac} and Lemma \ref{metrizable strata}.

\stop

\begin{cor}\label{compact equivalence} Any sequentially compact substack $\mathcal X$ of $\Msw$ is compact.
\end{cor}
\pf 

The proof proceeds by applying Lemma \ref{metrizable strata} to the deepest stratum of $\mathcal X$, then removing a neighbourhood of this stratum.  Each stratum of $\Msw$ has a tropical dimension. In particular if $f$ is a curve in $\Msw$, and $\tilde f$ is the family of curves topologically indistinguishable from $f$, then the dimension of $\totb{\ex F(\hat f)}$ depends only on the tropical part $\totb f$ of $f$, because $\tilde f$ is constructed by restricting a family with universal tropical structure.  This tropical dimension is therefore  the same for all curves in a stratum of $\Msw$. Moreover, if a sequence of curves $\hat f_i$ converges to to $f$, then either $f_i$ are eventually in the same stratum as $f$, or the dimension of $\totb{\ex F(\hat f_i)}$ is eventually strictly less than $\dim \totb{\ex F(\hat f)}$. 

As $\mathcal X$ is sequentially compact, it follows that the tropical dimension of strata of $\mathcal X$ is bounded above. Otherwise, there would exist a sequence of curves in strata with increasing tropical dimension, and this sequence would have no convergent subsequence. Let $\mathcal X_0$ be the union of all strata of $\mathcal X$ with maximal tropical dimension. It follows that $\mathcal X_0$ is sequentially compact. Moreover, Lemma \ref{metrizable strata} implies that  $\totl{(\mathcal X_0)_{top}}$ is metrizable, and therefore compact.  Given any open cover of $\mathcal X$, there therefore exists a finite subcover $\{\mathcal U_1,\dotsc,\mathcal U_n\}$ of $\mathcal X_0$. Moreover $\mathcal X\setminus \bigcup_{k=1}^n \mathcal U_k$ is sequentially compact, and has maximal tropical dimension strictly less than the maximal tropical dimension of $\mathcal X$. 

The corollary then follows by induction on the maximal tropical dimension of $\mathcal X$.

\stop  

We use the following technical lemma to know that we can shrink open substacks appropriately for defining extendible Kuranishi charts. 

\begin{lemma}\label{shrink} Suppose each connected component of $\Msw(\hat{\ex B})$ is proper over $\ex B_0$ in the following sense:  if a sequence of holomorphic curves in a connected component of $\Msw(\hat{ \ex B})$ has image in $\ex B_0$ that converges, then this sequence has a subsequence that converges in $\C\infty 1$.

If $\hat f/G$ is a core family for the open substack $\mathcal O$ of $\Msw(\hat{\ex B})$, and $f$ is a holomorphic curve in $\hat f$,  there exists a continuous function $\rho\co \mathcal O\longrightarrow \mathbb R$ such that $\rho(f)=1$, and such that any holomorphic curve in the closure (within $\Msw(\hat{\ex B})$) of  $\{\rho>0\}$ is contained in $\mathcal O$. 
\end{lemma}

\pf

Roughly speaking, this lemma holds  because  weak convergence of holomorphic curves to $f$ --- detected by continuous functions on $\mathcal O$ --- automatically implies $\C\infty1$ convergence.

Choose  a family of metrics on the fibers of  $\hat{\ex B}\longrightarrow \ex B_0$, and use $\text{dist}(x,y)$ to denote the distance between two points in $\hat{\ex B}$ using this metric.
Recall that Definition \ref{core family} tells us that  given any family of curves $\hat h\in\mathcal O$, there is a canonical fiberwise-holomorphic  map
\[\begin{tikzcd}\dar\ex C(\hat h)\rar{{\fun\Phi}^{+1}_{\hat h}}&\ex C(\hat f)/G\dar
\\\ex F(\hat h)\rar{{\fun\Phi}_{\hat h}}&\ex F(\hat f)/G \end{tikzcd}\]
and that composing this map with the map $\hat f/G\co \ex C(\hat f)/G\longrightarrow \hat{\ex B}$ gives another family of curves
\[\hat f/G\circ {\fun\Phi}^{+1}_{\hat h}\co \ex C(\hat h)\longrightarrow \hat{\ex B}\] which is a deformation of $\hat h$, so   $\text{dist}(\hat h(x),\hat f/G\circ {\fun\Phi}^{+1}_{\hat h}(x))$ is  bounded and  a continuous function on $\ex C(\hat h)$.

Choose a proper  $\C\infty1$  function $r_{0}\co \ex F(\hat f)/G\longrightarrow [0,\infty)$ equal to $0$ only on the $G$--orbit of $f$, (and curves with the same smooth part). So, given a sequence of points $p_i$ in $\ex F(\hat f)$,  $r_0(p_{i})\to 0$ implies $p_{i}$  converges to $f\in \ex F(\hat f)$ or some $G$--translate of $f$. For any family of curves $\hat h$ in $\mathcal O$ define 
\[r\co \ex F(\hat h)\longrightarrow [0,\infty)\]
as 
\[r(h): =\sup_{x\in\ex C(h)}\text{dist}(h(x),\hat f/G\circ {\fun\Phi}^{+1}_{h}(x))+r_{0}(\fun\Phi( h))\text{ for all curves $h$ in $\hat h$.}\]
Clearly, $r(h)$ depends only on $h$, and for each family $\hat h$ in $\mathcal O$, $r\co \ex F(\hat h)\longrightarrow [0,\infty)$ is continous, so $r$ defines a continuous map
 \[r\co \mathcal O\longrightarrow [0,\infty)\ .\]

\begin{claim}\label{stc} There exists an $\epsilon>0$  such that any holomorphic curve in the closure of $\{r<\epsilon\}\subset \Msw$  must be contained in $\mathcal O$.
 \end{claim}

To prove Claim \ref{stc}, suppose to the contrary that there exists a sequence of holomorphic curves $\{h_{i}\}$ not in $\mathcal O$ such that for all $\epsilon$, $\{h_{i}\}$ is eventually contained in the closure of the substack  where $r<\epsilon$. It follows that the image of $h_{i}$ in $\ex B_{0}$ converges, and $\{h_{i}\}$ is eventually contained in a connected component of $\Msw$. Our assumption on the properness of the map $\Mod\longrightarrow \ex B_{0}$ then implies that some subsequence of $\{h_{i}\}$ converges in $\C\infty1$ to a stable holomorphic curve $h$.

As $\mathcal O$ is open, our limiting holomorphic curve $h$ must not be in $\mathcal O$. On the other hand, for all $\epsilon$, the curve $h$ is in the closure of the substack of $\mathcal O$ where $r<\epsilon$. We shall achieve a contradiction by showing that this implies that $h$ and $f$ are topologically indistinguishable. 

Lemma \ref{topological sequential convergence} implies
that there is a sequence of curves $\{f_{i}\}$ in $\mathcal O$  converging in $\C\infty1$ to $h$, and such that $r(f_{i})$ converges to $0$. Lemma \ref{sequence family} implies that, by passing to a subsequence, we can assume that $f_{i}$ converge to $h$ within some $\C\infty1$ family $\hat h$ containing $h$.

 As $r_{0}\circ {\fun\Phi}_{f_{i}}$ converges to $0$, the images of $\ex C( f_{i})$ in $\ex C(\hat f)/G$ converge to $\ex C(f)$. By a judicious choice of resolution of the $G$--fold ambiguity of the map ${\fun\Phi}^{+1}_{f_{i}}$, we obtain inclusions $\ex C(f_{i})\longrightarrow \ex C(\hat f)$ converging to $\ex C(f)$ which we shall again call ${\fun\Phi}^{+1}_{f_{i}}$. Applying Lemma \ref{double family} as in the proof of Lemma \ref{hausdorff} then gives a stable holomorphic curve $\phi$ in $\ex C(\hat h)\times \ex C(\hat f)$ with image in $\ex C(h')\times \ex C(f')$ for curves $h'$ and $f'$ topologically indistinguishable from $h$ and $f$, and such that
 \[\dist{(\hat h\circ \phi,\hat f\circ \phi)}=0\ ,\]
  so $\hat h\circ \phi=\hat f\circ \phi$. As in the proof of Lemma \ref{hausdorff},  the stability of  $\phi$, $h'$ and $f'$ implies that $\hat h\circ\phi$ is isomorphic to $h'$ and $\hat f\circ \phi$ is isomorphic to $f'$, so $h$ and $f$ are topologically  indistinguishable. This completes the proof of Claim \ref{stc}.

 We now have that there exists some $\epsilon$ so that all holomorphic curves in the closure of $\{r>\epsilon\}$ are contained in $\mathcal O$. To complete the proof of our lemma, all we need to do is compose $r$ with a cut off function  equal to $1$ at $0$ and vanishing outside of an $\epsilon$--neighborhood of $0$. The resulting continuous function $\rho\co \mathcal O\longrightarrow [0,1]$ is $1$ at $f$ and has the property that every holomorphic curve in the closure of $\{\rho>0\}$ is contained in $\mathcal O$. This completes the proof of Lemma \ref{shrink}.
 
 \stop

\section{Locally representing the moduli  stack of solutions to $\dbar f\in V$}

\label{piV section}

The goal of this section is to prove Theorem \ref{V moduli stack}. In particular,   for a simply-generated subsheaf $V$ of $\Y$  transverse to $\dbar$ at a holomorphic curve $f$, the moduli stack  of curves $f'$ with $\dbar f'\in V(f')$,  is locally represented  by $\hat f/G$, where $\hat f$ is a family of curves containing $f$ and  $G$ is a finite group of automorphisms.

 Proposition \ref{model map to A} below is a way of locally representing the moduli stack of holomorphic curves  parametrized by the domain of a particular family of holomorphic curves. This proposition is then used in Lemma \ref{parametrize by curves} to prove that an arbitrary simply-generated subsheaf $V$ can be parametrized by a family of curves. With $V$ written in this special form, the results of \cite{reg} imply Theorem \ref{V moduli stack}.

\begin{prop}\label{model map to A}Given any family of holomorphic curves $\hat f$ in $\hat {\ex B}$ containing a given curve $f$, there exists a $\C\infty1$ family of curves $\hat f^\star$ in $\hat{\ex B}$  satisfying  the following properties:
\begin{enumerate}

\item

There is an inclusion $\iota\co \hat f\longrightarrow \hat f^\star$ and a map
\[\begin{tikzcd}\ex C(\hat f^\star)\arrow{r}{\psi}\arrow{d} &\ex C(\hat f)\arrow{d}
\\ \ex F(\hat f^\star)\arrow{r} &\ex F(\hat f)\end{tikzcd}\] that is a holomorphic isomorphism restricted to each fiber. Moreover $\ex C(\iota)$ and $\psi$ are related so that $\ex C(\hat f^\star)$ is isomorphic to a vector-bundle over $\ex C(\hat f)$ with projection $\psi$ and zero section $\ex C(\iota)$.
\item  \label{mmA2} There exists an open neighborhood $\mathcal O$ of $(\id, f)$ in the ambient moduli stack of $\C\infty1$ curves in $\ex C(\hat f)\times \hat{\ex B}$, such that 
every  family $(\psi', \hat h)$ of holomorphic curves in $\mathcal O$ has a unique map to $(\psi,\hat f^\star)$ in $\Ms(\ex C(\hat f)\times \hat{\ex B})$. In particular, $\hat h$ is a family of holomorphic curves in $\hat{\ex B}$ and $\psi'$ is fiberwise-holomorphic map $\psi'\co\ex C(\hat h)\longrightarrow \ex C(\hat f)$, and 
  there exists a unique map $\hat h\longrightarrow \hat f^\star$ such that the following diagram commutes. 

\[\begin{tikzcd}\ex C(\hat h)\arrow{r}{ \hat h }\arrow{d}{\psi'}\arrow{dr}{\exists!}&\hat{\ex B}
\\ \ex C(\hat f)&\arrow{l}{ \psi }\ex C(\hat f^\star)\arrow{u}[swap]{\hat f^\star}\end{tikzcd}\]
\end{enumerate}

 Suppose furthermore that $\hat {\ex B}$ has a group, $G_{0}$, of automorphisms and that there is a group $ G_{0}\times G$ of automorphisms of $\ex C(\hat f)\longrightarrow \ex F(\hat f)$ such that $\hat f$ is $G$--invariant and $G_{0}$--equivariant. 
Then $\hat f^\star$ can be constructed so that $\ex C(\hat f^\star)\longrightarrow \ex F(\hat f^\star)$ has a group $G_{0}\times G$ of automorphisms, $\psi$ and $\iota$ are $G_{0}\times G$--equivariant, and $\hat f^\star$ is $G$--invariant and  $G_{0}$--equivariant.
 
\end{prop}

By choosing $\mathcal O $ sufficiently small, we can also ensure that the map $\psi'\co\ex C(\hat h)\longrightarrow \ex C(\hat f)$ is an isomorphism restricted to each fiber. 

\begin{remark} The intersection of $\dbar\hat f^\star$ with $0$ should be regarded as representing the moduli space of holomorphic curves close to $f$ and parametrized by $\ex C(\hat f)\longrightarrow \ex F(\hat f)$. In the case when $G_{0}\times G$ is nontrivial, \[\{f'\in \hat f^\star \text{ such that }\dbar f'=0\}/G_{0}\times G\] should be regarded as representing the moduli stack of holomorphic curves in $\hat{\ex B}/G_{0}$ parametrized by $\ex C(\hat f)/(G_{0}\times G)$. The point of Proposition \ref{model map to A} is that this moduli stack is locally represented as a finite quotient of  a subset of a finite-dimensional $\C\infty1$ family of curves.
\end{remark}

\pf

We prove the equivariant case.  Using a $G_{0}$--invariant, smooth, $J$--preserving connection on $T_{vert}\hat {\ex B}$, we can construct a $G_{0}\times G$--invariant trivialization $(\mathcal F,\phi)$ to associate to $\hat f$. More precisely, define the map $\mathcal F\co f^{*}T_{vert}\hat{\ex B}\longrightarrow \hat{\ex B}$  by exponentiating using our invariant connection (and reparametrizing in a $G_{0}\times G$--equivariant way to ensure injectivity of $T\mathcal F$ restricted to any vertical tangent space), and define $\phi\co  \mathcal F^{*}T_{vert}\hat {\ex B}\longrightarrow \hat f^{*}T_{vert}\hat{\ex B}$  using parallel transport, along a straight line homotopy, using our connection. Such a $\mathcal F$ is $G$--invariant and $G_{0}$--equivariant, and $\phi$  is $G_{0}\times G$--equivariant.
Using such a trivialization, $\dbar\co X^{\infty,\underline 1}(\hat f)\longrightarrow \Y(\hat f)$ is $G_{0}\times G$--equivariant.

Consider a family  of curves $\hat h$ in $\hat {\ex B}$ and a fiberwise-holomorphic-isomorphism map $\psi'$.
\[\begin{tikzcd}\ex C(\hat h)\arrow{r}{\psi'}\arrow{d} &\ex C(\hat f)\arrow{d}
\\ \ex F(\hat h)\arrow{r} &\ex F(\hat f)\end{tikzcd}\]
 So long as  $\hat f\circ\psi'$ is sufficiently close to $\hat h$, the trivialization associated with $\hat f$ allows us to uniquely factor $\hat h$ as $\mathcal F$ composed with a map \begin{equation}\label{nuh}\nu_{\hat h}\co \ex C(\hat h)\longrightarrow f^{*}T_{vert}\hat{\ex B}\end{equation} so that the following diagram commutes. 
\[\begin{tikzcd}\ex C(\hat h)\arrow[bend left]{rr}{\hat h}\arrow{r}{\nu_{\hat h}}\arrow{dr}{\psi'}&f^{*}T_{vert}\hat{\ex B}\arrow{r}{\mathcal F}\arrow{d}&\hat{\ex B}
\\ &\ex C(\hat f)\end{tikzcd}\]

 Choose some $G_{0}\times G$--invariant collection  of marked-point sections $\{s_{i}\}$ of $\ex C(\hat f)\longrightarrow \ex F(\hat f)$ so that $D\dbar$ is injective when restricted to the subspace of sections of $f^{*}T_{vert}\hat{\ex B}$ vanishing at these extra marked points.  We shall use this injectivity with Theorem \ref{regularity theorem} to modify an extension, $\hat f^\star_{0}$,  of $\hat f$ to contain all relevant holomorphic curves.   Construct the domain of  $\hat f^\star_{0}$ as the following pullback.
\begin{equation}\label{f'def}\begin{tikzcd}\ex C(\hat f^\star_{0})\arrow{r}\arrow{d}&\ex C(\hat f)\arrow{d}
\\ \oplus_{i}(s_{i}\circ \hat f)^{*}T_{vert}\hat {\ex B}\arrow{r} &\ex F(\hat f)\end{tikzcd}\end{equation}
Note that $\ex C(\hat f^\star_{0})$ has a natural $G_{0}\times G$ action so that the above diagram is $G_{0}\times G$--equivariant.  Pull back the sections $s_{i}$ to give a $G_{0}\times G$--invariant collection of sections $\{s_{i}^\star\}$ of $\ex C(\hat f^\star_{0})\longrightarrow \ex F(\hat f^\star_{0})$.
\[\begin{tikzcd}\ex C(\hat f^\star_{0})\arrow{r}\arrow[transform canvas={xshift=-0.1cm}]{d}&\ex C(\hat f)\arrow[transform canvas={xshift=-0.1cm}]{d}  
\\ \oplus_{i}(s_{i}\circ \hat f)^{*}T_{vert}\hat {\ex B}\arrow[transform canvas={xshift=0.1cm}]{u}[swap]{s_{i}^\star}\arrow{r} &\ex F(\hat f)\ar[transform canvas={xshift=0.1cm}]{u}[swap]{s_{i}}\end{tikzcd}\]
Now, construct $\hat f^\star_{0}$ so that $\hat f^\star_{0}$ factors as a $G_{0}\times G$--equivariant map $\nu$ to $\hat f^{*}T_{vert}\ex B$ followed by $\mathcal F$, so that the $5$ inner loops in the following diagram commute,
\[\begin{tikzcd}
\\ \ex C(\hat f^\star_{0})\arrow{r}{\nu}\arrow[bend left]{rr}{\hat f^\star_{0}}\arrow[transform canvas={xshift=-0.1cm}]{d}\arrow{dr}&\hat f^*T_{vert}\hat{\ex B}\arrow{r}{\mathcal F}\arrow[transform canvas={xshift=-0.1cm}]{d}&\hat{\ex B}\arrow{d}
\\ \ex F(\hat f^\star_{0})\arrow{dr}\arrow[transform canvas={xshift=0.1cm}]{u}[swap]{s^\star_{i}} &\ex C(\hat f)\arrow[transform canvas={xshift=0.1cm}]{u}\arrow[transform canvas={xshift=0.1cm}]{d}\arrow{ur}& \ex G
\\ &\ex F(\hat f)\arrow[transform canvas={xshift=-0.1cm}]{u}{s_{i}}\arrow{ur}
\end{tikzcd}\]
and so that, if $\ex C(\hat f^\star_{0})$ is considered as a vector-bundle over $\ex C(\hat f)$, $\nu$ is a map of vector-bundles, and \[\nu\circ s^\star_{i}\co \oplus(s_{i}\circ \hat f)^{*}T_{vert}\hat {\ex B}\longrightarrow (s_{i}\circ \hat f)^{*}T_{vert}\hat{\ex B}\]
is projection onto the $i$th factor.

Our fiberwise-holomorphic-isomorphism map $\psi'\co\ex C(\hat h)\longrightarrow \ex C(\hat f)$ has a natural lift to a fiberwise-holmorphic-isomorphism map $\psi''\co \ex C(\hat h)\longrightarrow \ex C(\hat f^\star_{0})$, defined as follows. First define a map \[ x_{\hat h}\co \ex C(\hat h)\longrightarrow \oplus_{i}(s_{i}\circ \hat f)^{*}T_{vert}\hat{\ex B}\]
by evaluating the section $\nu_{\hat h}$ from equation (\ref{nuh}) at the image of each section $s_{i}$ using the following composition.  
\[ \begin{tikzcd}\ar[bend left]{rrr}{x_{\hat h}}\ex C(\hat h)\rar & \ex F(\hat h)\rar{(\psi')^{*}s_{i}} & \ex C(\hat h)\rar{\nu_{\hat h}} &(s_{i}\circ \hat f)^{*}T_{vert}\hat{\ex B}\end{tikzcd}\]
Then, the following diagram commutes.
\[\begin{tikzcd}\ex C(\hat h)\arrow{r}{\psi'}\arrow{d}{x_{\hat h}}& \ex C(\hat f)\arrow{d}
\\ \oplus_{i}(s_{i}\circ \hat f)^{*}T_{vert}\hat{\ex B}\arrow{r} &\ex F(\hat f)\end{tikzcd}\]
As $\ex C(\hat f^\star_{0})$ is defined by the pullback diagram (\ref{f'def}) there is an induced map 
\[\begin{tikzcd}\ex C(\hat h)\arrow{r}{\psi''}\arrow{d} &\ex C(\hat f_{0}^\star)\arrow{d}
\\\ex F(\hat h)\arrow{r} & \oplus_{i}(s_{i}\circ \hat f)^{*}T_{vert}\hat{\ex B}\end{tikzcd}\]
which is a holomorphic isomorphism on each fiber because $\psi'$ factorizes as $\psi''$ composed with the fiberwise-holomorphic-isomorphism map $\ex C(\hat f_0^\star)\longrightarrow \ex C(\hat f)$. Roughly speaking, this lift, $\psi''$, of $\psi'$ is determined by the condition that $ \hat h$ agrees with $ f_0^\star\circ \psi''$ when restricted to the pullback under $\psi'$ of the sections $s_{i}$.  The map $\psi''$ is also defined in the more general case that $\psi'$ is not necessarily fiberwise holomorphic, and in this case $\psi''$ is also not necessarily fiberwise holomorphic. 

Let us define $\nu_{\hat h}'$ analogously to  $\nu_{\hat h}$, using $\psi''$ in place of $\psi'$. Pull back our constructed trivialization (or construct another trivialization using the same connection) to give a $(G_{0}\times G)$--invariant trivialization  associated to our new family $\hat f^\star_{0}$. Again, call this trivialization $(\mathcal F,\phi)$. As before we can use this trivialization to factorize  $\hat h$ as follows:
 \[\begin{tikzcd}\ex C(\hat h)\arrow[bend left]{rr}{\hat h}\arrow{r}{ \nu'_{\hat h}}\arrow{dr}{\psi''}&(f^\star_{0})^{*}T_{vert}\hat{\ex B}\arrow{d}\arrow{r}{\mathcal F}&\hat{\ex B}
\\ &\ex C(\hat f^\star_{0})\end{tikzcd}\]
where $\nu_{\hat h}'$  is uniquely determined by the commutativity of the above diagram and the condition that it must vanish on the image of the sections $s_{i}$.

Let us modify $\hat f^\star_{0}$ using Theorem \ref{regularity theorem} in order to obtain a family containing all relevant holomorphic curves. At $f\in\hat f^\star_{0}$, we have that $D\dbar$ is injective  when restricted to the subspace sections of $f^{*}T_{vert}\hat {\ex B}$ that vanish on the image of the sections $s_{i}$. Theorem \ref{f replacement} implies that we can choose a $G_{0}\times G$--invariant,  finite-dimensional sub-vector-bundle $V$ of $\Y(\hat f^\star_{0})$ so that the pre-obstruction model $(\hat f^\star_{0},V)$ has $D\dbar\co  X^{\infty,\underline 1}(f)\longrightarrow \Y(f)$ complementary to $V(f)$. We can then apply Theorem \ref{regularity theorem} to $(\hat f^\star_{0},V)$ to obtain a unique section $\nu$ of $(\hat f^\star_{0})^{*}T_{vert}\hat{\ex B}$ defined near $f$ and vanishing on the image of all $s^\star_{i}$ so that $\dbar \nu$ is a section of $V$. (Without losing generality, we can assume that $\nu$ is defined on all of $\hat f^\star_{0}$.) Below, we shall show that the modified family $f^\star:=\mathcal F\circ \nu$ has the properties required by our proposition. 

\begin{claim}\label{nbhd} There exists a neighborhood $\mathcal O$ of $(\id, f)$ in $\Msw(\ex C(\hat f)\times \hat{\ex B})$ such that, if a family of holomorphic curves $(\psi',\hat h)$ is in $\mathcal O$, then
$\nu'_{\hat h}$ is the pullback of $\nu$ under the following diagram: 
 \[\begin{tikzcd}\ex C(\hat h)\arrow{r}{ \nu'_{\hat h}}\arrow{dr}{\psi''}&(f^\star_{0})^{*}T_{vert}\hat{\ex B}\arrow[transform canvas={xshift=0.1cm}]{d}\arrow{r}{\mathcal F}&\hat{\ex B}
\\ &\ex C(\hat f^\star_{0})\arrow[transform canvas={xshift=-0.1cm}]{u}{\nu}\end{tikzcd}\]
\end{claim}

To prove Claim \ref{nbhd}, we can choose  any neighborhood $\mathcal O$ of $(\id, f)$ so that 
\begin{enumerate}
\item the construction of $\nu'_{\hat h}$ makes sense for all  $(\psi',\hat h)$ in this neighborhood, and  

\item so that for any individual curve $(\psi',h)\in \mathcal O$, $\nu'_{\hat h}$ is in (the restriction of) the neighborhood $O$ of $0\in X^{\infty,\underline 1}(\hat f^\star)$  on which the uniqueness statement of Theorem \ref{regularity theorem}  applies. 
\end{enumerate}
As $\nu'_{\hat h}$ is natural and this second condition is an open condition on any family for which $\nu'_{\hat h}$ is defined, Lemma \ref{sequence family} implies that it is an open condition, so we can construct such an open neighborhood by choosing a neighborhood on which the first condition applies then restrict to a smaller open neighborhood so that the second condition holds.  For any family of holomorphic curves $(\psi',\hat h)$ in $\mathcal O$,  the uniqueness statement from Theorem \ref{regularity theorem} gives that $\nu_{\hat h}=\nu\circ \psi''$.   This  completes the proof of Claim \ref{nbhd}.

 Claim \ref{nbhd} implies that for $(\psi',\hat h)$ a family of holomorphic curves in $\mathcal O$,   \[\hat h=\mathcal F(\nu)\circ\psi''\] Letting $\hat f^\star:=\mathcal F(\nu)$, our map $\psi''$ therefore gives us the required unique map $\hat h\longrightarrow \hat f^\star$ so that the following diagram commutes.
 
 \[\begin{tikzcd}\ex C(\hat h)\arrow{r}{ \hat h}\arrow{d}{\psi'}\arrow{dr}{\psi''}&\hat{\ex B}
\\ \ex C(\hat f)&\arrow{l}{\psi}\ex C(\hat f^\star)\arrow{u}[swap]{\hat f^\star}\end{tikzcd}\]

Because the map $\dbar\co X^{\infty,\underline1}(\hat f^\star_{0})\longrightarrow \Y(\hat f^\star_{0})$ is defined using our $G_{0}\times G$--invariant trivialization, it  is $G_{0}\times G$--equivariant and $V$ is $G_{0}\times G$--invariant. Therefore, the unique solution $\nu$ to $\dbar\nu\in V$ must also be $G_{0}\times G$--equivariant. As the map $\mathcal F$ from our trivialization is $G_{0}$--equivariant and $G$--invariant, our family $\hat f^\star$ is $G_{0}$--equivariant and $G$--invariant as required.
 
 \stop

\

We wish to study the moduli stack of solutions to $\dbar \hat f\in V(\hat f)$ where $V$ is a simply-generated subsheaf of $\Y$ in the sense of Definition \ref{simply generated}.  The next lemma shows that any such $V$ can locally be parametrized by a family of curves $\hat f_{1}/G$ that is a vector-bundle over a core family $\hat f_{0}/G$. We prove  this in the slightly more general setting of a  moduli stack of decorated curves $\dmsw$; see Definition \ref{decorated}. We shall use this result to apply the analysis from \cite{reg}.

\begin{lemma}\label{parametrize by curves} If $V$ is a rank--$n$ simply-generated subsheaf of $\Y$, defined on a neighborhood (within $\dmsw$) of a holomorphic curve $f$ (with at least one smooth component), then there exists a neighborhood $\mathcal O$ of $f$ in $\dmsw$, a core family  $(({\fun\Phi}_0,{\fun\Phi}_0^{+1}),\hat f_{0}/G,\{s_{i}\})$ for $\mathcal O$, and a $G$--invariant family of curves,  $\hat f_{1}$ in $\mathcal O$, together with a fiberwise-holomorphic map 
\[\begin{tikzcd}\mathcal O^{+1}\arrow{r}{{\fun\Phi}_1^{+1}}\arrow{d} &\ex C(\hat f_{1})/G\arrow{d}
\\ \mathcal O\arrow{r}{{\fun\Phi}_1} &\ex F(\hat f_{1})/G\end{tikzcd}\]
 such that the following holds: \begin{enumerate}
 \item The family $\hat f_1$ parametrises $V$ in the following sense: There exists a locally-free
 $G$--invariant, rank--$n$  sheaf, $V_{1}$, of $\C\infty1(\ex F(\hat f_{1}))$--modules on $\ex F(\hat f_{1})$ which is a subsheaf of   $\Gamma^{(0,1)}(T^{*}_{vert}\ex C(\hat f_{1})\otimes T_{vert}\hat{\ex B})$ and which pulls back to give $V$ on $\mathcal O$. (See Definition \ref{Vpullback} for pulling  back  in this sense.)
 \item\label{pbc2} There is a canonical inclusion of $\hat f_{0}$ into $\hat f_1$, making $\hat f_0$  a $G$--invariant subfamily of $\hat f_{1}$, and there is a $G$--equivariant map  
 \[\psi\co \hat f_{1}\longrightarrow \hat f_{0}\]
  which is a projection in the sense that it is the identity on $\hat f_{0}\subset \hat f_1$, and  $\ex C(\hat f_{1})$ is isomorphic to a vector-bundle over $\ex C(\hat f_{0})$ with the given projection and zero section. Moreover, $\psi$ has the property that  the core family map $({\fun\Phi}_0,{\fun\Phi}_0^{+1})$ to $\ex C(\hat f_{0})/G$ factorizes as  
  \[\begin{tikzcd}\mathcal O^{+1}\arrow[swap]{r}{{\fun\Phi}_1^{+1}}\arrow{d}\arrow[bend left]{rr}{{\fun\Phi}_0^{+1}} &\ex C(\hat f_{1})/G\arrow{r}{\ex C(\psi)}\arrow{d}& \ex C(\hat f_{0})/G\arrow{d}
\\ \mathcal O\arrow{r} &\ex F(\hat f_{1})/G\arrow{r} &\ex F(\hat f_{0})/G\end{tikzcd}\]

  \item \label{pbc3}For all $t\in[0,1]$, the subspace $((1-t)D\dbar+tD\dbar^{\mathbb C})^{-1}(V(f))\subset{T_f\dmsw\ov{\hat{\ex B}}}$ does not contain the image\footnote{The image under the map from the short exact sequence (\ref{relative tangent sequence}) defining $T_f\dmsw\ov {\hat{\ex B}}$.} of any nonzero sections of $f^{*}T_{vert}\hat{\ex B}$  vanishing at the image of the marked point sections $\{s_{i}\}$.
  
 \end{enumerate}

\end{lemma}

\pf By the definition of a simply-generated subsheaf $V$ of $\Y$, Definition  \ref{simply generated}, we have that $V$ is pulled back from some sheaf over a family of exploded manifolds $\hat{\ex A}\longrightarrow \ex X$, with a group, $G_0$, of automorphisms. The purpose of this lemma is to replace $\hat{\ex A}\longrightarrow \ex X$ with the domain of a family of curves $\ex C(\hat f_1)\longrightarrow \ex F(\hat f_1)$. The main tool is Proposition \ref{model map to A}.
The following is a summary of the steps involved in the proof.
\begin{itemize}
\item[Step 1] Choose a core family $(({\fun\Phi},{\fun\Phi}^{+1}),\hat f/G_1,\{s_i\})$ containing $f$, for an open substack $\mathcal O$. Ensure condition (\ref{pbc3}) by choosing sufficiently many sections $s_i$. Moreover, choose $\mathcal O$ small enough that the simply generated subsheaf $V$ is defined over $\mathcal O$ using a fiberwise holmorphic map $(\fun\Theta,\fun\Theta^{+1})$.
\[\begin{tikzcd}\mathcal O^{+1}\arrow{r}{{\fun\Phi} ^{{+1}}}\arrow{d} &\ex C(\hat f)/G_{1} \arrow{d} && \mathcal O^{+1}\arrow{r}{\fun\Theta^{+1}}\arrow{d} &\hat{\ex A}/G_{0}\arrow{d}
\\ \mathcal O\arrow{r}{{\fun\Phi} }& \ex F(\hat f)/G_{1} && \mathcal O\arrow{r}{\fun\Theta} &\ex X/G_{0}\end{tikzcd}\]

\item[Step 2] Define the core family $(({\fun\Phi}_0,{\fun\Phi}_0^{+1}),\hat f_0/G,\{s_i\})$, with $G=G_0\times G_1$,  as follows. Applying $(\fun\Theta,\fun\Theta^{+1})$ to $\hat f$ provides a $G_0$--bundle $\hat f_0:=\hat f\bd{G_0}$ over $\hat f$, and a family $\fun\Theta^{+1}(\hat f)$ of holomorphic curves in $\hat{\ex A}$.
\[\fun\Theta^{+1}(\hat f)\co \ex C(\hat f_0)\longrightarrow \hat{\ex A}\]
 There is a canonical $G:=G_0\times G_1$ action on $\hat f_0$. Moreover  $\hat f/G_1$ is equivalent to $\hat f_0/G$. The sections $\{s_i\}$ of $\ex C(\hat f)$ induce sections of $\ex C(\hat f_0)$, which we again call $\{s_i\}$. Then  $(\hat f_0/G,\{s_i\})$ is the core family from the statement of this lemma. Define a fiberwise holomorphic map $({\fun\Phi}_0,{\fun\Phi}_0^{+1})$ using the following composition.
\[\begin{tikzcd} \mathcal O^{1}\rar{{\fun\Phi}^{+1}}\dar\ar[bend left]{rr}{{\fun\Phi}_0^{+1}} &\ex C(\hat f)/G_1 \rar{\equiv}\dar& \ex C(\hat f_0)/G\dar 
\\ \mathcal O \rar{{\fun\Phi}}\ar[bend right]{rr}{{\fun\Phi}_0} & \ex F(\hat f)/G_1 \rar{\equiv}& \ex F(\hat f_0)/G \end{tikzcd}\]
\item[Step 3] Use $({\fun\Phi}_0,{\fun\Phi}_0^{+1})$ and $(\fun\Theta,\fun\Theta^{+1})$ to construct a fiberwise holomorphic map $(\fun\Theta_1,\fun\Theta_1^{+1})$ with codomain $\hat{\ex A}/G$, where the $G_1$ factor in $G=G_0\times G_1$ acts trivially on $\hat{\ex A}$.  Applying $({\fun\Phi}_0,{\fun\Phi}_0^{+1})$ to a family $\hat h$ in $\mathcal O$ gives a $G$--bundle $\hat h\bd{G}$ over $\hat h$ with a fiberwise holomorphic map ${\fun\Phi}_0^{+1}(\hat h)\co \ex C(\hat h\bd G)\longrightarrow \ex C(\hat f)$. Moreover,   assuming $\mathcal O$ is small enough, this map will also be an isomorphism on each fiber. Applying the fiberwise holomorphic map $(\fun\Theta,\fun\Theta^{+1})$ to $\hat h$ gives a $G_0$--bundle $\hat h\bd {G_0}$ over $\hat h$ with a fiberwise holomorphic map $\fun\Theta^{+1}(\hat h)$ to $\hat{\ex A}$. After checking that the quotient of $\hat h\bd G$ by $G_1$ is $\hat h\bd{G_0}$, we then obtain a $G_1$--invariant, $G_0$--equivariant fiberwise holomorphic map $\fun\Theta^{+1}_1(\hat h)\co \ex C(\hat h\bd G)\longrightarrow \hat {\ex A}$ such that the following diagram commutes.
\[\begin{tikzcd} \ex C(\hat h\bd G)\dar \rar{\fun\Theta_1^{+1}(\hat h)} & \hat{\ex A}
\\ \ex C(\hat h\bd {G_0})\ar{ur}[swap]{\fun\Theta^{+1}(\hat h)}
\end{tikzcd}\]

\item[Step 4] Define the family of curves $\hat f_1$ that parametrises $V$.  Apply Proposition \ref{model map to A} to the family of holomorphic curves $\fun\Theta^{+1}(\hat f)\co \ex C(\hat f_0)\longrightarrow \hat{\ex A}$, to define a $G_1$--invariant and $G_0$--equivariant family of curves $(\fun\Theta^{+1}(\hat f))^\star$ in $\hat{\ex A}$, together with a canonical inclusion $\fun\Theta^{+1}(\hat f)\subset (\fun\Theta^{+1}(\hat f))^\star$ and a fiberwise holomorphic isomorphism map. 
\[\psi\co \ex C\left((\fun\Theta^{+1}(\hat f))^\star\right)\longrightarrow \ex C(\fun\Theta^{+1}(\hat f))=\ex C(\hat f_0)\]
Define the family of curves $\hat f_1$ in $\hat{\ex B}$ to have the same domain as $(\fun\Theta^{+1}(\hat f))^\star$, and define the map $\hat f_1:= \hat f_0\circ \psi \co \ex C((\fun\Theta^{+1}(\hat f))^\star)\longrightarrow \hat{\ex B}$. So, $\hat f_0$ is a $G$--equivariant subfamily of $\hat f_1$, and $\psi$ induces a $G$--equivariant morphism $\psi\co \hat f_1\longrightarrow \hat f_0$ within $\Msw(\hat{\ex B})$. 
\item[Step 5] Define the fiberwise holomorphic map $({\fun\Phi}_1,{\fun\Phi}_1^{+1})$ with codomain $\ex C(\hat f_1)/G$. Given $\hat h$ in $\mathcal O$, we have a family of holomorphic curves $({\fun\Phi}^{+1}_0(\hat h),\fun\Theta^{+1}_1(\hat h))$ in $\ex C(\hat f_0)\times \hat{\ex A}$, so Proposition \ref{model map to A}  part \ref{mmA2} implies that there exists a unique map $({\fun\Phi}^{+1}_0(\hat h),\fun\Theta^{+1}_1(\hat h))\longrightarrow (\psi,(\fun\Theta^{+1}(\hat f))^\star)$ in $\Ms(\ex C(\hat f_0)\times \hat{\ex A})$. That is, there exists a unique fiberwise holomorphic map  ${\fun\Phi}^{+1}_1(\hat h)\co \ex C(\hat h\bd G)\longrightarrow \ex C(\hat f_1)$ such that the following diagram commutes
\[\begin{tikzcd} \ex C(\hat h\bd G)\dar[swap]{{\fun\Phi}^{+1}_0(\hat h)} \ar{rr}{\fun\Theta^{+1}_1(\hat h)}\ar{drr}{{\fun\Phi}^{+1}_1(\hat h)} && \hat{\ex A}
\\ \ex C(\hat f_0) && \ex C(\hat f_1)\uar[swap]{(\fun\Theta^{+1}(\hat f))^\star}\ar{ll}{\psi} \end{tikzcd}\]
We check that, so long as $\mathcal O$ is small enough this suffices to define a fiberwise holomorphic map $({\fun\Phi}_1,{\fun\Phi}_1^{+1})$ 
\[\begin{tikzcd}\mathcal O^{+1}\dar \rar{{\fun\Phi}_1^{+1}} &\ex C(\hat f_1)/G\dar
\\ \mathcal O \rar & \ex F(\hat f_1)/G\end{tikzcd}\]
satisfying condition (\ref{pbc2}).
\item[Step 6] Verify that, on our suitably small open stack $\mathcal O$ the simply generated sheaf $V$ is pulled back from a $G$--equivariant sheaf on $\ex C(\hat f_1)$.

\end{itemize}

\noindent{\bf Step 1}:

Together with Remark \ref{decorated core family}, Proposition \ref{smooth model family} implies that there is some $\C\infty1$ core family $(\hat f/G_{1},\{s_{i}\})$ containing $f$.  As indicated by Proposition \ref{smooth model family} and Theorem \ref{f replacement}, we can construct $\hat f$ with enough marked point sections $s_{i}$ so that $((1-t)D\dbar+tD\dbar^{\mathbb C})^{-1}(V(f))$ contains no nonzero sections of $f^{*}T_{vert}\hat{\ex B}$ vanishing on the images of the $s_{i}$. This will ensure that condition \ref{pbc3} holds.
 
From the definition of a core family, there exists an open neighborhood $\mathcal O\subset\dmsw$ of $\hat f$ with a fiberwise-holomorphic map 
\[\begin{tikzcd}\mathcal O^{+1}\arrow{r}{{\fun\Phi} ^{{+1}}}\arrow{d} &\ex C(\hat f)/G_{1} \arrow{d}
\\ \mathcal O\arrow{r}{{\fun\Phi} }& \ex F(\hat f)/G_{1}\end{tikzcd}\]
in the sense of Definition \ref{fhm}. By choosing $\mathcal O$ small enough, we can ensure that for any family of curves $\hat h$ in $\mathcal O$, not only is ${\fun\Phi}^{+1}(\hat h)$  fiberwise holomorphic, this map is also an isomorphism on each fiber.   By the definition of a simply-generated subsheaf $V$ of $\Y$, Definition  \ref{simply generated}, we have that $V$ is pulled back from some sheaf over a family of exploded manifolds $\hat{\ex A}\longrightarrow \ex X$, with a group, $G_0$, of automorphisms; in particular,  if $\mathcal O$ is small enough, there is a fiberwise-holomorphic map
\[\begin{tikzcd}\mathcal O^{+1}\arrow{r}{\fun\Theta^{+1}}\arrow{d} &\hat{\ex A}/G_{0}\arrow{d}
\\ \mathcal O\arrow{r}{\fun\Theta} &\ex X/G_{0}\end{tikzcd}\]
and  sections $v_{1},\dotsc,v_{n}$ of $\Gamma^{(0,1)}(T_{vert}^{*}\hat {\ex A}\otimes T_{vert}\hat{\ex B})$  pulling back to  generate $V$. 
 
 \noindent{\bf Step 2:}
 
In particular, there is a $G_{0}$--bundle $\hat f\bd {G_0}\longrightarrow \hat f$ with a $G_{0}$--equivariant, fiberwise-holomorphic map 
\[\fun\Theta^{+1}(\hat f)\co \ex C(\hat f\bd {G_0})\longrightarrow \hat{\ex A}\ .\]
    Definition \ref{fhm} implies that the action of $G_{1}$ on $\hat f$ lifts to an action of $G_{1}$ on $\hat f\bd {G_0}$ and the family $\fun\Theta^{+1}(\hat f)$,  making the map $\fun\Theta^{+1}(\hat f)\co \ex C(\hat f\bd{G_1})\longrightarrow \hat{\ex A}$ a   $G_{1}$--invariant  and $G_{0}$--equivariant map. The group $G$ from the statement of this lemma is $G_0\times G_1$, and the core family $\hat f_{0}$ is $\hat f\bd {G_0}$, with the lift of the sections $\{s_i\}$ from $\hat f$. 
    \[f_0:=\hat f^{\bd {G_0}}\]
    Because the stacks $\ex C(\hat f)/G_1$ and $\ex C(\hat f\bd {G_0})/G$ are equivalent, $({\fun\Phi},{\fun\Phi}^{+1})$  is equivalent to a fiberwise holomorphic map $(\fun\Phi_0,\fun\Phi_0^{+1})$ as below.
\[\begin{tikzcd}\mathcal O^{+1}\ar[bend left]{rr}{{\fun\Phi}^{+1}}\arrow[swap]{r}{{\fun\Phi}_0 ^{{+1}}}\arrow{d} &\ex C(\hat f\bd {G_0})/G \arrow{d}\rar{\equiv} & \ex C(\hat f)/G_1\dar 
\\ \mathcal O\arrow{r}{{\fun\Phi}_0 }& \ex F(\hat f\bd {G_0})/G\rar{\equiv} & \ex F(\hat f)/G_1\end{tikzcd}\]
Unpacking definitions, $({\fun\Phi},\fun\Phi^{+1})$ applied to $\hat h$ is a $G_1$--bundle $\hat h\bd{G_1}\longrightarrow  \hat h$ and a $G_1$--equivariant, fiberwise-holomorphic-isomorphism map ${\fun\Phi}^{+1}(\hat h)\co\ex C(\hat h\bd {G_1})\longrightarrow \ex C(\hat f)$. Pulling back the $G_1$--equivariant  $G_0$--bundle $\ex C(\hat f\bd {G_0})\longrightarrow \ex C(\hat f)$ gives a $G_1$--equivariant $G_0$--bundle over $\ex C(\hat h\bd {G_1})$, and hence a $G_1\times G_0$ bundle $\hat h\bd {G_1\rtimes G_0}\longrightarrow \hat h$ such that the following is a $G_1\times G_0$--equivariant fiber-product diagram.
\[\begin{tikzcd}\ex C(\hat h\bd {G_1\rtimes G_0})\dar \rar{{\fun\Phi}_0^{+1}(\hat h)} & \ex C(\hat f\bd {G_0})\dar
\\ \ex C(\hat h\bd {G_1}) \rar{{\fun\Phi}^{+1}(\hat h)} &\ex C(\hat f) \end{tikzcd}\]
 
\noindent{\bf Step 3:}

We can also use ${\fun\Phi}$ to extend $(\fun\Theta,\fun\Theta^{+1})$ to a fiberwise holomorphic map  
\[\begin{tikzcd}\mathcal O^{+1}\arrow{r}{\fun\Theta_1^{+1}}\arrow{d} &\hat{\ex A}/(G_{0}\times G_1)\arrow{d}
\\ \mathcal O\arrow{r}{\fun\Theta_1} &\ex X/(G_{0}\times G_1)\end{tikzcd}\]
where $G_1$ acts trivially on $\hat{\ex A}$ and $\ex X$. In particular, for $\hat h$ a family in $\mathcal O$, let $\hat h\bd {G_1}\longrightarrow \hat h$ be the $G_1$--bundle from ${\fun\Phi}$, and let $\hat h\bd {G}\longrightarrow \hat h\bd {G_1}$ be the $G_1$--equivariant $G_0$--bundle produced by applying $\fun\Theta$ to $\hat h\bd {G_1}$. Then $\hat h\bd {G}\longrightarrow \hat h$ is a $G_0\times G_1$--bundle and hence a $G$--bundle. Given any map $\psi\co \hat g\longrightarrow \hat h$, we also get a map of $ G$--bundles $\psi\bd {G}\co \hat g\bd {G}\longrightarrow \hat h\bd {G}$ such that $(\phi\circ \psi)\bd {G}=\phi\bd {G}\circ \psi\bd {G}$. To complete the definition of $\fun\Theta_1$, the $G_0\times G_1$--invariant fiberwise-holomorphic map $\fun\Theta^{+1}_1(\hat h)$ is $\fun\Theta^{+1}(\hat h\bd {G_1})$.

We want to combine our fiberwise holomorphic maps to obtain a single fiberwise holomorphic map $({\fun\Phi}_0,\fun\Theta_1)$; for this, we must equate  the domains, $\ex C(\hat h\bd {G_1\rtimes G_0})$ of ${\fun\Phi}_0^{+1}(\hat h)$, and $\ex C(\hat h\bd {G})$ of $\fun\Theta^{+1}_1(\hat h)$.
\begin{claim}\label{equivalent bundles} Assuming $\mathcal O$ is small enough, (by replacing $\mathcal O$ with a smaller neighborhood of $\hat f_0/G$  if necessary), the following is true. For all $\hat h$ in $\mathcal O$, there is a canonical $G_0\times G_1$--equivariant map $\alpha_{\hat h}\co \hat h\bd {G}\longrightarrow \hat h\bd {G_1\rtimes G_0}$. Moreover, given any map $\beta\co \hat g\longrightarrow \hat h$, the following diagram commutes
\[\begin{tikzcd}\hat h\bd {G}\dar{\alpha_{\hat g}}\ar{rr}{\beta\bd {G}} && \hat h\bd {G}\dar{\alpha_{\hat h}}
\\ \hat g\bd {G_1\rtimes G_0}\ar{rr}{\beta\bd {G_1\rtimes G_0}}& & \hat h\bd {G_1\rtimes G_0}\end{tikzcd}\]
\end{claim}

To prove Claim \ref{equivalent bundles}, consider first the special  case of a family $\hat h$ with a map $\psi\co \hat h\bd {G_1}\longrightarrow \hat f$ such that  $\ex C(\psi)={\fun\Phi}^{+1}(\hat h)$. Then, recalling the definition of $\ex C(\hat h\bd {G_1\rtimes G_0})$ as a fiber product, we  define $\alpha_{\hat h}\co \hat h\bd {G_1\rtimes G_0}\longrightarrow \hat h\bd {G}$ to be the unique map such that the following diagram commutes.
\[\begin{tikzcd}\ex C(\hat h\bd {G})\rar{\ex C(\alpha_{\hat h})}\ar[bend right]{dr} \ar[bend left]{rr}{\ex C(\psi\bd{G_{0}})} & \ex C(\hat h\bd {G_1\rtimes G_0})\rar[swap]{{\fun\Phi}^{+1}_0(\hat h)}\dar   & \ex C(\hat f\bd{G_{0}})\dar
\\  & \ex C(\hat h\bd {G_1})\rar[swap]{\ex C(\psi)} & \ex C(\hat f)
  \end{tikzcd}\]
More generally, by choosing $\mathcal O$ small enough, we can assume that ${\fun\Phi}^{+1}(\hat h)$ is an isomorphism on each fiber. Remark \ref{deformation remark} implies that if $\mathcal O$ is small enough, there is a canonical homotopy $\hat h_t$  of $\hat h$ within $\mathcal O$  to a family of curves $\hat h_1$ with a map $\psi\co\hat h_1\bd{G_{1}}\longrightarrow \hat f$ such that $\ex C(\psi)={\fun\Phi}^{+1}(\hat h)$ and such that ${\fun\Phi}^{+1}(\hat h_t)$ is independent of $t$.  As argued above, we get a canoical map $\alpha_{\hat h_1}\co \hat h_1\bd {G_1\rtimes G_0}\longrightarrow \hat h_1\bd {G}$. This map extends uniquely along our homotopy $\hat h_t$ to define $\alpha_{\hat h}\co \hat h\bd {G_1\rtimes G_0}\longrightarrow \hat h\bd {G}$. Any map $\beta\co\hat g\longrightarrow \hat h$ extends uniquely to $\beta_t\co \hat g_t\longrightarrow \hat h_t$ along the homotopy from Remark \ref{deformation remark}, so it follows that the diagram from Claim \ref{equivalent bundles} commutes.


\

With Claim \ref{equivalent bundles} proved, we now have a fiberwise-holomorphic map $({\fun\Phi}_0,\fun\Theta_1)$ 
\[\begin{tikzcd}\mathcal O^{+1}\arrow{rr}{({\fun\Phi}_0 ,\fun\Theta_1)^{+1}}\arrow{d}& &(\ex C(\hat f\bd {G_1})\times \hat {\ex A})/G\arrow{d}
\\ \mathcal O\arrow{rr}{({\fun\Phi}_0,\fun\Theta_1)} &&(\ex F(\hat f\bd {G_1})\times \ex X)/G\end{tikzcd}\]
which associates to a family $\hat h$ in $\mathcal O$ the $G$--bundle $\hat h\bd {G}$ and the $G$--equivariant map $({\fun\Phi}_0,\fun\Theta_1)^{+1}(\hat h)\co \ex C(\hat h\bd {G})\longrightarrow \ex C(\hat f_0)\times \ex A$ defined as
\[({\fun\Phi}_0,\fun\Theta_1)^{+1}(\hat h)=\lrb{{\fun\Phi}_0^{+1}(\hat h)\circ \ex C(\alpha_{\hat h}),\fun\Theta_1^{+1}(\hat h)}\ .\]
 So, for any family $\hat h$ in $\mathcal O$, we have $({\fun\Phi}_0 ,\fun\Theta_1)^{+1}(\hat h)$ is a $G$--equivariant family of holomorphic  curves in $\ex C(\hat f_0)\times \hat{\ex A}\longrightarrow \ex F(\hat f_0)\times \ex X$, with domain $\ex C(\hat h\bd {G})$.

\noindent{\bf Step 4:}

Note that the map $\fun\Theta^{+1}(\hat f)\co \ex C(\hat f_0)\longrightarrow \hat{\ex A}$ is $G_1$ invariant and $G_0$ equivariant, so we refer to the family of curves $\fun\Theta^{+1}(\hat f)$ as being $G_1$ invariant and $G_0$ equivariant. Use Proposition \ref{model map to A} to extend the family of holomorphic curves $\fun\Theta^{+1}(\hat f)\co\ex C(\hat f_0)\longrightarrow \hat{\ex A}$ to a $G_1$--invariant and $G_0$--equivariant family of curves $\fun\Theta^{+1}(\hat f)^\star$ in $\hat {\ex A}$ with a map $\psi\co \ex C(\fun\Theta^{+1}(\hat f)^\star)\longrightarrow \ex C(\hat f_0)$ so that, given any family of curves $\hat h$ close enough to $ f$, the family of holomorphic curves  $({\fun\Phi}_0 ,\fun\Theta_1)^{+1}(\hat h)$ in $\ex C(\hat f_0)\times \hat{\ex A}$ has a unique map  \begin{equation}\label{psidef}\psi_{\hat h}\co ({\fun\Phi}_0 ,\fun\Theta_1)^{+1}(\hat h)\longrightarrow (\psi,\fun\Theta^{+1}(\hat f)^\star)\end{equation} in $\Msw(\ex C(\hat f_0)\times \hat{\ex A})$. As this map is unique and both families of curves are $G$--equivariant, $\psi_{\hat h}$ must also be $G$--equivariant.
 If necessary, replace $\mathcal O$ with a smaller open neighborhood of $f$ so that the above property holds for all families $\hat h$  in $\mathcal O$. 
 Note that  given any map $\beta\co \hat h\longrightarrow \hat g$, 
 \begin{equation}\label{beta compatibility}\ex C(\psi_{\hat h})=\ex C(\psi_{\hat g})\circ \ex C(\beta\bd {G})\ .\end{equation}

We are now ready to define the family of curves $\hat f_1$ that parametrises $V$.
\begin{construction}\label{f1def} Define the family of curves $\hat f_1$ in $\mathcal O$ to be the pullback  of $\hat f_0:=\hat f^{\bd {G_{0}}}$ under the map $\psi\co \ex C(\fun\Theta^{+1}(\hat f)^\star)\longrightarrow \ex C(\hat f_0)$. So, $\ex C(\hat f_1)=\ex C(\fun\Theta^{+1}(\hat f)^\star)$, and 
\[\hat f_1=\hat f_0\circ \psi\ .\]
\end{construction}
Note that $G$ still acts as a group of automorphisms of  $\hat f_1$.

\noindent{\bf Step 5}

Now, define our fiberwise holomorphic map $({\fun\Phi}_1,{\fun\Phi}_1^{+1})$ 
\[\begin{tikzcd}\mathcal O^{+1}\dar \rar{{\fun\Phi}_1^{+1}} & \ex C(\hat f_1)/G\dar
\\ \mathcal O\rar{{\fun\Phi}_1} & \ex F(\hat f_1)/G\end{tikzcd}\]
as follows: ${\fun\Phi}^{+1}_1$ applied to $\hat h$ is the $G$--bundle $\hat h\bd {G}\longrightarrow \hat h$, together with the $G$--equivariant fiberwise-holomorphic map
\[{\fun\Phi}_1^{+1}(\hat h):=\ex C(\psi_{\hat h})\co \ex C(\hat h\bd {G})\longrightarrow \ex C(\hat f_1)\ .\]
Equation (\ref{beta compatibility}) ensures that ${\fun\Phi}^{+1}_1$ satisfies Definition \ref{fhm}. Note also that $\psi\circ{\fun\Phi}^{+1}_0(\hat h)={\fun\Phi}^{+1}_1$, so we have verified condition \ref{pbc2} of our lemma.

Recall from Construction \ref{f1def}, that $\ex C(\hat f_1)=\ex C(\fun\Theta^{+1}(\hat f)^\star)$, and $\fun\Theta^{+1}(\hat f)^\star$ is a $G_1$--invariant, and $G_0$--equivariant family of curves in $\ex A$.
\[\begin{tikzcd}\ex C(\hat f_1)\dar\ar{rr}{\fun\Theta^{+1}(\hat f)^\star}& & \ex A\dar
\\ \ex F(\hat f_1)\ar{rr}& & \ex X\end{tikzcd}\]
Accordingly, composition with  $\fun\Theta^{+1}(\hat f)^\star$  defines a map of stacks. 
\[\begin{tikzcd}\ex C(\hat f_1)/G\dar \ar{rr}{\fun\Theta^{+1}(\hat f)^\star_G}& &\ex A/ G\dar
\\ \ex F(\hat f_1)/G \ar{rr} && \ex X/G \end{tikzcd}\]

\begin{claim}\label{pullback equivalence} $ \fun\Theta^{+1}(\hat f)^\star_G\circ {\fun\Phi}_1^{+1}=\fun\Theta_1^{+1}$.
\end{claim}

Proving Claim \ref{pullback equivalence} is a matter of tracing back through our definitions (which have been chosen  so that this is a strict equality.) Let us apply these maps to $\hat h$. First, applying ${\fun\Phi}_1^{+1}$, we get a $G$--bundle $\hat h\bd {G}\longrightarrow \hat h$ with a map $\ex C(\psi_{\hat h})\co \ex C(\hat h\bd {G})\longrightarrow \ex C(\hat f_1)=\ex C(\fun\Theta^{+1}(\hat f)^\star)$. The definition of $\psi_{\hat h}$ above equation (\ref{psidef}) then implies that  $\fun\Theta^{+1}(\hat f)^\star\circ {\fun\Phi}_1^{+1}(\hat h)=\fun\Theta_1^{+1}(\hat h)\co \ex C(\hat h\bd {G})\longrightarrow \ex A$. 

\noindent{\bf Step 6:}

With the help of Claim \ref{pullback equivalence}, we  examine how to parametrize our obstruction bundle $V$ using  $\hat f_{1}$. The map 
\[(\fun\Theta^{+1}(\hat f)^\star,\id)\co \ex C(\hat f_{1})\times \hat{\ex B}\longrightarrow \hat {\ex A}\times\hat{\ex B}\] pulls back the sections $v_{1},\dotsc,v_{n}$ of $\Gamma^{(0,1)}(T^{*}_{vert}\hat{\ex A}\otimes T_{vert}\hat{\ex B})$ to $G_{1}$--equivaraint sections
\[(\fun\Theta^{+1}(\hat f)^\star)^*v_{1},\dotsc, (\fun\Theta^{+1}(\hat f)^\star)^*v_{n}\] of $\Gamma (T^{*}_{vert}\ex C(\hat f_{1})\otimes T_{vert}\hat {\ex B})$. Because $\fun\Theta^{+1}(\hat f)^\star$ is not fiberwise holomorphic everywhere, these sections may not be in $\Gamma^{(0,1)}$, so we take their projection $((\fun\Theta^{+1}(\hat f)^\star)^*v_{i})^{(0,1)}$ to $\Gamma^{(0,1)}$. This projection does not affect our section on the image of ${\fun\Phi}_{1}^{+1}(\hat h)$, because $\fun\Theta^{+1}(\hat f)^\star$ is fiberwise holomorphic there. 

 Given any family $\hat h$ in $\mathcal O$, Claim \ref{pullback equivalence} implies that
\[({\fun\Phi}_1^{+1}(\hat h),\hat h\bd {G})^*((\fun\Theta^{+1}(\hat f)^\star)^*v_{i})^{(0,1)}=(\fun\Theta_{1}^{+1},\hat h\bd G)^{*}v_{i}\]
So, the map
\[({\fun\Phi}_1^{+1}(\hat h),\hat h\bd {G})\co \ex C(\hat h\bd {G})\longrightarrow \ex C(\hat f_{1})\times \hat{\ex B}\]
 pulls back these sections to $G_{1}$--equivariant sections 
 \[(\fun\Theta_1^{+1}(\hat h),\hat h\bd {G})^*v_{1},\dotsc, (\fun\Theta_1^{+1}(\hat h),\hat h\bd {G})^*v_{n}\] of $\Y(\hat h\bd {G})$. Because $\fun\Theta_1^{+1}(\hat h)=\fun\Theta^{+1}(\hat h\bd {G_1})$ these $G_1$--equivariant sections are the same as those obtained using $\fun\Theta^{+1}$ applied to $\hat h\bd {G_1}$.
 \[(\fun\Theta_1^{+1}(\hat h),\hat h\bd {G})^*v_i=(\fun\Theta^{+1}(\hat h\bd {G_1}),\hat h\bd {G})^*v_i\]
 So,  these sections of $\Y(\hat h\bd {G})$ are linearly independent at every curve in $\hat h\bd {G}$ and generate $V(\hat h\bd {G})$.  

In other words, the pullback of the sections $((\fun\Theta^{+1}(\hat f)^\star)^*v_{i})^{(0,1)}$ generate $V$ on $\mathcal O$. As the sheaf of $\C\infty1(\ex X)$--modules  over $\ex X$ generated by $v_{1},\dotsc,v_{n}$ is $G_{0}$--invariant and the map $\fun\Theta^{+1}(\hat f)^\star$ is $G_{1}$--invariant and $G_{0}$--equivariant, the sheaf of  $\C\infty1(\ex F)$--modules over $\ex F(\hat f_{1})$ generated by $((\fun\Theta^{+1}(\hat f)^\star)^*v_{i})^{(0,1)}$ is $G$--invariant. As the pullback of $((\fun\Theta^{+1}(\hat f)^\star)^*v_{i})^{(0,1)}$ to any curve in $\mathcal O$ are linearly independent, restricted to a neighborhood of the image of $\mathcal O$ in $\ex F(\hat f_{1})$, the sheaf $V_{1}$ generated by $((\fun\Theta^{+1}(\hat f)^\star)^*v_{i})^{(0,1)}$ is free and has rank $n$. Therefore, by restricting $\hat f_{1}$ to a $G$--invariant neighborhood of the image of $\mathcal O$, the sections $((\fun\Theta^{+1}(\hat f)^\star)^*v_{i})^{(0,1)}$ generate a rank--$n$ $G$--invariant vector-bundle $V_{1}$,   pulling back to give $V$ on $\mathcal O$.

\stop
 
   We are ready to prove our main theorem on the local structure of the moduli stack of holomorphic curves. The statement uses definitions \ref{Ydef}, \ref{decorated}, \ref{open substack}, \ref{simply generated} and  \ref{represented by family}.

\begin{thm}\label{V moduli stack} Let $V$ be a simply-generated subsheaf of $\Y$, defined on an open substack of $\dmsw$ such that, for some holomorphic curve $f$, we have that $V(f)\subset \Y(f)$  is transverse to $D\dbar\co T_{f}\dmsw(\hat {\ex B})\longrightarrow \Y(f)$. Then there exists an open neighborhood $\mathcal O\subset\dmsw$ of $f$ and a family of curves $\hat f$ in $\dmsw$ with automorphism group $G$ such that $\hat f/G$ represents the substack $\dbar^{-1} V\subset \mathcal O$ consisting of curves $h$ in $\mathcal O$ with the property that $\dbar h\in V(h)$.

Moreover, the map
\[T_f\hat f\co T_{f}\ex F(\hat f)\longrightarrow T_{f}\Msw(\hat{\ex B})\] 
is injective and has image  $D\dbar^{-1}(V(f))\subset T_{f}\Msw(\hat{\ex B})$. 

\end{thm}

\pf

See Lemma \ref{T stack} for the easy case,  in which the domain of $f$ is $\ex T$.  With this case excluded, we can assume that the domain of $f$ has at least one smooth component, so  Lemma \ref{parametrize by curves} applies.
The following is a summary of the steps in the proof.

\begin{itemize}
\item[Step 1] Using Lemma \ref{parametrize by curves}, construct a $G$--invariant family $\hat f_1$  that parametrises $V$. As part of Lemma \ref{parametrize by curves}, we obtain a fiberwise holomorphic map $(\fun\Phi_1,\fun\Phi_1^{+1})$ from a neighbourhood $\mathcal O$ of $f$ to the domain of $\hat f_1/G$, and a $G$--equivariant vector-bundle $V_1$ over $\ex F(\hat f_1)$. Construct a new family of curves $\hat f_2$ as the pullback of $\hat f_1$ using $V_1\longrightarrow \ex F(\hat f_1)$. 
\item[Step 2] On the substack $\dbar^{-1}V$, there is a natural lift of $(\fun\Phi_1,\fun\Phi_1^{+1})$ to a fiberwise holomorphic map 
\[\begin{tikzcd} (\dbar^{-1}V)^{+1}\dar \rar{{\fun\Phi}_2^{+1}} & \ex C(\hat f_2)/G\dar
\\ \dbar^{-1}V \rar{{\fun\Phi}_2} & \ex F(\hat f_2)/G\end{tikzcd}\]
defined as follows.

Given a family of curves $\hat h$ in $\mathcal O$,  a section $\theta_{\hat h}\in V(\hat h)$ corresponds to a section of $\fun\Phi_1(\hat h)^* V_1$. As $\hat f_2$ is defined as a pullback over $V_1\longrightarrow \ex F(\hat f_1)$, the  map $\fun\Phi_1^{+1}(\hat h)\co \ex C(\hat h\bd G)\longrightarrow \ex C(\hat f_1)$ induces a canonical map 
\[\psi_{\theta_{\hat h}}\co \ex C(\hat h \bd G)\longrightarrow \ex C(\hat f_2) \]
depending on $\hat h$ and $\theta_{\hat h}$. This fiberwise holomorphic map is related to  a tautological section $\theta$ of $\Gamma^{(0,1)}\lrb{T^{*}_{vert}\ex C(\hat f_{2})\otimes T_{vert}\hat{\ex B}}$ such that 
 \[(\psi_{\theta_{\hat h}},\hat h\bd {G})^*\theta=\theta_{\hat h}\]

For a family  $\hat h$ in the stack $\dbar^{-1}V$, we have $\dbar \hat h\in V(\hat h)$, so we obtain a canonical map $\psi_{\dbar\hat h}\co \ex C(\hat h \bd G)\longrightarrow \ex C(\hat f_2)$. This in turn defines the desired fiberwise holomorphic map $(\fun\Phi_2,\fun\Phi_2^{+1})$ by setting
 \[{\fun\Phi}^{+1}_2(\hat h):=\psi_{\dbar\hat h} \co \ex C(\hat h\bd {G})\longrightarrow \ex C(\hat f_2)\ .\]

 \item[Step 3] Modify the map $\hat f_2\co \ex C(\hat f_2)\longrightarrow \ex B$ to a map $\mathcal F(\nu)\co \ex C(\hat f_2)\longrightarrow \hat{\ex B}$ using Theorem \ref{regularity theorem}. In particular $\mathcal F(\nu)$ is a family of curves with the same domain as $\hat f_2$, but  satisfying an equation roughly of the form
 \[\dbar \mathcal F(\nu)-\theta\in V_2\]
 where $V_2$ is a obstruction bundle over $\hat f_2$ and $\theta$ is related to the tautological section above. (A precise description involves a choice of trivialisation $(\mathcal F,\phi)$ associated to $\hat f_2$.)
 \item[Step 4] Define a $G$--invariant family $\hat f_3\subset \mathcal F(\nu)$ as the transverse intersection of $\dbar \mathcal F(\nu)-\theta$ with the zero section. In particular, $\hat f_3$  is in $\dbar^{-1}V$, so $\dbar \hat f_3\in V(\hat f_3)$. For any family of curves $\hat h$ with $\dbar \hat h\in V(\hat h)$, the map $\fun\Phi_2^{+1}(\hat h)\co \ex C(\hat h\bd G)\longrightarrow \ex C(\hat f_2)$ has image within $\ex C(\hat f_3)\subset \ex C(\hat f_2)$, and moreover, this map is compatible with the maps $\hat h\bd G$ and $\hat f_3$ so we obtain a $G$--equivariant map of families $\hat h\bd G\longrightarrow \hat f_3$. This defines a map of stacks 
 \[\fun\Psi\co \dbar^{-1}V\longrightarrow \hat f_{3}/G\  \]
 with a left inverse
 \[\fun\Phi_{\mathcal X}\co \hat f_3/G\longrightarrow \dbar^{-1}V\]
 given by the  natural map $\fun\Phi_{\mathcal X}$ from Definition \ref{family quotient stack}.
 \item[Step 5] We have two maps from $\hat f_3$ to $\hat f_3/G$: the quotient map and $\fun\Psi(\hat f_3)$. These two maps do not coincide. Construct $\hat f\subset \hat f_3$ as the $G$--invariant subfamily on which these two maps coincide. Moreover,  prove that the map of stacks $\fun \Psi$ has image inside $\hat f/G\subset \hat f_3/G$. 
 \item[Step 6] Prove that $\fun \Psi\co \dbar^{-1}V\longrightarrow \hat f/G$ is an equivalence of stacks, with inverse $\fun\Phi_{\mathcal X}$. 
 \item[Step 7] Verify that $T_f\hat f\co T_f\ex F(\hat f)\longrightarrow D\dbar^{-1}V(f)$ is a bijection.
 
  \end{itemize}

\noindent{\bf Step 1:}

Use Lemma \ref{parametrize by curves} to construct a core family $(({\fun\Phi}_0,{\fun\Phi}_0^{+1}),\hat f_{0}/G,\{s_{i}\})$ for an open neighborhood $\mathcal O$ of $f$,  and a $G$--invariant extension, $\hat f_{1}$ of $\hat f_{0}$,  parametrizing $V$ in the sense of Lemma \ref{parametrize by curves}.  In particular, there is a $G$--invariant sub-bundle $V_{1}$ of $\Gamma^{(0,1)}\lrb{T^{*}_{vert}\ex C(\hat f_{1})\otimes T_{vert}\hat{\ex B}}$ so that, 
 given any family of curves $\hat h$ in $\mathcal O$, there is a $G$--bundle $\hat h\bd G \longrightarrow\hat h$ and a lift of the core-family   map ${\fun\Phi}_0^{+1}(\hat h)\co \ex C(\hat h\bd {G})\longrightarrow \ex C(\hat f_{0})$  to another $G$--equivariant map 
  \[\begin{tikzcd}\ex C(\hat h\bd {G})\arrow{d}\arrow{r}{{\fun\Phi}_{1}^{+1}(\hat h)}&\ex C(\hat f_{1}) \arrow{d}
 \\ \ex F(\hat h\bd {G})\arrow{r}{{\fun\Phi}_1(\hat h)} &\ex F(\hat f_{1})\end{tikzcd}\] 
such that the pullback of  $V_{1}$ using $({\fun\Phi}_{1}^{+1}(\hat h),\hat h\bd {G})$  is $V(\hat h\bd {G})$.

\begin{construction} \label{f2def}
Consider $V_{1}$ as a vector-bundle over $\ex F(\hat f_{1})$.  Let $\hat f_{2}$ be the family of curves defined using the following pullback diagram:
\[\begin{tikzcd}\ex C(\hat f_{2})\arrow[bend left]{rr}{\hat f_{2}}\arrow{r}{r}\arrow{d}&\ex C(\hat f_{1})\arrow{d}\arrow{r}{\hat f_{1}}&\hat{\ex B}
\\ \ex F(\hat f_{2})=V_{1}\arrow{r} &\ex F(\hat f_{1}) \end{tikzcd}\]
\end{construction}
We will regard $\hat f_{1}$  as the subfamily of $\hat f_{2}$ corresponding to the zero section of the bundle $V_{1}\longrightarrow \ex F(\hat f_{1})$.

\

\noindent{\bf Step 2:} 

We now define a fiberwise holomorphic map $(\fun\Phi_2,\fun \Phi_2^{+1})$ from $\dbar^{-1}(V)$ to the codomain of $\hat f_2/G$, and an associated tautological section $\theta$. 

There is a $G$--equivariant tautological section $\theta$ of $\Gamma^{(0,1)}\lrb{T^{*}_{vert}\ex C(\hat f_{2})\otimes T_{vert}\hat{\ex B}} $ defined as follows: A   point $(f',v)$ in $\ex F(\hat f_2)$ corresponds to a curve $f'$ in $\hat f_1$ and a vector $v\in V(f')$. By the definition of $V_1$, this $v$ corresponds to a section in $\Gamma^{(0,1)}\lrb{T^{*}\ex C( f')\otimes T_{vert}\hat{\ex B}}$, and we can pull this section back using the isomorphism $r\co\ex C(f',v)\longrightarrow \ex C(f')$. So, we define
\begin{equation}\label{tautological defn} \theta (f',v):=r^{*}v\ .\end{equation}
This tautological section $\theta$ has the following property: given any family of curves $\hat h$ in $\mathcal O$  and a section $\theta_{\hat h}$ of $V(\hat h)$, the canonical map $ {\fun\Phi}_{1}^{+1}(\hat h)\co \ex C(\hat h\bd {G})\longrightarrow \ex C(\hat f_{1})$ lifts  uniquely to a canonical $G$--equivariant map $\psi_{\theta_{\hat h}}$ such that 
such that the following diagram commutes
\[\begin{tikzcd}\ex C(\hat h\bd {G})\rar[swap]{\psi_{\theta_{\hat h}}}\dar\ar[bend left]{rr}{{\fun\Phi}_1^{+1}(\hat h)}&\ex C(\hat f_{2})\dar \rar{r} & \ex C(\hat f_1)\dar
 \\ \ex F(\hat h\bd {G})\rar &\ex F(\hat f_{2})\rar & \ex F(\hat f_1)\end{tikzcd}\] 
and such that 
\begin{equation}\label{psi2theta}(\psi_{\theta_{\hat h}},\hat h\bd {G})^*\theta=\theta_{\hat h} \end{equation}
where we have abused notation slightly on the right to identify $\theta_{\hat h}$ with its pullback to $\Y(\hat h\bd {G})$ using the $G$--bundle map $\hat h\bd {G}\longrightarrow \hat h$. Given a map $\alpha\co \hat g\longrightarrow \hat h$, the uniqueness of $\psi_{\theta_{\hat h}}$ implies that the following diagram commutes.
\begin{equation}\label{psi2theta commute}\begin{tikzcd}\ex C(\hat h\bd {G})\rar{\psi_{\theta_{\hat h}}} & \ex C(\hat f_2)
\\ \ex C(\hat g\bd {G}) \uar{\ex C(\alpha\bd {G})}\ar[swap]{ur}{\psi_{\alpha^*\theta_{\hat h}}}\end{tikzcd}\end{equation}

  A particular case of interest is when $\theta_{\hat h}=\dbar\hat h$, when we obtain a canonical $G$--equivariant map 
\begin{equation}\label{phi2e}{\fun\Phi}^{+1}_2(\hat h):=\psi_{\dbar\hat h} \co \ex C(\hat h\bd {G})\longrightarrow \ex C(\hat f_2) \ .\end{equation}
This map is only defined when $\dbar \hat h$ is a section of $V(\hat h)$, but for such $\hat h$ we get another fiberwise holomorphic map, $({\fun\Phi}_2,\fun\Phi_2^{+1})$.
\begin{construction}Define the fiberwise holomorphic map 
\[\begin{tikzcd}(\dbar^{-1}V)^{+1}\dar \rar{{\fun\Phi}_2^{+1}} & \ex C(\hat f_2)/G\dar
\\ \dbar^{-1}V \rar{{\fun\Phi}_2} & \ex F(\hat f_2)/G \end{tikzcd}\]
as follows. To $\hat h\in \dbar^{-1}V\subset \mathcal O$, associate the same $G$--bundle $\hat h\bd {G}\longrightarrow \hat h$ as for the maps ${\fun\Phi}_1$ and ${\fun\Phi}_0$ (along with the same $G$--bundle maps $\psi\bd {G}\co\hat g\bd {G}\longrightarrow \hat h\bd {G}$ for morphisms $\psi\co \hat g\longrightarrow \hat h$). Then define ${\fun\Phi}_2^{+1}(\hat h)$ using (\ref{phi2e}).
\end{construction}
Given any map $\alpha\co\hat g\longrightarrow \hat h$, we have that  $\alpha^*\dbar \hat h=\dbar\hat g$, so the commutative diagram (\ref{psi2theta commute}) implies that \[{\fun\Phi}_2^{+1}(\hat g)= {\fun\Phi}_2^{+1}(\hat h)\circ \ex C(\alpha\bd {G})\] so there exists a unique fiberwise holomorphic map $(\fun\Phi_2,\fun\Phi_2^{+1})$ satisfying  Definition \ref{fhm}, with $\fun\Phi_2^{+1}(\hat h)$ as defined above.

\

\noindent{\bf Step 3:}

Close to the image of $\ex C(f)$ within  $\ex C(\hat f_2)$,  the image of ${\fun\Phi}_2^{+1}$ will be a subfamily $\ex C(\hat f_3)\subset \ex C(\hat f_2)$; see Construction \ref{f_3} below. In what follows, we modify $\hat f_2$ to a family $\mathcal F(\nu)\co \ex C(\hat f_2)\longrightarrow \hat{\ex B}$ such that the restriction of $\mathcal F(\nu)$ to $\ex C(\hat f_3)$ is a family, $\hat f_{3}$,  in $\dbar^{-1}V$.

Recall that $\hat f_{1}$ is the pullback of $\hat f_{0}$, and $\hat f_{2}$ is the pullback of $\hat f_{1}$. 
Therefore, the sections $\{s_{i}\}$ of $\ex C(\hat f_{0})\longrightarrow \ex F(\hat f_{0})$ pull back to sections  $\{s'_{i}\}$ of $\ex C(\hat f_{2})\longrightarrow \ex F(\hat f_{2})$. As the set of sections $\{s_{i}\}$ is $G$--invariant and the maps $\hat f_{2}\longrightarrow \hat f_{1}\longrightarrow \hat f_{0}$ are $G$--equivariant, the set of sections $\{s'_{i}\}$  is $G$--invariant. 

Choose a $G$--invariant trivialization $(\mathcal F,\Phi)$  associated to $\hat f_{2}$, and use this trivialization and the tautological section $\theta$ from (\ref{tautological defn}) to define a simple perturbation of $\dbar$
\[\dbar'\co X^{\infty,\underline 1}(\hat f_{2})\longrightarrow \Y(\hat f_{2})\]
\[\dbar'(\nu):=\dbar \nu-\Phi((\id,\mathcal F(\nu))^{*}\theta)\]
satisfying Definition \ref{simple perturbation}. See  Example \ref{simple construction} for further details.   Note that $\dbar'$ is a $G$--equivariant map.

Recall, from Lemma \ref{parametrize by curves}, that we can choose $(\hat f_{0}/G,\{s_{i}\})$ so that $D\dbar$ is injective when restricted to the space of sections of $f^{*}T_{vert}\hat{\ex B}$  vanishing at the image of $\{s_{i}\}$. As $\theta$ vanishes on $\ex C(f)\subset \ex C(\hat f_{2})$, we also get that $D\dbar'$ is  injective when restricted to the space of sections of $f^{*}T_{vert}\hat{\ex B}$ vanishing at the image of $\{s_{i}\}$. The space $X^{\infty,\underline 1}(f)$ is the space of such sections. We can therefore choose 
  a $G$--invariant obstruction bundle $V_{2}$ over $\hat f_{2}$ such that 
  \[D\dbar'\co X^{\infty,\underline 1}(f)\longrightarrow \Y(f)\]
  has image complementary  to $V_{2}$. Note that $V_{2}$ need not be related to $V$.
  
 Apply Theorem \ref{regularity theorem} to $(\hat f_{2},V_{2})$ to obtain a unique section $\nu$ in $X^{\infty,\underline 1}(\hat f_{2})$ defined in a neighborhood of $f$ in $\hat f_{2}$ such that $\dbar' \nu$ is a section of $V_{2}$. Using the uniqueness statement in Theorem \ref{regularity theorem} and the fact that $\dbar'$ is $G$--equivariant, we have that $\nu$ is $G$--equivariant.

If necessary, reduce the size of $\mathcal O$ so that every curve in $\mathcal O$ is isomorphic to some curve in $\mathcal F(O)$ where $O$ is the neighborhood referred to in Theorem \ref{regularity theorem}. Remark \ref{deformation remark} implies that this is possible because $(\hat f_{0}/G,\{s_{i}\})$ is a core family.
 Then, we have the following. 
\begin{claim}\label{hFnu}Suppose that $\hat h$ is a family in $\dbar^{-1}V\subset \mathcal O$. Then 
\[\hat h\bd {G}=\mathcal F(\nu)\circ {\fun\Phi}_{2}^{+1}(\hat h)\]
and $\dbar' \nu=0$ on the image of ${\fun\Phi}_{2}(\hat h)$.
\end{claim}

To prove Claim \ref{hFnu}, it  suffices to consider the case of a single curve $h$. We have that  $h\bd {G}=\mathcal F(\nu_{h})\circ {\fun\Phi}_{2}^{+1}(h)$ for some section $\nu_{h}$ in the neighborhood $O$ from Theorem \ref{regularity theorem}. Using (\ref{phi2e}), we deduce that  $({\fun\Phi}_{2}^{+1})^{*}\theta=\dbar h\bd {G}$, so  ${\fun\Phi}_{2}(h)^{*}\dbar'\nu_{h}=0$, and we conclude that $\dbar'\nu_{h}$ vanishes on the image of ${\fun\Phi}_{2}( h)$. The uniqueness statement from Theorem \ref{regularity theorem} then implies that $\nu_{h}$ and $\nu$ must coincide on the image of ${\fun\Phi}_{2}( h)$. Therefore, $ h\bd {G}=\mathcal F(\nu)\circ {\fun\Phi}_{2}^{+1}( h)$, as required for Claim \ref{hFnu}.

\

\noindent{\bf Step 4:}

We now define a subfamily $\hat f_3\subset \mathcal F(\nu)$ using the transverse intersection of $\dbar'\nu$ with the $0$--section, and upgrade the fiberwise holomorphic map $(\fun\Phi_2,\fun\Phi_2^{+1})$ to a map of stacks $\fun\Psi\co \dbar^{-1}(V)\longrightarrow \hat f_3/G$.

\begin{claim}\label{dbar'nu} The section $\dbar'\nu$ is transverse to  the $0$--section of $V_{2}$ at the curve $f$ in $\hat f_{2}$.
\end{claim}
To prove Claim \ref{dbar'nu} above, we  need  the condition that \[D\dbar\co T_{f}\dmsw(\hat{\ex B})\longrightarrow \Y(f)\] is transverse to $V$. This transversality and Lemma \ref{tangent1} imply that, given any element  $v_{0}$  of $V_{2}(f)$, there exists a family $\hat h$ parametrized by $\mathbb R$ and containing $f$ at $0\in\mathbb R$, and a section $\theta_{\hat h}$ of $V(\hat h)$ so that the derivative at $0$ of $\dbar \hat h-\theta_{\hat h}$ is $v_{0}$. Recall that the tautological section $\theta$ of $\Gamma^{(0,1)}(T^{*}_{vert}\ex C(\hat f_{2})\otimes T_{vert}\hat{\ex B})$ is defined so that there exists a canonical lift of the $G$--equivariant map $ {\fun\Phi}_{1}^{+1}(\hat h)\co \ex C(\hat h\bd {G})\longrightarrow \ex C(\hat f_{1})$ to a map
%
\[\psi_{\theta_{\hat h}}\co \ex C(\hat h\bd {G})\longrightarrow \ex C(\hat f_{2})\]
such that \[\psi_{\theta_{\hat h}}^{*}\theta=\theta_{\hat h}\ ,\]  or more accurately, $\psi_{\theta_{\hat h}}^{*}\theta$ is equal to the lift of $\theta_{\hat h}$ to $V(\hat h\bd {G})$.
We can then express $\hat h\bd {G}$ using the trivialization $(\mathcal F, \Phi)$  associated to $\hat f_{2}$, so there is a section $\nu_{\hat h\bd {G}}$ of $\psi_{\theta_{\hat h}}^{*}\hat f_{2}^{*}T_{vert}\hat{\ex B}$ such that $\hat h\bd {G}=\mathcal F(\psi_{\theta_{\hat h}},\nu_{\hat h\bd {G}})$. 
The map $\psi_{\theta_{\hat h}}$ is constructed so that 
\[\dbar'\nu_{\hat h\bd {G}}=\Phi(\dbar \hat h-\theta_{\hat h})\ .\] 
In particular, $\dbar'\nu_{\hat h\bd {G}}$ is zero at $0$ and has first derivative equal to $v_{0}$. Theorem \ref{regularity theorem} implies that there is a unique section $\nu'_{\hat h\bd {G}}$ of $\psi_{\theta_{\hat h}}^{*}\hat f_{2}^{*}T_{vert}\hat{\ex B}$ such that $\dbar'\nu'_{\hat h\bd {G}}\in V_{2}$, and Corollary \ref{trt} implies that that $\nu'_{\hat h\bd {G}}$ and $\nu_{\hat h\bd {G}}$ are equal to first order at $0$.  It follows that the derivative at $0$ of $\dbar'\nu'_{\hat h\bd {G}}$ is equal to $v_{0}$. The uniqueness part of Theorem \ref{regularity theorem} implies that $\nu'_{\hat h\bd {G}}=\psi_{\theta_{\hat h}}^{*}\nu$, so  $v_{0}$ must be in the image of the derivative of $\dbar'\nu$ at the curve $f$ in $\hat f_{2}$. As this argument holds for any $v_{0}$ in $V_{2}(f)$, it follows that $\dbar'\nu$ is transverse to $0$ at $f=0$, and the proof of Claim \ref{dbar'nu} is complete.

\

With Claim \ref{dbar'nu} complete, we can shrink $\hat f_{2}$ to an open subfamily (and shrink $\mathcal O$ to a smaller open substack accordingly) until we can assume that  $\dbar'\nu$ is transverse to $0$, and $D\dbar'$ is transverse to $V_{2}$ at $\nu$. 

 \begin{construction}\label{f_3}Let  $\hat f_{3}$ be the subfamily of $\mathcal F(\nu)$ given by the intersection of $\dbar'\nu$ with $0$. 
\end{construction}

So $\iota\co \hat f_3\longrightarrow \mathcal F(\nu)$ is a subfamily with $\dbar\hat f_3\subset V(\hat f_3)$, and $\ex C(\iota)\co \ex C(\hat f_3)\longrightarrow \ex C(\hat f_2)=\ex C(\mathcal F(\nu))$ is also a subfamily. Moreover,  
\begin{equation}\label{dbarf3}\dbar\hat f_3=(\ex C(\iota),\hat f_3)^*\theta\end{equation}
where $\theta$ is the tautological section from (\ref{tautological defn}).  
\begin{claim}\label{f3uniqueness} For any family of curves $\hat h$, and map $r'\co \ex C(\hat h)\longrightarrow \ex C(\hat f_1)$, there exists at most one map $\psi\co \hat h\longrightarrow \hat f_3$ such that 
\[r\circ \ex C(\psi)=r'\]
where $r\co \ex C(\hat f_3)\longrightarrow \ex C(\hat f_1)$ is the restriction of the map $r$ in Construction \ref{f2def} to the subfamily $\ex C(\hat f_3)\subset \ex C(\hat f_2)$. 
\end{claim}
In fact, equation (\ref{dbarf3}) implies that curves over the fiber of $\ex F(\hat f_3)\longrightarrow \ex F(\hat f_1)$ over different points in $
\ex F(\hat f_3)\subset \ex F(\hat f_2)=V_1$ can never be isomorphic, because $\dbar$ of these curves is determined by the map to $V_1$. As any morphism $\psi\co \hat h\longrightarrow \hat f_3$ is determined by $r\circ\ex C(\psi)=r'$ and $\ex F(\psi)\co \ex F(\hat h)\longrightarrow \ex F(\hat f_3)$, Claim \ref{f3uniqueness} follows.

The uniqueness property of $\nu$ implies that  $\hat f_{3}$ is a $G$--invariant family of curves. Claim \ref{hFnu} gives that our fiberwise holomorphic map ${\fun\Phi}^{+1}_{2}$ has image in $\ex C(\hat f_{3})\subset \ex C(\hat f_{2})$, and also implies that this fiberwise holomorphic map can be regarded as a map of stacks 
\begin{equation}\fun\Psi\co \dbar^{-1}V\longrightarrow \hat f_{3}/G\end{equation}
so that $\fun\Psi(\hat h)$ is the $G$--bundle $h\bd {G}\longrightarrow \hat h$ and  the $G$--equivariant map $\fun\Psi(\hat h)\co h\bd {G}\longrightarrow \hat f_{3}$ such that $\ex C(\fun\Psi(\hat h))={\fun\Phi}^{+1}_{2}(\hat h)$. As with ${\fun\Phi}_{2}$ and ${\fun\Phi}_{1}$, $\fun\Psi$ sends the morphism $\alpha\co \hat g\longrightarrow \hat h$ to the map $\alpha\bd {G}\co \hat g\bd {G}\longrightarrow \hat h\bd {G}$ of $G$--bundles.

Moreover, the natural map $\fun\Phi_{\mathcal X}$ from Definition \ref{family quotient stack} (with $\mathcal X=\dbar^{-1}(V)$) is a left inverse to $\fun\Psi$. This is because we obtain $\hat h$ when we apply $\fun\Phi_{\mathcal X}$ to the $G$--bundle $h\bd G\longrightarrow \hat h$ with the $G$--equivariant map $\fun\Psi(\hat h)\co h\bd {G}\longrightarrow \hat f_{3}$. Note, however, that $\fun\Phi_{\mathcal X}$ is not a right inverse to $\fun\Psi$, because the $G$--equivariant map $\fun\Psi(\hat f_3)\co \hat f_3\bd G\longrightarrow \hat f_3$ is not equivalent to  the quotient map $\hat f_3\longrightarrow \hat f_3/G$.

\

\noindent{\bf Step 5:}

We construct a $G$--invariant subfamily $\hat f\subset \hat f_3$ on which the map $\fun\Psi(\hat f_3)\co \hat f_3\bd G\longrightarrow \hat f_3$ coincides with the quotient map $\hat f_3\longrightarrow \hat f_3/G$. 

  Claim \ref{f3uniqueness} implies that the map   $\fun\Psi(\hat h)\co \hat h\bd {G}\longrightarrow \hat f_{3}$ is also the unique map  with the property that $ {\fun\Phi}_{1}^{+1}(\hat h)$ factorizes as follows: \begin{equation}\label{uniqueness property}\begin{tikzcd}\ex C(\hat h\bd {G})\arrow[bend right]{rrr}{ {\fun\Phi}_{1}^{+1}(\hat h)}\arrow{r}{\ex C(\fun\Psi(\hat h))}&\ex C(\hat f_{3})\arrow[hook]{r}{\ex C(\iota)}& \ex C(\hat f_{2})\arrow{r}{r}  &\ex C(\hat f_{1})\end{tikzcd}\end{equation}
so ${\fun\Phi}_{1}^{+1}(\hat h)=r\circ \ex C(\fun\Psi(\hat h))$, where $r$ appears in Construction \ref{f2def}.

We have not proved --- and it is not true --- that $\hat f_{3}/G$ represents $\dbar^{-1}V\subset \mathcal O$. We must restrict $\hat f_{3}$ to a subfamily to achieve this. Moreover, the maps  ${\fun\Phi}_1^{+1}(\hat f_3)$ and $r$ do not coincide on $\ex C(\hat f_{3})$. Roughly speaking, the family $\hat f\subset \hat f_{3}$ we are looking for is the subfamily where $r$ and ${\fun\Phi}_1^{+1}(\hat f_{3})$ coincide.

  The fiberwise holomorphic map ${\fun\Phi}_{1}$ applied to $\hat f_{3}$ gives a $G$--fold cover $\hat f_{3}\bd {G}$ of $\hat f_{3}$ and a map 
\[{\fun\Phi}_1^{+1}(\hat f_3)\co \ex C(\hat f_{3}\bd {G})\longrightarrow \ex C(\hat f_{1})\ .\]
For $\mathcal O$ small enough, and $\hat f_{i}$ restricted accordingly,  we can assume that $\hat f_{3}\bd {G}$ consists of $\abs G$ copies of $\hat f_{3}$, so ${\fun\Phi}_1^{+1}(\hat f_3)$ above consists of $\abs G$ maps $\ex C(\hat f_{3})\longrightarrow \ex C(\hat f_{1})$. Both $\hat f_{3}$ and $\hat f_{1}$ have a canonical inclusion of the curve $f$. From Lemma \ref{parametrize by curves}, we know that  the map $\fun\Phi^{+1}_0(\hat f_3)\co \ex C(\hat f_3\bd G)\longrightarrow \ex C(\hat f_{0})$ factorizes through ${\fun\Phi}_{1}^{+1}(\hat f_{3})$; therefore,  exactly one of these maps $\ex C(\hat f_{3})\longrightarrow \ex C(\hat f_{1})$ must be the identity on $\ex C(f)$. There is therefore a unique lift
\begin{equation}\label{lmap} l\co \hat f_{3}\longrightarrow \hat f_{3}\bd {G}\end{equation} (or section of the $G$--bundle $\hat f_{3}\bd {G}\longrightarrow \hat f_{3}$) such that the map
\begin{equation}\label{phi1def}\phi_1:={\fun\Phi}_1^{+1}(\hat f_3)\circ l\co \ex C(\hat f_{3})\longrightarrow \ex C(\hat f_{1})\end{equation}
 is the identity on $\ex C(f)$.

 As $\hat f_{3}\bd {G}\longrightarrow \hat f_{3}$ is a $G$ bundle, and the $G$--action on $\hat f_{3}$ lifts to an action on this $G$--bundle,  there is an action of $G\times G$ on $\hat f\bd {G}_{3}$ such that ${\fun\Phi}_1^{+1}(\hat f_3)$ is invariant under the first factor and equivariant under the second factor, and the bundle map $\hat f\bd {G}_{3}\longrightarrow \hat f_{3}$ is equivariant under the first factor and invariant under the second factor. The lift $l\co \hat f_{3}\longrightarrow \hat f_{3}\bd {G}$ is $G$--equivariant when the diagonal action of $G$ is used on $\hat f\bd {G}$, so the map $\phi_{1}$ is $G$--equivariant.

 The map $\phi_1$ has the following important property:
 \begin{equation}\label{factorization2}\phi_{1}\circ \ex C(\fun\Psi(\hat h))={\fun\Phi}^{+1}_{1}(\hat h) \end{equation}

The two factorizations of $ {\fun\Phi}_{1}^{+1}(\hat h)\co \ex C(\hat h\bd {G})\longrightarrow \ex C(\hat f_{1})$ from  (\ref{uniqueness property}) and (\ref{factorization2}) imply that the map $\fun\Psi(\hat h)\co \hat h\bd {G}\longrightarrow \hat f_{3}$ has image in the subset of $\ex C(\hat f_{3})$ where $r$ and $\phi_{1}$ coincide. These two maps $r$ and $\phi_1$ are not the same, but they do coincide after projecting to $\ex C(\hat f_{0})$. Recall from Lemma \ref{parametrize by curves} that $ {\fun\Phi}_{1}^{+1}(\hat h)$ followed by the projection $\ex C(\hat f_{1})\longrightarrow \ex C(\hat f_{0})$ is the core-family map ${\fun\Phi}_0^{+1}(\hat h)$ to $\hat f_{0}$. The map $r$ followed by projection to $\ex C(\hat f_{0})$ and quotiented by $G$ is also the stack version of the core-family map  $({\fun\Phi}^{+1}_{0})_{\hat f_{2}}\co \ex C(\hat f_{2})\longrightarrow \ex C(\hat f_{0})/G$. Definition \ref{core family} part \ref{cc3} implies that this core-family map is the same for $\hat f_{2}$ and its deformation $\mathcal F(\nu)$, so the restriction of $r$ to $\ex C(\hat f_{3})$ followed by a quotient by $G$ is also $({\fun\Phi}_{0}^{+1})_{\hat f_{3}}$. As both $r$ and $\phi_1$ are the identity on the canonical inclusion of $f$ in $\hat f_{3}$,  the step of taking a quotient by $G$ is not necessary,   and  the following diagram commutes:

\[\begin{tikzcd}\ex C(\hat f_{3})\arrow{r}{\phi_1}\arrow{d}{r}&\ex C(\hat f_{1})\arrow{d}
\\ \ex C(\hat f_{1})\arrow{r} &\ex C(\hat f_{0})\end{tikzcd}\]

\begin{claim}\label{mnb} At the curve $f$ in $\hat f_{3}$, the map 
\[(r,\phi_1)\co \ex C(\hat f_{3})\longrightarrow \ex C(\hat f_{1})\times_{\ex C(\hat f_{0})}\ex C(\hat f_{1})\]
is transverse to the diagonal section of $\ex C(\hat f_{1})\times_{\ex C(\hat f_{0})}\ex C(\hat f_{1})\longrightarrow \ex C(\hat f_{0})$.
\end{claim}

To prove Claim \ref{mnb}, recall that there is an inclusion of $\hat f_{1}$ into $\hat f_{2}$ on which the tautological section $\theta$ from (\ref{tautological defn}) is $0$. The inverse  image, $\hat g$, of $f$ under the map $\hat f_{1}\longrightarrow \hat f_{0}$ is holomorphic. Therefore, $\nu$ is $0$ on the image of $\hat g$ in $\hat f_{2}$, so there is an inclusion $\hat g\longrightarrow \hat f_{3}\subset \mathcal F(\nu)$ lifting the inclusion $\hat g\longrightarrow \hat f_{1}$. As every curve in $\hat g$ is isomorphic to $f$, $\phi_1$ restricted to $\ex C(\hat g)$ is constant. On the other hand, $r$ restricted to $\ex C(\hat g)$ is an isomorphism onto the fiber of $\ex C(\hat f_{1})\longrightarrow \ex C(\hat f_{0})$ over $\ex C(f)$, therefore $(r,\phi_1)$ is transverse to the diagonal at $f$ as required.

\

By restricting $\hat f_{3}$ to a possibly smaller neighborhood $\mathcal O$ of $f$, Claim \ref{mnb} implies that we can assume that  $(r,\phi_1)$ is transverse to the diagonal.

\begin{construction}\label{fldef} Let $\hat f$ be the restriction of $\hat f_{3}$  to the inverse image of the diagonal under $(r,\phi_1)\co \ex C(\hat f_{3})\longrightarrow \ex C(\hat f_{1})\times_{\ex C(\hat f_{0})}\ex C(\hat f_{1})$.
\end{construction} As the diagonal  and $(r,\phi_1)$ are $G$--equivariant, $\hat f$ is a $G$--invariant family. 

\

\noindent{\bf Step 6:}

As noted above, (\ref{uniqueness property}) and (\ref{factorization2}) imply that $\fun\Psi\co \dbar^{-1}V\longrightarrow \hat f_{3}/G$ has image in the substack where $r$ and $\phi_1$ coincide, so we have a map $\fun\Psi\co \dbar^{-1}V\longrightarrow\hat f/G$. In the following claims, we  prove that $\fun\Psi$ is an equivalence if stacks, so $\hat f/G$ represents the stack $\dbar^{-1}V$.

\begin{claim}\label{qm} The fiber product  $\hat f\times_{\mathcal O}\hat f$ is represented by the family of curves $\hat f\bd G$ together with the map $\fun\Psi(\hat f)\co \hat f^{\bd G}\longrightarrow \hat f\subset \hat f_3$ and the bundle map $\pi:\hat f\bd G\longrightarrow \hat f$; see Definition \ref{family fiber product}. In particular,  given a family $\hat h$ in $\mathcal O$ with two maps $\alpha,\beta\co\hat h\longrightarrow \hat f$, there exists a unique map 
$\psi_{\alpha,\beta}\co\hat h\longrightarrow \hat f\bd {G}$ such that the following diagram commutes.
\[\begin{tikzcd}\hat h\ar{dr}{\psi_{\alpha,\beta}}\rar{\alpha}\dar{\beta}& \hat f
\\ \hat f & \hat f\bd {G}\lar{\pi}\uar[swap]{\fun\Psi(\hat f)}\end{tikzcd}\]

\end{claim}

To prove Claim \ref{qm}, we will construct an equivariant lift $l'\co \hat f\longrightarrow \hat f\bd {G}$ such that  both $\pi\circ l'$ and $\fun\Psi(\hat f)\circ l'$ are the identity. In particular, $l'$ is the pullback of the lift $l\co \hat f_3\longrightarrow \hat f_3\bd {G}$ from (\ref{lmap}) using the inclusion $\iota\co \hat f\hookrightarrow \hat f_3$.
\[\begin{tikzcd}\hat f\bd {G}\rar{\iota\bd {G}}\dar[shift right,swap]{\pi} &\hat f_3\bd {G}\dar[shift right]
\\ \hat f \rar{\iota}\uar[shift right,swap]{l'} & \hat f_3\uar[shift right,swap]{l}
\end{tikzcd}\]
From  Definition \ref{fhm}, we have
\[{\fun\Phi}^{+1}_1(\hat f_3)\circ \ex C(\iota\bd {G})={\fun\Phi}^{+1}_1(\hat f)\]
so (\ref{phi1def}) implies that
\[{\fun\Phi}^{+1}_1(\hat f)\circ \ex C(l')=\phi_1 \ .\]
From  (\ref{uniqueness property}), we also have that  $r\circ \ex C(\fun\Psi(\hat f))={\fun\Phi}^{+1}_1(\hat f)$, and Construction \ref{fldef} implies that $\phi_1$ and $r$ coincide restricted to $\ex C(\hat f_3)$, so 
\[r\circ \ex C(\fun\Psi(\hat f))\circ\ex C(l')=r\]
and Claim \ref{f3uniqueness} then implies that $\fun\Psi(\hat f)\circ l'$ is the identity. By constriction, $l'$ is also a section of the $G$--bundle $\pi\co \hat f\bd {G}\longrightarrow \hat f$, so $l'$ allows us to identify $\hat f\bd {G}$ with $\hat f\times G$, where $l'(f)=(f,\text{id})$, $\pi( f,g)= f$, and $\fun\Psi(\hat f)(f,g)=g*f$.

Assume without loss of generality that $\hat h$ is  a connected family of curves in $\mathcal O$. Given maps $\alpha$ and $\beta$ from $\hat h$ into $\hat f$, we get two sections $\alpha^*l'$ and $\beta^*l'$ of the $G$--bundle $\hat h\bd {G}\longrightarrow \hat h$, and  hence a unique element $g_{\alpha,\beta}\in G$ such that $g_{\alpha,\beta}*(\beta^*l')=\alpha^*l'$. Therefore,  there is a map $\psi_{\alpha,\beta}\co \hat h\longrightarrow \hat f\bd {G}=\hat f\times G$ such that the required diagram commutes. In particular, 
\[\psi_{\alpha,\beta}=(\beta,g_{\alpha,\beta})\]
As $G$ does not act trivially on any curve in $\hat f$, it follows that $\psi_{\alpha,\beta}$ is the unique map such that the required diagram commutes. This completes the proof of Claim \ref{qm}.

\begin{claim}\label{hfp}For $\hat h$ a family in $\dbar^{-1}V\cap \mathcal O$, the fiber product $\hat h\times_\mathcal O\hat f$ is represented by $\hat h\longleftarrow \hat h\bd {G}\xrightarrow{\fun\Psi(\hat h)}\hat f$; see Definition \ref{family fiber product}.
\end{claim}

So, given $\alpha\co \hat g\longrightarrow\hat f$ and $\beta\co \hat g\longrightarrow \hat h$, we must prove that there exists a unique map $\psi_{\alpha,\beta}\co \hat g\longrightarrow \hat h\bd {G}$ such that the following diagram commutes.
\[\begin{tikzcd}\hat g\ar{dr}{\psi_{\alpha,\beta}}\rar{\alpha}\dar{\beta}& \hat f
\\ \hat h & \hat h\bd {G}\lar{\pi}\uar[swap]{\fun\Psi(\hat h)}\end{tikzcd}\]
Without losing generality, we can assume that $\hat g$ is connected and small enough that it maps into a subfamily of $\hat h$ over which the $G$--bundle $\hat h\bd {G}\longrightarrow \hat h$ is trivial, so there are exactly $\abs{G}$ lifts of $\beta$, permuted by the action of $G$ on this $G$--bundle. Moreover, Claim \ref{qm} implies that there are exactly $\abs{G}$ maps $\hat g\longrightarrow \hat f$, permuted by the action of $G$ on $\hat f$. As $\fun\Psi(\hat h)$ is $G$--equivariant, it follows that there is exactly one lift $\psi_{\alpha,\beta}$ of $\beta$ such that $\fun\Psi(\hat h)\circ \psi_{\alpha,\beta}=\alpha$. Claim \ref{hfp} follows.

\begin{claim}\label{rbs} The stack $\dbar^{-1}V\subset \mathcal O$ is represented by $\hat f/G$ in the sense of Definition \ref{represented by family}. Moreover, the map ${\fun\Phi}_{\mathcal X}$  from Definition \ref{family quotient stack} (with $\dbar^{-1}V=\mathcal X$) is an inverse to the  map $\fun\Psi\co \dbar^{-1}V\longrightarrow \hat f/G$.

\end{claim}

By definition, $\fun\Psi(\hat h)$ is a $G$--bundle $\hat h\bd {G}\longrightarrow \hat h$ with a $G$--equivariant map to $\hat f$, so ${\fun\Phi}_{\mathcal X}(\fun\Psi(\hat h))=\hat h$. Similarly for a morphism, $\alpha$ in $\dbar^{-1}V$, we have  ${\fun\Phi}_{\mathcal X}(\fun\Psi(\alpha))=\alpha$. So, ${\fun\Phi}_{\mathcal X}\circ\fun\Psi$ is the identity. To prove Claim \ref{rbs}, it remains to construct a $2$--isomorphism  $\eta$ from the identity to $\fun\Psi\circ {\fun\Phi}_{\mathcal X}$. 

Consider a family in $\hat f/G$; this is a $G$--bundle $\hat h'\longrightarrow \hat h$ internal to $\mathcal O$ and a $G$--equivariant map $\alpha\co \hat h'\longrightarrow \hat f$; call this family $\alpha$. We have ${\fun\Phi}_{\mathcal X}(\alpha)=\hat h$, so $\fun\Psi({\fun\Phi}_{\mathcal X}(\alpha))$ is the $G$--bundle $\hat h\bd {G}\longrightarrow \hat h$ with the $G$--equivariant map $\fun\Psi(\hat h)\co \hat h\bd {G}\longrightarrow \hat f$. We want an isomorphism $\eta_\alpha$ from  $\alpha$ to $\fun\Psi(\hat h)$, or equivalently,  an isomorphism of $G$--bundles $\eta_{\alpha}\co\hat h'\longrightarrow \hat h\bd {G}$ such that the following diagram commutes.
\[\begin{tikzcd}\hat h'\dar \ar{dr}{\eta_\alpha}\rar{\alpha}& \hat f
\\ \hat h & \lar\hat h\bd {G}\uar[swap]{\fun\Psi(\hat h)}\end{tikzcd}\] 
In fact,  Claim \ref{hfp} implies that there is a unique map $\eta_\alpha$ such that the above diagram commutes. Claim \ref{qm} implies that $\eta_\alpha$ must be $G$--equivariant, so $\eta_\alpha$ is the required isomorphism of $G$--bundles. Moreover, Claim \ref{hfp} implies that, given any map $\psi'$ of $G$--bundles, the following diagram commutes:
\[\begin{tikzcd}\hat g\bd {G}\rar{\psi\bd {G}}&\hat h\bd {G}
\\ \hat g'\uar{\eta_{\alpha\circ \psi'}}\rar{\psi'}\dar &\uar{\eta_\alpha}\dar\hat h'\rar{\alpha}& \hat f
\\ \hat g\rar{\psi} & \hat h\end{tikzcd}\]
so $\eta$ defines a natural transformation. This completes the proof of Claim \ref{rbs}.

\

\noindent{\bf Step 7:}

To complete the proof of our theorem, it remains to verify that $T_f\hat f\co T_f\ex F(\hat f)\longrightarrow D\dbar^{-1}V(f)$ is a bijection.

The following claim allows us to think of $\ex C(\hat f)$ as embedded in both $\ex C(\hat f_{1})$ and $\ex C(\hat f_{0})$. It implies that any $G$--invariant section of $V(\hat f)$ can be extended to a global section of $V$ defined on a neighborhood of $\hat f$ in $\dmsw$, and similarly any $G$--invariant, $\C\infty1$ function of $\ex F(\hat f)$ can be extended to a $G$--invariant, $\C\infty1$ function defined on a neighborhood of $\hat f$ in $\dmsw$.

\begin{claim}\label{Vembedding}The maps
\[\phi_1\co \ex C(\hat f)\longrightarrow \ex C(\hat f_{1})\]
and 
\[\ex C(\hat f)\xrightarrow{\phi_1}\ex C(\hat f_{1})\longrightarrow \ex C(\hat f_{0}) \]
are embeddings in a neighborhood of $f$.
\end{claim}

To prove Claim \ref{Vembedding}, note that $\ex F(\hat f_{3})\longrightarrow \ex F(\hat f_{2})$ is an embedding, with image locally the zero set of some  $\C\infty1$ functions which are transverse to $0$. (This was  proved in Claim \ref{dbar'nu}.) $\ex C(\hat f)\subset \ex C(\hat f_{3})$ is defined as the subset  where $r$ and $\phi_1$ agree. These maps are proved to be transverse in Claim \ref{mnb}. As required by Lemma \ref{parametrize by curves} part \ref{pbc2}, $\ex C(\hat f_{1})\longrightarrow \ex C(\hat f_{0})$ is isomorphic to a vector-bundle. The maps $r$ and $\phi_1$ agree after composition with the map $\ex C(\hat f_{1})\longrightarrow \ex C(\hat f_{0})$, therefore $\ex F(\hat f)\subset \ex F(\hat f_{2})$ is also an embedding  locally defined by the transverse vanishing of some $\C\infty1$ functions.

As $\ex F(\hat f_{2}) $ is isomorphic to a vector-bundle over $\ex F(\hat f_{1})$, and $\ex F(\hat f_{1})$ is isomorphic to a vector-bundle over $\ex F(\hat f_{0})$, to finish Claim \ref{Vembedding} it now suffices to prove that the derivative of the map $\ex F(\hat f)\longrightarrow \ex F(\hat f_{0})$ is injective at $f$. Let $v$ be a vector in $T_{f}\ex F(\hat f)$  sent to $0$ in $T_{f}\ex F(\hat f_{0})$. As this map $\ex F(\hat f)\longrightarrow \ex F(\hat f_{0})$ corresponds to the core-family map, it follows that $v$ must  equal  a section of $f^{*}T_{vert}\hat{\ex B}$ vanishing on the extra marked points $\{s_i\}$ in the definition of the core family $\hat f_{0}$. Lemma \ref{parametrize by curves} part \ref{pbc3} specifies that $D\dbar^{-1}V(f)$  does not contain any nonzero such vector, but $D\dbar(v)$ must be in $V(f)$ because $\dbar\hat f$ is a section of $V(\hat f)$. It follows that the image of $v$ in $T_{f}\dmsw$ is $0$. In particular, $D\dbar(v)=0$, therefore $v$ is tangent to $\ex F(\hat f_{1})\subset \ex F(\hat f_{2})$.  The projection of $T_{f}\ex F(\hat f_{2})$ onto $T_{f}\ex F(\hat f_{1})$ comes from the map $r$, which coincides on $T_{f}\ex F(\hat f)\subset T_{f}\ex F(\hat f_{2})$ with  $\phi_1$. As $\phi_1$ comes from a fiberwise-holomorphic map ${\fun\Phi}_1$ of a neighborhood of $f$ in $\dmsw$, Lemma \ref{tangent2} implies that the map $T_{f}\ex F(\hat f)\longrightarrow T_{f}\ex F(\hat f_{1})$ factors through  $T_{f}
\dmsw$, therefore the image of $v$ in $T_{f}\ex F(\hat f_{1})$ must be zero. As we have already established that $v$ is tangent to $T_{f}\ex F(\hat f_{1})$, it follows that $v$ is the zero vector in $T_{f}\ex F(\hat f)$ and the map $T_{f}\ex F(\hat f)\longrightarrow T_{f}\ex F(\hat f_{0})$ is injective.  

It follows that, on some neighborhood of $f$ in  $\ex F(\hat f)$, the map $\ex F(\hat f)\longrightarrow \ex F(\hat f_{0})$ is an embedding, locally defined by the transverse vanishing of some $\C\infty1$ functions. As $\ex F(\hat f_{1})\longrightarrow\ex F(\hat f_{0})$ is isomorphic to a vector-bundle, the same holds for the map $\ex F(\hat f)\longrightarrow \ex F(\hat f_{1})$. This completes the proof of Claim \ref{Vembedding}

\

To complete the proof of Theorem \ref{V moduli stack}, it remains to prove the following: 

\begin{claim}\label{TV moduli} The map 
\[T_f\hat f\co T_{f}\ex F(\hat f)\longrightarrow T_{f}\dmsw(\hat{\ex B})\]
 is injective, and has image
\[D\dbar^{-1}(V(f))\subset T_{f}\dmsw(\hat{\ex B})\ .\] 
\end{claim}

To prove Claim \ref{TV moduli},  note that Claim \ref{Vembedding} gives that the map $T_{f}\ex F(\hat f)\longrightarrow  T_{f}\ex F(\hat f_{0})$ is injective. As this map comes from the core family map ${\fun\Phi}_0$ to  $\ex C(\hat f_{0})/G$, Lemma \ref{tangent2} implies that it factors through $T_{f}\dmsw$, therefore $T_{f}\ex F(\hat f)$ injects into $T_{f}\dmsw$. 

So far, we have seen that $T_f\hat f\co T_{f}\ex F(\hat f)\longrightarrow T_{f}\dmsw$ is injective. Obviously, $T_{f}\ex F(\hat f)$ has image contained inside $D\dbar^{-1}(V(f))$, so it remains to show that the image of $T_{f}\ex F(\hat f)$ contains $D\dbar^{-1}(V(f))$. Given any vector $v$ in $D\dbar^{-1}(V(f))$, Lemma \ref{tangent1} implies that there exists a family $\hat h$ of curves in $\mathcal O$ parametrized by $\mathbb R$, such that $f$ is the curve over $0$, and the derivative of $\hat h$ at $0$ is  $v$.  There therefore exists a section $\theta_{\hat h}$ of $V(\hat h)$ so that  $\dbar \hat h-\theta_{\hat h}$ is tangent to the zero-section at $0$. Then there exists a map $\psi'_{\theta_{\hat h}}\co \ex C(\hat h)\longrightarrow \ex C(\hat f_{2})$ and a section $\nu_{\hat h}$ of $(\hat f_{2}\circ \psi'_{\theta_{\hat h}})^{*}T_{vert}\hat{\ex B}$ vanishing on the inverse image of the marked point sections $\{s_{i}\}$ so that $\hat h=\mathcal F(\nu_{\hat h})$ and $\dbar'\nu_{\hat h}$ is tangent to the zero-section at $0$. Theorem \ref{regularity theorem} then implies that, close to $0$, there exists a section $\nu'_{\hat h}$ for which $\dbar'\nu'_{\hat h}=0$, and Corollary \ref{trt} implies that $\nu'_{\hat h}$ is tangent to $\nu_{\hat h}$ at $0$. The uniqueness part of Theorem \ref{regularity theorem} implies that $\nu'_{\hat h}=(\psi'_{\theta_{\hat h}})^{*}\nu$, so $\hat h$ is tangent at $0$ to the family $\hat f_{3}$ at $f$.  This proves that the image of $T_f\hat f_3\co T_{f}\ex F(\hat f_{3})\longrightarrow T_{f}\dmsw$ contains $D\dbar^{-1}(V(f))$, so the image of $T_{f}\hat f$ also contains $D\dbar^{-1}(V(f))$, because $\fun\Psi(\hat f_3)$ maps a $G$--bundle over $\hat f_3$ into $\hat f$, and represents the $G$--quotient map when restricted to $\hat f\subset \hat f_3$. It follows that $T_f\hat f\co T_{f}\ex F(
\hat f)\longrightarrow D\dbar^{-1}(V(f))$ is a bijection. 

This completes the proof of Claim \ref{TV moduli} and Theorem \ref{V moduli stack}.

\stop

\

\section{Construction of an embedded Kuranishi structure}\label{K construction}

The goal of this section is to construct compatible Kuranishi charts covering  the moduli stack of stable holomorphic curves $\mathcal M(\hat{\ex B})\subset \Msw(\hat{\ex B})$ in some family $\hat {\ex B}\longrightarrow \ex B_0$ of exploded manifolds. 
Throughout this section, we  assume the map $\Mod\longrightarrow  \ex B_{0}$ is proper when restricted to any connected component of $\Msw$ in the sense\footnote{Note that lemmas \ref{topological sequential convergence} and \ref{open pullback} imply that a substack of $\Msw$ is compact if and only if it is sequentially compact. The existence of core families make $\Msw$ relatively tame, but  more general stacks  are non-Hausdorff in nastier ways than exploded manifolds, so the correct definition of a proper map of stacks is probably stronger than the inverse image of compact substacks being compact.} that we can extract a convergent subsequence from any sequence of holomorphic curves in $\Mod(\hat{\ex B})$ contained in a connected component of $\Msw(\hat{\ex B})$ and over a compact subset of $\ex B_{0}$; moreover, this convergent  subsequence  converges in $\C\infty1$.  This compactness property for the moduli stack of holomorphic curves is proved for many targets $\hat{\ex B}$ in \cite{cem}, and Definition \ref{decorated} implies that the same property holds for $\dmod\subset\dmsw$. We need this assumption  to use Lemma \ref{shrink}, and  to construct our Kuranishi charts to give a locally finite cover of $\mathcal M$ or a moduli stack $\dmod$ of decorated holomorphic curves.

Before beginning the construction, we need an analogue of Lemma \ref{shrink} to apply to $\dmsw$ so that we can shrink open substacks of $\dmsw$ appropriately.

\begin{lemma}\label{dshrink} Given any holomorphic  curve $f$ in $\dmsw$, and an open neighborhood $\mathcal U\subset\dmsw$ of $f$,  there exists an open neighborhood $\mathcal O\subset\mathcal U$ of $f$ and a $\C\infty1$ function $\rho\co \mathcal O\longrightarrow [0,1]$ so that
$\rho(f)=1$, and every holomorphic curve in the closure (within $\dmsw$) of $ \{\rho>0\}$ is contained in $\mathcal O$.  
\end{lemma}

\pf 

Use $\fun\pi\co\dmsw\longrightarrow \Msw$ for the map forgetting decorations. Theorem \ref{V moduli stack} implies that we can choose a neighborhood $\mathcal O_{1}$ of the undecorated $\fun\pi(f)$ in $\Msw$ so that the moduli stack of holomorphic curves in $\mathcal O_{1}$ is represented by some substack of $\hat f/G$, and Claim \ref{Vembedding} implies that $\hat f/G$ is embedded  in a core family $\hat f_{0}/G$ for $\mathcal O_{1}$.  Definition \ref{decorated} implies that the lift $(\fun\pi^{*}\hat f)/G$ of $\hat f/G$ to $\dmsw$ contains the moduli stack of holomorphic curves in the inverse image $\fun\pi^{-1}\mathcal O_{1}$ of $\mathcal O_{1}$ within $\dmsw$. Remark \ref{decorated core family} also  implies  that the lift $(\fun\pi^{*}\hat f_{0})/G$ of $\hat f_{0}/G$ is a core family for $\fun\pi^{-1}\mathcal O_{1}$. Note that this means that any  $G$--invariant $\C\infty1$ function on $\fun\pi^{*}\hat f$ can be extended to a $G$--invariant $\C\infty1$ function on $\fun\pi^{-1}\mathcal O_{1}$. Lemma \ref{shrink} implies that there exists an open neighborhood $\mathcal O_{2}\subset\mathcal O_{1}$ of the undecorated curve $\fun\pi(f)$ such that every holomorphic curve in the closure of $\mathcal O_{2}$ within $\Msw$ is contained in $\mathcal O_{1}$.   Choose a $\C\infty1$ function $\rho\co \fun\pi^{-1}\mathcal O_{1}\longrightarrow [0,1]$ to be equal to $1$ at $f$, and to have compact support when restricted to the intersection of $\fun\pi^{-1}\hat f$ with $\mathcal U\cap\fun\pi^{-1}\mathcal O_{2}$. The restriction of $\rho$ to $\mathcal O: =\mathcal U\cap \fun\pi^{-1}\mathcal O_{2}$ satisfies the required property.     

\stop

\begin{lemma}\label{transverse sheaf} Let  $\mathcal O$ be an open neighborhood of a holomorphic curve  $f\in\dmsw$ with a $\C\infty1$ submersion
\[{\fun\Phi}\co \mathcal O\longrightarrow\ex X\]
to an exploded manifold or  orbifold $\ex X$. Then, there exists an open neighborhood $\mathcal U$ of $f\in \dmsw$, and, on $\mathcal U$, a simply-generated complex subsheaf $V$ of $\Y$ so that
\begin{enumerate}\item $V$  is strongly transverse to $\dbar$ at all holomorphic curves in $\mathcal U$, in the sense of Definition \ref{strongly transverse}, 
\item and, at any holomorphic curve $f'$ in $\mathcal U$, $D\dbar$ restricted to the kernel of $T_{f'}{\fun\Phi}$ is strongly transverse to $V(f)$, so  for any $t\in[0,1]$
\[T_{f'}{\fun\Phi}\lrb{((1-t)D\dbar+tD\dbar^{\mathbb C})^{-1}(V(f'))}=T_{{\fun\Phi}(f')}\ex X\ .\] 
\end{enumerate}

Suppose further that there is a finite collection of open substacks $\mathcal U_{i}$ of $\dmsw$ on which  simply-generated subsheaves $V_{i}$ of $\Y$ are defined, and suppose that for each $\mathcal U_{i}$, there is a chosen substack $C_{i}\subset \mathcal U_{i}$, closed within $\dmsw$. Then $V$ can be modified so that, in addition to the above conditions, for any holomorphic  curve $f'$ in $C_{i}\cap\mathcal U$, the intersection of $V_{i}(f')$ and $V(f')$ is $0$.
\end{lemma}

\pf

If $f$ has domain $\ex T$, this lemma follows from Lemma \ref{T stack} and the observations that precede it on page \pageref{T stack}. Going forward, we assume that the domain of $f$ is not $\ex T$.

Remark \ref{decorated core family} together with Proposition \ref{smooth model family} gives that there exists an open neighborhood $\mathcal U$ of $f$ in $\dmsw$ and a core family $\hat f/G$ for $\mathcal U$  containing $f$. Theorem \ref{f replacement}  implies that there exists a finite-dimensional subspace of $\Y(f)$  strongly transverse to $\dbar$ in the sense of Definition \ref{strongly transverse}. As the codimension of $\ker T_{f}{\fun\Phi}\subset T_{f}\dmsw$ is finite, we can also construct our finite dimensional subspace of $\Y(f)$ to be strongly transverse to $D\dbar$ restricted to $\ker T_{f}{\fun\Phi}$.  We can also assume that   $G$ preserves the smooth part, $\totl f$, of $f$ within $\totl{\ex F(\hat f)}$. Then,   as $\Y(f)$ only depends on the smooth part of $f$, the action of $G$ on $\totl{f}$ gives a $G$--action on $\Y(f)$. As this  $G$--action on $\Y(f)$ is complex, we  can choose a finite-dimensional, complex,  $G $--invariant subspace $V(f)$ of $\Y(f)$,  strongly transverse to $\dbar$ when restricted to $\ker T_{f}{\fun\Phi}$.

There exists a complex  basis $\{v_{1}(f),\dotsc,v_{n}(f)\}$ for $V(f)$ so that  the action of $g\in G $ in this basis is given by a $n\times n$  complex matrix $A_{g}$. Let us extend $v_{i}$ and this action.
As the inclusion $\ex C( f)\longrightarrow \ex C(\hat f )$ is an isomorphism onto a fiber of $\ex C(\hat f )$, Claim \ref{section extension} implies that there exist sections $v_{i}'$ of $\Gamma^{(0,1)}(T^{*}_{vert}\ex C(\hat f )\otimes T_{vert}\hat{\ex B})$, considered as a sheaf of $\C\infty1(\ex F(\hat f ),\mathbb C)$--modules, such that the pullback of $v_{i}'$ to $\Y(f)$ is $v_{i}(f)$.  Let $v'$ indicate the vector with components $v_{i}'$, then define 
\[ v:=\sum_{g\in G } A_{g}^{-1}g*v'\ . \]
Note that 
\[g*v=A_{g}v\]
so  the sheaf of $\C\infty1(\ex F(\hat f ),\mathbb C)$--modules generated by the components of $v$ is $G$--invariant. Restricted to $\ex C(f)$, we have  $g*v'=A_{g}v'$, so the pullback of the $i$th component of $v$ to $\Y(f)$ is  $\abs {G }v_{i}(f)$. It follows that the complex subsheaf, $V\subset \Y$,  generated by the pullback of this sheaf of $\C\infty1(\ex F(\hat f ),\mathbb C)$--modules is equal to $V(f)$ at $f$, and  is simply generated on some neighborhood $\mathcal U$ of $f$.

We have chosen $V$ so that $V(f)$ is strongly  transverse to $D\dbar$ restricted to $\ker T_{f}{\fun\Phi}$. Such strong transversality will also hold at all holomorphic curves in a neighborhood of $f$:
Theorem \ref{V moduli stack}  states that, if $\mathcal U$ is small enough, the moduli stack $\dbar^{-1}V\subset \mathcal U$  is represented by the quotient of some family $\hat f'$ of curves by an automorphism group $G'$. Theorem \ref{regularity theorem} then  implies that, at all holomorphic curves $f'$ in a neighborhood of $f$ in $\hat f'$,  $V$ is strongly transverse to $D\dbar$ when restricted to $\ker T_{f'}{\fun\Phi}$; see also Lemma \ref{k extend}. In particular, as $T_{f'}{\fun\Phi}$ is surjective, this is equivalent to requiring that $V(f')$ is strongly transverse to $D\dbar$, and that  
\[T_{f'}{\fun\Phi}(((1-t)D\dbar+tD\dbar^{\mathbb C})^{-1}V(f'))=T_{{\fun\Phi}(f')}\ex X\ .\]

%
%

It remains to prove that $V$ can be modified so that, for any holomorphic curve $h$ in $C_{i}$,  $V(h)\cap V_{i}(h)=0$. For applications of this lemma, it is important that the domain $\mathcal U$ of definition of this modified $V$ does not depend on $\mathcal U_{i}$ and $V_{i}$. Choose our $\mathcal U$ so that the stack of holomorphic curves in the closure, $\bar{\mathcal U}$, of $\mathcal U$, is compact, and so that  $V$ is still defined and satisfies the required transversality conditions  on a larger neighborhood $\mathcal U'$ containing all such holomorphic curves --- Lemma \ref{dshrink}  implies that such a reduction of the size of $\mathcal U'$ is possible. We shall use  $\hat f /G$ to indicate  the core family for this larger $\mathcal U'$.

 Recall that $V$ is pulled back from a $G$--invariant complex subsheaf of $\Gamma^{(0,1)}(T^{*}_{vert}\ex C(\hat f )\otimes T_{vert}\hat{\ex B})$. This  subsheaf is  a sheaf of $\C\infty1(\ex F(\hat f ),\mathbb C)$--modules, but there is also an action of $\C\infty1(\ex C(\hat f )\times \hat {\ex B},\mathbb C)$ by multiplication on $\Gamma^{(0,1)}(T^{*}_{vert}\ex C(\hat f )\otimes T_{vert}\hat{\ex B})$. Multiplication by any $G$--invariant, $\mathbb C^{*}$--valued, function $m$ on $\ex C(\hat f )\times \hat{\ex B}$ sends $V$ to some other complex, simply-generated  subsheaf $m V$ of $\Y$. 
 
 Consider a family $m_{t}$ of such $G$--invariant $\mathbb C^{*}$--valued functions on $\ex C(\hat f )\times \hat{\ex B}$ parametrized by $\mathbb R$ such that $m_{0} =1$. Then we can consider $(\hat f\times \mathbb R)/G$ to be a core family for $\mathcal U'\times \mathbb R$ in $\dmsw(\hat {\ex B}\times \mathbb R)$, and $m_{t}V$ is a simply-generated complex subsheaf of $\Y$ on $\mathcal U'\times \mathbb R$. If  $m_{t}V$ satisfies the required transversality conditions at any holomorphic curve $h$, it satisfies these transversality conditions at all holomorphic curves in a neighborhood of $h$. As the required transversality conditions hold for all holomorphic curves in $\mathcal U'$, and the set of holomorphic curves within $ \bar{\mathcal  U}$ is compact and contained in $\mathcal U'$, it follows that for some neighborhood $O$ of $0$ in $\mathbb R$, the required transversality conditions hold for all holomorphic curves in $\mathcal U\times O$.

Let  $h$ be a holomorphic curve in $\bar{\mathcal U}\cap \Mod$.   
Given any nonzero $v\in V(h)$, there exists some $G$--invariant function $m_{v}$ on $\ex C(\hat f)\times \hat{\ex B}$ so that, if our holomorphic curve $h\in\mathcal U_{i}$, then
\[m_{v}v\notin V(h)\oplus V_{i}(h)\ .\]  
It follows that for all curves $h'$ and $v'\in V(h')$ within a neighborhood of $(h,v)$ in $V(\hat f')$
\[m_{v}v'\notin V(h')\oplus_{i}V_{i}(h')\ .\]  
 
The compactness of the set of holomorphic curves within $C_{i}\cap\bar{\mathcal U}$  and the fact that $V$ is finitely generated then imply that there exists some  $m$ so that for all holomorphic curves $h\in C_{i}$, and nonzero $v\in V(h)$, 
\[mv\notin V(h)\oplus V_{i}(h)\ .\]
Then for $t>0$ small enough,   \[e^{tm}v\notin \oplus V_{i}(h)\]
so $e^{tm}V(h)\cap\oplus V_{i}(h)=0$ for all $h\in \Mod\cap \mathcal U$. As our transversality conditions also hold for $t$ small enough, it follows that $e^{tm}V$ is a modification of $V$ with the required properties for any $t>0$ small enough.

\stop

 We are now ready to construct the embedded Kuranishi structures from Definition \ref{K def}. We want each chart to be extendible in the sense of Definition \ref{K extend def}, and given a submersion $\fun\Phi$ from the ambient moduli stack $\Msw$ to some evaluation space $\ex X$, we can also want each chart to be $\fun\Phi$--submersive; see Definition \ref{submersive def}. More importantly, we also require all our Kuranishi charts to be compatible in the sense of Definition \ref{K compatible}.
 
 \begin{thm}\label{K existence} Suppose that the map $\Mod(\hat {\ex B})\longrightarrow \ex B_{0}$ is proper when restricted to any connected component of $\Msw(\hat{\ex B})$. Then there exists an embedded Kuranishi structure  on $\Mod\subset\Msw$.  More generally,  there exists an embedded Kuranishi structure on the decorated moduli stack $\dmod\subset\dmsw$. 
 
 Moreover, given any submersion ${\fun\Phi}\co \dmsw\longrightarrow \ex X$ where $\ex X$ is an exploded manifold or orbifold, all Kuranishi charts can be chosen  ${\fun\Phi}$--submersive. 
 
 This embedded Kuranishi structure can be chosen to include any countable, locally-finite, compatible collection of extendible  ${\fun\Phi}$--submersive Kuranishi charts $\{(\mathcal U_{k},V_{k},\hat f_{k}/G_{k})\}$.
 \end{thm}
 
 \pf

 We prove the more general case, and construct an embedded Kuranishi structure on a decorated moduli stack $\dmod\subset \dmsw$.
 We are given a locally-finite, countable collection $\{(\mathcal U_{k},V_{k},\hat f_{k}/G_{k})\}$ of extendible, ${\fun\Phi}$--submersive Kuranishi charts. Index these charts by negative integers $k$, leaving positive integers free for the rest of our Kuranishi charts. As specified in definitions \ref{K extend def} and \ref{K compatible collection}, there are extensions 
 \[(\mathcal U\x _{k},V_{k},\hat f\x _{k}/G_{k})\]
 of $(\mathcal U_{k},V_{k},\hat f_{k}/G_{k})$ so that each pair of these extended Kuranishi charts are compatible, and $\{\mathcal U\x _{k}\}$ is a locally-finite collection of substacks of $\Msw$. In particular, as specified by Definition \ref{K compatible collection}, each holomorphic curve $f$ has a neighborhood  intersecting only finitely many $\mathcal U\x_{i}$. Definition \ref{K extend def} gives that there is a continuous function $\rho_{k}\co \mathcal U_{k}\x\longrightarrow (0,1]$ such that $\mathcal U_{k}=\{\rho_{k}>0.5\}\subset \mathcal U_{k}\x$ and for any $t>0$, any holomorphic curve in  the closure (in $\dmsw$) of $\{\rho_{k}>t\}$ is contained in $\{\rho_{k}\geq t\}\subset \mathcal U_{k}\x$. We shall use the restriction of $(\mathcal U_{k}\x,V_{k},\hat f\x_{k}/G_{k})$ to $\{\rho_{k}>0.4\}$ instead of our original extension. The above local finiteness implies that, for any holomorphic curve $f$,  there exists a neighborhood of $f$ intersecting only finitely many $\mathcal U_{k}\x$, and  not intersecting $\{\rho_{k}>0.1\}$ if $f$ is not in $\mathcal U_{k}\x$.

 Lemmas \ref{transverse sheaf} and  \ref{dshrink} imply that each holomorphic curve $f$ has a neighborhood $\mathcal O$ with a $\C\infty1$ function $\rho\co \mathcal O\longrightarrow [0,1]$ so that
 \begin{itemize}\item all holomorphic curves in the closure of $\{\rho>0\}$ are contained in $\mathcal O$,
 \item $\rho(f)=1$,
 \item $\mathcal O$  satisfies the conditions on $\mathcal U$ within Lemma \ref{transverse sheaf},
 \item $\mathcal O$ intersects only finitely many $\mathcal U_{k}\x$, and intersects $\{\rho_{k}>0.1\}$ only if $f$ is in $\mathcal U_{k}\x$.
 \end{itemize}
Our assumption that $\Mod\longrightarrow \ex B_{0}$ is proper restricted to each connected component of $\Msw$ and Corollary \ref{compact equivalence} implies $\Mod$ admits an exhaustion by compact substacks and   that these $(\mathcal O,\rho)$ can be chosen so that there is a countable collection $\{(\mathcal O_{i},\rho_{i})\}$ of them indexed by natural numbers so that the sets $\{\rho_{i}>0.5\}\subset \mathcal O_{i}$ cover $\dmod$, and each $\mathcal O_{i}$ intersects only finitely many $\mathcal O_{j}$.

Construct $V_i $ on $\mathcal O_i$ by induction on $i$. Given $V_j$ for $j<k$,   the second part of Lemma \ref{transverse sheaf} implies that on $\mathcal O_{i}$, we can choose a simply-generated complex subsheaf $V_{i}$ of $\Y$ satisfying the transversality conditions of Lemma \ref{transverse sheaf}, so that, given any holomorphic curve $f$ in $\mathcal O_{i}$, any $V_{k}(f)$  with $\rho_{k}\geq 0.1$ and the $V_{j}(f)$ with $0<j\leq i$ and $\rho_{j}\geq 0.1$ are linearly independent --- so, 
    $V_{k}(f)\oplus_{j} V_{j}(f)\subset \Y(f)$ has dimension $\dim V_k(f)+\sum_{k,j}\dim V_j(f)$.
    
    For any set $A$ of negative integers, and nonempty $I\subset \mathbb N$,  define the sheaf
    \[V_{A,I}:=\oplus_{k\in A}V_{k}\oplus_{i\in I}V_{i}\] 
on $\cap_{k\in A}\mathcal U\x_{k}\cap_{i\in I}\mathcal O_{i} $. (Recall that negative integers are indexing our already-constructed Kuranishi charts.) Because the $V_{k}$ are subsheaves of each other, the maximum dimension that $V_{A,I}$ can be is
\[\dim V_{A,I}=\max_{k\in A}\dim V_{k} + \sum_{i\in I}\dim V_{i}\ .\]
 The substack on which the dimension of $V_{A,I}$ is maximal is open; see Definition \ref{open substack}. Let $\mathcal O_{A,I}$ denote this open substack. As noted above, $\mathcal O_{A,I}$ is an open neighborhood of the holomorphic curves where $\rho_{j}\geq0.1$ for $j\in A\cup I$.

 As we want to define compatible Kuranishi charts, we must determine where to use $V_{A,I}$ carefully. In particular, we shall use $V_{A,I}$ on an open substack $\mathcal O'_{A\cup I}\subset \mathcal O_{A,I}$ with the following properties:
 \begin{enumerate}
 \item\label{o1} If $f\in \mathcal O'_{S}$ and $\rho_{j}(f)>0.4$, then $j\in S$.
 \item\label{o2} If $f\in \mathcal O'_{S}$ and $\rho_{j}$ is not defined at $f$, or $\rho_{j}\leq0.1$, then $j\notin S$. 
 \item\label{o3} $\mathcal O'_{S}$ intersects $\mathcal O'_{S'}$ nontrivially  only if $S\subset S'$ or $S'\subset S$.
 \end{enumerate} 
We need to define $\mathcal O'_{S}$ satisfying the above for finite subsets $S\subset \mathbb Z$ containing at least one natural number. Let $n_{S}$ be the number of $j'$ with $\rho_{j'}(f)\geq 0.1$ for some $f$  so that $\rho_{j}(f)\geq 0.5$ for some $j\in S$. We constructed $\mathcal O_{i}$ so that $n_{S}$ is finite. Now define $\mathcal O'_{S}$ to be the interior of the following stack:
\[O''_{S}:=\left\{\min\lrb{0.4,\min _{j\in S} \rho_{j}} -\max\lrb{0.1,\max_{j'\notin S}\rho_{j'}}>\frac{0.1}{n_{S}}\right\}\subset \mathcal O_{S}\]
 In the above, we set $\rho_{j}$ to be $0$ where it is not yet defined. As the functions $\rho_{j}$ do not necessarily extend to be continuous functions on $\dmsw$, the above inequality does not necessarily define an open substack, and we must take $\mathcal O'_{S}$ to be its interior. Explicitly, we need to remove the closure of $\{\rho_{j}'>0.1\}$ from $\mathcal O''_{S}$ for all $j'\notin S$. Each of the above required properties of $\mathcal O'_{S}$ are easy to check   for $\mathcal O''_{S}$, and therefore they also hold for $\mathcal O'_{S}\subset \mathcal O''_{S}$. 
 
 \begin{claim}\label{oc} $\{\mathcal O'_{S}\}$ is an open cover of the substack of holomorphic curves within $ \dmsw$. \end{claim}
 To prove Claim \ref{oc}, we must show that each stable holomorphic curve $f$  is in $\mathcal O'_{S}$ for some $S$. We already know that for some $i\in\mathbb N$, $\rho_{i}(f)>0.5$. There are at most $n_{i}$ indices $j$ such that $\rho_{j}(f)\geq0.1$, therefore, there exists some set $S$ containing $i$ so that 
 \[\min\lrb{0.4,\min _{j\in S} \rho_{j}(f)} -\max\lrb{0.1,\max_{j'\notin S}\rho_{j'}(f)}\geq 0.3/n_{i}\ .\]
 As by definition, $n_{S}\geq n_{i}$ when $i\in S$, $f$ is in $\mathcal O''_{S}$. It remains to verify that $f$ is in the interior of $\mathcal O''_{S}$. If $f$ was in the boundary of $\mathcal O''_{S}$, then there would be some $j'$ so that $\rho_{j'}$ was not defined at $f$, but $f$ was in the closure of $\{\rho_{j'}>0.1\}$. One of our conditions on $\rho_{j'}$ is that this is not possible for holomorphic curves $f$, therefore no holomorphic curve $f$ is in the boundary of $\mathcal O''_{S}$. Therefore $f\in \mathcal O'_{S}$ and Claim \ref{oc} is proved. 
 
 \
 
 We can now construct our Kuranishi charts.
 Restricted to $\mathcal O'_{A\cup I}$, the sheaf $V_{A,I}$ is a simply-generated  complex subsheaf of $\Y$  strongly transverse to $D\dbar$ at all holomorphic curves. Theorem \ref{V moduli stack} then implies that each holomorphic curve $f$ in $\mathcal O'_{A\cup I}$ has an open neighborhood $\mathcal U$ on which  $\dbar^{-1}V_{A,I}$
 is locally represented by some $\hat f/G$. For $\mathcal U$ small enough,  $(\mathcal U,V_{A,I},\hat f/G)$ is then a Kuranishi chart containing $f$.
 
One of the transversality conditions from Lemma \ref{transverse sheaf} is that, at holomorphic curves $f$, $D\dbar$ restricted to $\ker T_{f}{\fun\Phi}$ is strongly transverse to $V_{i}$ --- this implies that, by choosing $\mathcal U$ small enough if necessary, we can assume that $\hat f$ is ${\fun\Phi}$--submersive in the sense of Definition \ref{submersive def}.

Claim \ref{Vembedding} implies that,  if $\mathcal U$ is chosen small enough, $\hat f/G$ is embedded in a core family for $\mathcal U$. It follows that there exists a continuous function, $\rho\co \mathcal U\longrightarrow [0,1]$, equal to $1$ at $f$, and with compact support on $\hat f/G$.  Lemma \ref{dshrink} implies furthermore that $\rho$ can be chosen so that any holomorphic curve in the closure (within $\dmsw$) of $\{\rho>0\}$ is contained in $\mathcal U$. Restricting $(\mathcal U,V_{A,I},\hat f/G)$ to $\{\rho>0.5\}$ gives a Kuranishi chart with an extension to $\{\rho>0\}$. Condition \ref{o3} on $\mathcal O'_{A\cup I}$ implies that all such charts are compatible, and condition \ref{o1} on $\mathcal O'_{A\cup I}$ implies that any such chart is compatible with any of our original $(\mathcal U_{k},V_{k},\hat f_{k}/G_{k})$ on the extension where $\rho_{k}>0.4$.

Our properness assumption on $\Mod\longrightarrow\ex B_{0}$ together with Corollary \ref{compact equivalence} implies  that the moduli stack of homomorphic curves within $\dmsw$ has an exhaustion by compact substacks, therefore we can choose a countable, locally-finite collection of extendible Kuranishi charts of the  above type  covering $\dmod$. This collection of Kuranishi charts together with our original collection of Kuranishi charts is our required embedded Kuranishi structure

\stop

\

 \begin{cor}\label{K homotopy}If $\Mod\longrightarrow \hat{\ex B_{0}}$ is proper restricted to every connected component of $\Msw$, then any two embedded Kuranishi structures on $\dmod(\hat {\ex B})$ are homotopic:  there exists an embedded Kuranishi structure on $\dmod(\hat {\ex B}\times \mathbb R)$  pulling back to the given embedded Kuranishi structures via the inclusions  of $\hat {\ex B}$ over the points $0$ and $1$ in $\mathbb R$. 
 
 Moreover, if there is a submersion 
 \[{\fun\Phi}\co \dmsw(\hat {\ex B})\longrightarrow \ex X\]
 and the two original embedded Kuranishi structures are ${\fun\Phi}$--submersive, then the homotopy can be chosen  ${\fun\Phi}'$--submersive where ${\fun\Phi}'$ is  the composition
 \[\begin{tikzcd}\dmsw(\hat {\ex B}\times \mathbb R)\rar \ar[bend left]{rr}{{\fun\Phi}'}&\dmsw(\hat {\ex B})\rar{{\fun\Phi}}&\ex X \end{tikzcd}\]
\end{cor}

\pf 

Pull back the first embedded Kuranishi structure to a collection of Kuranishi charts over $(-\infty,\frac 13)$ and pull back the second embedded Kuranishi structure to a collection of Kuranishi charts over $(\frac 23,\infty)$. Together, these charts give a countable, locally finite, extendible collection of ${\fun\Phi}'$--submersive, compatible Kuranishi charts on $\dmsw(\hat{\ex B}\times\mathbb R)$. 

Theorem \ref{K existence} implies that we can expand this collection of Kuranishi charts to a ${\fun\Phi}'$--submersive embedded Kurnanishi structure.

\stop

\subsection{The case of curves with domain $\ex T$}\label{T section}

\

\

Observe the following:
\begin{itemize}
\item If any curve in a connected family of curves $\hat f$ has domain $\ex T$, then all curves in $\hat f$ have domain $\ex T$.
\item If any curve in a connected component of $\Msw$ has domain $\ex T$, then all curves in that connected component have domain $\ex T$.
\item If a curve $f$ has domain $\ex T$,  then $f$ is holomorphic and $ \Y(f)$ is trivial. 
\item As the smooth part of $\ex T$ is a single point, a curve $f$ in $\hat{\ex B}$ with domain $\ex T$ is always contained in a single coordinate chart, and any curve in a neighborhood of $f$ within  $\Msw$ is also contained in the same coordinate chart.
\item Any coordinate chart on $\hat{\ex B}$ containing a stable curve $f$ with domain $\ex T$ can be written in the form $U\times \ex T$ so that 
\[f(\tilde z)=(u,c\e a\tilde z^{n})\]
where $u\in U$ and  $c\e a\in\mathbb C^{*}\e{\mathbb R}$ are constant, $n$ is a positive integer, and $\tilde z$ is the standard coordinate on $\ex T$; see \cite[Example 3.4]{iec}. All nearby curves $f'$ will be in the same form, with different constants, but the same integer $n$.
\end{itemize}

\begin{lemma}\label{T stack} If $n$ is a positive integer, and $U$ is a connected exploded manifold, the moduli stack $\mathcal X$ of  curves $f$ in $U\times \ex T$ in the form
\[f(\tilde z)=(u,c\e a\tilde z^{n})\]
is represented\footnote{Definition \ref{represented by family}.} by the quotient of the family of curves 
\[\begin{array}{ccc}U\times \ex T&\xrightarrow{(\id,\tilde z^{n})} &U\times \ex T
\\\downarrow 
\\ U\end{array}\]
by its automorphism group, $\mathbb Z_{n}$.
\end{lemma}

\pf
Call this family of curves $\hat f_{0}$. The automorphism group of $\hat f_{0}$ is $\mathbb Z_{n}$,  acting trivially on $U=\ex F(\hat f_{0})$, and acting by multiplying the $\ex T$ fibers of $\ex C(\hat f_{0})$ by $n$th roots of unity. To satisfy Definition \ref{represented by family} we must construct an inverse $\fun\Psi\co \mathcal X\longrightarrow \hat f_0/\mathbb Z_n$ to the map ${\fun\Phi}_\mathcal X\co \hat f_0/\mathbb Z_n\longrightarrow \mathcal X$ from Definition \ref{family quotient stack}. 

Given any family $\hat f$ of curves in $\mathcal X$,  we need to define ${\fun\Psi}(f)$ as some $\mathbb Z_n$--bundle $\hat f\bd{\mathbb Z_{n}}\longrightarrow \hat f$  with a $\mathbb Z_n$--equivariant map $\hat f\bd{\mathbb Z_{n}}\longrightarrow \hat f_0$. Define $\ex C(\hat f\bd {\mathbb Z_n})$ to be the $n$--fold cover of $\ex C(\hat f)$ given by the fiber product
\[\ex C(\hat f\bd {\mathbb Z_n}):=\ex C(\hat f)\fp{\hat f}{\hat f_{0}}(U\times \ex T)\]
and define $\ex F(\hat f\bd {\mathbb Z_n})$ to be the fiber product of $\hat f$ with the inclusion of $U$ into $U\times \ex T$ sending $u$ to $(u,1\e 0)$.
\[\ex F\bd {\mathbb Z_n}(\hat f):=\ex C(\hat f)\fp {\hat f}{(\id,1\e 0)}U\]
The projection $U\times \ex T\longrightarrow U$ induces a map  $\ex C(\hat f\bd {\mathbb Z_n})\longrightarrow \ex F(\hat f\bd {\mathbb Z_n})$  making $\ex C(\hat f\bd {\mathbb Z_n})\longrightarrow \ex F(\hat f\bd {\mathbb Z_n})$ the family of $\ex T$s pulled back from the following diagram:
\[\begin{array}{ccc}\ex C(\hat f\bd {\mathbb Z_n})&\longrightarrow& \ex C(\hat f_{0})
\\\downarrow &&\downarrow
\\\ex F(\hat f\bd {\mathbb Z_n})&\longrightarrow &\ex F(\hat f_{0})\end{array}\]
Pull back $\hat f_{0}$ via the above map $\ex C(\hat f\bd {\mathbb Z_n})\longrightarrow \ex C(\hat f_{0})$ to define a family of curves $\hat f\bd {\mathbb Z_n}$ with maps to $\hat f_{0}$ and $\hat f$.
\[\begin{tikzcd}\hat f\bd {\mathbb Z_n}\dar\rar& \hat f_{0}
\\ \hat f\end{tikzcd}\]
There are two actions of $\mathbb Z_{n}$ on $\hat f\bd {\mathbb Z_n}$, induced by multiplying fibers of $\ex C(\hat f)$ or $\ex C(\hat f_{0})$ by $n$th roots of unity. The map $\hat f\bd {\mathbb Z_n}\longrightarrow \hat f_{0}$ is invariant under the first action and equivariant under the second action, whereas the map $\hat f\bd {\mathbb Z_n}\longrightarrow \hat f$ is equivariant under the first action and invariant under the second action. This second action makes $\hat f\bd {\mathbb Z_n}\longrightarrow \hat f$ a $\mathbb Z_{n}$--fold cover, so the above is equivalent to a family $\fun\Psi(\hat f)$ in the stack $\hat f_0/\mathbb Z_n$; see Definition \ref{family quotient}.

 Now suppose that $\hat h\bd {\mathbb Z_n}\longrightarrow \hat h$  is a  $\mathbb Z_{n}$--bundle over $\hat h$ with a  map $\psi\co\hat h\longrightarrow \hat f$ and a $\mathbb Z_n$--equivariant map $\hat h\bd {\mathbb Z_n}\longrightarrow \hat f_{0}$. The construction of $\ex C(\hat f\bd {\mathbb Z_n})$ as a fiber product gives a unique map 
 \[\psi\bd {\mathbb Z_n}\co \hat h\bd {\mathbb Z_n}\longrightarrow \hat f\bd {\mathbb Z_n}\]
  commuting with the maps to $\hat f$ and $\hat f_{0}$. This unique map $\psi\bd {\mathbb Z_n}$ must be $\mathbb Z_{n}$--equivariant because the action of $\mathbb Z_{n}$ on $\ex C(\hat f_{0})$ is free. It follows that $\psi\bd {\mathbb Z_n}$ defines a unique lift of $\psi$ to a morphism in the stack $\hat f_0/\mathbb Z_n$. Defining $\fun\Psi(\psi)=\psi\bd {\mathbb Z_n}$ gives a functor $\fun\Psi\co \mathcal X\longrightarrow \hat f_0/\mathbb Z_n$.
  
   Clearly, ${\fun\Phi}_{\mathcal X}\circ \fun\Psi$ is the identity. Moreover if  $\hat f'\longrightarrow \hat f$ is any family in $\hat f/\mathbb Z_{n}$, $\text{id}_{\hat f}\bd {\mathbb Z_n}$ is  a canonical isomorphism between $\hat f'\longrightarrow \hat f$ and $\fun\Psi\circ {\fun\Phi}_{\mathcal X}(\hat f'\longrightarrow \hat f)$, defining a $2$--isomorphism between $\fun\Psi\circ {\fun\Phi}_{\mathcal X}$ and the identity. 
 Therefore our moduli stack of curves is represented by $\hat f_{0}/\mathbb Z_{n}$ as required.
 
 \stop

\section{Relative complex structure and orientation of Kuranishi charts}\label{rcs section}

Suppose that $(\mathcal U,V,\hat f/G)$ is a Kuranishi chart on the moduli stack of holomorphic curves in $\hat{\ex B}\longrightarrow \ex B_0$. In the case that $\ex B_0$ is a point,  at each holomorphic curve $f$ in $\hat f$, we will construct a canonical homotopy class of $G$--invariant complex structure on $T_{f}\ex F(\hat f)$. Or, in the more general case, we will construct a complex structure on $T_{f}\ex F(\hat f)\ov{\ex B_{0}}$. On a neighborhood of the holomorphic curves in $\hat f$, this defines a canonical orientation of $\ex F(\hat f)$ relative to $\ex B_{0}$.

Let us define this homotopy class of complex structures in a single Kuranishi chart. If $f$ is holomorphic, $T_{f}\dmsw\ov{\ex B_{0}}$ is complex, so $D\dbar\co T_{f}\dmsw\ov{\ex B_{0}}\longrightarrow \Y(f)$ has a complex-linear part, $D\dbar^{\mathbb C}$. For a Kuranishi chart, $(\mathcal U_{i},V_{i},\hat f_{i}/G_{i})$, containing $f$, define
\[K_{i,t}(f):=\lrb{(1-t)D\dbar+tD\dbar^{\mathbb C}}^{-1}(V_{i}(f))\subset T_{f}\dmsw\ov{\ex B_{0}}\ .\]
This $K_{i,t}(f)$ for $t\in[0,1]$ is a smooth family of sub-vector-spaces of $T_{f}\dmsw$. Because  $\dbar$ is strongly transverse to $V_{i}$, we can prove this using  Theorem \ref{f replacement}. At the start of our homotopy,   Theorem \ref{V moduli stack} implies that $K_{i,0}(f)$ is equal to $T_{f}\ex F(\hat f)\ov{\ex B_{0}}$. At the other end,  $K_{i,1}(f)$ is the inverse image of the complex vector-space $V_{i}(f)$ under a complex map, so it is a complex subspace of $T_{f}\dmsw\ov{\ex B_{0}}$. Therefore, there is a canonical homotopy class of complex structure on $K_{i,0}(f)=T_{f}\ex F(\hat f)\ov{\ex B_{0}}$, namely, those complex structures  that are   homotopic to the complex structure on $K_{i,1}(f)$. To check that the resulting relative orientation of Kuranishi charts is compatible, we must globalize this construction.

To globalize  our construction,   we will  show that $K_{i,t}$ can be regarded as a family of $G_{i}$--invariant, $\C\infty1$ sub-vector-bundles of $T_{\hat f_{i}}\dmsw\ov{\ex B_{0}}$ restricted to $\{\dbar\hat f_{i}=0\}$. Moreover, in Proposition \ref{t-trivialization}, we show we can trivialise the bundle $K_{i,t}$ in the $t$--direction, and that $K_{i,t}$ can be identified for all $t$ in a way which is $G_{i}$--invariant, compatible with all inclusions $K_{i,t}\subset K_{j,t}$, constant on $\mathbb R$--nil vectors, and compatible with any chosen submersion ${\fun\Phi}\co \dmsw\longrightarrow \ex X$.

 Let us first verify that we can keep the canonical complex structure on $\mathbb R$--nil vectors.
 The $\mathbb R$--nil vectors in $T_{f}\dmsw\ov{\ex B_{0}}$ --- the vectors acting trivially as derivations on $\C\infty1$ functions --- form a complex-linear subspace of $T_{f}\dmsw\ov{\ex B_{0}}$,  in the kernel of $D\dbar$, and therefore  always contained in $K_{i,t}(f)$.

\begin{lemma}Suppose that $f$ is a holomorphic curve contained in $\hat f$ and $T_{f}\ex F(\hat f)\longrightarrow T_{f}\dmsw$ is injective. Then  the $\mathbb R$--nil vectors from $T_{f}\ex F(\hat f)\ov{\ex B_{0}}$ form a complex-linear subspace of $T_{f}\dmsw\ov{\ex B_{0}}$  contained in the kernel of $D\dbar\co T_{f}\dmsw\ov{\ex B_{0}}\longrightarrow \Y(f)$. 
\end{lemma}

\pf

 Without losing generality, we can restrict to the case that $\totl{\ex F(\hat f)}$ is a single point. As in this case each curve in $\hat f$ is holomorphic, it follows that $T_{f}\ex F(\hat f)\ov{\ex B_{0}}$ is in the kernel of $D\dbar$. As we are dealing with tangent spaces relative to $\ex B_{0}$ in this lemma, it suffices to consider the case when $\ex B_{0}$ is a point, so we can talk of $T\dmsw(\ex B)$ instead of $T\dmsw(\hat{\ex B})\ov{\ex B_{0}}$.

There exists a unique complex structure $\bold j$ on $\ex C(\hat f)$ extending the given fiberwise almost-complex structure. In particular, the (locally defined) fiberwise-holomorphic exploded functions on $\ex C(\hat f)$ are always equal, in local coordinates, to some monomial times  a holomorphic  $\mathbb C^{*}$--valued function on the smooth part of the coordinate chart. These fiberwise-holomorphic exploded functions define a sheaf of holomorphic exploded functions on $\ex C(\hat f)$ and give the canonical complex structure $\bold j$ on $\ex C(\hat f)$. Also,   there is a unique complex structure on $\ex F(\hat f)$, because the smooth part of $\ex F(\hat f)$ is a single point. 

The map $\ex C(\hat f)\longrightarrow \ex F(\hat f)$ is holomorphic. 
In fact, any fiberwise-holomorphic $\C\infty1$ map from $\ex C(\hat f)$ to a complex exploded manifold must also be holomorphic with respect to this canonical complex structure $\bold j$, because any $\C\infty1$ map must be holomorphic restricted to  $\mathbb R$--nil vectors, and the tangent space of $\ex C(\hat f)$ is spanned by the vertical tangent space and $\mathbb R$--nil vectors. In particular, the map $\hat f$ is also $\bold j$--holomorphic. 

We have established that $\hat f$ is a holomorphic family of curves. We shall now check that the map $T_f\hat f\co T_{f}\ex F(\hat f) \longrightarrow  T_{f}\dmsw$ is complex linear. Let $v$ be the lift of any vector-field   on $\ex F(\hat f)$ to vector-field  on $\ex C(\hat f)$. As $\bold j$ is integrable, calculation in $\mathbb C^{n}$ gives that \[\bold j \circ L_{v}\bold j=L_{\bold j v}\bold j\ .\]
As $\hat f$ is holomorphic, $J d\hat f(v)=d\hat f(\bold j v)$. In particular, the map $v\mapsto ( L_{v}\bold j,df(v))$ is complex, so  the map $T_{f}\ex F(\hat f)\longrightarrow T_{f}\dmsw$ is complex linear.

\stop
%

\

The canonical homotopy class of complex structure on $T_{f}\ex F(\hat f)\ov{\ex B_{0}}$  gives a canonical homotopy class of complex structure on $T\ex F(\hat f)\ov{\ex B_{0}}$ in a neighborhood of $f$. To verify that this agrees with the canonical homotopy class of complex structure at  other holomorphic curves in this neighborhood, we shall extend the definition of $K_{i,t}$ to all curves in a neighborhood of $f$ in $\hat f$. As a first step, we modify $V_{i}$ off $\hat f_{i}$ to make later application of Theorem \ref{f replacement} easier.

\begin{lemma}\label{Vfembedd}After restricting  $\hat f_{i}$ to a neighborhood of $f$ if necessary, there exists a core family $(\hat f/G,\{s_{l}\})$ and
\begin{itemize}
\item an identification of $G$ with $G_{i}$ and a $G$--equivariant embedding
\[\hat f_{i}\longrightarrow \hat f\]
so that, for $t\in[0,1]$, there is no nonzero section of $f^{*}T_{vert}\hat{\ex B}$ in $((1-t)D\dbar+tD\dbar^{\mathbb C})^{-1}V_{i}(f)$  also vanishing on the image of all the sections $\{s_{l}\}$,
\item and a  locally free,  $G$--invariant subsheaf $V'$ of $\Gamma^{(0,1)}(T^{*}_{vert}\ex C(\hat f)\otimes T_{vert}\hat{\ex B})$ with the same dimension as $V_{i}$, and  pulling back to $\hat f_{i}$ (in the sense of Definition \ref{Vpullback}) to give $V_{i}(\hat f_{i})$.
\end{itemize}
\end{lemma}

\pf

After restricting $\hat f_{i}$ to a neighborhood of any holomorphic curve $f$ in $\hat f_{i}$, Claim \ref{Vembedding}  provides an equivariant, fiberwise-holomorphic embedding
\[\label{emb}\begin{tikzcd}\ex C(\hat f_{i})\arrow{r}{\phi}\arrow{d}&\ex C(\hat f_{0})\arrow{d}
\\\ex F(\hat f_{i})\arrow{r}&\ex F(\hat f_{0})\end{tikzcd}\]
where, as specified by  Lemma \ref{parametrize by curves} part (\ref{pbc3}), $(\hat f_{0}/G,\{s_{l}\})$ is a core family containing $f$ with enough sections $\{s_{l}\}$ so that $((1-t)D\dbar+tD\dbar^{\mathbb C})^{-1}V(f)$ contains no nonzero sections of $f^{*}T_{vert}\hat{\ex B}$  vanishing on the image of all $\{s_{l}\}$. The above map $\phi$ removes the $G$--fold ambiguity of the core family map $\ex C(\hat f_{i})\longrightarrow \ex C(\hat f_{0})/G$, so $\hat f_{i}$ is equal to $\phi$ composed with $\hat f_{0}$ followed by exponentiation of some section of $\hat f_{0}^{*}T_{vert}\hat{\ex B}$ vanishing on the image of the sections $\{s_{l}\}$. We can choose this section to be $G$--equivariant, and define $\hat f$ to be $f_{0}$ followed by exponentiation of this section. Now the above map $\phi$ corresponds to a $G$--equivariant embedding
\[\hat f_{i}\longrightarrow \hat f\]
and $(\hat f/G,\{s_{l}\})$ is a core family satisfying the requirements of this lemma. In particular, as $\hat f_{i}\longrightarrow \hat f$ is an embedding, we can choose a locally free subsheaf $V'$ of $\Gamma^{(0,1)}(T^{*}_{vert}\ex C(\hat f)\otimes T_{vert}\hat{\ex B})$  pulling back to  give $V_{i}(\hat f_{i})$, and with the same rank as $V_{i}$, when considered as a sheaf of $\C\infty1(\ex F(\hat f))$--modules. We can also choose $V'$ to be  $G$--invariant.

\stop

On a neighborhood of $\hat f_{i}$, the pullback of $V'$ defines some subsheaf $V$ of $\Y$. Let us verify that $\dbar^{-1}V=\dbar^{-1}V_{i}$.  For any family of curves $\hat h$ in this neighborhood, consider the subset comprising  curves $h$ in  $\hat h$ for which the dimension of $V(h)$ is  equal to the rank of $V'$; this subset is open, and contains any curve isomorphic to $f$. Therefore, on a neighborhood of  $f$,  $V$ is simply generated in the sense of Definition \ref{simply generated}. 
As $V(f)=V_{i}(f)$, $\dbar$ is strongly transverse to $V(f)$, so we can apply Theorem \ref{V moduli stack} to see that, in a neighborhood of $f$,  the moduli stack $\dbar^{-1}V$ is represented by the quotient  of a family of curves by a group of automorphisms. Theorem \ref{V moduli stack} also implies that this family of curves has the same dimension as $\hat f_{i}$, therefore this moduli stack  $\dbar^{-1}V$ equals our original moduli stack $\dbar^{-1}V_{i}$,  represented by $\hat f_{i}/G$ (restricted to a neighborhood of $f$ if necessary).
We can therefore continue to define an extension of $K_{i,t}(f)$ using $V$ instead of $V_{i}$.

To linearize $\dbar$ at a non-holomorphic curve, we would need a connection on $\Y$. Instead, we  linearize $\pi_{V}\dbar$, where $\pi_{V}$ is the projection
\[\pi_{V}\co \Y\longrightarrow \Y/V\ .\]
 \begin{construction} For any curve $f$ so that $\dbar f\in V(f)$, define the linear map 
 \[D\pi_{V}\dbar\co T_{f}\dmsw\longrightarrow \Y(f)/V(f)\]
as follows: In light of Lemma \ref{tangent2}, we need only construct $D\pi_{V}\dbar$ on $T_{f}\ex F(\hat f)$ for an arbitrary $\C\infty1$ family $\hat f$ containing $f$. So long as this construction commutes with maps of families and gives a linear map, Lemma \ref{tangent2} implies that it defines a linear map from $T_{f}\dmsw$.
 Choose a section $v$ of $V(\hat f)$ equal to $\dbar f$ at $f$ in $\hat f$. Then $\dbar f_{t}-v$ is a section of $\Y(\hat f)$ vanishing at $f$, therefore its derivative at $f$ is a linear map
\[  T_{f}\ex F(\hat f)\longrightarrow \Y(f)\ .\]
 Another choice of $v$ would change the above map by a linear map to $V(f)$. This derivative  therefore gives a well defined linear map 
 \[D\pi_{V}\dbar\co T_{f}\ex F(\hat f)\longrightarrow \Y(f)/V(f)\ .\]
  As the above construction is compatible with maps of families, Lemma \ref{tangent2} implies that this construction gives a well defined linear map  $T_{f}\dmsw\longrightarrow \Y(f)/V(f)$.
\end{construction}

\begin{lemma}\label{DpiV} On a neighborhood $U$ of all holomorphic curves in $\hat f_{i}$, there exists a map, of sheaves of $\C\infty1(\ex F(\hat f_{i}))$--modules over $\ex F(\hat f_{i})$,
\[D\pi_{V}\dbar\co T_{\hat f_{i}}\dmsw\longrightarrow \Y(\hat f_{i})/V(\hat f_{i})\] 
 which becomes $D\pi_{V}\dbar\co T_{ f'}\dmsw\longrightarrow \Y(f')/V(f')$ when restricted to each curve $f'$ in $U$.
\end{lemma}

\pf 

If such a map exists it is uniquely determined by its restriction to $T_{f'}\dmsw$ for each $f'$ in $U$. We therefore need only construct such a map in a neighborhood of a holomorphic curve $f$ in $\hat f_{i}$. 

To globalize the definition of $D\pi_{V}\dbar$, note that there exists a neighborhood of $f$ in $\dmsw$ with a  section $\theta$ of the sheaf $V$ such that $\theta(\hat f_{i})=\dbar\hat f_{i}$. In particular, as $\hat f_{i}$ embeds into the core family $\hat f$, the section $\dbar\hat f_{i}$ is the pullback of some section $\theta'$ of $V'$ over $\ex F(\hat f)$. As $\dbar$ is $G$--invariant and the inclusion $\hat f_{i}\longrightarrow \hat f$ is $G$--equivariant, averaging allows us to construct $\theta'$ to be $G$--invariant. The pullback of $\theta'$ to a neighborhood of $f$ defines a   section $\theta$ of $V\subset \Y$.

Given a section $w$ of $T_{\hat f_{i}}\dmsw$, Lemma \ref{tangent3} implies that there is a one-dimensional deformation $\hat f_{t}$ of $\hat f_{i}$ such that the $t$--derivative of $\hat f_{t}$ at $t=0$ is $w$.  The derivative of $\dbar \hat f_{t}-\theta(\hat f_{t})$ at $t=0$ is a section of $\Y(\hat f_{i})$ that, when restricted to any curve $f'$ in $\hat f_{i}$ and followed by projection to $\Y(f')/V(f')$,  is equal to $D\pi_{V}\dbar(w(f'))$. So,  given any $\C\infty1$ section $w$ of $T_{\hat f_{i}}\dmsw$, there exists a $\C\infty1$ section $D\dbar\pi_{V}w$ of $\Y(\hat f_{i})/V(\hat f_{i})$  restricting to each $\Y(f')/V(f')$ to equal  $D\pi_{V}\dbar(w(f'))$. As this characterization uniquely determines $D \pi_{V}\dbar w$, it follows that 
\[D\pi_{V}\dbar\co T_{\hat f_{i}}\dmsw\longrightarrow \Y(\hat f_{i})/V(\hat f_{i})\]
is a well-defined map of sheaves of $\C\infty1$--modules over $\ex F(\hat f_{i})$.

\stop 

\

The quotient $\Y(\hat f_{i})/V(\hat f_{i})$ has a canonical complex structure because, on $\hat f_{i}$, $V$  coincides with the complex sub-bundle $V_{i}$. In order to talk of the complex-linear part $(D\pi_{V}\dbar)^{\mathbb C}$ of $D\pi_{V}\dbar$ restricted to $T_{\hat f_{i}}\dmsw\ov{\ex B_{0}}$, we  must choose a complex structure on $T_{\hat f_{i}}\dmsw\ov{\ex B_{0}}$.

\begin{claim}\label{Cstr}There exists a $G_{i}$--invariant complex structure $J'$ on $T_{\hat f_{i}}\dmsw\ov{\ex B_{0}}$  defined within a neighborhood $U$ of $f$ so that
\begin{itemize}\item\label{Cstr1} the restriction of this complex structure $J'$ to any holomorphic curve $f'\in U$ is the canonical complex structure on $T_{f'}\dmsw\ov{\ex B_{0}}$,
\item\label{Cstr2} if $X(\hat f_{i})$ indicates the sheaf of $\C\infty1$ sections of $\hat f_{i}^{*}T_{vert}\hat{\ex B}$ vanishing on the image of the core-family sections $\{s_{l}\}$ from  Lemma \ref{Vfembedd}, the inclusion 
\[X(\hat f_{i})\longrightarrow T_{\hat f_{i}}\dmsw\ov{\ex B_{0}}\]
is complex.
\end{itemize}\end{claim}

\pf

Recall that the sheaf $T_{\hat f_{i}}\dmsw\ov{\ex B_{0}}$ is defined using the following short exact sequence:
\[\begin{tikzcd} \Gamma(T_{vert}\ex C(\hat f_{i}))\rar[hook]& \Gamma^{0,1}(T_{vert}\ex C(\hat f_{i})\otimes T_{vert}^*\ex C(\hat f_{i}))\times \Gamma( \hat f_{i}^{*} T_{vert}\hat{\ex B})\dar
\\ & T_{\hat f_{i}}\dmsw(\hat{\ex B})\ov{\ex B_{0}}\end{tikzcd}\]
The sheaf $\Gamma^{0,1}(T_{vert}\ex C(\hat f_{i})\otimes T_{vert}^*\ex C(\hat f_{i}))\times \Gamma( \hat f_{i}^{*} T_{vert}\hat{\ex B})$ has a complex structure, but $T_{\hat f_{i}}\dmsw\ov{\ex B_{0}}$ fails to have a canonical complex structure because the lefthand map above need   not be complex. Below, we construct a complex structure on $T_{\hat f_{i}}\dmsw\ov{\ex B_{0}}$ by splitting the above exact sequence with an inclusion of $T_{\hat f_{i}}\dmsw\ov{\ex B_{0}}$ as a complex subsheaf of $\Gamma^{0,1}(T_{vert}\ex C(\hat f_{i})\otimes T_{vert}^*\ex C(\hat f_{i}))\times \Gamma( \hat f_{i}^{*} T_{vert}\hat{\ex B})$.

Recall that we have $\hat f_{i}$ embedded in a core family $(\hat f,\{s_l\})$, and that $X(\hat f_{i})$ indicates the sheaf of $\C\infty1$ sections of $\hat f_{i}^{*}T_{vert}\hat{\ex B}$ vanishing on the image of each of the core-family sections $s_{l}\co \ex F(\hat f)\longrightarrow \ex C(\hat f)$. There is a canonical injective map $X(\hat f_{i})\longrightarrow T_{\hat f_{i}}\dmsw\ov{\ex B_{0}}$ with cokernel isomorphic to the pullback of $T\ex F(\hat f)\ov{\ex B_{0}}$ to $\ex F(\hat f_{i})$ (because $(\hat f/G,\{s_{l}\})$ is a core family.) This $X(\hat f_{i})$ is a complex subsheaf of $\Gamma^{0,1}(T_{vert}\ex C(\hat f_{i})\otimes T_{vert}^*\ex C(\hat f_{i}))\times \Gamma( \hat f_{i}^{*} T_{vert}\hat{\ex B})$. Moreover,  at our holomorphic curve $f$, the inclusion $X(f)\longrightarrow T_{f}\dmsw\ov{\ex B_{0}}$ is complex, and has finite-dimensional cokernel. So,  there exists a $G_{i}$--invariant, locally-free, finite-rank, complex  subsheaf $W$ of $\Gamma^{0,1}(T_{vert}\ex C(\hat f_{i})\otimes T_{vert}^*\ex C(\hat f_{i}))\times \Gamma( \hat f_{i}^{*} T_{vert}\hat{\ex B})$ which, restricted to $f$, is a finite-dimensional vector-space  complementary to the  direct sum of the image of $\Gamma(T\ex C(f))$ with $X( f)$. 

Restricted to any other curve $f'$, this complex subsheaf $W$ is complementary to the direct sum of the image of $\Gamma(T\ex C(f'))$ with $X( f')$ if and only if its image in  $T_{f'}\ex F(\hat f)\ov{\ex B_{0}}$ is surjective; moreover this condition holds for all $f'$ in a neighborhood of $f$ in $\hat f_{i}$. On this neighborhood the map $X(\hat f_{i})\oplus W\longrightarrow T_{\hat f_{i}}\dmsw\ov{\ex B_{0}}$ is an isomorphism, so we can use the complex structure from $X(\hat f_{i})\oplus W$ to give  a $G_{i}$--invariant complex structure on $T_{\hat f_{i}}\dmsw\ov{\ex B_{0}}$ in a neighborhood of $f$. For $f'$ holomorphic, this complex structure on $T_{f'}\dmsw\ov{\ex B_{0}}$ agrees with the canonical one.

\stop

\begin{remark}\label{tangent metric} In the proof of Claim \ref{Cstr} above, we embedded $T_{\hat f_{i}}\dmsw\ov{\ex B_{0}}$  inside $\Gamma^{0,1}(T_{vert}\ex C(\hat f_{i})\otimes T_{vert}^*\ex C(\hat f_{i}))\times \Gamma( \hat f_{i}^{*} T_{vert}\hat{\ex B})$. We can use this embedding to construct a metric on $T_{\hat f_{i}}\dmsw\ov{\ex B_{0}}$.
 Our embedding represents sections of $T_{\hat f_{i}}\dmsw\ov{\ex B_{0}}$ as sections of some vector-bundle over $\ex C(\hat f_{i})$. So,  we can choose a $\C\infty1$ inner product $<\cdot,\cdot>_{0}$ on this vector-bundle, and a $\C\infty1$ fiberwise volume form $\theta$ on $\pi\co \ex C(\hat f_{i})\longrightarrow\ex F(\hat f_{i})$ to define an inner product  $<\cdot,\cdot>$ on $T_{\hat f_{i}}\dmsw\ov{\ex B_{0}}$ as
\[<v,w>:=\pi_{!}<v,w>_{0}\theta\]
so, $<v,w>$ is a $\C\infty1$ function on $\ex F(\hat f_{i})$ with value at a point $p$ the integral of $<v,w>_{0}\theta$ over the fiber of $\ex C(\hat f_{i})$ over $p$. The fact that $<v,w>$ is $\C\infty1$ follows from Theorem 6.1 of  \cite{dre}.
\end{remark}

  With our chosen complex structure on $T_{\hat f_{i}}\dmsw\ov{\ex B_{0}}$, we can define the complex-linear part, $(D\pi_{V}\dbar)^{\mathbb C}$ of $D\pi_{V}\dbar$. Use the following notation:
 \[A_{t}:= (1-t)D\pi_{V}\dbar+t(D\pi_{V}\dbar)^{\mathbb C}\co T_{\hat f_{i}}\dmsw\ov{\ex B_{0}}\longrightarrow \Y(\hat f_{i})/V(\hat f_{i})\]

\begin{lemma}\label{k extend} On a neighborhood of $f$ in $\hat f_{i}$, 
\[\ker A_{t}\subset T_{\hat f_{i}}\dmsw\ov{\ex B_{0}}\]
is a one-dimensional family of finite-dimensional sub-vector-bundles of $T_{\hat f_{i}}\dmsw\ov{\ex B_{0}}$ --- so,  $ker A_{t}$ is a one-dimensional family of locally-free, finite-rank subsheaves of $T_{\hat f_{i}}\dmsw$.
\end{lemma}

\pf 

We  use  notation from the proof of Lemma \ref{DpiV}.  Within the proof of Lemma \ref{DpiV}, we constructed $D\pi_{V}\dbar$  as the linearization of $(\dbar-\theta)$ followed by projection to $\Y/V$.
 The section $\theta$ of $V$ is the pullback of a section $\theta'$ of $V'$ over $\ex F(\hat f)$ for some core family $(\hat f/G,\{s_l\} )$. As $\hat f_{i}$ is a subfamily of $\hat f$,  we can restrict $\theta'$ to $\ex F(\hat f_{i})$ (and again call it $\theta'$).
 
Let $X(\hat f_{i})$ denote the sheaf of $\C\infty1$ sections of $\hat f_{i}^{*}T_{vert}\hat{\ex B}$ vanishing on the pullback of the core-family sections $\{s_{l}\}$. After choosing a $G$--invariant trivialization for $\hat f_{i}$ in the sense of Definition \ref{trivialization def}, we can use $\theta'$ to  get a simple perturbation  $\dbar':=\dbar-\theta'\co  X(\hat f_{i})\longrightarrow \Y(\hat f_{i})$ as in Example \ref{simple construction}.
Now apply Theorem \ref{f replacement} to the linearization of $\dbar'$, and its linear homotopy to a complex operator,   as allowed by Remark \ref{homotopy application}. Theorem \ref{f replacement} implies that there exists a vector sub-bundle $W$ of $\Y(\hat f_{i})$ containing $V$ so that, on some neighborhood of $f$, we have that $((1-t)D\dbar'+tD\dbar'^{\mathbb C})$ restricted to $X(\hat f_{i})$ surjects onto $\Y(\hat f)/W$ for all $t$, and on this neighborhood,  the inverse image of $W$ is a $1$--dimensional family of  vector sub-bundles of $X(\hat f_{i})$. As $T_{\hat f_{i}}\dmsw \ov {\ex B_{0}}/ X(\hat f_{i})$ is a finite-dimensional $\C\infty1$ vector-bundle, the same holds for $((1-t)D\dbar'+tD\dbar'^{\mathbb C})$ with domain extended from $X(\hat f_{i})$ to $T_{\hat f_{i}}\dmsw \ov {\ex B_{0}}$. By our strong transversality assumption, this map restricted to $T_{f}\dmsw\ov{\ex B_{0}}$ is transverse to $V(f)$ for all $t$, therefore it follows that on some neighborhood of $f$, the inverse image of $V$ is a $1$--dimensional family of finite-dimensional  sub-vector-bundles of $T_{\hat f_{i}}\dmsw\ov{\ex B_{0}}$, as required. 

\stop 

For any holomorphic curve $f'$ in $\hat f_{i}$, the complex structure on $T_{f'}\dmsw\ov{\ex B_{0}}$, chosen in Claim \ref{Cstr}, agrees with the canonical complex structure, so the restriction of  $\ker A_{t}$ to $f'$ is  $K_{i,t}(f')$. Therefore, the kernel of $ A_{t}$ is our extension of $K_{i,t}$ to a family of vector-bundles defined on a neighborhood of holomorphic curves within $\hat f_{i}$. Moreover, as $A_{1}$ is complex, $\ker A_{1}$ has a complex structure  agreeing  with the complex structure on $K_{i,1}(f')$. As our choice of complex structure on $T_{\hat f_{i}}\dmsw\ov{\ex B_{0}}$ was $G_{i}$--invariant, and $D\pi_{V}\dbar$ is intrinsically defined, $\ker A_{t}$ is a $G_{i}$--invariant family of sub-vector-bundles of $T_{\hat f_{i}}\dmsw\ov{\ex B_{0}}$. We therefore have a canonical homotopy class of $G_{i}$--invariant complex structures on $\ker A_{0}$, restricting to equal our previously constructed one at any holomorphic curve $f'$.

\begin{lemma}In a neighborhood of the holomorphic curves in $\hat f_{i}$,  $\ker A_{0}$ is  $T\ex F(\hat f_{i})\ov{\ex B_{0}}$.
\end{lemma}

\pf At holomorphic curves $f'$, Theorem \ref{V moduli stack} implies that $K_{i,0}(f')$ and $\ker A_{0}(f')$  both equal  $T_{f'}\ex F(\hat f_{i})\ov{\ex B_{0}}$. 

More generally, the fact that $\dbar\hat  f_{i}$ is a section of $V(\hat f_{i})$ implies that the  map $T\ex F(\hat f_{i})\ov{\ex B_{0}}\longrightarrow T_{\hat f_{i}}\dmsw\ov{\ex B_{0}}$ has image inside $\ker A_{0}=\ker D\pi_{V}\dbar$. As noted above, this $\C\infty1$ map of finite-dimensional vector-bundles $T\ex F(\hat f_{i})\ov G\longrightarrow \ker A_{0}$ is an isomorphism at holomorphic curves $f'$. Therefore, it is an isomorphism in a neighborhood of these holomorphic curves.

\stop

We have now locally constructed an extension of $K_{i,t}$ to a family of $G_{i}$--equivariant  sub vector-bundles $\ker A_{t}\subset T_{\hat f_{i}}\dmsw\ov{\ex B_{0}}$.  More generally, say $K'_{i,t}$ is an extension of $K_{i,t}$ if it is a family of $G_{i}$--equivariant sub vector-bundles of $T_{\hat f_{i}}\dmsw\ov{\ex B_{0}}$ so that, for any holomorphic curve $f$ in the domain of definition, $K'_{i,t}(f)=K_{i,t}(f)$. With these extensions, we can think of $K_{i,t}$ as a family of $G_{i}$--equivariant vector-bundles over $\{\dbar\hat f_{i}=0\}$. In particular,  define a $\C\infty1$ section of $K_{i,t}$ to be a section  extending to a $\C\infty1$ section of $T_{\hat f_{i}}\dmsw\ov{\ex B_{0}}$.  The following lemma implies that any such $\C\infty1$ section of $K_{i,t}$ extends to a $\C\infty1$ section of any extension $K'_{i,t}$ of $K_{i,t}$.

\begin{lemma}\label{kei}Any section of $K_{i,t}$  extending to a $\C\infty1$ section of $T_{\hat f_{i}}\dmsw\ov{\ex B_{0}}$ also extends to a $\C\infty1$ section of any extension  $K'_{i,t}$ of $K_{i,t}$.\end{lemma}

\pf Let us use the metric on $T_{\hat f_{i}}\dmsw\ov{\ex B_{0}}$ constructed in Remark \ref{tangent metric}. Given any $\C\infty1$ section of $T_{\hat f_{i}}\dmsw\ov{\ex B_{0}}$, its orthogonal projection to $K'_{i,t}$ is a $\C\infty1$ section of $K'_{i,t}$. In particular, any $\C\infty1$ extension of a section of $K_{i,t}$ to $T_{\hat f_{i}}\dmsw\ov{\ex B_{0}}$ projects to a $\C\infty1$ extension contained within $K'_{i,t}$.

\stop

Lemma \ref{kei} implies that the $\C\infty1$ vector-bundle structure on $K_{i,t}$ induced by including $K_{i,t}$ inside $K'_{i,t}$ does not depend on the choice of extension $K'_{i,t}$ of $K_{i,t}$.

\begin{defn} A {\bf$t$--trivialization} of $K_{i,t}$ is a choice of $G_{i}$--equivariant identification of $K_{i,t}$ with $K_{i,0}$ for all $t$ which extends to a $\C\infty1$ family of isomorphisms $K'_{i,t}\longrightarrow K'_{i,0}$ and which is the identity on the subspace  of $\mathbb R$--nil vectors within $K_{i,t}$. 

Given a submersion ${\fun\Phi}\co \dmsw\longrightarrow \ex X$, say that a $t$--trivialization is {\bf ${\fun\Phi}$--submersive} if the following diagram commutes
\[\begin{tikzcd}K_{i,t}(f)\dar\rar{T_{f}{\fun\Phi}} &T\ex X
\\ K_{i,0}(f)\ar{ur}[swap]{T_{f}{\fun\Phi}}\end{tikzcd}\]

A  {\bf $t$--trivialization} for an embedded Kuranishi structure, $\{(U_{i},V_{i},\hat f_{i}/G_{i})\}$,  is a choice of $t$--trivialization for $K_{i,t}$ on $(U\x_{i},V_{i},\hat f\x_{i}/G_{i})$ for all $i$ such that, whenever there is an inclusion $K_{i,t}(f)\longrightarrow K_{j,t}(f)$, the following diagram commutes:
\[\begin{tikzcd}[column sep=large]K_{i,t}(f)\dar \rar&K_{j,t}(f)\ar{rrr}{(1-t)D\pi_{V_{i}}\dbar+t(D\pi_{V_{i}}\dbar)^{\mathbb C}}\dar&&& V_{j}/V_{i}
\\ K_{i,0}(f)\rar&K_{j,0}(f)\ar{urrr}{D\pi_{V_{i}}\dbar}\end{tikzcd}\]
\[\begin{tikzcd}\end{tikzcd}\]

\end{defn} 

Lemma \ref{kei} implies that the definition of a $t$--trivialization is  independent of  the choice of extension $K'_{i,t}$ of $K_{i,t}$.

The most useful case of  ${\fun\Phi}$--submersive $t$--trivializations  is when ${\fun\Phi}$ is holomorphic in the sense of Definition \ref{holomorphic submersion}. The following proposition constructs a ${\fun\Phi}$--submersive $t$--trivialization. 

\begin{prop}\label{t-trivialization}Given a holomorphic submersion ${\fun\Phi}\co \dmsw\longrightarrow \ex X$, and a  ${\fun\Phi}$--submersive embedded Kuranishi structure $\{(\mathcal U_{i},V_{i},\hat f_{i}/G_{i})\}$ on $\dmod\subset \dmsw$, there exists a  ${\fun\Phi}$--submersive $t$--trivialization of $K_{i,t}$.

Such a $t$--trivialization can be chosen compatibly with a given $t$--trivialization defined on a neighborhood of a closed substack $\Mod'\subset \dmod$, so that our new $t$--trivialization  coincides with the given one when restricted to a (possibly smaller) neighborhood of $\Mod'$.
\end{prop}

\pf 
The construction of a $t$--trivialization shall proceed by transfinite induction. In particular, we shall choose a well ordering of our Kuranishi charts, then construct our $t$--trivialization in this order. At each step we shall need to shrink the domain of definition a little, so we shall also need to specify more than one extension of each Kuranishi chart. We shall use the notation $\hat f\exte \hat f^{\sharp}$ to indicate that $\hat f^{\sharp}$ is an extension of $\hat f$.

\begin{claim}\label{well order} There is a well-ordering $\prec$ of the Kuranishi charts  such that
$j\prec i$ if $\dim V_{j}<\dim V_{i}$,  and such that, for any $i$,  there are only finitely many $j$ with $j\prec i$ and $\dim V_{j}=\dim V_{i}$.

There are compatible extensions, $(\mathcal U_{i,k},V_{i},\hat f_{i,k}/G_{i})$  and $(\mathcal U_{i,\underline k},V_{i},\hat f_{i,\underline k}/G_{i})$, of $(\mathcal U_{i},V_{i},\hat f_{i}/G_{i})$ for all $ k\succeq i$  so that
 \[\hat f_{i,k}\exte \hat f_{i,\underline k}\exte \hat f_{i,k'}\ \ \ \text{ whenever }k\succ k'\ ,\]
  and the intersection of  $\hat f_{i,k}$ for all $k\succeq i$ contains an extension of $\hat f_{i}$.  

\end{claim}
 
As there are only a countable number of Kuranishi charts, they can be well-ordered as above.  
Embed our well-ordered set of indices into $(0,3/8)$ as follows:
\[x_{k}:=3/8-2^{-\dim V_k-2}(1+1/(n_k+2))\]
where $n_k$ is the number of indices $j\preceq k$ with $\dim V_j=\dim V_k$. 
The only important property of $x_{k}$ is that   
\[x_{k}>\sup_{j\prec k}x_{j}\ .\] 
By definition, $(\mathcal U_i,V_i,\hat f_i/G_i)$ is extendible, so there exists an extension $(\mathcal U_{i}\x, V_{i},\hat f_{i}\x/G_{i})$ and a $\C\infty1$ map $\rho\co \mathcal U_{i}\x\longrightarrow(0,1]$ satisfying the requirements of Definition \ref{K extend def} so that $\mathcal U_i$ is the substack where $\rho>1/2$. Then, define
\begin{equation}\label{uik}\mathcal U_{i,k}:=\{\rho>x_{k}\}\subset \mathcal U_{i}\x\end{equation}
and 
\begin{equation}\label{uikbar}\mathcal U_{i,\underline k}:=\{\rho> x_{k}/2+\sup_{j\prec k}x_{j}/2\}\subset \mathcal U_{i}\x\ .\end{equation}
These open substacks $\mathcal U_{i,k}$ satisfy all the requirements of Claim \ref{well order} above.

\

For $j\prec i$, we only require compatibility between $K_{j,t}$ and $K_{i,t}$  on the intersection of $\mathcal U_{j,i}$ with $\mathcal U_{i,i}$. Similarly, when matching a given $t$--trivialization on a neighborhood of $\Mod'$, our new $t$--trivialization need only agree with the given $t$--trivialization on  some fixed open neighborhood $O$ of $\Mod'\subset \dmod$ whose closure is  contained in the domain of definition of our given $t$--trivialization.

   Let $f$ be a holomorphic curve in $\hat f_{i,i}$. Suppose that, for $j\preceq k\prec i$, we have a compatible choice of $t$--trivialization for all  $K_{j,t}$ on $\hat f_{j,k}$. We shall construct a $t$--trivialization of $K_{i,t}$ in a neighborhood of $f$  using the following three methods.
\begin{enumerate}
\item \label{case1}If $f$ is contained in $\Mod'$, then there is already a given $t$--trivialization of $K_{i,t}$ compatible with our trivializations of $K_{j,t}$ for $j\prec i$ (by inductive assumption).
\item\label{case2} If $f$ is contained in $\mathcal U_{j,\underline i}$ for some $j\prec i$ but $f$ is not in $\Mod'$, proceed as follows: Without losing generality, assume that $\dim V_{j}\geq \dim V_{j'}$ for all $j'\prec i$ such that  $f$ is in $\mathcal U_{j',\underline i}$. We shall construct a $t$--trivialization of $K_{i,t}$ on some open subset of $\hat f_{i,i}\cap \mathcal U_{j,\underline i}$ where we already have a $t$--trivialization of $K_{j,t}\subset K_{i,t}$.

Recall the complex structure on $T_{\hat f_{i,i}}\dmsw\ov{\ex B_{0}}$ from Claim \ref{Cstr}. On a neighborhood of $f$, construct  an extension $K'_{i,t}$ of $K_{i,t}$ by using this complex structure  to construct a homotopy of $D\pi_{V_{i}}\dbar$ to its complex-linear part, and letting $K'_{i,t}$ be the kernel. Because this complex structure is $G_{i}$--invariant, and there is a unique map $\hat f_{j,\underline i}\longrightarrow \hat f_{i,i}/G_{i}$ (defined in a neighborhood of $f$), we can pull this complex structure back to a complex structure on $T_{\hat f_{j,\underline i}}\dmsw\ov{\ex B_{0}}$. Using this pulled back complex structure, construct an extension $K'_{j,t}$ of $K_{j,t}$ on a neighborhood of $f$ in $\hat f_{j,\underline i}$ with the property that $K'_{j,t}\subset K'_{i,t}$. In particular, construct $K'_{j,t}$ as the kernel of the following operator:

\[A_{t}:=(1-t)D\pi_{V_{j}}\dbar+t (D\pi_{V_{j}}\dbar)^{\mathbb C}\co T_{\hat f_{j,\underline i}}\dmsw\ov{\ex B_{0}}\longrightarrow \Y(\hat f_{j,\underline i})/V_{j}\]
Note that $K'_{i,t}$ is defined using the analogous operator with $V_{i}$ in place of $V_{j}$. As $V_{j}\subset V_{i}$, and the same complex structure is used to define both $K'_{i,t}$ and $K'_{j,t}$,  we have that 
$K'_{i,t}(f')=A_{t}^{-1}(V_{i}/V_{j})$, and $K_{j,t}'\subset K'_{i,t}$.

As our embedded Kuranishi structure is ${\fun\Phi}$--submersive, $A_{t}$ restricted to the kernel of $T_{f}{\fun\Phi}$ is surjective, $T_{f}{\fun\Phi}\co K'_{j,t}(f)\longrightarrow T_{{\fun\Phi}(f)}\ex X\ov{\ex X_{0}}$ is surjective,  and the same holds for all $f'$ in some neighborhood of $f$. It follows that 
\[K'_{i,t}(f')/K'_{j,t}(f')\xrightarrow{A_{t}} V_{i}/V_{j}\]
and
\[\lrb{\ker T_{f'}{\fun\Phi}\cap K'_{i,t}(f')}/\lrb{\ker T_{f'}{\fun\Phi}\cap K'_{j,t}(f')}\xrightarrow{A_{t}} V_{i}/V_{j}\]
are both isomorphisms for $f'$ in a neighborhood of $f$.

Use the notation  $(K'_{j,t})^{\perp}$ to denote the orthogonal complement of  $\ker T{\fun\Phi}\cap K'_{j,t}$ within $\ker T{\fun\Phi}\cap K'_{i,t}$. (For this, use the equivariant metric from Remark \ref{tangent metric}.)   We can split $K'_{i,t}$    into $K'_{j,t}\oplus (K'_{j,t})^{\perp}$. 
 On a neighborhood of $f$ in $\hat f_{j,\underline i}$, the map $A_{t}$ defines an isomorphism of vector-bundles as follows:
\[A_{t}\co (K'_{j,t})^{\perp}\longrightarrow V_{i}/V_{j}\] 
Our $t$--trivialization of $K_{j,t}$ extends by definition to a $t$--trivialization of $K'_{j,t}$, which we can take to be $G_{j}$--equivariant. There is a canonical $G_{j}$--equivariant $t$--trivialization of $(K'_{j,t})^{\perp}$ so that the diagram 
\[\begin{tikzcd} (K'_{j,t})^{\perp}\rar{A_{t}}& V_{i}/V_{j}
\\ \uar (K'_{j,0})^{\perp} \ar{ur}{A_{0}} \end{tikzcd}\]
commutes. The corresponding $G_{j}$--equivariant $t$--trivialization of $K'_{j,t}\oplus (K'_{j,t})^{\perp}$ provides a locally defined,  $G_{i}$--equivariant $t$--trivialization of the restriction of $K'_{i,t}$ to $(\dbar\hat f_{i,i})^{-1} V_{j}$ with the property that the following diagram commutes
\[\begin{tikzcd}K'_{j,t}\rar & K'_{i,t}\rar{A_{t}}&V_{i}/V_{j}
\\ K'_{j,0}\uar \rar &K'_{i,0}\uar \ar{ur}{A_{0}} \end{tikzcd}\]
As $\dbar\hat f_{i,i}$ is transverse to $V_{j}$, we can extend the above $t$--trivialization  to a $t$--trivialization of $K'_{i,t}$ in a neighborhood of $f$ within $\hat f_{i,i}$. The diagram above commutes for all holomorphic curves $f'$. As all $\mathbb R$--nil vectors in $T_{f'}\dmsw\ov{\ex B_{0}}$ are contained in $K_{j,t}(f')$ and the $t$--trivialization of $K_{j,t}$ is constant on $\mathbb R$--nil vectors,  the resulting $t$--trivialization of $K_{i,t}$ is constant on all $\mathbb R$--nil vectors.

 As $(K'_{j,t})^{\perp}(f')$ is contained in the kernel of $T_{f'}{\fun\Phi}$, the diagram
\[\begin{tikzcd}K_{i,t}(f')\ar{rr}{T_{f'}{\fun\Phi}}& &T_{{\fun\Phi}(f')}\ex X
\\ \uar K_{i,0}(f')\ar{urr}[swap]{T_{f'}{\fun\Phi}}
\end{tikzcd}\]
 commutes for holomorphic curves $f'$. 
 
 So far we have constructed a ${\fun\Phi}$--submersive $t$--trivialization of $K_{i,t}$ on a neighborhood of $f$ within $\hat f_{i,i}$, compatible with the $t$--trivialization of $K_{j,t}$.  Given any $j'\prec i$ so that $\dim V_{j'}\leq \dim V_{j}$,  the diagram
\[\begin{tikzcd}[column sep=large]K_{j',t}(f')\dar \rar&K_{j,t}(f')\ar{rrr}{(1-t)D\pi_{V_{j'}}\dbar+t(D\pi_{V_{j'}}\dbar)^{\mathbb C}}\dar&&& V_{j}/V_{j'}
\\ K_{j',0}(f')\rar&K_{j,0}(f')\ar{urrr}[swap]{D\pi_{V_{i}}\dbar}\end{tikzcd}\]
\[\begin{tikzcd}\end{tikzcd}\]
commutes whenever $K_{j',t}(f')$ and $K_{j,t}(f')$ are both defined. It follows that our locally constructed $t$--trivialization of $K_{i,t}$ is automatically compatible with the $t$--trivialization of $K_{j',t}(f')$. On the other hand, we have no reason to expect that our constructed $t$--trivialization is compatible with a $t$--trivialization of $K_{j',t}$ if $\dim V_{j'}>\dim V_{j}$, and we also have no reason to expect compatibility with the already defined $t$--trivialization  on a neighborhood of $\Mod'$.

\item If $f$ is not contained in $\Mod'$ or $\mathcal U_{j,\underline i}$ for any $j\prec i$, then proceed as follows:
On a neighborhood of $f$, choose an extension $K'_{i,t}$ of $K_{i,t}$. Choose a $G_{i}$--equivariant metric on $T_{\hat f_{i,i}}\dmsw\ov{\ex B_{0}}$ as in Remark \ref{tangent metric}.  For all $f'$ in a neighborhood of $f$ in $\hat f_{i,i}$, $T_{f'}{\fun\Phi} \co K'_{i,t}(f')\longrightarrow T_{{\fun\Phi}(f')}\ex X\ov{\ex X_{0}}$ is surjective. We can therefore locally  choose a 
$G_{i}$--equivariant splitting of $K'_{i,t}$
\[K'_{i,t}=\lrb{\ker T{\fun\Phi}\cap K'_{i,t}}\oplus W\]
and a $G_{i}$--equivariant splitting of $T_{\hat f_{i,i}}\dmsw\ov{\ex B_{0}}$ as follows: 
\[T_{\hat f_{i,i}}\dmsw\ov{\ex B_{0}}=\lrb{\ker T{\fun\Phi}\cap T_{\hat f_{i,i}}\dmsw\ov{\ex B_{0}}}\oplus W\] 
 Denote by \[\pi\co T_{\hat f_{i,i}}\dmsw\ov{\ex B_{0}}\longrightarrow K'_{i,t}\] the projection which in the above splittings is the orthogonal projection of  $\ker T{\fun\Phi}\cap T_{\hat f_{i,i}}\dmsw\ov{\ex B_{0}}$ to $\ker T{\fun\Phi}\cap K'_{i,t}$ and the identity on $W$. As our splittings and our  metric  are $G_{i}$--equivariant, $\pi$ is a $G_{i}$--equivariant projection.
 
To construct a $t$--trivialization of $K'_{i,t}$, define a connection in the $t$--direction of $K_{i,t}$ as follows.  A section $\sigma_{t}$ of $K'_{i,t}$ for all $t$ can be viewed as a family of sections of $T_{\hat f_{i,i}}\dmsw\ov{\ex B_{0}}$. The derivative $\frac {d\sigma_{t}}{dt}$ is again a section of $T_{\hat f_{i,i}}\dmsw\ov{\ex B_{0}}$. Define 
 \[\nabla_{t}\sigma_{t}:=\pi\lrb{\frac {d\sigma_{t}}{dt}}\ .\] 

Note that $\nabla_{t}f\sigma_{t}=\frac{\partial f}{\partial t}\sigma+f\nabla_{t}\sigma$, so $\nabla_{t}$  defines a $G_{i}$--invariant connection in the $t$--direction on $K'_{i,t}$. Therefore parallel transport in the $t$ direction gives $G_{i}$--equivariant $t$--trivialization maps $K'_{i,0}\longrightarrow K'_{i,t}$. For holomorphic curves $f'$, $K_{i,t}(f')$ contains all $\mathbb R$--nil vectors, so our $t$--trivialization is the identity on these $\mathbb R$--nil vectors, as required. Notice too that if $\sigma_{t}$ is a section for which $T{\fun\Phi}(\sigma_{t})$ is independent of $t$, then $T{\fun\Phi}(\nabla_{t}\sigma)$ is the zero section. It follows that our $t$--trivialization commutes with $T{\fun\Phi}$, so the diagram
 \[\begin{tikzcd}K_{i,t}(f')\dar\rar{T_{f}{\fun\Phi}} &T\ex X
\\ K_{i,0}(f')\ar{ur}[swap]{T_{f}{\fun\Phi}}\end{tikzcd}\]
commutes.

\end{enumerate}

The above three methods give  a locally defined, ${\fun\Phi}$--submersive $t$--trivialization of $K_{i,t}$ around every holomorphic curve in $\hat f_{i,i}$. As these $t$--trivializations have no reason to match up, we shall average them using a $G_{i}$--invariant partition of unity. In particular, given any finite collection of $t$--trivializations of $K_{i,t}$ defined in a neighborhood of $f$, we can extend them all to locally defined $t$--trivializations of some extension $K'_{i,t}$ of $K_{i,t}$. These $t$--trivializations correspond to connections on $K'_{i,t}$ in the $t$ direction, which can be glued using a partition of unity to create another $t$--trivialization of $K'_{i,t}$. Note the following.
\begin{itemize}
\item The resulting glued $t$--trivialization of $K_{i,t}$ does not depend on the choice of extension $K'_{i,t}$ used.
\item As gluing $G_{i}$--invariant connections using a $G_{i}$--invariant partition of unity gives a $G_{i}$--invariant connection, the resulting $t$--trivialization of $K_{i,t}$ is $G_{i}$--invariant.
\item A $t$--trivialization is ${\fun\Phi}$--submersive if and only if, for each holomorphic curve $f$, the corresponding connection $\nabla_{t}$ satisfies the following property: if $\sigma_{t}$ is a section of $K_{i,t}(f)$ such that  $T_{f}{\fun\Phi}(\sigma_{t})$ is constant, then  $T{\fun\Phi}(\nabla_{t}\sigma_{t})=0$. This property is preserved when we glue connections satisfying it, so the resulting $t$--trivialization of $K_{i,t}$ is also ${\fun\Phi}$--submersive.
\item Gluing connections equal on a subspace produces a connection unchanged on the given subspace. It follows that our glued $t$--trivialization is constant on $\mathbb R$--nil vectors.
\item If all the original $t$--trivializations were compatible with a given $t$--trivialization of $K_{j,t}\subset K_{i,t}$, then the glued $t$--trivialization is also compatible with the given $t$--trivialization. This follows because compatibility with the inclusion $K_{j,t}\subset K_{i,t}$ specifies what our connections must be when restricted to the subspace $K_{j,t}(f)\subset K_{i,t}(f)$, and gluing preserves  the property that the isomorphism  $K_{i,t}(f)/K_{j,t}(f)\longrightarrow V_{i}/V_{j}$ is constant in our $t$--trivialization.  

\end{itemize}

Recall the definition of $\mathcal U_{j,i}$ and $\mathcal U_{j,\underline i}$ from (\ref{uik}) and (\ref{uikbar}). In particular $U_{j,i}\subset U_{j,\bar i}$ are both substacks  in the form $\{\rho>c\}$, where $c>0$ and $\rho$ obeys Definition \ref{K extend def}. It follows that each holomorphic curve in the closure of $\mathcal U_{j, i}$ is in $\mathcal U_{j,\underline i}$. As embedded Kuranishi structures are by definition locally finite, each holomorphic curve $f$ in $\hat f_{i,i}$ has an open neighborhood  in $\hat f_{i,i}$ intersecting $\mathcal U_{j,i}$ for $j\prec i$ only if $f\in \mathcal U_{j,\underline i}$.  Choose this open neighborhood so that it either contains  the domain of definition of our given $t$--trivialization, or  does not intersect $O$, our chosen neighborhood of $\mathcal M'$.  We can also choose this neighborhood of $f$ small enough so that the relevant method above for constructing a $t$--trivialization applies, and choose such a $t$--trivialization. Choose a $G_{i}$--equivariant partition of unity on $\ex F(\hat f_{i,i})$ subordinate to the corresponding open cover of the holomorphic curves in $\hat f_{i,i}$, and average our $t$--trivializations using this partition of unity. 

As noted in the bullet points above, the corresponding $t$--trivialization is $G_{i}$--equivariant, ${\fun\Phi}$--submersive, agrees with the previously chosen trivialization on $O$,  and for all $j\prec i$ is compatible with the $t$--trivialization of $K_{j,t}$ on $\mathcal U_{j,i}$. By   (transfinite) induction, we can choose such compatible $t$--trivializations on $\hat f_{i,i}$ for all $i$. These define the required $t$--trivialization on our Kuranishi structure. 
\stop

\bibliographystyle{plain}

\bibliography{../ref}

\begin{thebibliography}{10}

\bibitem{AMS2021}
Mohammed Abouzaid, Mark McLean, and Ivan Smith.
\newblock Complex cobordism, hamiltonian loops and global kuranishi charts,
  2021.

\bibitem{acgw}
Dan Abramovich and Qile Chen.
\newblock Stable logarithmic maps to {D}eligne-{F}altings pairs {II}.
\newblock {\em The Asian Journal of Mathematics}, 18(3):465--488, 2014.

\bibitem{Bai2022}
Shaoyun Bai and Guangbo Xu.
\newblock Arnold conjecture over integers, 2022.

\bibitem{Bai2025}
Shaoyun Bai and Guangbo Xu.
\newblock A new transversality condition on orbifolds and integer-valued
  gromov-witten type invariants, 2025.

\bibitem{diffstacks}
Kai Behrend and Ping Xu.
\newblock Differentiable stacks and gerbes.
\newblock {\em J. Symplectic Geom.}, 2011.

\bibitem{CLW}
Bohui Chen, An-Min Li, and Bai-Ling Wang.
\newblock Virtual neighborhood technique for pseudo-holomorphic spheres.
\newblock arXiv:1306.3276.

\bibitem{DeligneMumford}
P.~Deligne and D.~Mumford.
\newblock The irreducibility of the space of curves of given genus.
\newblock {\em Inst. Hautes \'Etudes Sci. Publ. Math.}, (36):75--109, 1969.

\bibitem{stacks}
Barbara Fantechi.
\newblock Stacks for everybody.
\newblock In {\em European Congress of Mathematics, Vol. I (Barcelona, 2000)},
  volume 201 of {\em Progr. Math.}, pages 349--359. Birkh\"auser, Basel, 2001.

\bibitem{KFOOO}
Kenji Fukaya, Yong-Geun Oh, Hiroshi Ohta, and Kaoru Ono.
\newblock Technical details on {K}uranishi structure and virtual fundamental
  chain.
\newblock \href{http://arxiv.org/abs/1209.4410}{arXiv:12094410}.

\bibitem{FO}
Kenji Fukaya and Kaoru Ono.
\newblock Arnold conjecture and {G}romov-{W}itten invariant.
\newblock {\em Topology}, 38(5):933--1048, 1999.

\bibitem{FOinteger}
Kenji Fukaya and Kaoru Ono.
\newblock Floer homology and {G}romov-{W}itten invariant over integer of
  general symplectic manifolds---summary.
\newblock In {\em Taniguchi {C}onference on {M}athematics {N}ara '98},
  volume~31 of {\em Adv. Stud. Pure Math.}, pages 75--91. Math. Soc. Japan,
  Tokyo, 2001.

\bibitem{GSlogGW}
Mark Gross and Bernd Siebert.
\newblock Logarithmic {G}romov-{W}itten invariants.
\newblock {\em J. Amer. Math. Soc.}, 26(2):451--510, 2013.

\bibitem{polyfold0}
H.~Hofer.
\newblock A general {F}redholm theory and applications.
\newblock In {\em Current developments in mathematics, 2004}, pages 1--71. Int.
  Press, Somerville, MA, 2006.

\bibitem{polyfold1}
H.~Hofer, K.~Wysocki, and E.~Zehnder.
\newblock A general {F}redholm theory. {I}. {A} splicing-based differential
  geometry.
\newblock {\em J. Eur. Math. Soc. (JEMS)}, 9(4):841--876, 2007.

\bibitem{polyfoldint}
H.~Hofer, K.~Wysocki, and E.~Zehnder.
\newblock Integration theory on the zero sets of polyfold {F}redholm sections.
\newblock {\em Math. Ann.}, 346(1):139--198, 2010.

\bibitem{polyfoldgw}
Helmut Hofer, Kris Wysocki, and Eduard Zehnder.
\newblock A general {F}redholm theory. {III}. {F}redholm functors and
  polyfolds.
\newblock {\em Geom. Topol.}, 13(4):2279--2387, 2009.

\bibitem{polyfoldsc}
Helmut Hofer, Kris Wysocki, and Eduard Zehnder.
\newblock sc-smoothness, retractions and new models for smooth spaces.
\newblock {\em Discrete Contin. Dyn. Syst.}, 28(2):665--788, 2010.

\bibitem{polyfold2}
Helmut Hofer, Krzysztof Wysocki, and Eduard Zehnder.
\newblock A general {F}redholm theory. {II}. {I}mplicit function theorems.
\newblock {\em Geom. Funct. Anal.}, 19(1):206--293, 2009.

\bibitem{joycebook}
Dominic Joyce.
\newblock D-manifolds and d-orbifolds: a theory of derived diferential
  geometry.
\newblock Preliminary version of book available here:
  \href{http://people.maths.ox.ac.uk/~joyce/dmanifolds.html}{http://people.maths.ox.ac.uk/~joyce/dmanifolds.html}.

\bibitem{kuranishihomology}
Dominic Joyce.
\newblock Kuranishi bordism and kuranishi homology.
\newblock math.SG/0707.3572v4, 2008.

\bibitem{kim}
Bumsig Kim.
\newblock Logarithmic stable maps.
\newblock In {\em New developments in algebraic geometry, integrable systems
  and mirror symmetry ({RIMS}, {K}yoto, 2008)}, volume~59 of {\em Adv. Stud.
  Pure Math.}, pages 167--200. Math. Soc. Japan, Tokyo, 2010.

\bibitem{KnudsenDM}
Finn~F. Knudsen.
\newblock The projectivity of the moduli space of stable curves. {II}. {T}he
  stacks {$M_{g,n}$}.
\newblock {\em Math. Scand.}, 52(2):161--199, 1983.

\bibitem{orbifoldstack}
Eugene Lerman.
\newblock Orbifolds as stacks?
\newblock {\em Enseign. Math. (2)}, 56(3-4):315--363, 2010.

\bibitem{Tian-Li}
Jun Li and Gang Tian.
\newblock Virtual moduli cycles and {G}romov-{W}itten invariants of general
  symplectic manifolds.
\newblock In {\em Topics in symplectic {$4$}-manifolds ({I}rvine, {CA}, 1996)},
  First Int. Press Lect. Ser., I, pages 47--83. Int. Press, Cambridge, MA,
  1998.

\bibitem{MW2}
Dusa McDuff and Katrin Wehrheim.
\newblock The fundamental class of smooth {K}uranishi atlases with trivial
  isotropy.
\newblock arXiv:1508.01560.

\bibitem{MW3}
Dusa McDuff and Katrin Wehrheim.
\newblock Smooth {K}uranishi atlases with isotropy.
\newblock arXiv:1508.01556.

\bibitem{KMW}
Dusa McDuff and Katrin Wehrheim.
\newblock Smooth {K}uranishi structures with trivial isotropy.
\newblock \href{http://arxiv.org/abs/1208.1340}{arXiv:1208.1340}.

\bibitem{pardon}
John Pardon.
\newblock An algebraic approach to virtual fundamental cycles on moduli spaces
  of pseudo-holomorphic curves.
\newblock {\em Geom. Topol.}, 20(2):779--1034, 2016.

\bibitem{scgp}
Brett Parker.
\newblock Notes on exploded manifolds and a tropical gluing formula for
  {G}romov-{W}itten invariants.
\newblock \href{http://arxiv.org/abs/1605.00577}{arXiv:1605.00577}.

\bibitem{iec}
Brett Parker.
\newblock Exploded manifolds.
\newblock {\em Adv. Math.}, 229:3256--3319, 2012.
\newblock \href{http://arxiv.org/abs/0910.4201}{arXiv:0910.4201}.

\bibitem{elc}
Brett Parker.
\newblock Log geometry and exploded manifolds.
\newblock {\em Abh. Math. Sem. Hamburg}, 82:43--81, 2012.
\newblock \href{http://arxiv.org/abs/1108.3713}{arxiv:1108.3713}.

\bibitem{icc}
Brett Parker.
\newblock Integral counts of pseudo-holomorphic curves.
\newblock \href{http://arxiv.org/abs/1309.0585}{arXiv:1309.0585}, 2013.

\bibitem{uts}
Brett Parker.
\newblock Universal tropical structures for curves in exploded manifolds.
\newblock \href{http://arxiv.org/abs/1301.4745}{arXiv:1301.4745}, 2013.

\bibitem{cem}
Brett Parker.
\newblock Holomorphic curves in exploded manifolds: compactness.
\newblock {\em Adv. Math.}, 283:377--457, 2015.
\newblock \href{http://arxiv.org/abs/0911.2241}{arXiv:0911.2241}.

\bibitem{gfgw}
Brett Parker.
\newblock Tropical gluing formulae for {G}romov-{W}itten invariants.
\newblock \href{http://arxiv.org/abs/1703.05433}{arXiv:1703.05433}, 2017.

\bibitem{dre}
Brett Parker.
\newblock De {R}ham theory of exploded manifolds.
\newblock {\em Geometry and Topology}, 22(1):1--54, 2018.
\newblock \href{http://arxiv.org/abs/1003.1977}{arXiv:1003.1977}.

\bibitem{reg}
Brett Parker.
\newblock Holomorphic curves in exploded manifolds: regularity.
\newblock {\em Geom. Topol.}, 23(4):1621--1690, 2019.

\bibitem{vfc}
Brett Parker.
\newblock Holomorphic curves in exploded manifolds: virtual fundamental class.
\newblock {\em Geom. Topol.}, 23(4):1877--1960, 2019.

\bibitem{SalamonDM}
Joel~W. Robbin and Dietmar~A. Salamon.
\newblock A construction of the {D}eligne-{M}umford orbifold.
\newblock {\em J. Eur. Math. Soc. (JEMS)}, 8(4):611--699, 2006.

\end{thebibliography}
\end{document}